\documentclass[12pt,twoside,reqno,openany]{amsart}

\usepackage{amsbsy,amscd,amsfonts,amsmath,amssymb,amsthm,color,
fancybox,fancyhdr,footmisc,graphics,graphicx,ifthen,mathrsfs,
multicol,pdfpages,rotating,times,wasysym}
\usepackage[dvipsnames,svgnames,x11names]{xcolor}

\usepackage{pifont}

\usepackage[all]{xy}
\usepackage[utf8]{inputenc}
\usepackage[T1]{fontenc}
\usepackage{mathtools}
\newtagform{EngelLie}[\scriptstyle]{$}{$}
\makeatletter\newcommand{\leqnomode}{\tagsleft@true}
\newcommand{\reqnomode}{\tagsleft@false}\makeatother

\newtheorem{Theorem}[equation]{Theorem}

\newtheorem{Proposition}[equation]{Proposition}

\newtheorem{Lemma}[equation]{Lemma}

\newtheorem{Assertion}[equation]{Assertion}

\newtheorem{Observation}[equation]{Observation}

\theoremstyle{definition} 

\newtheorem{Hypothesis}[equation]{Hypothesis}

\newtheorem{Notation}[equation]{Notation}

\newtheorem{Convention}[equation]{Convention}

\newtheorem{Question}[equation]{Question}

\newcommand{\C}{\mathbb{C}}

\newcommand{\R}{\mathbb{R}}

\newcommand{\FF}{\text{\sc f}}

\newcommand{\LL}{\text{\sc l}}

\newcommand{\kaux}{{\text{\usefont{T1}{qcs}{m}{sl}k}}}

\newcommand{\paux}{{\text{\usefont{T1}{qcs}{m}{sl}p}}}

\newcommand{\Aaux}{{\text{\usefont{T1}{qcs}{m}{sl}A}}}
\newcommand{\Baux}{{\text{\usefont{T1}{qcs}{m}{sl}B}}}
\newcommand{\Caux}{{\text{\usefont{T1}{qcs}{m}{sl}C}}}

\newcommand{\Faux}{{\text{\usefont{T1}{qcs}{m}{sl}F}}}
\newcommand{\Gaux}{{\text{\usefont{T1}{qcs}{m}{sl}G}}}
\newcommand{\Haux}{{\text{\usefont{T1}{qcs}{m}{sl}H}}}

\newcommand{\Jaux}{{\text{\usefont{T1}{qcs}{m}{sl}J}}}
\newcommand{\Kaux}{{\text{\usefont{T1}{qcs}{m}{sl}K}}}

\newcommand{\Paux}{{\text{\usefont{T1}{qcs}{m}{sl}P}}}

\newcommand{\Raux}{{\text{\usefont{T1}{qcs}{m}{sl}R}}}

\newcommand{\Uaux}{{\text{\usefont{T1}{qcs}{m}{sl}U}}}
\newcommand{\Vaux}{{\text{\usefont{T1}{qcs}{m}{sl}V}}}
\newcommand{\Waux}{{\text{\usefont{T1}{qcs}{m}{sl}W}}}
\newcommand{\Xaux}{{\text{\usefont{T1}{qcs}{m}{sl}X}}}
\newcommand{\Yaux}{{\text{\usefont{T1}{qcs}{m}{sl}Y}}}
\newcommand{\Zaux}{{\text{\usefont{T1}{qcs}{m}{sl}Z}}}

\definecolor{blue}{cmyk}{1.,1.,0.,0.63}
\definecolor{red}{cmyk}{0.,1.,1.,0.63}
\definecolor{green}{cmyk}{1.,0.,1.,0.63}
\definecolor{black}{cmyk}{1.,1.,1.,1.}

\newcommand{\blue}{\textcolor{blue}}

\newcommand{\red}{\textcolor{red}}

\makeatletter
\renewcommand{\@fnsymbol}[1]
{\ensuremath{\ifcase#1\or $*$\or $**$\or $***$\or $****$\or $*****$
\else\@ctrerr\fi}}
\makeatother

\newcommand{\HEAD}[2]{%
\pagestyle{fancy}
\fancyhead[RO]{\tiny\sf\thepage}
\fancyhead[CO]{{\tiny\sf #1}}
\fancyhead[LE]{\tiny\sf\thepage}
\fancyhead[CE]{{\tiny\sf #2}}
\fancyfoot{}}

\numberwithin{equation}{section}

\newcommand{\Section}[1]{
\renewcommand{\thesection}{\bf\arabic{section}}
\section{#1}
\renewcommand{\thesection}{\arabic{section}}}

\newcommand{\style}[1]{\text{\footnotesize{\sf #1}}}

\newcommand{\conjugates}{\style{conjugates}}

\renewcommand{\Im}{\style{Im}}

\newcommand{\Ker}{\style{Ker}}

\renewcommand{\lim}{\style{lim}}

\renewcommand{\mod}{\style{mod}}

\renewcommand{\Re}{\style{Re}}

\newcommand{\something}{\style{something}}

\newcommand{\Torsion}{\style{Torsion}}

\newcommand{\centersmallbullet}{{}_{{}^{{}^{
\scriptscriptstyle{\bullet\!}}}}}

\newcommand{\isqrt}{i}

\newcommand{\LF}{\LL\FF}

\newcommand{\linestop}{\medskip\centerline{\bf 
-----------------}\medskip}

\newcommand{\smallbullet}{{\scriptscriptstyle{\bullet}}}

\newcommand{\smallsum}[1]{
\underset{#1}{\raisebox{1pt}{$\sum$\,}}
}

\newcommand{\vf}{\vfill\end{document}}

\newcommand{\zero}[1]{\underline{#1}_{\red{\circ}}}

\newcommand{\zerozero}[1]{\underline{#1}_{\red{\circ\circ}}}

\begin{document}

\setcounter{section}{0}

$\:$

\bigskip\bigskip\bigskip\bigskip\bigskip

\begin{center}

{\large\bf Differential $\{e\}$-structures for equivalences}
\label{e-structure-M-5-C-3}

\medskip

{\large\bf of $2$-nondegenerate Levi rank $1$ 
hypersurfaces $M^5 \subset \C^3$}

\bigskip\bigskip

Wei Guo {\sc Foo}\footnote{Hua Loo-Keng Center for Mathematical
Sciences, AMSS, Chinese Academy of Sciences, Beijing, China.} 
and Joël~{\sc Merker}\footnote{Laboratoire de Mathématiques d'Orsay,
Université Paris-Sud, CNRS, Université Paris-Saclay, 91405 Orsay Cedex,
France.}

\end{center}\bigskip

\begin{center}
\begin{minipage}[t]{12.5cm}
\parindent 0.53cm
\footnotesize
\noindent
{\sc Abstract}.
The class $\text{\sf IV}_2$ of $2$-nondegenerate constant Levi rank
$1$ hypersurfaces $M^5 \subset \C^3$ is governed by Pocchiola's two
primary invariants $\Waux_0$ and $\Jaux_0$. Their vanishing
characterizes equivalence of 
such a hypersurface $M^5$ to the tube $M_{\sf LC}^5$ over the
real light cone in $\R^3$. When either $\Waux_0 \not\equiv 0$ or
$\Jaux_0 \not\equiv 0$, by normalization of certain two group
parameters ${\sf c}$ and ${\sf e}$, an invariant coframe can be built
on $M^5$, showing that the dimension of the CR automorphism
group drops from $10$ to $5$.

This paper constructs an explicit $\{e\}$-structure in case $\Waux_0$
and $\Jaux_0$ do not necessarily vanish.  Furthermore, Pocchiola's
calculations hidden on a computer now appear in details, especially
the determination of a secondary invariant $R$, expressed in terms of
the first jet of $\Waux_0$. All other secondary invariants
of the $\{e\}$-structure are also expressed explicitly
in terms of $\Waux_0$ and $\Jaux_0$.
\end{minipage}
\end{center}

\Section{\bf Introduction}
\label{introduction-e-structure-Pocchiola}
\HEAD{{\ref{introduction-e-structure-Pocchiola}}.~{\sf Introduction}
}{
Wei Guo {\sc Foo} (Beijing) and Joël {\sc Merker} (Orsay)}

We study the equivalence problem under biholomorphisms of real
hypersurfaces $M^{5}\subset\C^3$\,\,---\,\,hence of CR dimension
$2$\,\,---\,\,whose Levi form is degenerate of constant rank $1$, and
whose Freeman form is nowhere zero, or equivalently, which are
$2$-nondegenerate. There are previous approaches to this problem, and
we refer our readers to the article of Isaev et al.
\cite{Isaev-Zaitsev-2013}, and also to the article of Medori-Spiro
\cite{Medori-Spiro-2014, Medori-Spiro-2015}, in which a {\em Cartan
connection} was constructed.

In a recently published article
\cite{Merker-Pocchiola-2018}, the authors exhibited two important
primary invariants, $\Waux_{0}$ and $\Jaux_{0}$, whose existence was
not previously discovered prior to Pocchiola's
prepublication~{\cite{Pocchiola-2013}}, and which, in depth, required
the help of a computer algebra system. These invariants have useful
applications, such as in Isaev's study~{\cite{Isaev-2018}} of tube
hypersurfaces in $\C^3$ that are $2$-nondegenerate and uniformly Levi
degenerate of rank 1.

Our first objective here is to reconstruct $\Waux_{0}$ and
$\Jaux_{0}$, by presenting fully detailed computations, only by hand,
without the help of any computer.  In contrast
to~{\cite{Pocchiola-2013, Merker-Pocchiola-2018}}, the present text
has the ambition of exhibiting all calculations, without requiring any
extra work from the readers: `{\sl no pen needed, no computer
needed}'. Within the Cartan theory, this sounds quite like a challenge
opposite to a certain tradition of hiding a lot of computations.
But we believe that fully detailed articles can be read,
checked and studied more rapidly.

As a second objective, we construct an {\em explicit}
$\{e\}$-structure which characterizes equivalences under
biholomorphisms of these types of hypersurfaces $M^5 \subset \C^3$.
This way, we give a theoretical proof which will provide a definitive
confirmation of the existence of {\em exactly $2$} primary invariants,
$\Waux_{0}$ and $\Jaux_{0}$. Unlike the approach
of~{\cite{Pocchiola-2013, Merker-Pocchiola-2018}} which proceeded at
each step with {\em systematic} and {\em explicit} calculations of
{\em all} torsion coefficients, we will bypass some of these steps,
thereby economizing some computations.  On the way, we will closely
observe the evolution of the modified
Maurer-Cartan $1$-forms during the Cartan process.

The basic principle of Cartan's approach is to create a collection of
$1$-forms (a coframe), by absorbing as many as possible torsion terms,
in order that the structure of this coframe be as close as possible to
the structure of the Maurer-Cartan coframe on the 
(prolongation of the) model $M_{\sf LC}^5
\subset \C^3$, the tube over the real light cone 
$\big\{ x_1^2 + x_2^2 = x_3^2\}$ in $\R^3$: 
\[
M_{\sf LC}^5
\,:=\,
\big\{
(z_1,z_2,z_3)
\in\C^3
\colon\,
(\Re\,z_1)^2
+
(\Re\,z_2)^2
=
(\Re\,z_3)^2
\big\},
\]
whose local CR automorphism group is known to be isomorphic to ${\sf
SO}_{3,2}(\R)$.

Recall that a Maurer-Cartan form $\omega$
valued in some Lie algebra $\mathfrak{g}$
satisfies the structure equation with no curvature:
\[
d\omega
+
{\textstyle{\frac{1}{2}}}
\big[
\omega
\wedge
\omega
\big]
\,=\,
0.
\]
In practice, as in our current case, the right-hand
side of the equation is
not always zero, and this constitutes the {\sl default} of $\omega$
being a Maurer-Cartan form. This happens when an invariant is written
as a linear combination of torsion terms, and such a linear combination
fails to follow the structure equations, thus obstructing the
absorption process.

We now give a summary of our results. Recall that if $J$ denotes the
complex structure of $T\C^3$, then the tangent bundle $TM^5$ has a
distribution $T^cM^5 := TM^5\cap JTM^5\subseteq TM^5$ of codimension $1$
which is invariant under $J$ at each point of $M^5$. Let $\rho$ be a
real $1$-form with $\Ker\, \rho=T^cM^5$. The {\sl Levi form} is a
bilinear map on $T^cM^5$ defined as $(X,Y)\mapsto d\rho(X,JY)$ for any
two sections $X$, $Y$ of $T^cM^5$.

Letting $\C TM^5 := \C \otimes_\R TM^5$ be the complexification of the
tangent bundle of $M^5$, by defining $T^{1,0}M^5 := \C TM^5\cap
T^{1,0}\C^3$ together with its complex conjugate $T^{0,1}M^5 :=
\overline{T^{1,0}M^5}$, we have the (classical) direct sum
decomposition $\C T^cM^5=T^{1,0}M^5\oplus T^{0,1}M^5$. Let
$\{\mathcal{L}_1,\mathcal{L}_2\}$ be two
local generators of $T^{1,0}M^5$,
{\em i.e.} a frame for $T^{1,0}M^5$.

Section~{\ref{local-geometry-hypersurfaces-C-3}} 
provides more information, while 
complete background may be found 
in~{\cite{Merker-Pocchiola-Sabzevari-2013-5-CR-II}}.

By the assumption that the Levi form is uniformly of rank $1$
at each point of $M$, there exists 
by~{\cite{Merker-Pocchiola-Sabzevari-2013-5-CR-II}} 
a uniquely determined
{\sl slant function} $\kaux \colon M \longrightarrow \C$ 
such that the vector field:
\[
\mathcal{K}
\,:=\,
\kaux\,
\mathcal{L}_{1}
+
\mathcal{L}_{2}
\]
generates the kernel of the Levi form, of constant
rank $2 - 1 = 1$.  If we let $\mathcal{T}$ denote
a vector field with $\rho (\mathcal{T}) \equiv 1$, we may consider the
coframe $\big\{ \rho,\kappa_{0},\zeta_{0} \big\}$ dual to $\big\{
\mathcal{T},\mathcal{L}_{1},\mathcal{K} \big\}$. 
In fact, the conjugates $\overline{\kappa}_0$, $\overline{\zeta}_0$
and $\overline{\mathcal{L}}_1$, $\overline{\mathcal{K}}$
also come into play in order to really
make up a (co)frame on $\C TM^5$, while $\overline{\rho} = \rho$
and $\overline{\mathcal{T}} = \mathcal{T}$ are real.
A certain appropriate real $1$-form $\rho$ will be chosen,
and denoted $\rho_0$.

Performing the Cartan
process, we will make a series of changes to these 1-forms:
\[
(\rho_0,\kappa_{0},\zeta_{0})
\,\leadsto\,
(\rho_0,\kappa_{0}',\zeta_{0}''),
\]
and after (really a lot of)
computations, we will
obtain a $4$-dimensional ${\sf G}$-structure
whose lifted $1$-forms write up as:
\[
\left(
\begin{matrix}
\rho\\
\kappa\\
\zeta
\end{matrix}
\right)
:=
\left(
\begin{matrix}
{\sf c}\overline{\sf c} & 0 & 0\\
-i\,\overline{\sf c}\sf{e} & {\sf c} & 0\\
-\frac{i}{2}\frac{\overline{\sf c}\sf{e}\sf{e}}{\sf c} & \sf{e} & \frac{\sf c}{\overline{\sf c}}
\end{matrix}
\right)
\left(
\begin{matrix}
\rho_0\\
\kappa_{0}'\\
\zeta_{0}''
\end{matrix}
\right).
\]
Also, after a long process, we will construct 
modified Maurer-Cartan forms:
\[
\aligned
\pi^1
&
\,:=\,
\alpha
-
\Big(
{\sf t}
-
\frac{\isqrt}{2}\,
\Im\,Z^2
\Big)\,
\rho
-
\Big(
R^1
-
\overline{K}^6
\Big)\,
\kappa
-
R^2\,
\zeta
-
K^6\,\overline{\kappa}
-
0,
\\
\pi^2
&
\,:=\,
\beta
-
\isqrt\,Z^1\,\rho
-
\Big(
{\sf t}
-
\frac{\isqrt}{2}\,
\Im\,Z^2
+
K^1
\Big)\,
\kappa
-
K^2\,\zeta
-
K^3\,\overline{\kappa}
-
K^4\,\overline{\zeta},
\endaligned
\]
with $R^i$, $K^i$, $Z^i$ being some explicit functions on $M^5 \times
G^4$, where ${\sf t}$ is a new real variable, and then, after
meticulous absorption work, we will obtain as is stated below in
Theorem~{\ref{Theorem-finalization-absorption-d-rho-d-kappa-d-zeta}}
on p.~{\pageref{Theorem-finalization-absorption-d-rho-d-kappa-d-zeta}},
three finalized structure equations of the neat shape:
\[
\aligned
d\rho
&
\,=\,
\big(
\pi^1
+
\overline{\pi}^1
\big)
\wedge\rho
+
\isqrt\,\kappa\wedge\overline{\kappa},
\notag
\\
d\kappa
&
\,=\,
\pi^2\wedge\rho
+
\pi^1\wedge\kappa
+
\zeta\wedge\overline{\kappa},
\notag
\\
d\zeta
&
\,=\,
\big(\pi^1-\overline{\pi}^1\big)
\wedge\zeta
+
\isqrt\,\pi^2\wedge\kappa
\,+
\notag
\\
&
\ \ \ \ \
+
R\,
\rho\wedge\zeta
+
J\,\rho\wedge\overline{\kappa}
+
W\,
\kappa\wedge\zeta,
\endaligned
\]
in which are present Pocchiola's two primary invariants:
\[
W
\,=\,
\frac{1}{\sf c}\,\Waux_{0}
\ \ \ \ \ \ \ \ \ \ \ \ \ \ \ \ \ \
\text{and}
\ \ \ \ \ \ \ \ \ \ \ \ \ \ \ \ \ \
J
\,=\,
\frac{i}{{\overline{\sf c}}^3}\,
\overline{\Jaux}_{0},
\]
together with a single secondary (derived) invariant:
\[
R
\,=\,
\Re\,
\left[
\isqrt\,
\frac{{\sf e}}{{\sf c}{\sf c}}\,
\Waux_0
+
\frac{1}{{\sf c}\overline{\sf c}}
\bigg(
-\,\frac{\isqrt}{2}\,
\overline{\mathcal{L}}_1\big(\Waux_0\big)
+
\frac{\isqrt}{2}\,
\bigg(
-\,\frac{1}{3}\,
\frac{\overline{\mathcal{L}}_1\big(
\overline{\mathcal{L}}_1(\kaux)\big)}{
\overline{\mathcal{L}}_1(\kaux)}
+
\frac{1}{3}\,
\overline{\Paux}
\bigg)\,
\Waux_0
\bigg)
\right].
\]

We would like to mention that the two invariants that Pocchiola
denoted $W$ and $J$ are now 
denoted in our paper $\Waux_0$ and $\Jaux_0$,
with the subscript $(\centersmallbullet)_0$ designating functions
defined on $M^5$ alone, independently of any extra group variable.

The expression of $R$ was discovered by Pocchiola
in~{\cite{Pocchiola-2013, Merker-Pocchiola-2018}} 
thanks to intensive
computer explorations, but no details of proof appeared in print at
all.  In Section~{\ref{computation-R}} of this paper, a complete, {\em
detailed}, {\em hand-done} proof, will be provided, thus verifying
that $R$ is indeed a function of the first jet of $\Waux_{0}$, hence a
{\sl secondary} invariant.

We will also construct a certain real $1$-form $\Lambda = d{\sf t} +
\cdots$, and in Section~{\ref{final-e-structure}}, the final
$\{e\}$-structure that we obtain will take the following form
(conjugate equations are unwritten):
\begin{equation*}
\begin{aligned}
d\rho &= 
\pi^{1}\wedge\rho+\overline{\pi}^{1}\wedge\rho+i\kappa\wedge\overline{\kappa},\\
d\kappa &= \pi^{1}\wedge\kappa+\pi^{2}\wedge\rho+\zeta\wedge\overline{\kappa},\\
d\zeta &= i\pi^{2}\wedge\kappa+\pi^{1}\wedge\zeta-\overline{\pi}^{1}\wedge\zeta
+W\kappa\wedge\zeta+R\rho\wedge\zeta+J\rho\wedge\overline{\kappa},\\
d\pi^{1} &= \Lambda\wedge\rho-i\overline{\pi}^{2}\wedge\kappa+\zeta\wedge\overline{\zeta}+\widehat{\Omega}_{1},\\
d\pi^{2} &= \Lambda\wedge\kappa+\pi^{2}\wedge\overline{\pi}^{1}-\overline{\pi}^{2}\wedge\zeta+\widehat{\Omega}_{2}+h\rho\wedge\kappa,\\
d\Lambda &= \Lambda\wedge\pi^{1}+\Lambda\wedge\overline{\pi}^{1}+i\pi^{2}\wedge\overline{\pi}^{2}+\Phi,
\end{aligned}
\end{equation*}
with:
\begin{equation*}
\begin{aligned}
\widehat{\Omega}_{1} &=
-{\textstyle{\frac{1}{4}}}W\pi^{2}\wedge\rho
+
{\textstyle{\frac{1}{4}}}\overline{W}\overline{\pi}^{2}\wedge\rho
-
{\textstyle{\frac{1}{2}}}(R_{\kappa}-\overline{J_{\zeta}})\rho\wedge\kappa
-
{\textstyle{\frac{1}{2}}}R_{\zeta}\rho\wedge\zeta\\
&\hspace{0.5cm}
+
{\textstyle{\frac{1}{2}}}(R_{\overline{\kappa}}-J_{\zeta})\rho\wedge\overline{\kappa}
+
{\textstyle{\frac{1}{2}}}R_{\overline{\zeta}}\rho\wedge\overline{\zeta}
+
\bigg({\textstyle{\frac{1}{2}}}W_{\overline{\kappa}}-iR\bigg)\kappa\wedge\overline{\kappa}
-\overline{W}\kappa\wedge\overline{\zeta}
-W\zeta\wedge\overline{\kappa},\\
\widehat{\Omega}_{2}
&=
-R\pi^{2}\wedge\rho
-{\textstyle{\frac{1}{4}}}W\pi^{2}\wedge\kappa
+{\textstyle{\frac{1}{4}}}\overline{W}\overline{\pi}^{2}\wedge\kappa
-i(W_{\rho}-2R_{\kappa}+\overline{J_{\zeta}})\rho\wedge\zeta\\
&\hspace{0.5cm}
-i(WJ-J_{\kappa})\rho\wedge\overline{\kappa}
-iJ\rho\wedge\overline{\zeta}
-{\textstyle{\frac{1}{2}}}R_{\zeta}\kappa\wedge\zeta
+{\textstyle{\frac{1}{2}}}(R_{\overline{\kappa}}-J_{\zeta})\kappa\wedge\overline{\kappa}
+{\textstyle{\frac{1}{2}}}R_{\overline{\zeta}}\kappa\wedge\overline{\zeta}\\
&\hspace{0.5cm}
-R\zeta\wedge\overline{\kappa}.
\end{aligned}
\end{equation*}

Furthermore, we will show that $h$ and $\Phi$ can be expressed in
terms of $\widehat{\Omega}_{1}$, of $\widehat{\Omega}_{2}$ and of
their first-order derivatives. Thus, this demonstrates that there are
exactly $2$ primary invariants.

Clearly, when $W \equiv J \equiv 0$, the $\{e\}$-structure collapses
to:
\begin{equation*}
\begin{aligned}
d\rho &= 
\pi^{1}\wedge\rho+\overline{\pi}^{1}\wedge\rho+i\kappa\wedge\overline{\kappa},\\
d\kappa &= \pi^{1}\wedge\kappa+\pi^{2}\wedge\rho+\zeta\wedge\overline{\kappa},\\
d\zeta &= i\pi^{2}\wedge\kappa+\pi^{1}\wedge\zeta-\overline{\pi}^{1}\wedge\zeta,\\
d\pi^{1} &= \Lambda\wedge\rho-i\overline{\pi}^{2}\wedge\kappa+\zeta\wedge\overline{\zeta},\\
d\pi^{2} &= \Lambda\wedge\kappa+\pi^{2}\wedge\overline{\pi}^{1}-\overline{\pi}^{2}\wedge\zeta,\\
d\Lambda &= \Lambda\wedge\pi^{1}+\Lambda\wedge\overline{\pi}^{1}+i\pi^{2}\wedge\overline{\pi}^{2},
\end{aligned}
\end{equation*}
and these constant coefficients equations correspond to the structure
equations of the tube $M_{\sf LC}^5$ over the light cone, which
is the reference model for this equivalence problem.

We would like to mention that, strictly speaking, Cartan's equivalence
method of producing homogeneous models {\em requires} to normalize any
group variable which occurs in some essential torsion term, and this is
what Pocchiola did in Section~7 of~{\cite{Pocchiola-2013}} for ${\sf
c} := (\Jaux_0)^{1/3}$ and in Section~8 for ${\sf c} := \Waux_0$,
showing afterwards that ${\sf e}$ can also be normalized in both
cases.

For this deep reason, Pocchiola then {\em disregarded}
the\,\,---\,\,essentially useless\,\,---\,\,task of constructing a
general $\{e\}$-structure, since, when $\Jaux_0 \equiv \Waux_0 \equiv
0$, the final Section~9 of~{\cite{Pocchiola-2013}} shows that one
comes uniquely to the structure equations of the model $M_{\sf LC}^5$,
{\em without any further nonzero essential torsion appearing}.  And
this was really a discovery, because most of the times in CR geometry,
primary invariants appear {\em after} a first prolongation.

However, because there is a tradition of setting up $\{e\}$-structures,
even in absence of explicit computations, even without discovering
invariants at all, and because the needs for {\em verifiable}
computations has been expressed by some experts, we decided to set up
the present article. While re-building 
this chapter~{\cite{Pocchiola-2013}}
of Pocchiola's Ph.D. 
(Orsay University, September 2014), 
we found a few
copying mistakes 
in some intermediate formulas of~{\cite{Pocchiola-2013,
Merker-Pocchiola-2018}}, but no error in either statements or
final formulas, {\em e.g.} $\Waux_0$ and $\Jaux_0$
are correct.

For a more informative exposition of introductory aspects, the reader
should read now the brief and complementary {\sl Introduction} to the
Addendum to~{\cite{Merker-Pocchiola-2018}}, reproduced as an Appendix,
after the end, on p.~{\pageref{addendum-Merker-Pocchiola-JGA}}.

\smallskip

This paper is organized as follows. In 
Section~{\ref{local-geometry-hypersurfaces-C-3}}, 
we recall the local
geometry of $2$-nondegenerate Levi rank 1 real hypersurfaces $M^5$ in
$\C^3$. In 
Section~{\ref{initial-G-1-structure-biholomorphic-equivalences}}, 
we give a description of the
$G_{1}$-structure of the biholomorphic equivalences of such real
hypersurfaces. 
Section~{\ref{labyrinthmap-Pocchiola-calculations}} 
gives a quick glimpse of a series of
normalizations of parameters, which will be detailed in 
Sections~{\ref{first-loop-reduction-f}} 
to~{\ref{absorption-pi-1-pi-2}},
with the first appearance of $\Waux_{0}$ in 
Section~{\ref{third-loop-reduction-d}}. 
The explicit expression of the invariant $\Jaux_{0}$ is given in 
Section~{\ref{computation-Jaux-0}}, 
and a complete proof of the above formula for $R$ 
is detailed in 
Section~{\ref{computation-R}}. 
Section~{\ref{summarized-structure-equations}}
gives a short summary of the things that have been done in the previous
sections, and finally 
Section~{\ref{final-e-structure}}
gives a proposed $\{e\}$-structure for the equivalence problem.

\medskip\noindent
{\bf Acknowledgments.}
Both authors benefited from enlightening exchanges with 
Pawe{\l} Nu\-row\-ski.

\Section{\bf Local Geometry of $2$-nondegenerate Levi rank $1$
hypersurfaces $M^5 \subset \C^3$}
\label{local-geometry-hypersurfaces-C-3}
\HEAD{{\ref{local-geometry-hypersurfaces-C-3}}.~{\sf Local Geometry 
of $2$-nondegenerate Levi rank $1$ hypersurfaces $M^5 \subset \C^3$}
}{
Wei Guo {\sc Foo} (Beijing) and Joël {\sc Merker} (Orsay)}

This section only summarizes what has been presented and
detailed in~{\cite{Merker-Pocchiola-Sabzevari-2013-5-CR-II,
Merker-2013-5-CR-V, Merker-Pocchiola-2018}}.
Let $M^5 \subset \C^3$ be a $\mathcal{C}^\omega$ (real-analytic)
smooth, local or global, real hypersurface and let $p_0 \in M$. In any
affine holomorphic coordinate system:
\[
\big(z_1,z_2,\,w\big)
\,\in\,
\C^3
\ \ \ \ \ \ \ \ \ \ \ \ \ \ \ \ \ \
\text{with}
\ \ \ \ \ \ \ \ \ \ \ \ \ \ \ \ \ \
w
\,=\,
u
+
\isqrt\,v,
\]
centered at $p_0 = (0,0,0) = 0$ 
in which $\frac{\partial}{\partial u}\big\vert_0 \not\in T_0 M$,
there is a local $\mathcal{C}^\omega$ graphing function $\Faux = 
\Faux \big(z_1, z_2, \overline{z}_1, \overline{z}_2, v \big)$ 
with $\Faux (0) = 0$ such that $M$ is represented, in some (possibly
small) open neighborhood of the origin $0$ by:
\[
u
\,=\,
\Faux
\big(
z_1,z_2,\overline{z}_1,\overline{z}_2,v
\big).
\]

\begin{Convention}
From now on, the hypersurface will be identified with its
localization in some small open neighborhood of the origin,
and it will always be denoted by $M$.
\end{Convention}

As is known ({\em
see}~{\cite{Merker-Pocchiola-Sabzevari-2013-5-CR-II}} for detailed
background), the complexified tangent bundle $\C TM := \C \otimes_\R
TM$ inherits from $\C T\C := \C \otimes_\R T\C^3$ two
biholomorphically invariant complex rank $2$
vector subbundles:
\[
T^{1,0}M
\,:=\,
T^{1,0}\C^3
\cap 
\C TM
\ \ \ \ \ \ \ \ \ \ \ \ \ \ \ \ \ \
\text{and}
\ \ \ \ \ \ \ \ \ \ \ \ \ \ \ \ \ \
T^{0,1}M
\,:=\,
T^{0,1}\C^3
\cap 
\C TM
\,=\,
\overline{T^{1,0}M},
\]
which are conjugate one to another. Then a check shows that the two
vector fields written in the intrinsic coordinates $(z_1, z_2,
\overline{z}_1, \overline{z}_2, v)$ on $M$:
\[
\mathcal{L}_1
\,:=\,
\frac{\partial}{\partial z_1}
+
\Aaux^1\,
\frac{\partial}{\partial v}
\ \ \ \ \ \ \ \ \ \ \ \ \ \ \ \ \ \
\text{and}
\ \ \ \ \ \ \ \ \ \ \ \ \ \ \ \ \ \
\mathcal{L}_2
\,:=\,
\frac{\partial}{\partial z_2}
+
\Aaux^2\,
\frac{\partial}{\partial v},
\]
whose coefficients are defined by:
\[
\Aaux^i
\,:=\,
-\,\isqrt\,
\frac{\Faux_{z_i}}{
1+\isqrt\,\Faux_v}
\eqno
{\scriptstyle{(i\,=\,1,\,2)}},
\]
generate $T^{1,0}M$, locally. Hence their two conjugates 
$\overline{\mathcal{L}}_1$, $\overline{\mathcal{L}}_2$
generate the bundle $T^{0,1}M$, also of complex rank $2$.

Then visibly the differential $1$-form:
\[
\varrho_0
\,:=\,
dv
-
\Aaux^1\,dz_1
-
\Aaux^2\,dz_2
-
\overline{\Aaux}^1\,
d\overline{z}_1
-
\overline{\Aaux}^2\,
d\overline{z}_2
\]
has kernel:
\[
\big\{
\varrho_0
=
0
\big\}
\,=\,
T^{1,0}M
\oplus
T^{0,1}M.
\]
There are various (equivalent) aspects of the concept of {\sl Levi
form} of $M$, but they will not be recalled here, since several
sources treat that.  Here, the Levi form of $M$ can be represented as
a function of the points:
\[
p
\,=\,
\big(z_1, z_2, \overline{z}_1, \overline{z}_2, v\big)
\,\in\,
M,
\]
valued in the space of Hermitian $2 \times 2$
matrices, and in terms of $\varrho_0$ and of
the Lie brackets of the above vector fields, it writes as:
\[
\LF_M(p)
\,:=\,
\left(\!
\begin{array}{cc}
\varrho_0\big(\isqrt\,[\mathcal{L}_1,\overline{\mathcal{L}}_1]\big)
&
\varrho_0\big(\isqrt\,[\mathcal{L}_2,\overline{\mathcal{L}}_1]\big)
\\
\varrho_0\big(\isqrt\,[\mathcal{L}_1,\overline{\mathcal{L}}_2]\big)
&
\varrho_0\big(\isqrt\,[\mathcal{L}_2,\overline{\mathcal{L}}_2]\big)
\end{array}
\!\right)
(p).
\]

As is known, the biholomorphic invariance of the Levi form 
legitimates our current

\begin{Hypothesis}
\label{Hypothesis-Levi-nondegeneracy}
{\bf [Uniform Levi rank $1$]}
At all points $p \in M$, the Levi matrix (form)
$\LF_M(p)$ has constant rank $1$.\hfill$\bigtriangleup$
\end{Hypothesis}

After a linear change of coordinates in the $(z_1, z_2)$ space,
we may assume that its $(1,1)$-entry vanishes nowhere on $M$:
\[
\varrho_0\big(\isqrt\,[\mathcal{L}_1,\overline{\mathcal{L}}_1]\big)(p)
\,\neq\,
0
\eqno
{\scriptstyle{(\forall\,p\,\in\,M)}}.
\]
This means that the {\em real} vector field:
\[
\mathcal{T}
\,:=\,
\isqrt\,
\big[\mathcal{L}_1,\overline{\mathcal{L}}_1\big]
\,=\,
\isqrt\,
\Big(
\mathcal{L}_1\big(\overline{\Aaux}^1\big)
-
\overline{\mathcal{L}}_1\big(\Aaux^1\big)
\Big)
\frac{\partial}{\partial v}
\,=:\,
\ell\,
\frac{\partial}{\partial v},
\]
has nowhere vanishing real coefficient that will be abbreviated as:
\[
\ell
\,:=\,
\isqrt\,
\Big(
\overline{\Aaux}_{z_1}^1
+
\Aaux^1\,\overline{\Aaux}_v^1
-
\Aaux_{\overline{z}_1}^1
-
\overline{\Aaux}^1\,
\Aaux_v^1
\Big)
\,\,\neq\,\,
0.
\]

Furthermore, since the $2 \times 2$ Levi matrix has constant rank $1$,
the collection of its $1$-dimensional kernels at all points $p \in M$ 
spans a $\mathcal{C}^\omega$
smooth subdistribution $K^{1,0} M \subset T^{1,0} M$ which
satisfies (\cite{ Merker-Pocchiola-Sabzevari-2013-5-CR-II}, pp.~72--73):
\reqnomode\usetagform{EngelLie}
\begin{align}
\big[K^{1,0}M,\,K^{1,0}M\big]
&
\,\subset\,
K^{1,0}M,
\notag
\\
\big[K^{0,1}M,\,K^{0,1}M\big]
&
\,\subset\,
K^{0,1}M
\tag{(K^{0,1}M\,:=\,\overline{K^{1,0}M}),}
\\
\big[K^{1,0}M,\,K^{0,1}M\big]
&
\,\subset\,
K^{1,0}M
\oplus
K^{0,1}M.
\notag
\end{align}

With this, 
a vector field generator $\mathcal{K}$ of $K^{1,0}M$ writes uniquely
as:
\[
\mathcal{K}
\,:=\,
\kaux\,\mathcal{L}_1
+
\mathcal{L}_2,
\]
where the function $\kaux$\,\,---\,\,very important in the 
theory\,\,---\,\,is the negative of the quotient 
of two entries of the Levi matrix:
\[
\kaux
\,:=\,
-\,
\frac{
\mathcal{L}_2\big(\overline{\Aaux}^1\big)
-
\overline{\mathcal{L}}_1\big(\Aaux^2\big)}{
\mathcal{L}_1\big(\overline{\Aaux}^1\big)
-
\overline{\mathcal{L}}_1\big(\Aaux^1\big)}.\,
\]

\begin{Hypothesis}
\label{Hypothesis-2-nondegeneracy}
{\bf [$2$-nondegeneracy]}
At all points $p \in M$, the {\sl Freeman form}
has constant (maximal possible) rank $1$.\hfill$\bigtriangleup$
\end{Hypothesis}

For a detailed presentation of this
second concept of form, also biholomorphically
invariant, 
{\em see}~{\cite{Merker-Pocchiola-Sabzevari-2013-5-CR-II}}.

\begin{Proposition}
{\rm ({\cite{Merker-Pocchiola-Sabzevari-2013-5-CR-II}})}
In this formalism, $M$ is $2$-nondegenerate if and only if:
\[
\overline{\mathcal{L}}_1(\kaux)
\,\neq\,
0
\eqno
{\scriptstyle{(\text{\rm everywhere on}\,\,M)}}.
\]
\end{Proposition}

In summary, {\em two} functions will be assumed to be nowhere
vanishing on $M$, corresponding to the two
Hypotheses~{\ref{Hypothesis-Levi-nondegeneracy}}
and~{\ref{Hypothesis-2-nondegeneracy}}:
\[
\ell(p)
\,\neq\,
0
\ \ \ \ \ \ \ \ \ \ \ \ \ \ \ \ \ \
\text{and}
\ \ \ \ \ \ \ \ \ \ \ \ \ \ \ \ \ \
\overline{\mathcal{L}}_1(\kaux)(p)
\,\neq\,
0
\eqno
{\scriptstyle{(\forall\,p\,\in\,M)}}.
\]

Next, along with $\kaux$, 
introduce a second and last fundamental function:
\[
\Paux
:=
\frac{\ell_{z_1}+\Aaux^1\,\ell_v-\ell\,\Aaux_v^1}{\ell}.
\] 
All invariants and semi-invariants in this paper will express
in terms of $\kaux$ and $\Paux$.

Next, according to~{\cite{Merker-2013-5-CR-V, Pocchiola-2013,
Merker-Pocchiola-2018}}, 
there are $10$ Lie bracket identities:
\[
\aligned
\big[\mathcal{T},\mathcal{L}_1\big]
&
=
-\,\Paux\cdot\mathcal{T},
\\
\big[\mathcal{T},\mathcal{K}\big]
&
=
\mathcal{L}_1(\kaux)\cdot\mathcal{T}
+
\mathcal{T}(\kaux)\cdot\mathcal{L}_1,
\\
\big[\mathcal{T},\overline{\mathcal{L}}_1\big]
&
=
-\,\overline{\Paux}\cdot\mathcal{T},
\\
\big[\mathcal{T},\overline{\mathcal{K}}\big]
&
=
\overline{\mathcal{L}}_1\big(\overline{\kaux}\big)
\cdot
\mathcal{T}
+
\mathcal{T}
\big(\overline{\kaux}\big)
\cdot
\overline{\mathcal{L}}_1,
\\
\big[\mathcal{L}_1,\mathcal{K}\big]
&
=
\mathcal{L}_1(\kaux)\cdot\mathcal{L}_1,
\\
\big[\mathcal{L}_1,\overline{\mathcal{L}}_1\big]
&
=
-\,\isqrt\,\mathcal{T},
\\
\big[\mathcal{L}_1,\overline{\mathcal{K}}\big]
&
=
\mathcal{L}_1\big(\overline{\kaux}\big)
\cdot\overline{\mathcal{L}}_1,
\\
\big[\mathcal{K},\overline{\mathcal{L}}_1\big]
&
=
-\,
\overline{\mathcal{L}}_1(\kaux)
\cdot\mathcal{L}_1,
\\
\big[\mathcal{K},\overline{\mathcal{K}}\big]
&
=
0,
\\
\big[\overline{\mathcal{L}}_1,\overline{\mathcal{K}}\big]
&
=
\overline{\mathcal{L}}_1\big(\overline{\kaux}\big)\cdot
\overline{\mathcal{L}}_1.
\endaligned
\]

\begin{Lemma}
\label{Lemma-K-bar-k-K-bar-P}
{\rm ({\cite{Merker-Pocchiola-Sabzevari-2013-5-CR-II,
Merker-2013-5-CR-V}})}
The following $3$ functional identities hold identically on $M$:
\reqnomode\usetagform{EngelLie}
\begin{align}
\mathcal{K}\big(\overline{\kaux}\big)
&
\,\equiv\,
0,
\notag
\\
\mathcal{K}(\Paux)
&
\,\equiv\,
-\,\Paux\,
\mathcal{L}_1(\kaux)
-
\mathcal{L}_1\big(\mathcal{L}_1(\kaux)\big),
\notag
\\
\mathcal{K}
\big(\overline{\Paux}\big)
&
\,\equiv\,
-\,\Paux\,
\overline{\mathcal{L}}_1(\kaux)
-
\overline{\mathcal{L}}_1\big(\mathcal{L}_1(\kaux)\big)
-
\isqrt\,\mathcal{T}(\kaux).
\tag{\qed}
\end{align}
\end{Lemma}

Then the coframe:
\[
\big\{
\rho_0,\,
\kappa_0,\,
\zeta_0,\,
\overline{\kappa}_0,\,
\overline{\zeta}_0
\big\}
\]
dual to the frame:
\[
\big\{
\mathcal{T},\,
\mathcal{L}_1,\,
\mathcal{K},\,\,
\overline{\mathcal{L}}_1,\,
\overline{\mathcal{K}}
\big\},
\]
{\em i.e.} which satisfies by definition:
\[
\begin{array}{ccccc}
\rho_0(\mathcal{T})=1, \ \ \ & \ \ \
\rho_0(\mathcal{L}_1)=0, \ \ \ & \ \ \
\rho_0(\mathcal{K})=0, \ \ \ & \ \ \
\rho_0(\overline{\mathcal{L}}_1)=0, \ \ \ & \ \ \
\rho_0(\overline{\mathcal{K}})=0,
\\
\kappa_0(\mathcal{T})=0, \ \ \ & \ \ \
\kappa_0(\mathcal{L}_1)=1, \ \ \ & \ \ \
\kappa_0(\mathcal{K})=0, \ \ \ & \ \ \
\kappa_0(\overline{\mathcal{L}}_1)=0, \ \ \ & \ \ \
\kappa_0(\overline{\mathcal{K}})=0, 
\\
\zeta_0(\mathcal{T})=0, \ \ \ & \ \ \
\zeta_0(\mathcal{L}_1)=0, \ \ \ & \ \ \
\zeta_0(\mathcal{K})=1, \ \ \ & \ \ \
\zeta_0(\overline{\mathcal{L}}_1)=0, \ \ \ & \ \ \
\zeta_0(\overline{\mathcal{K}})=0,
\\
\overline{\kappa}_0(\mathcal{T})=0, \ \ \ & \ \ \
\overline{\kappa}_0(\mathcal{L}_1)=0, \ \ \ & \ \ \
\overline{\kappa}_0(\mathcal{K})=0, \ \ \ & \ \ \
\overline{\kappa}_0(\overline{\mathcal{L}}_1)=1, \ \ \ & \ \ \
\overline{\kappa}_0(\overline{\mathcal{K}})=0,
\\
\overline{\zeta}_0(\mathcal{T})=0, \ \ \ & \ \ \
\overline{\zeta}_0(\mathcal{L}_1)=0, \ \ \ & \ \ \
\overline{\zeta}_0(\mathcal{K})=0, \ \ \ & \ \ \
\overline{\zeta}_0(\overline{\mathcal{L}}_1)=0, \ \ \ & \ \ \
\overline{\zeta}_0(\overline{\mathcal{K}})=1,
\end{array}
\]
has its $5$ component $1$-forms given explicitly by:
\[
\aligned
\rho_0
&
=
\frac{dv-\Aaux^1dz_1-\Aaux^2dz_2
-\overline{\Aaux}^1d\overline{z}_1
-\overline{\Aaux}^2d\overline{z}_2}{\ell},
\\
\kappa_0
&
=
dz_1-\kaux\,dz_2,
\\
\zeta_0
&
=
dz_2,
\\
\overline{\kappa}_0
&
=
d\overline{z}_1
-
\overline{\kaux}\,d\overline{z}_2,
\\
\overline{\zeta}_0
&
=
d\overline{z}_2.
\endaligned
\]
Notice that a different notation $\rho_0 \neq \varrho_0$ has
been employed just now.
Hence using a classical formula which goes back at least
to Lie (\cite[Chap.~5]{Lie-Merker-2015}) which holds for 
two arbitrary vector fields $X$ and $Y$  
and for any differential 1-form $\omega$:
\[
d\omega(X,Y)
\,=\,
X\,\big(\omega(Y)\big) 
- 
Y\big(\omega(X)\big) 
- 
\omega\big( 
\big[X,Y\big] 
\big),
\]
by representing the $10$ Lie brackets in some appropriate array:
\[
\footnotesize
\begin{array}{cccccccccccc}
& & \mathcal{T} & & \mathcal{L}_1 & &
\mathcal{K} & &
\overline{\mathcal{L}}_1 & & \overline{\mathcal{K}}
\\
& & \boxed{d\rho_0} & & \boxed{d\kappa_0} & &
\boxed{d\zeta_0} & & 
\boxed{d\overline{\kappa_0}} & & \boxed{d\overline{\zeta}_0}\medskip
\\
\big[\mathcal{T},\,\mathcal{L}_1\big] & = &
-\,\Paux\cdot\mathcal{T} & + &
0 & + & 0
& + & 0 & + & 0 & \boxed{\rho_0\wedge\kappa_0}
\\
\big[\mathcal{T},\,\mathcal{K}\big] & = &
\mathcal{L}_1(\kaux)\cdot\mathcal{T} & + &
\mathcal{T}(\kaux)\cdot\mathcal{L}_1 & + & 0 & + & 0 & + & 0 &
\boxed{\rho_0\wedge\zeta_0}
\\
\big[\mathcal{T},\,\overline{\mathcal{L}}_1\big] & = &
-\,\overline{\Paux}\cdot\mathcal{T} & + &
0 & + & 0 & + & 0 &
+ & 0 & \boxed{\rho_0\wedge\overline{\kappa}_0}
\\
\big[\mathcal{T},\,\overline{\mathcal{K}}\big] & = &
\overline{\mathcal{L}}_1\big(\overline{\kaux}\big)\cdot\mathcal{T} & + &
0 & + &
0 & + & \mathcal{T}\big(\overline{\kaux}\big)\cdot\overline{\mathcal{L}}_1 & + & 0 &
\boxed{\rho_0\wedge\overline{\zeta}_0}
\\
\big[\mathcal{L}_1,\,\mathcal{K}\big] & = & 
0 & + & \mathcal{L}_1(\kaux)\cdot\mathcal{L}_1 & + &
0 & + & 0 & + & 0 &
\boxed{\kappa_0\wedge\zeta_0}
\\
\big[\mathcal{L}_1,\,\overline{\mathcal{L}}_1\big] & = &
-\,\isqrt\,\cdot\mathcal{T} & + & 0 & + & 0 & + & 0 & + & 0 &
\boxed{\kappa_0\wedge\overline{\kappa}_0}
\\
\big[\mathcal{L}_1,\,\overline{\mathcal{K}}\big] & = &
0 & + & 0 & + & 0 & + & 
\mathcal{L}_1(\overline{\kaux})\cdot\overline{\mathcal{L}}_1 & + & 0 &
\boxed{\kappa_0\wedge\overline{\zeta}_0}
\\
\big[\mathcal{K},\,\overline{\mathcal{L}}_1\big] & = &
0 & + & -\,\overline{\mathcal{L}}_1(\kaux)\cdot\mathcal{L}_1 & + &
0 & + & 0 & + & 0 & 
\boxed{\zeta_0\wedge\overline{\kappa}_0}
\\
\big[\mathcal{K},\,\overline{\mathcal{K}}\big] & = & 0 & + & 
0 & +
& 0 & + & 0 & + & 0 & 
\boxed{\zeta_0\wedge\overline{\zeta}_0}
\\
\big[\overline{\mathcal{L}}_1,\,\overline{\mathcal{K}}\big] & = &
0 & + & 
0
\cdot\overline{\mathcal{L}}_1 & +
& 0 & + & \overline{\mathcal{L}}_1\big(\overline{\kaux}\big) & + & 0 &
\boxed{\overline{\kappa}_0\wedge\overline{\zeta}_0}
\end{array}
\]
and by reading this array {\em vertically}, 
we obtain the {\sl initial Darboux-Cartan structure:}
\[
\aligned 
d\rho_0
&
=
\Paux
\cdot
\rho_0\wedge\kappa_0
-
\mathcal{L}_1(\kaux)
\cdot
\rho_0\wedge\zeta_0
+
\overline{\Paux}
\cdot
\rho_0\wedge\overline{\kappa}_0
-
\overline{\mathcal{L}}_1\big(\overline{\kaux}\big)
\cdot
\rho_0\wedge\overline{\zeta}_0
+
\isqrt\,
\kappa_0\wedge\overline{\kappa}_0,
\\
d\kappa_0
&
=
-\,\mathcal{T}(\kaux)
\cdot
\rho_0\wedge\zeta_0
-
\mathcal{L}_1(\kaux)
\cdot
\kappa_0\wedge\zeta_0
+
\overline{\mathcal{L}}_1(\kaux)
\cdot
\zeta_0\wedge\overline{\kappa}_0,
\\
d\zeta_0
&
=
0,
\\
d\overline{\kappa}_0
&
=
-\,\mathcal{T}\big(\overline{\kaux}\big)
\cdot
\rho_0\wedge\overline{\zeta}_0
-
\mathcal{L}_1\big(\overline{\kaux}\big)
\cdot
\kappa_0
\wedge
\overline{\zeta}_0
-
\overline{\mathcal{L}}_1\big(\overline{\kaux}\big)
\cdot
\overline{\kappa}_0
\wedge
\overline{\zeta}_0,
\\
d\overline{\zeta}_0
&
=
0.
\endaligned
\]

The fact that the frame $\big\{ \mathcal{T}, \mathcal{L}_1,
\mathcal{K}, \overline{\mathcal{L}}_1, \overline{\mathcal{K}} \big\}$
is dual to the coframe $\big\{ \rho_0, \kappa_0, \zeta_0,
\overline{\kappa}_0, \overline{\zeta}_0 \big\}$ yields a formula
that shall be used several times later.

\begin{Lemma}
\label{Lemma-d-G-0}
The exterior differential of 
any function $\Gaux = \Gaux \big(z_1, z_2, \overline{z}_1,
\overline{z}_2, v \big)$ on $M$ expresses as:
\[
d\Gaux
\,=\,
\mathcal{T}\big(\Gaux\big)\,
\rho_0
+
\mathcal{L}_1\big(\Gaux\big)\,
\kappa_0
+
\mathcal{K}\big(\Gaux\big)\,
\zeta_0
+
\overline{\mathcal{L}}_1\big(\Gaux\big)\,
\overline{\kappa}_0
+
\overline{\mathcal{K}}\big(\Gaux\big)\,
\overline{\zeta}_0.
\]
\end{Lemma}

\proof
Indeed, starting from the definition:
\[
d\Gaux
\,=\,
\frac{\partial\Gaux}{\partial v}\,dv
+
\frac{\partial\Gaux}{\partial z_1}\,dz_1
+
\frac{\partial\Gaux}{\partial z_2}\,dz_2
+
\frac{\partial\Gaux}{\partial\overline{z}_1}\,d\overline{z}_1
+
\frac{\partial\Gaux}{\partial\overline{z}_2}\,d\overline{z}_2,
\]
and inverting the above coframe:
\[
\aligned
dz_2
&
\,=\,
\zeta_0,
\\
dz_1
&
\,=\,
\kappa_0
+
\kaux\,\zeta_0,
\\
dv
&
\,=\,
\ell\,\rho_0
+
\Aaux^1\,
\big(
\kappa_0
+
\kaux\,\zeta_0
\big)
+
\Aaux^2\,
\zeta_0
+
\overline{\Aaux}^1\,
\big(
\overline{\kappa}_0
+
\overline{\kaux}\,
\overline{\zeta}_0
\big)
+
\overline{\Aaux}^2\,
\overline{\zeta}_0
\\
&
\,=\,
\ell\,\rho_0
+
\Aaux^1\,\kappa_0
+
\big(
\Aaux^2
+
\kaux\,\Aaux^1
\big)\,\zeta_0
+
\conjugates
\endaligned
\]
we can replace, reorganize\,\,---\,\,unwritting
the redundant conjugates\,\,---\,\,and reach the formula:
\begin{footnotesize}
\begin{align}
d\Gaux
&
\,\equiv\,
\frac{\partial\Gaux}{\partial v}\,
\Big(
\ell\,\rho_0
+
\Aaux^1\,\kappa_0
+
\big(
\Aaux^2
+
\kaux\,\Aaux^1
\big)\,
\zeta_0
\Big)
+
\frac{\partial\Gaux}{\partial z_1}\,
\big(
\kappa_0
+
\kaux\,\zeta_0
\big)
+
\frac{\partial\Gaux}{\partial z_2}\,
\zeta_0
\notag
\\
&
\,\equiv\,
\bigg(
\ell\,
\frac{\partial}{\partial v}
\bigg)
\big(\Gaux\big)
\cdot
\rho_0
+
\bigg(
\frac{\partial}{\partial z_1}
+
\Aaux^1\,\frac{\partial}{\partial v}
\bigg)
\big(\Gaux\big)
\cdot
\kappa_0
+
\bigg(
\frac{\partial}{\partial z_1}
+
\Aaux^2\,\frac{\partial}{\partial v}
+
\kaux\,\frac{\partial}{\partial z_2}
+
\kaux\,\Aaux^1\,\frac{\partial}{\partial v}
\bigg)
\big(\Gaux\big)
\cdot
\zeta_0.
\qedhere
\end{align}
\end{footnotesize}
\endproof

For later much deeper computations, we need strong notational
conventions. The order succession for our five $1$-forms
which we will constantly use:
\[
\big\{
\rho_0,\,
\kappa_0,\,
\zeta_0,\,
\overline{\kappa}_0,
\overline{\zeta}_0,
\big\},
\]
induces an order succession for
the ten generated $2$-forms on the $5$-dimensional CR manifold $M$:
\[
\aligned
\underset{{\text{\bf 1}}}{
\rho_0\wedge\kappa_0}
\ \ \ \ \ \ \ \ \ \
\underset{{\text{\bf 2}}}{
\rho_0\wedge\zeta_0}
\ \ \ \ \ \ \ \ \ \
\underset{{\text{\bf 3}}}{
\rho_0\wedge\overline{\kappa}_0}
\ \ \ \ \ \ \ \ \ \
\underset{{\text{\bf 4}}}{
\rho_0\wedge\overline{\zeta}_0}
\\
\underset{{\text{\bf 5}}}{
\kappa_0\wedge\zeta_0}
\ \ \ \ \ \ \ \ \ \
\underset{{\text{\bf 6}}}{
\kappa_0\wedge\overline{\kappa}_0}
\ \ \ \ \ \ \ \ \ \
\underset{{\text{\bf 7}}}{
\kappa_0\wedge\overline{\zeta}_0}
\\
\underset{{\text{\bf 8}}}{
\zeta_0\wedge\overline{\kappa}_0}
\ \ \ \ \ \ \ \ \ \
\underset{{\text{\bf 9}}}{
\zeta_0\wedge\overline{\zeta}_0}
\\
\underset{{\text{\bf 10}}}{
\overline{\kappa}_0\wedge\overline{\zeta}_0}.
\endaligned
\]
With such a numbering, we can abreviate the structure equations 
as\,\,---\,\,dropping their conjugates\,\,---:
\[
\aligned
d\rho_0
&
\,=\,
\Raux_0^1\,
\rho_0\wedge\kappa_0
+
\Raux_0^2\,
\rho_0\wedge\zeta_0
+
\Raux_0^3\,
\rho_0\wedge\overline{\kappa}_0
+
\Raux_0^4\,
\rho_0\wedge\overline{\zeta}_0
+
\isqrt\,
\kappa_0\wedge\overline{\kappa}_0,
\\
d\kappa_0
&
\,=\,
\Kaux_0^2\,
\rho_0\wedge\zeta_0
+
\Kaux_0^5\,
\kappa_0\wedge\zeta_0
+
\Kaux_0^8\,
\zeta_0\wedge\overline{\kappa}_0,
\\
d\zeta_0
&
\,=\,
0.
\endaligned
\]

\begin{Convention}
All functions of $p = (z_1, z_2, \overline{z}_1, \overline{z}_2, v)
\in M$
will be denoted with a lower index 
$(\centersmallbullet)_0$, always employing
the special auxiliary font characters
$\Aaux, \Baux, \Caux, \dots$.\hfill$\bigtriangleup$
\end{Convention}

After some transformations in the next
sections, this initial coframe will change and become
more complicated (unwriting the conjugates):
\[
\big\{
\rho_0,\,
\kappa_0,\,
\zeta_0
\big\}
\ \ \ \ \
\leadsto
\ \ \ \ \
\big\{
\rho_0,\,
\kappa_0,\,
\zeta_0'
\big\}
\ \ \ \ \
\leadsto
\ \ \ \ \
\big\{
\rho_0,\,
\kappa_0',\,
\zeta_0'
\big\}
\ \ \ \ \
\leadsto
\ \ \ \ \
\big\{
\rho_0,\,
\kappa_0',\,
\zeta_0''
\big\},
\]
and new structure function $\Raux_0^{i\prime}$, $\Kaux_0^{i\prime}$,
$\Zaux_0^{i\prime}$, \dots will appear.

\smallskip

We end up this section by stating some 
technical commutation relations
that shall be constantly necessary to {\em normalize} 
incoming (complicated) expressions in order to avoid
ambiguities. In fact, we can take advantage of 
$\overline{\mathcal{K}} (\kaux) = 0$
from Lemma~{\ref{Lemma-K-bar-k-K-bar-P}}, 
to make $\overline{\mathcal{K}}$ `jump' 
above iterated derivatives like {\em e.g.} in:\bigskip
\[
\xymatrix{ 
\overline{\mathcal{K}}
\ar@/^2pc/[r] 
&
\!\!\!\!\!\!\!\!\!\!\!\!\!\!\!
\big(\overline{\mathcal{L}}_1
(\kaux)\big)},
\ \ \ \ \ \ \ \ \ \ \ \ \ \ \ \ \ \ \ \ \ \ \ \ \ \
\xymatrix{ 
\overline{\mathcal{K}}
\ar@/^2pc/[r] 
&
\!\!\!\!\!\!\!\!\!\!\!\!\!\!\!
\big(\overline{\mathcal{L}}_1
\big(\overline{\mathcal{L}}_1
(\kaux))\big)}.
\]
Precisely, the last, 10\textsuperscript{\,th} Lie bracket relation
preceding Lemma~{\ref{Lemma-K-bar-k-K-bar-P}}:
\leqnomode\usetagform{default}
\begin{align}
\label{K-bar-L-1-bar-bracket}
-\,\overline{\mathcal{L}}_1
\big(\overline{\kaux}\big)
\cdot
\overline{\mathcal{L}}_1
(\centersmallbullet)
\,=\,
\big[
\overline{\mathcal{K}},
\overline{\mathcal{L}}_1
\big]
(\centersmallbullet),
\end{align}
when applied to the function $\centersmallbullet := \kaux$ yields:
\[
\aligned
-\,\overline{\mathcal{L}}_1
\big(\overline{\kaux}\big)\,
\overline{\mathcal{L}}_1
(\kaux)
\,=\,
\big[
\overline{\mathcal{K}},
\overline{\mathcal{L}}_1
\big]
(\kaux)
&
\,=\,
\overline{\mathcal{K}}
\big(
\overline{\mathcal{L}}_1(\kaux)
\big)
-
\overline{\mathcal{L}}_1
\big(
\zero{
\overline{\mathcal{K}}(\kaux)}
\big)
\\
&
\,=\,
\overline{\mathcal{K}}
\big(
\overline{\mathcal{L}}_1(\kaux)
\big).
\endaligned
\]

\begin{Lemma}
\label{Lemma-K-bar-penetrates-L1-bar-k}
One has the $3$ relations:
\leqnomode\usetagform{default}
\begin{align}
\overline{\mathcal{K}}
\big(\overline{\mathcal{L}}_1(\kaux)\big)
&
\,=\,
-\,
\overline{\mathcal{L}}_1\big(\overline{\kaux}\big)\,
\overline{\mathcal{L}}_1(\kaux),
\tag{\bf 1}
\\
\overline{\mathcal{K}}
\big(\overline{\mathcal{L}}_1
\big(\overline{\mathcal{L}}_1(\kaux)\big)\big)
&
\,=\,
-\,2\,
\overline{\mathcal{L}}_1
\big(\overline{\kaux}\big)\,\,
\overline{\mathcal{L}}_1\big(
\overline{\mathcal{L}}_1(\kaux)\big)
-
\overline{\mathcal{L}}_1
\big(\overline{\mathcal{L}}_1
\big(\overline{\kaux}\big)\big)\,\,
\overline{\mathcal{L}}_1(\kaux),
\tag{\bf 2}
\\
\ \ \ \ \ \ \ \ \ \
\overline{\mathcal{K}}
\big(\overline{\mathcal{L}}_1
\big(\overline{\mathcal{L}}_1
\big(\overline{\mathcal{L}}_1(\kaux)\big)\big)\big)
&
\,=\,
-\,3\,
\overline{\mathcal{L}}_1
\big(\overline{\kaux}\big)\,\,
\overline{\mathcal{L}}_1\big(
\overline{\mathcal{L}}_1\big(
\overline{\mathcal{L}}_1(\kaux)\big)\big)
\,-
\tag{\bf 3}
\\
&
\ \ \ \ \
-\,
3\,
\overline{\mathcal{L}}_1
\big(\overline{\mathcal{L}}_1
\big(\overline{\kaux}\big)\big)\,\,
\overline{\mathcal{L}}_1\big(
\overline{\mathcal{L}}_1(\kaux)\big)
-
\overline{\mathcal{L}}_1\big(
\overline{\mathcal{L}}_1\big(
\overline{\mathcal{L}}_1(\overline{\kaux})\big)\big)\,\,
\overline{\mathcal{L}}_1(\kaux).
\notag
\end{align}
\end{Lemma}

\proof
As {\small\bf (1)} is done, we can apply $\overline{\mathcal{L}}_1
(\centersmallbullet)$ to it, reversing sides:
\[
-\,
\overline{\mathcal{L}}_1\big(
\overline{\mathcal{L}}_1\big(
\overline{\kaux}\big)\big)\,\,
\overline{\mathcal{L}}_1(\kaux)
-
\overline{\mathcal{L}}_1\big(\overline{\kaux}\big)\,\,
\overline{\mathcal{L}}_1\big(
\overline{\mathcal{L}}_1
(\kaux)\big)
\,\,=\,\,
\overline{\mathcal{L}}_1\big(
\overline{\mathcal{K}}\big(
\overline{\mathcal{L}}_1
(\kaux)\big)\big).
\]
Similarly, we apply~({\ref{K-bar-L-1-bar-bracket}}) to
$\centersmallbullet := \overline{\mathcal{L}}_1 (\kaux)$
and we reach {\small\bf (2)} after a replacement:
\[
\aligned
-\,
\overline{\mathcal{L}}_1
\big(\overline{\kaux}\big)\,\,
\overline{\mathcal{L}}_1\big(
\overline{\mathcal{L}}_1
(\kaux)\big)
\,\,=\,\,
\big[
\overline{\mathcal{K}},
\overline{\mathcal{L}}_1
\big]
\big(
\overline{\mathcal{L}}_1(\kaux)
\big)
&
\,\,=\,\,
\overline{\mathcal{K}}\big(
\overline{\mathcal{L}}_1\big(
\overline{\mathcal{L}}_1
(\kaux)\big)\big)
-
\underbrace{
\overline{\mathcal{L}}_1\big(
\overline{\mathcal{K}}\big(
\overline{\mathcal{L}}_1
(\kaux)\big)\big)}_{\sf replace}.
\endaligned
\]

Now, as {\small\bf (2)} is done, we can apply $\overline{\mathcal{L}}_1
(\centersmallbullet)$ to it, and get after reorganization:
\[
\footnotesize
\aligned
\overline{\mathcal{L}}_1\big(
\overline{\mathcal{K}}\big(
\overline{\mathcal{L}}_1\big(
\overline{\mathcal{L}}_1
(\kaux)\big)\big)\big)
\,\,=\,\,
-\,
2\,
\overline{\mathcal{L}}_1\big(\overline{\kaux}\big)\,\,
\overline{\mathcal{L}}_1\big(
\overline{\mathcal{L}}_1\big(
\overline{\mathcal{L}}_1
(\kaux)\big)\big)
-
3\,
\overline{\mathcal{L}}_1\big(
\overline{\mathcal{L}}_1\big(
\overline{\kaux}\big)\big)\,\,
\overline{\mathcal{L}}_1\big(
\overline{\mathcal{L}}_1
(\kaux)\big)
-
\overline{\mathcal{L}}_1\big(
\overline{\mathcal{L}}_1\big(
\overline{\mathcal{L}}_1\big(
\overline{\kaux}\big)\big)\big)\,\,
\overline{\mathcal{L}}_1
(\kaux).
\endaligned
\]
Lastly, we apply~({\ref{K-bar-L-1-bar-bracket}}) to
$\centersmallbullet := \overline{\mathcal{L}}_1
\big( \overline{\mathcal{L}}_1 (\kaux) \big)$
and we reach {\small\bf (3)} after a replacement:
\begin{align}
-\,
\overline{\mathcal{L}}_1
\big(\overline{\kaux}\big)\,\,
\overline{\mathcal{L}}_1\big(
\overline{\mathcal{L}}_1\big(
\overline{\mathcal{L}}_1
(\kaux)\big)\big)
&
\,\,=\,\,
\big[
\overline{\mathcal{K}},
\overline{\mathcal{L}}_1
\big]
\Big(
\overline{\mathcal{L}}_1\big(
\overline{\mathcal{L}}_1
(\kaux)\big)
\Big)
\notag
\\
&
\,\,=\,\,
\overline{\mathcal{K}}\big(
\overline{\mathcal{L}}_1\big(
\overline{\mathcal{L}}_1\big(
\overline{\mathcal{L}}_1
(\kaux)\big)\big)\big)
-
\underbrace{
\overline{\mathcal{L}}_1\big(
\overline{\mathcal{K}}\big(
\overline{\mathcal{L}}_1\big(
\overline{\mathcal{L}}_1
(\kaux)\big)\big)\big)}_{\sf replace}.
\qedhere
\end{align}
\endproof

\Section{\bf Initial $G_1$-structure for local biholomorphic
equivalences $h \colon M \overset{\sim}{\longrightarrow} M'$}
\label{initial-G-1-structure-biholomorphic-equivalences}
\HEAD{{\ref{initial-G-1-structure-biholomorphic-equivalences}}.~{\sf
Initial $G_1$-structure
for local biholomorphic
equivalences $h \colon M \overset{\sim}{\longrightarrow} M'$}
}{
Wei Guo {\sc Foo} (Beijing) and Joël {\sc Merker} (Orsay)}

Now, let $h \colon U \overset{\sim}{\longrightarrow} U' \subset \C^3$
be a (local) biholomorphism from an open set $U \subset \C^3$
containing $U\ni 0$ the origin onto its image:
\[
h(U)
\,=:\,
U'
\,\ni\;
0'
\,=\,
h(0),
\]
which is
also an open set $U' \subset {\C'}^3$ containing the origin $0'$
in another target complex
Euclidean space ${\C'}^3$ having the same dimension.

\begin{center}
\begin{picture}(0,0)%
\includegraphics{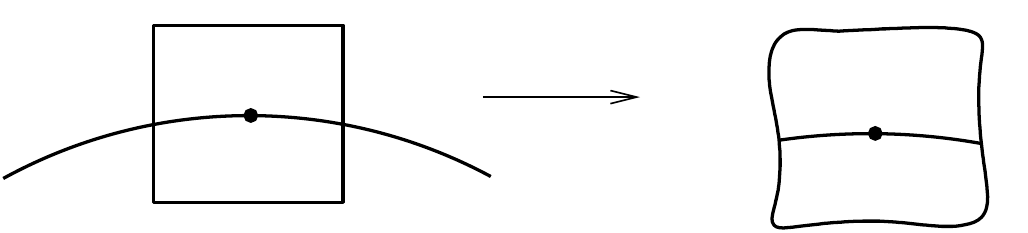}%
\end{picture}%
\setlength{\unitlength}{4144sp}%
\begingroup\makeatletter\ifx\SetFigFont\undefined%
\gdef\SetFigFont#1#2#3#4#5{%
  \reset@font\fontsize{#1}{#2pt}%
  \fontfamily{#3}\fontseries{#4}\fontshape{#5}%
  \selectfont}%
\fi\endgroup%
\begin{picture}(4693,1052)(874,-2319)
\put(943,-1939){\makebox(0,0)[lb]{\smash{{\SetFigFont{9}{10.8}{\familydefault}{\mddefault}{\updefault}{\color[rgb]{0,0,0}$M$}%
}}}}
\put(3339,-1656){\makebox(0,0)[lb]{\smash{{\SetFigFont{9}{10.8}{\familydefault}{\mddefault}{\updefault}{\color[rgb]{0,0,0}$h$}%
}}}}
\put(5117,-1833){\makebox(0,0)[lb]{\smash{{\SetFigFont{9}{10.8}{\familydefault}{\mddefault}{\updefault}{\color[rgb]{0,0,0}$M'$}%
}}}}
\put(902,-1414){\makebox(0,0)[lb]{\smash{{\SetFigFont{9}{10.8}{\familydefault}{\mddefault}{\updefault}{\color[rgb]{0,0,0}$\C^{3}$}%
}}}}
\put(5552,-1470){\makebox(0,0)[lb]{\smash{{\SetFigFont{9}{10.8}{\familydefault}{\mddefault}{\updefault}{\color[rgb]{0,0,0}${\C'}^{3}$}%
}}}}
\put(4844,-2009){\makebox(0,0)[lb]{\smash{{\SetFigFont{9}{10.8}{\familydefault}{\mddefault}{\updefault}{\color[rgb]{0,0,0}$0'$}%
}}}}
\put(4459,-1546){\makebox(0,0)[lb]{\smash{{\SetFigFont{9}{10.8}{\familydefault}{\mddefault}{\updefault}{\color[rgb]{0,0,0}$U'$}%
}}}}
\put(1612,-1507){\makebox(0,0)[lb]{\smash{{\SetFigFont{9}{10.8}{\familydefault}{\mddefault}{\updefault}{\color[rgb]{0,0,0}$U$}%
}}}}
\put(1990,-1927){\makebox(0,0)[lb]{\smash{{\SetFigFont{9}{10.8}{\familydefault}{\mddefault}{\updefault}{\color[rgb]{0,0,0}$0$}%
}}}}
\end{picture}%

\end{center}

As in Cartan's equivalence theory,
assume that $h \big(M \cap U \big) \subset M'$
is contained in another real hypersurface $M' \subset {\C'}^3$,
also passing through the origin $0' \in M'$, represented
in holomorphic coordinates $\big(z_1', z_2', w' = u' + \isqrt\,
v' \big)$ by a similar $\mathcal{C}^\omega$ graphed equation:
\[
u'
\,=\,
\Faux'
\big(
z_1',z_2',\overline{z}_1',\overline{z}_2',v'
\big).
\]
We now make the convention of not mentioning the open sets
that must sometimes be shrunk, so that 
we think of $h \colon M
\overset{\sim}{\longrightarrow} M'$ as being
a CR equivalence between hypersurfaces 
$M \subset \C^3$ and $M' \subset {\C'}^3$.

In the target space, 
introduce similar generators $\mathcal{L}_1'$, $\mathcal{L}_2'$
for $T^{1,0}M'$.
Since $h$ is holomorphic, its differential
$h_\ast \colon \C T\C^3 \longrightarrow \C T{\C'}^3$
stabilizes holomorphic $(1,0)$ and holomorphic
$(0,1)$ vector fields:
\[
h_\ast\big(T^{1,0}\C^3\big)
\,=\,
T^{1,0}{\C'}^3
\ \ \ \ \ \ \ \ \ \ \ \ \ \ \ \ \ \
\text{and}
\ \ \ \ \ \ \ \ \ \ \ \ \ \ \ \ \ \
h_\ast
\big(
T^{0,1}M
\big)
\,=\,
T^{0,1}M'.
\]
Furthermore, by invariancy of the Freeman form,
$h$ respects the Levi-kernel distributions:
\[
h_\ast
\big(
K^{1,0}M
\big)
\,=\,
K^{1,0}M'.
\]
Consequently, there exist functions $f'$, $c'$, $e'$ on $M'$ such
that:
\[
\aligned
h_\ast(\mathcal{K})
&
\,=\,
f'\,\mathcal{K}',
\\
h_\ast\big(\mathcal{L}_1\big)
&
\,=\,
c'\,\mathcal{L}_1'
+
e'\,\mathcal{K}',
\endaligned
\]
whence by conjugation:
\[
\aligned
h_\ast\big(\overline{\mathcal{K}}\big)
&
\,=\,
\overline{f}'\,
\overline{\mathcal{K}}',
\\
h_\ast\big(\overline{\mathcal{L}}_1\big)
&
\,=\,
\overline{c}'\,
\overline{\mathcal{L}}_1'
+
\overline{e}'\,
\overline{\mathcal{K}}'.
\endaligned
\]

On the other hand, there is {\em a priori} no special condition
that shall be satisfied by $h_\ast (\mathcal{T})$,
except that it be a real vector field, because $\mathcal{T}$
is real. Thus, there are a real-valued function $a'$ and
two complex-valued $b'$ and $d'$ on $M'$ such that:
\[
h_\ast(\mathcal{T})
\,=\,
a'\,\mathcal{T}'
+
b'\,\mathcal{L}_1'
+
d'\,\mathcal{K}'
+
\overline{b}'\,
\overline{\mathcal{L}}_1'
+
\overline{d}'\,
\overline{\mathcal{K}}'.
\]
In fact, the function $a'$ is determined, because:
\[
\aligned
h_\ast(\mathcal{T})
\,=\,
h_\ast
\big(
\isqrt\,
\big[\mathcal{L}_1,\overline{\mathcal{L}}_1
\big]
\big)
&
\,=\,
\isqrt\,
\big[
h_\ast\big(\mathcal{L}_1\big),\,
h_\ast\big(\overline{\mathcal{L}}_1\big)
\big]
\\
&
\,=\,
\isqrt\,
\big[
c'\mathcal{L}_1'
+
e'\mathcal{K}',\,\,
\overline{c}'\overline{\mathcal{L}}_1'
+
\overline{e}'\overline{\mathcal{K}}'
\big]
\\
&
\,\equiv\,
c'\overline{c}'\,
\isqrt\,
\big[
\mathcal{L}_1',\,
\overline{\mathcal{L}}_1'
\big]
\ \ \ \ \
\mod\,\,
\big(
T^{1,0}M'
\oplus
T^{0,1}M'
\big),
\endaligned
\]
whence necessarily:
\[
a'
\,=\,
c'\overline{c}'.
\]

Summarizing, we have the following matrix relations:
\[
h_\ast\,
\left(\!\!
\begin{array}{c}
\mathcal{T}
\\
\mathcal{L}_1
\\
\mathcal{K}
\\
\overline{\mathcal{L}}_1
\\
\overline{\mathcal{K}}
\end{array}
\!\!\right)
\,\,=\,\,
\left(\!\!
\begin{array}{ccccc}
c'\overline{c}' & b' & d' & \overline{b}' & \overline{d}'
\\
0 & c' & e' & 0 & 0
\\
0 & 0 & f' & 0 & 0
\\
0 & 0 & 0 & \overline{c}' & \overline{e}'
\\
0 & 0 & 0 & 0 & \overline{f}'
\end{array}
\!\!\right)\,
\left(\!\!
\begin{array}{c}
\mathcal{T}'
\\
\mathcal{L}_1'
\\
\mathcal{K}'
\\
\overline{\mathcal{L}}_1'
\\
\overline{\mathcal{K}}'
\end{array}
\!\!\right).
\]
As $h_\ast$ is invertible, the function $f'$, and then the
function $c'$ too, must be nowhere vanishing. The
relation between the coframe $\big\{ \rho_0, \kappa_0, \zeta_0, 
\overline{\kappa}_0,
\overline{\zeta}_0 \big\}$
in the source space and  
the coframe $\big\{ \rho_0', \kappa_0', \zeta_0', \overline{\kappa}_0',
\overline{\zeta}_0' \big\}$
in the target space is therefore given by
a plain transposition:
\[
h^\ast\,
\left(\!\!
\begin{array}{c}
\rho_0'
\\
\kappa_0'
\\
\zeta_0'
\\
\overline{\kappa}_0'
\\
\overline{\zeta}_0'
\end{array}
\!\!\right)
\,\,=\,\,
\left(\!\!
\begin{array}{ccccc}
c'\overline{c}' & 0 & 0 & 0 & 0
\\
b' & c' & 0 & 0 & 0
\\
d' & e' & f' & 0 & 0
\\
\overline{b}' & 0 & 0 & \overline{c}' & 0
\\
\overline{d}' & 0 & 0 & \overline{e}' & \overline{f}'
\end{array}
\!\!\right)\,
\left(\!\!
\begin{array}{c}
\rho_0
\\
\kappa_0
\\
\zeta_0
\\
\overline{\kappa}_0
\\
\overline{\zeta}_0
\end{array}
\!\!\right).
\]

These preliminaries, 
also explained in~{\cite{Merker-2013-5-CR-IV, Pocchiola-2013,
Merker-Pocchiola-2018}}, 
justify that the initial $G$-structure
for such equivalences of CR manifolds 
is the matrix ambiguity group $G_1$ is constituted of $5 \times 5$ 
matrices of the form:
\[
\left(\!\!
\begin{array}{ccccc}
{\sf c}\overline{\sf c} & 0 & 0 & 0 & 0
\\
{\sf b} & {\sf c} & 0 & 0 & 0
\\
{\sf d} & {\sf e} & {\sf f} & 0 & 0
\\
\overline{\sf b} & 0 & 0 & \overline{\sf c} & 0
\\
\overline{\sf d} & 0 & 0 & \overline{\sf e} & \overline{\sf f}
\end{array}
\!\!\right),
\]
with free variable complex entries:
\[
{\sf c},\,{\sf f}
\,\in\,
\C
\backslash
\{0\}
\ \ \ \ \ \ \ \ \ \ \ \ \ \ \ \ \ \
\text{and}
\ \ \ \ \ \ \ \ \ \ \ \ \ \ \ \ \ \
{\sf b},\,{\sf d},\,{\sf e}
\,\in\,\C,
\]
namely:
\[
\left(\!\!
\begin{array}{c}
\rho
\\
\kappa
\\
\zeta
\\
\overline{\kappa}
\\
\overline{\zeta}
\end{array}
\!\!\right)
\,:=\,
\left(\!\!
\begin{array}{ccccc}
{\sf c}\overline{\sf c} & 0 & 0 & 0 & 0
\\
{\sf b} & {\sf c} & 0 & 0 & 0
\\
{\sf d} & {\sf e} & {\sf f} & 0 & 0
\\
\overline{\sf b} & 0 & 0 & \overline{\sf c} & 0
\\
\overline{\sf d} & 0 & 0 & \overline{\sf e} & \overline{\sf f}
\end{array}
\!\!\right)\,\,
\left(\!\!
\begin{array}{c}
\rho_0
\\
\kappa_0
\\
\zeta_0
\\
\overline{\kappa}_0
\\
\overline{\zeta}_0
\end{array}
\!\!\right).
\]

Eliminating the conjugate $1$-forms $\overline{\kappa}$,
$\overline{\zeta}$ for which the structure
equations are redundant, this can be abbreviated as:
\[
\left(\!\!
\begin{array}{c}
\rho
\\
\kappa
\\
\zeta
\end{array}
\!\!\right)
\,:=\,
\left(\!\!
\begin{array}{ccc}
{\sf c}\overline{\sf c} & 0 & 0 
\\
{\sf b} & {\sf c} & 0
\\
{\sf d} & {\sf e} & {\sf f}
\end{array}
\!\!\right)\,\,
\left(\!\!
\begin{array}{c}
\rho_0
\\
\kappa_0
\\
\zeta_0
\end{array}
\!\!\right).
\]

\Section{\bf A Labyrinthmap to Pocchiola's Calculations}
\label{labyrinthmap-Pocchiola-calculations}
\HEAD{{\ref{labyrinthmap-Pocchiola-calculations}}.~{\sf A Labyrinthmap 
to Pocchiola's Calculations}
}{
Wei Guo {\sc Foo} (Beijing) and Joël {\sc Merker} (Orsay)}

The successive reductions of this $G_1$ structure will look as:
\[
\aligned
g
\,:=\,
\left(\!\!
\begin{array}{ccc}
{\sf c}\overline{\sf c} & 0 & 0 
\\
{\sf b} & {\sf c} & 0
\\
{\sf d} & {\sf e} & {\sf f}
\end{array}
\!\!\right)
\ \ \ \ \
\leadsto
\ \ \ \ \
g
\,:=\,
\left(\!\!
\begin{array}{ccc}
{\sf c}\overline{\sf c} & 0 & 0 
\\
{\sf b} & {\sf c} & 0
\\
{\sf d} & {\sf e} & \frac{{\sf c}}{\overline{\sf c}}
\end{array}
\!\!\right)
\ \ \ \ \
\leadsto
\ \ \ \ \
g
&
\,:=\,
\left(\!\!
\begin{array}{ccc}
{\sf c}\overline{\sf c} & 0 & 0 
\\
-\isqrt\,\overline{\sf c}{\sf e} & {\sf c} & 0
\\
{\sf d} & {\sf e} & \frac{{\sf c}}{\overline{\sf c}}
\end{array}
\!\!\right)
\\
\ \ \ \ \
\leadsto
\ \ \ \ \
g
&
\,:=\,
\left(\!\!
\begin{array}{ccc}
{\sf c}\overline{\sf c} & 0 & 0 
\\
-\isqrt\,\overline{\sf c}{\sf e} & {\sf c} & 0
\\
-\frac{\isqrt}{2}\,\frac{\overline{\sf c}{\sf e}^2}{{\sf c}} & 
{\sf e} & \frac{{\sf c}}{\overline{\sf c}}
\end{array}
\!\!\right),
\endaligned
\]
thanks to successive normalization of some group parameters
(offered by some essential torsion coefficients yielding 
invariants that are deeper than Levi and Freeman forms):
\[
\aligned
{\sf f}
\,:=\,
\frac{{\sf c}}{\overline{\sf c}}\,
\overline{\mathcal{L}}_1(\kaux),
\ \ \ \ \ \ \ \ \ \ \ \ \ \ \ \ \ \
{\sf b}
&
\,:=\,
-\,\isqrt\,
\overline{\sf c}\,{\sf e}
+
\frac{\isqrt}{3}\,
{\sf c}\,
\Baux_0,
\\
{\sf d}
&
\,:=\,
-\,\frac{\isqrt}{2}\,
\frac{\overline{\sf c}\,{\sf e}\,{\sf e}}{{\sf c}}
+
\isqrt\,
\frac{{\sf c}}{\overline{\sf c}}\,
\Haux_0,
\endaligned
\]
in terms of the following two function on $M$:
\[
\aligned
\Baux_0
&
\,:=\,
\frac{\overline{\mathcal{L}}_1\big(\overline{\mathcal{L}}_1(\kaux)\big)}{
\overline{\mathcal{L}}_1(\kaux)}
-
\overline{\Paux},
\\
\Haux_0
&
\,:=\,
\,-
\frac{1}{6}\,
\frac{\overline{\mathcal{L}}_1\big(
\overline{\mathcal{L}}_1\big(\overline{\mathcal{L}}_1(\kaux)\big)
\big)
}{
\overline{\mathcal{L}}_1(\kaux)}
+
\frac{2}{9}\,
\frac{
\overline{\mathcal{L}}_1\big(\overline{\mathcal{L}}_1(\kaux)\big)^2
}{
\overline{\mathcal{L}}_1(\kaux)^2}
+
\frac{1}{18}\,
\frac{\overline{\mathcal{L}}_1\big(\overline{\mathcal{L}}_1(\kaux)\big)\,
\overline{\Paux}}{
\overline{\mathcal{L}}_1(\kaux)}
+
\frac{1}{6}\,
\overline{\mathcal{L}_1}\big(\overline{\Paux}\big)
-
\frac{1}{9}\,
\overline{\Paux}^2.
\endaligned
\]
This function $\Haux_0$ coincides with Pocchiola's function $H$.

\smallskip

The next sections will present in details these successive 
reductions of $G$-structures, by these normalizations 
of the group parameters ${\sf f}$, ${\sf b}$, ${\sf d}$.
Contrary to~{\cite{Merker-Pocchiola-2018, Pocchiola-2013}}, 
all computations will be progressive, simple, detailed, readable,
clear, without needing any help of either a computer or
a pen. A great care will be devoted to readability.

\Section{\bf First Loop: Reduction of the Group Parameter ${\sf f}$}
\label{first-loop-reduction-f}
\HEAD{{\ref{first-loop-reduction-f}}.~{\sf First Loop: 
Reduction of the Group Parameter ${\sf f}$}
}{
Wei Guo {\sc Foo} (Beijing) and Joël {\sc Merker} (Orsay)}

We recall that the initial Darboux-Cartan structure of the coframe
$\big\{ \rho_0, \kappa_0, \zeta_0, \overline{\kappa}_0,
\overline{\zeta}_0 \big\}$ is, without writing conjugate
equations\,\,---\,\,remind $\overline{\rho}_0 = \rho_0$\,\,---:
\leqnomode\usetagform{default}
\begin{align}
\label{d-rho-0-kappa-0-zeta-0}
d\rho_0
&
\,=\,
\Paux\,
\rho_0\wedge\kappa_0
-
\mathcal{L}_1(\kaux)\,
\rho_0\wedge\zeta_0
+
\overline{\Paux}\,
\rho_0\wedge\overline{\kappa}_0
-
\overline{\mathcal{L}}_1\big(\overline{\kaux}\big)\,
\rho_0\wedge\overline{\zeta}_0
+
\isqrt\,
\kappa_0\wedge\overline{\kappa}_0,
\notag
\\
\ \ \ \ \ 
d\kappa_0
&
\,=\,
-\,\mathcal{T}(\kaux)\,
\rho_0\wedge\zeta_0
-
\mathcal{L}_1(\kaux)\,
\kappa_0\wedge\zeta_0
+
\overline{\mathcal{L}}_1(\kaux)\,
\zeta_0\wedge\overline{\kappa}_0,
\\
d\zeta_0
&
\,=\,
0.
\notag
\end{align}

With the first $G$-structure exhibited above, 
introduce the {\sl lifted differential
forms}, defined by:
\[
\left(\!
\begin{array}{c}
\rho
\\
\kappa
\\
\zeta
\end{array}
\!\right)
\,:=\,
\left(\!
\begin{array}{ccc}
{\sf c}\overline{\sf c} & 0 & 0
\\
{\sf b} & {\sf c} & 0
\\
{\sf d} & {\sf e} & {\sf f}
\end{array}
\!\right)
\left(\!
\begin{array}{c}
\rho_0
\\
\kappa_0
\\
\zeta_0
\end{array}
\!\right),
\]
{\em id est:}
\[
\aligned
\rho
&
\,:=\,
{\sf c}\overline{\sf c}\,\rho_0,
\\
\kappa
&
\,:=\,
{\sf b}\,\rho_0
+
{\sf c}\,\kappa_0,
\\
\zeta
&
\,:=\,
{\sf d}\,\rho_0
+
{\sf e}\,\kappa_0
+
{\sf f}\,\zeta_0.
\endaligned
\]
Here ${\sf c}, {\sf f} \in \C^\ast$ and 
${\sf b}, {\sf e}, {\sf d} \in \C$.
Mind that conjugate equations giving $\overline{\kappa}$ and
$\overline{\zeta}$ are not written, but will
be used.

An inversion yields:
\leqnomode\usetagform{default}
\begin{align}
\label{inverse-lifted-base-loop-1}
\rho_0
&
\,=\,
\frac{1}{{\sf c}\overline{\sf c}}\,
\rho,
\notag
\\
\kappa_0
&
\,=\,
\frac{1}{{\sf c}}\,
\kappa
-
\frac{{\sf b}}{{\sf c}{\sf c}\overline{\sf c}}\,
\rho,
\\
\zeta_0
&
\,=\,
\frac{{\sf b}{\sf e}-{\sf c}{\sf d}}{{\sf c}{\sf c}\overline{\sf c}
{\sf f}}\,
\rho
-
\frac{{\sf e}}{{\sf c}{\sf f}}\,
\kappa
+
\frac{1}{{\sf f}}\,
\zeta.
\notag
\end{align}

With the above $3 \times 3$ matrix $g$ representing
the general element of a $10$-dimensional (real)
group $G^{10} \subset {\sf GL}_3(\C)$, 
the Maurer-Cartan matrix is:
\[
\aligned
dg\cdot g^{-1}
&
\,=\,
\left(\!
\begin{array}{ccc}
\overline{\sf c}\,d{\sf c}+{\sf c}d\overline{\sf c} & 0 & 0
\\
d{\sf b} & d{\sf c} & 0
\\
d{\sf d} & d{\sf e} & d{\sf f}
\end{array}
\!\right)
\left(\!
\begin{array}{ccc}
\frac{1}{{\sf c}\overline{\sf c}} & 0 & 0
\\
-\frac{{\sf b}}{{\sf c}{\sf c}\overline{\sf c}} & \frac{1}{{\sf c}} & 0
\\
\frac{{\sf b}{\sf e}-{\sf c}{\sf d}}{{\sf c}{\sf c}\overline{\sf c}
{\sf f}}\, & -\frac{{\sf e}}{{\sf c}{\sf f}} & \frac{1}{{\sf f}}
\end{array}
\!\right)
\\
&
\,=:\,
\left(\!
\begin{array}{ccc}
\alpha+\overline{\alpha} & 0 & 0
\\
\beta & \alpha & 0
\\
\gamma & \delta & \varepsilon
\end{array}
\!\right),
\endaligned
\]
in terms of the group-invariant $1$-forms:
\[
\aligned
\alpha
&
\,:=\,
\frac{d{\sf c}}{{\sf c}},
\\
\beta
&
\,:=\,
\frac{d{\sf b}}{{\sf c}\overline{\sf c}}
-
\frac{{\sf b}\overline{\sf c}\,d{\sf c}}{{\sf c}{\sf c}},
\\
\gamma
&
\,:=\,
\frac{d{\sf d}}{{\sf c}\overline{\sf c}}
-
\frac{{\sf b}\,d{\sf e}}{{\sf c}{\sf c}\overline{\sf c}}
+
\frac{{\sf b}{\sf e}-{\sf c}{\sf d}}{{\sf c}{\sf c}\overline{\sf c}
{\sf f}}\,
d{\sf f},
\\
\delta
&
\,:=\,
\frac{d{\sf e}}{{\sf c}}
-
\frac{{\sf e}\,d{\sf f}}{{\sf c}{\sf f}},
\\
\varepsilon
&
\,:=\,
\frac{d{\sf f}}{\sf f}.
\endaligned
\]

As is known, after painful computations whose outcomes
are presented extensively 
in~{\cite{Pocchiola-2013, Merker-Pocchiola-2018}},
one can re-express, 
using~({\ref{d-rho-0-kappa-0-zeta-0}}) 
and~({\ref{inverse-lifted-base-loop-1}}),
the exterior differentials of the $3$ lifted $1$-forms
$\rho$, $\zeta$, $\kappa$ as:
\[
\aligned
d\rho
&
\,=\,
\alpha\wedge\rho
+
\overline{\alpha}\wedge\rho
\,+
\\
&
\ \ \ \ \
+
R^1\,
\rho\wedge\kappa
+
R^2\,
\rho\wedge\zeta
+
R^3\,
\rho\wedge\overline{\kappa}
+
R^4\,
\rho\wedge\overline{\zeta}
+
\isqrt\,\kappa\wedge\overline{\kappa},
\\
d\kappa
&
\,=\,
\beta\wedge\rho
+
\alpha\wedge\kappa
\,+
\\
&
\ \ \ \ \ 
K^1\,
\rho\wedge\kappa
+
K^2\,
\rho\wedge\zeta
+
K^3\,
\rho\wedge\overline{\kappa}
+
K^4\,
\rho\wedge\overline{\zeta}
\,+
\\
&
\ \ \ \ \ 
+
K^5\,
\kappa\wedge\zeta
+
K^5\,
\kappa\wedge\overline{\kappa}
+
\boxed{K^8}\,
\zeta\wedge\overline{\kappa},
\\
d\zeta
&
\,=\,
\gamma\wedge\rho
+
\delta\wedge\kappa
+
\varepsilon\wedge\zeta
\,+
\\
&
\ \ \ \ \ 
+
Z^1\,
\rho\wedge\kappa
+
Z^2\,
\rho\wedge\zeta
+
Z^3\,
\rho\wedge\overline{\kappa}
+
Z^4\,
\rho\wedge\overline{\zeta}
\,+
\\
&
\ \ \ \ \
+
Z^5\,
\kappa\wedge\zeta
+
Z^6\,
\kappa\wedge\overline{\kappa}
+
Z^8\,
\zeta\wedge\overline{\kappa},
\endaligned
\]
in terms of certain complicated functions $R^i$, $K^i$, $Z^i$ of the
horizontal variables and of the group parameters as well:
\[
\big( 
z_1,z_2,\overline{z}_1,\overline{z}_2,v
\big)
\times
\big(
{\sf c},\overline{\sf c},{\sf f},\overline{\sf f},
{\sf b},\overline{\sf b},
{\sf d},\overline{\sf d},{\sf e},\overline{\sf e}
\big)
\,\,\in\,\,
M^5
\times
G^{10},
\]
but we shall not need the expressions of all these functions,
and focus only on the boxed one, $K^8$, since it will
bring an interesting normalization for the
diagonal group parameter
${\sf f}$. 

\begin{Notation}
Given a differential $2$-form $\Omega \in 
\Gamma(M, \Lambda^2 T^\ast M)$ on an $n$-dimensional manifold
$M$ equipped with a coframe $\big\{ \omega^1, \dots, \omega^n \big\}$
for its cotangent bundle $T^\ast M$, 
which is expanded as:
\[
\Omega
\,=\,
\sum_{1\leqslant i<j\leqslant n}\,
A_{i,j}\,
\omega^i\wedge\omega^j,
\]
with uniquely determined coefficients-functions
$A_{\smallbullet, \smallbullet}$, 
for fixed $i < j$, the coefficient $A_{i,j}$ 
of $\omega^i \wedge \omega^j$ will be denoted by:
\[
\big[
\omega^i
\wedge
\omega^j
\big]
\big\{
\Omega
\big\}
\,:=\,
A_{i,j}.
\]
\end{Notation}

To capture $K^8$ without pain, the computation\big/re-expression
of $d\kappa$ starts from $\kappa = {\sf b}\, \rho_0 + 
{\sf c}\, \kappa_0$ as follows to see how Maurer-Cartan forms
enter the play:
\[
\aligned
d\kappa
&
\,=\,
d{\sf b}
\wedge
\rho_0
+
d{\sf c}\wedge\kappa_0
+
{\sf b}\,d\rho_0
+
{\sf c}\,d\kappa_0
\\
&
\,=\,
d{\sf b}
\wedge
\big(
{\textstyle{\frac{1}{{\sf c}\overline{\sf c}}}}\,
\rho
\big)
+
d{\sf c}
\wedge
\big(
{\textstyle{\frac{1}{{\sf c}}}}\,\kappa
-
{\textstyle{\frac{{\sf b}}{{\sf c}{\sf c}\overline{\sf c}}}}\,
\rho
\big)
+
\Torsion
\\
&
\,=\,
\big(
d{\sf b}
-
{\textstyle{\frac{{\sf b}\,d{\sf c}}{
{\sf c}{\sf c}\overline{\sf c}}}}\,
\big)
\wedge
\rho
+
\big(
{\textstyle{\frac{d{\sf c}}{{\sf c}}}}
\big)
\wedge
\kappa
+
\Torsion
\\
&
\,=\,
\beta\wedge\rho
+
\alpha\wedge\kappa
+
\Torsion.
\endaligned
\]

Certainly, $K^8$ belongs to the torsion remainder, and we want
to determine only:
\[
K^8
\,:=\,
\big[
\zeta\wedge\overline{\kappa}
\big]
\big\{
d\kappa
\big\}
\,=\,
\big[
\zeta\wedge\overline{\kappa}
\big]
\big\{
{\sf b}\,d\rho_0
+
{\sf c}\,d\kappa_0
\big\}.
\]
For the first term ${\sf b}\, d\rho_0$, we look 
at~({\ref{d-rho-0-kappa-0-zeta-0}})
in which we replace {\em visually} 
$\rho_0$, $\zeta_0$, $\kappa_0$ by 
$\rho$, $\zeta$, $\kappa$ watching 
simultaneously~({\ref{inverse-lifted-base-loop-1}})\,\,---\,\,no pen
needed! computers shut down!\,\,---\,\,and we get:
\[
{\sf b}\,
\big[
\zeta\wedge\overline{\kappa}
\big]
\big\{
d\rho_0
\big\}
\,=\,
0
+
0
+
0
+
0
+
0
\,=\,
0.
\]
Proceeding similarly, just with eyes:
\[
\aligned
{\sf c}\,
\big[
\zeta\wedge\overline{\kappa}
\big]
\big\{
d\kappa_0
\big\}
&
\,=\,
0
+
0
+
{\sf c}\,
\overline{\mathcal{L}}_1(\kaux)\,
\big[
\zeta\wedge\overline{\kappa}
\big]
\bigg\{
\Big(
\frac{{\sf b}{\sf e}-{\sf b}{\sf d}}{{\sf c}{\sf c}\overline{\sf c}
{\sf f}}\,\rho
-
\frac{{\sf e}}{{\sf c}{\sf f}}\,\kappa
+
\frac{1}{{\sf f}}\,\zeta
\Big)
\wedge
\Big(
-\,\frac{\overline{\sf b}}{\overline{\sf c}\overline{\sf c}{\sf c}}\,
\rho
+
\frac{1}{\overline{\sf c}}\,
\overline{\kappa}
\Big)
\\
&
\,=\,
{\sf c}\,
\overline{\mathcal{L}}_1(\kaux)\,
\big(
{\textstyle{\frac{1}{{\sf f}}}}
\big)\,
\big(
{\textstyle{\frac{1}{\overline{\sf c}}}}
\big),
\endaligned
\]
whence adding:
\[
K^8
\,=\,
\frac{{\sf c}}{\overline{\sf c}{\sf f}}\,
\overline{\mathcal{L}}_1(\kaux).
\]

Furthermore, without computation, we see that $K^8$ is {\em not}
absorbable in the Maurer-Cartan part $\beta \wedge \rho + 
\alpha \wedge \kappa$ by means of any replacement:
\[
\aligned
\alpha
&
\,=\,
\alpha'
+
a_1\,\rho
+
a_2\,\kappa
+
a_3\,\zeta
+
a_4\,\overline{\kappa}
+
a_5\,\overline{\zeta},
\\
\beta
&
\,=\,
\beta'
+
b_1\,\rho
+
b_2\,\kappa
+
b_3\,\zeta
+
b_4\,\overline{\kappa}
+
b_5\,\overline{\zeta},
\endaligned
\]
because the result will always be:
\[
\something
\wedge
\rho
+
\something
\wedge
\kappa,
\]
whereas $K^8\, \zeta \wedge \overline{\kappa}$ is not 
$\wedge$-divisible by either $\wedge \rho$ or $\wedge \kappa$.

Consequently, $K^8$ is an essential torsion coefficient,
and by general Cartan theory, $K^8$ may bring a
group parameter normalization.

In fact, 
since the diagonal coefficients ${\sf c} \neq 0 \neq {\sf f}$
of the invertible triangular matrix
must be nonvanishing, and since $\overline{\mathcal{L}}_1(\kaux) \neq 
0$ is nowhere
vanishing by our assumption of $2$-nondegeneracy, it is natural, then,
to normalize $K^8$
to be constant nonzero, {\em e.g.} $K^8 := 1$,
and this yields a reduction of the $G^{10}$-structure to an 
eight-dimensional $G^8$-structure by setting:
\[
\boxed{\,
{\sf f}
\,:=\,
\frac{{\sf c}}{\overline{\sf c}}\,
\overline{\mathcal{L}}_1(\kaux).\,
}
\]

Inserting this in the lifted coframe:

\[
\left(\!
\begin{array}{c}
\rho
\\
\kappa
\\
\zeta
\end{array}
\!\right)
\,:=\,
\left(\!
\begin{array}{ccc}
{\sf c}\overline{\sf c} & 0 & 0
\\
{\sf b} & {\sf c} & 0
\\
{\sf d} & {\sf e} & \frac{{\sf c}}{\overline{\sf c}}\,
\overline{\mathcal{L}}_1(\kaux)
\end{array}
\!\right)
\left(\!
\begin{array}{c}
\rho_0
\\
\kappa_0
\\
\zeta_0
\end{array}
\!\right),
\]
we are conducted to change the initial coframe by introducing
the new {\sl horizontal}\,\,---\,\,{\em i.e.}
defined on $M$\,\,---\,\,$1$-form:
\leqnomode\usetagform{default}
\begin{align}
\label{def-zeta-prime-0}
\zeta_0'
\,:=\,
\overline{\mathcal{L}}_1(\kaux)\,
\zeta_0.
\end{align}
As anticipated in a summary {\em supra}, 
we are thus changing of horizontal coframe:
\[
\big\{
\rho_0,\,
\kappa_0,\,
\zeta_0,\,
\overline{\kappa}_0,\,
\overline{\zeta}_0
\big\}
\ \ \ \ \
\leadsto
\ \ \ \ \
\big\{
\rho_0,\,
\kappa_0,\,
\zeta_0',\,
\overline{\kappa}_0,\,
\overline{\zeta}_0'
\big\},
\]
and unavoidably, we have to set up its Darboux-Cartan structure.

Thanks to Lemma~{\ref{Lemma-d-G-0}}, we can compute:
\[
\aligned
d\zeta_0'
&
\,=\,
d\big(
\overline{\mathcal{L}}_1(\kaux)
\big)
\wedge
\zeta_0
+
\overline{\mathcal{L}}_1(\kaux)
\wedge
\zero{d\zeta_0}
\\
&
\,=\,
\mathcal{T}
\big(
\overline{\mathcal{L}}_1(\kaux)
\big)\,
\rho_0\wedge\zeta_0
+
\mathcal{L}_1
\big(
\overline{\mathcal{L}}_1(\kaux)
\big)\,
\kappa_0\wedge\zeta_0
+
\mathcal{K}\big(\overline{\mathcal{L}}_1(\kaux)\big)\,
\zero{\zeta_0\wedge\zeta_0}
+
\overline{\mathcal{L}}_1\big(
\overline{\mathcal{L}}_1
(\kaux)\big)\,
\overline{\kappa}_0
\wedge
\zeta_0
\\
&
\ \ \ \ \ \ \ \ \ \ \ \ \ \ \ \ \ \ \ \ \ \ \ \ \ \ \ \ \ \ \ \ \ \
\ \ \ \ \ \ \ \ \ \ \ \ \ \ \ \ \ \ \ \ \ \ \ \ \ \ \ \ \ \ \ \ \ \
\ \ \ \ \ \ \ \ \ \ \ \ \ \ \ 
+
\blue{
\overline{\mathcal{K}}\big(
\overline{\mathcal{L}}_1
(\kaux)\big)}\,
\overline{\zeta}_0\wedge\zeta_0
+
0,
\endaligned
\]
and next, replacing everywhere $\zeta_0 = \frac{\zeta_0'}{
\overline{\mathcal{L}}_1(\kaux)}$, reorganizing, 
and transforming the last term above in application of
Lemma~{\ref{Lemma-K-bar-penetrates-L1-bar-k}}~{\small\bf (1)},
we obtain the structure equations enjoyed by
this new initial base coframe:
\leqnomode\usetagform{default}
\begin{align}
\label{2-loop-structure-0-coframe}
d\rho_0
&
\,=\,
\Paux\,
\rho_0\wedge\kappa_0
-
\frac{\mathcal{L}_1(\kaux)}{\overline{\mathcal{L}}_1(\kaux)}\,
\rho_0\wedge\zeta_0'
+
\overline{\Paux}\,
\rho_0\wedge\overline{\kappa}_0
-
\frac{\overline{\mathcal{L}}_1(\overline{\kaux})}{
\mathcal{L}_1(\overline{\kaux})}\,
\rho_0\wedge\overline{\zeta}_0'
+
\isqrt\,
\kappa_0\wedge\overline{\kappa}_0,
\notag
\\
\ \ \ \ \ \ \ \ \ \
d\kappa_0
&
\,=\,
-\,
\frac{\mathcal{T}(\kaux)}{
\overline{\mathcal{L}}_1(\kaux)}\,
\rho_0\wedge\zeta_0'
-
\frac{\mathcal{L}_1(\kaux)}{\overline{\mathcal{L}}_1(\kaux)}\,
\kappa_0\wedge\zeta_0'
+
\zeta_0'\wedge\overline{\kappa}_0,
\\
d\zeta_0'
&
\,=\,
\frac{\mathcal{T}\big(\overline{\mathcal{L}}_1(\kaux)\big)}{
\overline{\mathcal{L}}_1(\kaux)}\,
\rho_0\wedge\zeta_0'
+
\frac{\mathcal{L}_1\big(\overline{\mathcal{L}}_1(\kaux)\big)}{
\overline{\mathcal{L}}_1(\kaux)}\,
\kappa_0\wedge\zeta_0'
-
\frac{\overline{\mathcal{L}}_1\big(
\overline{\mathcal{L}}_1(\kaux)\big)}{
\overline{\mathcal{L}}_1(\kaux)}\,
\zeta_0'\wedge\overline{\kappa}_0
+
\frac{\overline{\mathcal{L}}_1(\overline{\kaux})}{
\mathcal{L}_1(\overline{\kaux})}\,
\zeta_0'\wedge\overline{\zeta}_0'.
\notag
\end{align}

Sometimes, it can be useful to abbreviate these formulas as:
\[
\aligned
d\rho_0
&
\,=\,
\Raux_0^1\,
\rho_0\wedge\kappa_0
+
\Raux_0^2\,
\rho_0\wedge\zeta_0'
+
\overline{\Raux}_0^1\,
\rho_0\wedge\overline{\kappa}_0
+
\overline{\Raux}_2^0\,
\rho_0\wedge\overline{\zeta}_0'
+
\isqrt\,
\kappa_0\wedge\overline{\kappa}_0,
\\
d\kappa_0
&
\,=\,
\Kaux_0^2\,
\rho_0\wedge\zeta_0'
+
\Kaux_0^5\,
\kappa_0\wedge\zeta_0'
+
\zeta_0'\wedge\overline{\kappa}_0,
\\
d\zeta_0'
&
\,=\,
\Zaux_0^2\,
\rho_0\wedge\zeta_0'
+
\Zaux_0^5\,
\kappa_0\wedge\zeta_0'
+
\Zaux_0^8\,
\zeta_0'\wedge\overline{\kappa}_0
+
\Zaux_0^9\,
\zeta_0'\wedge\overline{\zeta}_0',
\endaligned
\]
and no {\sl primes} will be appended to these coefficients-functions,
for the reason that {\em exactly two} further changes of initial
base coframes:
\[
\big\{
\rho_0,\,
\kappa_0,\,
\zeta_0',\,
\overline{\kappa}_0,\,
\overline{\zeta}_0'
\big\}
\ \ \ \ \
\leadsto
\ \ \ \ \
\big\{
\rho_0,\,
\kappa_0',\,
\zeta_0',\,
\overline{\kappa}_0',\,
\overline{\zeta}_0'
\big\}
\ \ \ \ \
\leadsto
\ \ \ \ \
\big\{
\rho_0,\,
\kappa_0',\,
\zeta_0'',\,
\overline{\kappa}_0',\,
\overline{\zeta}_0''
\big\}
\]
will force us to introduce {\em e.g.} $\Zaux_0^{i\prime}$ and
$\Zaux_0^{i\prime\prime}$, so that we will avoid
to use primes trice.

\Section{\bf Second Loop: Reduction of the Group Parameter ${\sf b}$}
\label{second-loop-reduction-b}
\HEAD{{\ref{second-loop-reduction-b}}.~{\sf Second Loop: 
Reduction of the Group Parameter ${\sf b}$}
}{
Wei Guo {\sc Foo} (Beijing) and Joël {\sc Merker} (Orsay)}

With this new reduced (real) eight-dimensional group $G^8$,
the lifted coframe, in which for simplicity we use the same 
letters $\rho$, $\kappa$, $\zeta$ as before, becomes:
\[
\left(\!
\begin{array}{c}
\rho
\\
\kappa
\\
\zeta
\end{array}
\!\right)
\,:=\,
\left(\!
\begin{array}{ccc}
{\sf c}\overline{\sf c} & 0 & 0
\\
{\sf b} & {\sf c} & 0
\\
{\sf d} & {\sf e} & \frac{{\sf c}}{\overline{\sf c}}
\end{array}
\!\right)
\left(\!
\begin{array}{c}
\rho_0
\\
\kappa_0
\\
\zeta_0'
\end{array}
\!\right)
\ \ \ \ \ \ \ \ \ \ \ \ \ \ \ \ \ \
\Longleftrightarrow
\ \ \ \ \ \ \ \ \ \ \ \ \ \ \ \ \ \
\left\{
\aligned
\rho
&
\,:=\,
{\sf c}\overline{\sf c}\,
\rho_0,
\\
\kappa
&
\,:=\,
{\sf b}\,\rho_0
+
{\sf c}\,\kappa_0,
\\
\zeta
&
\,:=\,
{\sf d}\,\rho_0
+
{\sf e}\,\kappa_0
+
\frac{{\sf c}}{\overline{\sf c}}\,
\zeta_0',
\endaligned\right.
\]
and inverse formulas are:
\leqnomode\usetagform{default}
\begin{align}
\label{2-loop-initial-inverted}
\rho_0
&
\,=\,
\frac{1}{{\sf c}\overline{\sf c}}\,
\rho,
\notag
\\
\kappa_0
&
\,=\,
-\,\frac{{\sf b}}{{\sf c}{\sf c}\overline{\sf c}}\,
\rho
+
\frac{1}{{\sf c}}\,
\kappa,
\\
\zeta_0'
&
\,=\,
\frac{{\sf b}{\sf e}-{\sf c}{\sf d}}{{\sf c}{\sf c}{\sf c}}\,
\rho
-
\frac{\overline{\sf c}{\sf e}}{{\sf c}{\sf c}}\,
\kappa
+
\frac{\overline{\sf c}}{{\sf c}}\,
\zeta.
\notag
\end{align}

The Maurer-Cartan matrix becomes:
\[
\aligned
dg\cdot g^{-1}
&
\,=\,
\left(\!
\begin{array}{ccc}
\overline{\sf c}\,d{\sf c}+{\sf c}d\overline{\sf c} & 0 & 0
\\
d{\sf b} & d{\sf c} & 0
\\
d{\sf d} & d{\sf e} & \frac{d{\sf c}}{\overline{\sf c}}
- \frac{{\sf c}\,d\overline{\sf c}}{\overline{\sf c}\overline{\sf c}}
\end{array}
\!\right)
\left(\!
\begin{array}{ccc}
\frac{1}{{\sf c}\overline{\sf c}} & 0 & 0
\\
-\frac{{\sf b}}{{\sf c}{\sf c}\overline{\sf c}} & \frac{1}{{\sf c}} & 0
\\
\frac{{\sf b}{\sf e}-{\sf c}{\sf d}}{{\sf c}{\sf c}{\sf c}}\, 
& -\frac{\overline{\sf c}{\sf e}}{{\sf c}{\sf c}} & 
\frac{\overline{\sf c}}{{\sf c}}
\end{array}
\!\right)
\\
&
\,=:\,
\left(\!
\begin{array}{ccc}
\alpha+\overline{\alpha} & 0 & 0
\\
\beta & \alpha & 0
\\
\gamma & \delta & \alpha-\overline{\alpha}
\end{array}
\!\right),
\endaligned
\]
in terms of the group-invariant $1$-forms:
\[
\aligned
\alpha
&
\,:=\,
\frac{d{\sf c}}{{\sf c}},
\\
\beta
&
\,:=\,
\frac{d{\sf b}}{{\sf c}\overline{\sf c}}
-
\frac{{\sf b}\,d{\sf c}}{{\sf c}{\sf c}\overline{\sf c}},
\\
\gamma
&
\,:=\,
\frac{d{\sf d}}{{\sf c}\overline{\sf c}}
-
\frac{{\sf b}\,d{\sf e}}{{\sf c}{\sf c}\overline{\sf c}}
+
\frac{{\sf b}{\sf e}-{\sf c}{\sf d}}{{\sf c}{\sf c}{\sf c}
\overline{\sf c}}\,
d{\sf c}
-
\frac{{\sf b}{\sf e}-{\sf c}{\sf d}}{{\sf c}{\sf c}\overline{\sf c}
\overline{\sf c}}\,
d\overline{\sf c},
\\
\delta
&
\,:=\,
\frac{d{\sf e}}{{\sf c}}
-
\frac{{\sf e}\,d{\sf c}}{{\sf c}{\sf c}}
+
\frac{{\sf e}\,d\overline{\sf c}}{{\sf c}\overline{\sf c}}.
\endaligned
\]

Now, let us exterior-differentiate the lifted coframe on the product
manifold equipped with coordinates:
\[
\big( 
z_1,z_2,\overline{z}_1,\overline{z}_2,v
\big)
\times
\big(
{\sf c},\overline{\sf c},
{\sf b},\overline{\sf b},
{\sf d},\overline{\sf d},{\sf e},\overline{\sf e}
\big)
\,\,\in\,\,
M^5
\times
G^8.
\]
The computation starts as:
\leqnomode\usetagform{default}
\begin{align}
\label{2-loop-before-replacement}
d\rho
&
\,=\,
\big(
\overline{\sf c}\,d{\sf c}
+
{\sf c}\,d\overline{\sf c}
\big)
\wedge
\rho_0
\,+
\notag
\\
&
\ \ \ \ \ \ \ \ \ \ \ \ \ \ \ \ \ \ \ \ \ \ \ \ \
+
{\sf c}\overline{\sf c}\,
d\rho_0,
\notag
\\
d\kappa
&
\,=\,
d{\sf b}\wedge\rho_0
+
d{\sf c}\wedge\kappa_0
\,+
\\
&
\ \ \ \ \ \ \ \ \ \ \ \ \ \ \ \ \ \ \ \ \ \ \ \ \
+
{\sf b}\,d\rho_0
+
{\sf c}\,d\kappa_0,
\notag
\\
d\zeta
&
\,=\,
d{\sf d}\wedge\rho_0
+
d{\sf e}\wedge\kappa_0
+
\Big(
\frac{d{\sf c}}{{\sf c}}
-
\frac{{\sf c}\,d\overline{\sf c}}{\overline{\sf c}\overline{\sf c}}
\Big)
\wedge\zeta_0'
\,+
\notag
\\
&
\ \ \ \ \ \ \ \ \ \ \ \ \ \ \ \ \ \ \ \ \ \ \ \ \
+
{\sf d}\,d\rho_0
+
{\sf e}\,d\kappa_0
+
\frac{{\sf c}}{\overline{\sf c}}\,
d\zeta_0'.
\notag
\end{align}

As is known, one must replace in
second lines $d\rho_0$, $d\kappa_0$,
$d\zeta_0'$ by the structure 
equations~({\ref{2-loop-structure-0-coframe}}),
and after, replace everywhere $\rho_0$, $\kappa_0$, $\zeta_0'$,
using the inversion 
formulas~({\ref{2-loop-initial-inverted}}).

However, contrary to Pocchiola's systematic approach, we will {\em
not} perform these calculations completely, but select only meaningful
terms.

At least, at the level of Maurer-Cartan forms, 
after replacements of $\rho_0$, $\kappa_0$, $\zeta_0'$
in the first lines of~({\ref{2-loop-before-replacement}}) above
using~({\ref{2-loop-initial-inverted}}), 
we have as usual:
\[
\aligned
d\rho
&
\,=\,
\big(\alpha+\overline{\alpha}\big)
\wedge
\rho
+
\Torsion,
\\
d\kappa
&
\,=\,
\beta\wedge\rho
+
\alpha\wedge\kappa
+
\Torsion,
\\
d\zeta
&
\,=\,
\gamma\wedge\rho
+
\delta\wedge\kappa
+
\big(
\alpha
-
\overline{\alpha}
\big)
\wedge\zeta
+
\Torsion.
\endaligned
\]

\begin{Question}
{\sl Without computing everything, 
what are the shapes of the three Torsion remainders?}
\end{Question}

Consider for instance what happens of the last term
$\frac{\sf c}{\overline{\sf c}}\, d\zeta_0'$ in $d\zeta$,
when peforming the required replacements, and
restrict attention even to the last term of 
$\frac{\sf c}{\overline{\sf c}}\,d\zeta_0'$
in~({\ref{2-loop-structure-0-coframe}}), which becomes:
\[
\frac{\sf c}{\overline{\sf c}}\,
\frac{\overline{\mathcal{L}}_1(\overline{\kaux})}{
\mathcal{L}_1(\overline{\kaux})}\,
\zeta_0'\wedge\overline{\zeta}_0'
\,=\,
\frac{\sf c}{\overline{\sf c}}\,
\frac{\overline{\mathcal{L}}_1(\overline{\kaux})}{
\mathcal{L}_1(\overline{\kaux})}\,
\bigg(
\frac{{\sf b}{\sf e}-{\sf c}{\sf d}}{{\sf c}{\sf c}{\sf c}}\,\rho
-
\frac{\overline{\sf c}{\sf e}}{{\sf c}{\sf c}}\,\kappa
+
\frac{\overline{\sf c}}{\sf c}\,\zeta
\bigg)
\wedge
\bigg(
\frac{\overline{\sf b}\overline{\sf e}-
\overline{\sf c}\overline{\sf d}}{
\overline{\sf c}\overline{\sf c}}\,
\rho
-
\frac{{\sf c}\overline{\sf e}}{\overline{\sf c}\overline{\sf c}}\,
\overline{\kappa}
+
\frac{{\sf c}}{\overline{\sf c}}\,
\overline{\zeta}
\bigg).
\]
After expansion, we see that are present the eight $2$-forms:
\[
\aligned
(\centersmallbullet)\,
\rho\wedge\kappa,
\ \ \ \ \ \ \ \ \ \ \ \ \ \ \ 
&
(\centersmallbullet)\,
\rho\wedge\zeta,
\ \ \ \ \ \ \ \ \ \ \ \ \ \ \ 
(\centersmallbullet)\,
\rho\wedge\overline{\kappa},
\ \ \ \ \ \ \ \ \ \ \ \ \ \ \ 
(\centersmallbullet)\,
\rho\wedge\overline{\zeta},
\\
&
(\centersmallbullet)\,
\kappa\wedge\overline{\kappa},
\ \ \ \ \ \ \ \ \ \ \ \ \ \ \ 
(\centersmallbullet)\,
\kappa\wedge\overline{\zeta},
\ \ \ \ \ \ \ \ \ \ \ \ \ \ \ 
(\centersmallbullet)\,
\zeta\wedge\overline{\kappa},
\ \ \ \ \ \ \ \ \ \ \ \ \ \ \ 
(\centersmallbullet)\,
\zeta\wedge\overline{\zeta}.
\endaligned
\]
Doing the same for all torsion terms, we may
realize\,\,---\,\,although it is not necessary to check this for what
follows\,\,---\,\,with almost no computation that the {\em
nonexplicit} shape of the structure equations of the lifted coframe
is:
\[
\aligned
d\rho
&
\,=\,
\big(\alpha+\overline{\alpha}\big)
\wedge
\rho
\,+
\\
&
\ \ \ \ \
+
R^1\,
\rho\wedge\kappa
+
R^2\,
\rho\wedge\zeta
+
\boxed{\overline{R}^1}\,
\rho\wedge\overline{\kappa}
+
\overline{R}^2\,
\rho\wedge\overline{\zeta}
+
\isqrt\,
\kappa\wedge\overline{\kappa},
\\
d\kappa
&
\,=\,
\beta\wedge\rho
+
\alpha\wedge\kappa
\,+
\\
&
\ \ \ \ \ 
+
K^1\,
\rho\wedge\kappa
+
K^2\,
\rho\wedge\zeta
+
K^3\,
\rho\wedge\overline{\kappa}
+
K^4\,
\rho\wedge\overline{\zeta}
\,+
\\
&
\ \ \ \ \ \ \ \ \ \ \ \ \ \ \ \ \ \ \ \ \ \ \ \ \ \ \ \ \ \ \ \ \ \ 
\ \ \ \ \ \ \ \ \ \   
+
K^5\,
\kappa\wedge\zeta
+
\boxed{K^6}\,
\kappa\wedge\overline{\kappa}
+
{\bf 1}
\cdot
\zeta\wedge\overline{\kappa},
\\
d\zeta
&
\,=\,
\gamma\wedge\rho
+
\delta\wedge\kappa
+
\big(\alpha-\overline{\alpha}\big)\wedge\zeta
\,+
\\
&
\ \ \ \ \ 
+
Z^1\,
\rho\wedge\kappa
+
Z^2\,
\rho\wedge\zeta
+
Z^3\,
\rho\wedge\overline{\kappa}
+
Z^4\,
\rho\wedge\overline{\zeta}
\,+
\\
&
\ \ \ \ \ \ \ \ \ \ \ \ \ \ \ \ \ \ \ \ \ \ \ \ \ \ \ \ \ \ \ \ \ \ 
\ \ \ \ \ \ \ \ \ \   
+
Z^5\,
\kappa\wedge\zeta
+
Z^6\,
\kappa\wedge\overline{\kappa}
+
Z^7\,
\kappa\wedge\overline{\zeta}
+
\boxed{Z^8}\,
\zeta\wedge\overline{\kappa}
+
Z^9\,
\zeta\wedge\overline{\zeta}.
\endaligned
\]
Of course, the preceding normalization ${\sf f} := \frac{{\sf c}}{
\overline{\sf c}}\, \overline{\mathcal{L}}_1(\kaux)$ forces:
\[
{\bf 1}
\,=\,
\big[
\zeta\wedge\overline{\kappa}
\big]
\big\{
d\kappa
\big\},
\]
a fact that can also be confirmed by a direct computation of this
torsion coefficient (exercise).

So we do not compute all torsion coefficients like Pocchiola did,
but we determine {\em before} some essential torsions,
so that we may focus on just the useful torsion terms.
In advance, we have boxed above 
the $3$ useful ones, shown by Pocchiola.
The subtle thing is that all three structure equations are needed.

\begin{Lemma}
\label{Lemma-R-1-2-K-6-Z-8}
Here is an essential linear combination of torsion terms:
\[
\overline{R}^1
-
2\,K^6
+
Z^8.
\]
\end{Lemma}

\proof
In order to '{\sl absorb}' as many torsion coefficients as possible,
let us substitute:
\[
\aligned
\alpha
&
\,=:\,
\alpha'
+
a_1\,\rho
+
a_2\,\kappa
+
a_3\,\zeta
+
a_4\,\overline{\kappa}
+
a_5\,\overline{\zeta},
\\
\beta
&
\,=:\,
\beta'
+
b_1\,\rho
+
b_2\,\kappa
+
b_3\,\zeta
+
b_4\,\overline{\kappa}
+
b_5\,\overline{\zeta},
\\
\gamma
&
\,=:\,
\gamma'
+
c_1\,\rho
+
c_2\,\kappa
+
c_3\,\zeta
+
c_4\,\overline{\kappa}
+
c_5\,\overline{\zeta},
\\
\delta
&
\,=:\,
\delta'
+
d_1\,\rho
+
d_2\,\kappa
+
d_3\,\zeta
+
d_4\,\overline{\kappa}
+
d_5\,\overline{\zeta}.
\endaligned
\]
At first, we have to transform the structure equations after 
such a substitution, the task is easy, and we write out the details
so that the reader needs no pen and no computer.

Substituting, the Maurer-Cartan part of $d\rho$ becomes:
\[
\aligned
\big(\alpha+\overline{\alpha}\big)
\wedge
\rho
\,=\,
\big(\alpha'+\overline{\alpha}'\big)
\wedge
\rho
&
+
0
+
a_2\,
\kappa\wedge\rho
+
a_3\,
\zeta\wedge\rho
+
a_4\,
\overline{\kappa}\wedge\rho
+
a_5\,
\overline{\zeta}\wedge\rho
\,+
\\
&
+
0
+
\overline{a}_2\,
\overline{\kappa}\wedge\rho
+
\overline{a}_3\,
\overline{\zeta}\wedge\rho
+
\overline{a}_4\,
\kappa\wedge\rho
+
\overline{a}_5\,
\zeta\wedge\rho,
\endaligned
\]
hence adding and reorganizing visually, we get:
\[
\aligned
d\rho
&
\,=\,
\big(\alpha'+\overline{\alpha}'\big)
\wedge
\rho
\,+
\\
&
\ \ \ \ \
+
\rho\wedge\kappa
\Big(
R^1
-
a_2
-
\overline{a}_4
\Big)
+
\rho\wedge\zeta
\Big(
R^2
-
a_3
-
\overline{a}_5
\Big)
+
\rho\wedge\overline{\kappa}
\Big(
\boxed{
\overline{R}^1
-
a_4
-
\overline{a}_2}
\Big)
\,+
\\
&
\ \ \ \ \ 
+
\rho\wedge\overline{\zeta}
\Big(
\overline{R}^2
-
a_5
-
\overline{a}_3
\Big)
+
\isqrt\,
\kappa\wedge\overline{\kappa}.
\endaligned
\]

Next:
\[
\aligned
\beta\wedge\rho
+
\alpha\wedge\kappa
&
\,=\,
\beta'\wedge\rho
+
0
+
b_2\,\kappa\wedge\rho
+
b_3\,\zeta\wedge\rho
+
b_4\,\overline{\kappa}\wedge\rho
+
b_5\,\overline{\zeta}\wedge\rho
\,+
\\
&
\ \ \ \ \
+
\alpha'\wedge\kappa
+
a_1\,\rho\wedge\kappa
+
0
+
a_3\,\zeta\wedge\kappa
+
a_4\,\overline{\kappa}\wedge\kappa
+
a_5\,\overline{\zeta}\wedge\kappa,
\endaligned
\]
hence:
\[
\aligned
d\kappa
&
\,=\,
\beta'\wedge\rho
+
\alpha'\wedge\kappa
\,+
\\
&
\ \ \ \ \ 
+
\rho\wedge\kappa\,
\Big(
K^1
+
a_1
-
b_2
\Big)
+
\rho\wedge\zeta\,
\Big(
K^2
-
b_3
\Big)
+
\rho\wedge\overline{\kappa}\,
\Big(
K^3
-
b_4
\Big)
+
\rho\wedge\overline{\zeta}\,
\Big(
K^4
-
b_5
\Big)
\,+
\\
&
\ \ \ \ \
+
\kappa\wedge\zeta\,
\Big(
K^5
-
a_3
\Big)
+
\kappa\wedge\overline{\kappa}\,
\Big(
\boxed{
K^6
-
a_4}
\Big)
+
\kappa\wedge\overline{\zeta}\,
\big(
-a_5
\big)
+
\zeta\wedge\overline{\kappa}.
\endaligned
\]

Lastly:
\[
\aligned
\gamma\wedge\rho
+
\delta\wedge\kappa
+
\big(\alpha-\overline{\alpha}\big)
\wedge
\zeta
&
\,=\,
\gamma'\wedge\rho
+
0
+
c_2\,
\kappa\wedge\rho
+
c_3\,
\zeta\wedge\rho
+
c_4\,
\overline{\kappa}\wedge\rho
+
c_5\,
\overline{\zeta}\wedge\rho
\,+
\\
&
\ \ \ \ \
+
\delta'\wedge\kappa
+
d_1\,
\rho\wedge\kappa
+
0
+
d_3\,
\zeta\wedge\kappa
+
d_4\,
\overline{\kappa}\wedge\kappa
+
d_5\,
\overline{\zeta}\wedge\kappa
\,+
\\
&
\ \ \ \ \ 
+
\alpha'\wedge\zeta
+
a_1\,
\rho\wedge\zeta
+
a_2\,
\kappa\wedge\zeta
+
0
+
a_4\,
\overline{\kappa}\wedge\zeta
+
a_5\,
\overline{\zeta}\wedge\zeta
\,-
\\
&
\ \ \ \ \
-\,
\overline{\alpha}'\wedge\zeta
-
\overline{a}_1\,
\rho\wedge\zeta
-
\overline{a}_2\,
\overline{\kappa}\wedge\zeta
-
\overline{a}_3\,
\overline{\zeta}\wedge\zeta
-
\overline{a}_4\,
\kappa\wedge\zeta
-
0,
\endaligned
\]
hence:
\[
\aligned
d\zeta
&
\,=\,
\gamma'\wedge\rho
+
\delta'\wedge\kappa
+
\big(
\alpha'-\overline{\alpha}'
\big)
\wedge
\zeta
\,+
\\
&
\ \ \ \ \
+
\rho\wedge\kappa
\Big(
Z^1
-
c_2
+
d_1
\Big)
+
\rho\wedge\zeta\,
\Big(
Z^2
-
c_3
+
a_1
-
\overline{a}_1
\Big)
+
\rho\wedge\overline{\kappa}\,
\Big(
Z^3
-
c_4
\Big)
+
\rho\wedge\overline{\zeta}\,
\Big(
Z^4
-
c_5
\Big)
\,+
\\
&
\ \ \ \ \
+
\kappa\wedge\zeta\,
\Big(
Z^5
-
d_3
+
a_2
-
\overline{a}_4
\Big)
+
\kappa\wedge\overline{\kappa}\,
\Big(
Z^6
-
d_4
\Big)
+
\kappa\wedge\overline{\zeta}\,
\Big(
Z^7
-
d_5
\Big)
\,+
\\
&
\ \ \ \ \
+
\zeta\wedge\overline{\kappa}\,
\Big(
\boxed{
Z^8
-
a_4
+
\overline{a}_2}
\Big)
+
\zeta\wedge\overline{\zeta}\,
\Big(
Z^9
-
a_5
+
\overline{a}_3
\Big).
\endaligned
\]
Extracting the boxed three new torsion coefficients:
\[
\aligned
\overline{R}^{1\prime}
&
\,=\,
\overline{R}^1
-
a_4
-
\overline{a}_2,
\\
K^{6\prime}
&
\,=\,
K^6
-
a_4,
\\
Z^{8\prime}
&
\,=\,
Z^8
-
a_4
+
\overline{a}_2,
\endaligned
\]
we see well the announced essentiality\big/invariancy of this
torsion combination:
\[
\overline{R}^{1\prime}
-
2\,K^{6\prime}
+
Z^{8\prime}
\,=\,
\overline{R}^1
-
2\,K^6
+
Z^8.
\qedhere
\]
\endproof

Consequently, we may restrict ourselves to computing only these
three torsion coefficients.

\begin{Lemma}
\label{Lemma-bar-R-1-K-6-Z-8}
Their explicit expressions are:
\[
\aligned
\overline{R}^1
&
\,=\,
\frac{\overline{\Paux}}{\overline{\sf c}}
+
\frac{{\sf c}\overline{\sf e}}{\overline{\sf c}\overline{\sf c}}\,
\frac{\overline{\mathcal{L}}_1(\overline{\kaux})}{
\mathcal{L}_1(\overline{\kaux})}
-
\isqrt\,
\frac{{\sf b}}{{\sf c}\overline{\sf c}},
\\
K^6
&
\,=\,
\isqrt\,
\frac{{\sf b}}{{\sf c}\overline{\sf c}}
-
\frac{{\sf e}}{{\sf c}},
\\
Z^8
&
\,=\,
\frac{{\sf e}}{{\sf c}}
-
\frac{1}{\overline{\sf c}}\,
\frac{\overline{\mathcal{L}}_1\big(
\overline{\mathcal{L}}_1(\kaux)\big)}{
\overline{\mathcal{L}}_1(\kaux)}
-
\frac{{\sf c}\overline{\sf e}}{\overline{\sf c}\overline{\sf c}}\,
\frac{\overline{\mathcal{L}}_1(\overline{\kaux})}{
\mathcal{L}_1(\overline{\kaux})}.
\endaligned
\]
\end{Lemma}

\proof
We proceed by chasing coefficients. Let us treat $\overline{R}^1$.
From~({\ref{2-loop-before-replacement}}),  
replacing in~({\ref{2-loop-structure-0-coframe}}) by means
of~({\ref{2-loop-initial-inverted}}), we reach its expression:
\[
\aligned
\overline{R}^1
\,=\,
\big[
\rho\wedge\overline{\kappa}
\big]
\big\{
{\sf c}\overline{\sf c}\,d\rho_0
\big\}
&
\,=\,
0
+
0
+
\big[
\rho\wedge\overline{\kappa}
\big]
\bigg\{
{\sf c}\overline{\sf c}\,
\overline{\Paux}\,
\Big(
\frac{1}{{\sf c}\overline{\sf c}}\,
\rho
\Big)
\wedge
\Big(
-
\frac{\overline{\sf b}}{\overline{\sf c}\overline{\sf c}{\sf c}}\,
\rho
+
\frac{1}{\overline{\sf c}}\,
\overline{\kappa}
\Big)
\,-
\\
&
\ \ \ \ \ \ \ \ \ \ \ \ \ \ \ \ \ \ \ \ \ \ \ \ \ \ \ \ \ \ \ \ \ \ \ \
\,-
{\sf c}\overline{\sf c}\,
\frac{\overline{\mathcal{L}}_1(\overline{\kaux})}{
\mathcal{L}_1(\overline{\kaux})}\,
\Big(
\frac{1}{{\sf c}\overline{\sf c}}\,
\rho
\Big)
\wedge
\Big(
\frac{\overline{\sf b}\overline{\sf e}-
\overline{\sf c}\overline{\sf d}}{
\overline{\sf c}\overline{\sf c}\overline{\sf c}}\,
\rho
-
\frac{{\sf c}\overline{\sf e}}{
\overline{\sf c}\overline{\sf c}}\,
\overline{\kappa}
+
\frac{{\sf c}}{\overline{\sf c}}\,
\overline{\zeta}
\Big)
\,+
\\
&
\ \ \ \ \ \ \ \ \ \ \ \ \ \ \ \ \ \ \ \ \ \ \ \ \ \ \ \ \ \ \ \ \ \ \ \
+
{\sf c}\overline{\sf c}\,\isqrt\,
\Big(
-
\frac{{\sf b}}{{\sf c}{\sf c}\overline{\sf c}}\,
\rho
+
\frac{1}{{\sf c}}\,
\kappa
\Big)
\wedge
\Big(
-
\frac{\overline{\sf b}}{{\sf c}\overline{\sf c}\overline{\sf c}}\,
\rho
+
\frac{1}{\overline{\sf c}}\,
\overline{\kappa}
\Big)
\bigg\}
\\
&
\,=\,
\zero{{\sf c}\overline{\sf c}}\,
\overline{\Paux}\,
\frac{1}{\zero{{\sf c}\overline{\sf c}}}\,
\frac{1}{\overline{\sf c}}
+
\zero{{\sf c}\overline{\sf c}}\,
\frac{\overline{\mathcal{L}}_1(\overline{\kaux})}{
\mathcal{L}_1(\overline{\kaux})}\,
\frac{1}{\zero{{\sf c}\overline{\sf c}}}\,
\frac{{\sf c}\overline{\sf e}}{\overline{\sf c}\overline{\sf c}}
-
\isqrt\,
\zero{{\sf c}\overline{\sf c}}\,
\frac{{\sf b}}{{\sf c}\zero{{\sf c}\overline{\sf c}}}\,
\frac{1}{\overline{\sf c}}.
\endaligned
\]

Next, from~({\ref{2-loop-before-replacement}}), let us treat:
\[
K^6
\,=\,
\big[
\kappa\wedge\overline{\kappa}
\big]
\big\{
{\sf b}\,d\rho_0
+
{\sf c}\,d\kappa_0
\big\}.
\]
In ${\sf b}\, d\rho_0$, the first four terms 
in~({\ref{2-loop-structure-0-coframe}})
have zero contribution, since they are multiples of $\rho_0$, hence
of $\rho$, whence:
\[
\aligned
\big[
\kappa\wedge\overline{\kappa}
\big]
\big\{
{\sf b}\,d\rho_0
\big\}
&
\,=\,
0
+
0
+
0
+
0
+
\big[
\kappa\wedge\overline{\kappa}
\big]
\big\{
{\sf b}\,
\isqrt\,
\kappa_0
\wedge
\overline{\kappa}_0
\big\}
\\
&
\,=\,
\big[
\kappa\wedge\overline{\kappa}
\big]
\bigg\{
\isqrt\,{\sf b}\,
\Big(
-
\frac{{\sf b}}{{\sf c}{\sf c}\overline{\sf c}}\,
\rho
+
\frac{1}{{\sf c}}\,
\kappa
\Big)
\wedge
\Big(
-
\frac{\overline{\sf b}}{{\sf c}\overline{\sf c}\overline{\sf c}}\,
\rho
+
\frac{1}{\overline{\sf c}}\,
\overline{\kappa}
\Big)
\bigg\}
\\
&
\,=\,
\isqrt\,
\frac{{\sf b}}{{\sf c}\overline{\sf c}}.
\endaligned
\]
Also, in ${\sf c}\, d\kappa_0$, the first two terms contribute
$0$, and it remains:
\[
\aligned
\big[
\kappa\wedge\overline{\kappa}
\big]
\big\{
{\sf c}\,d\kappa_0
\big\}
&
\,=\,
0
+
0
+
\big[
\kappa\wedge\overline{\kappa}
\big]
\big\{
{\sf c}\,
\zeta_0'\wedge\overline{\kappa}_0
\big\}
\\
&
\,=\,
\big[
\kappa\wedge\overline{\kappa}
\big]
\bigg\{
{\sf c}\,
\Big(
-
\frac{\overline{\sf c}{\sf e}}{{\sf c}{\sf c}}\,
\kappa
\Big)
\wedge
\Big(
\frac{1}{\overline{\sf c}}\,
\overline{\kappa}
\Big)
\bigg\}
\\
&
\,=\,
-\,
\frac{{\sf e}}{{\sf c}}.
\endaligned
\]

Lastly:
\[
Z^8
\,=\,
\big[
\zeta\wedge\overline{\kappa}
\big]
\Big\{
{\sf d}\,d\rho_0
+
{\sf e}\,d\kappa_0
+
\frac{{\sf c}}{\overline{\sf c}}\,
d\zeta_0'
\Big\}.
\]
Here, ${\sf d}\, d\rho_0$ contributes $0$. Next, the first two terms
in ${\sf e}\, d\kappa_0$ contribute $0$, and it remains:
\[
\aligned
\big[
\zeta\wedge\overline{\kappa}
\big]
\big\{
{\sf e}\,
d\kappa_0
\big\}
&
\,=\,
\big[
\zeta\wedge\overline{\kappa}
\big]
\big\{
{\sf e}\,
\zeta_0'\wedge\kappa_0
\big\}
\\
&
\,=\,
\big[
\zeta\wedge\overline{\kappa}
\big]
\bigg\{
{\sf e}\,
\Big(
\frac{\overline{\sf c}}{{\sf c}}\,
\zeta
\Big)
\wedge
\Big(
\frac{1}{\overline{\sf c}}\,
\overline{\kappa}
\Big)
\bigg\}
\\
&
\,=\,
\frac{{\sf e}}{{\sf c}}.
\endaligned
\]
Also, in $\frac{{\sf c}}{\overline{\sf c}}\, d\zeta_0'$, the first
two terms contribute $0$, and the last two terms are:
\[
\aligned
\big[
\zeta\wedge\overline{\kappa}
\big]
\Big\{
\frac{{\sf c}}{\overline{\sf c}}\,
d\zeta_0'
\Big\}
&
\,=\,
-\,
\frac{{\sf c}}{\overline{\sf c}}\,
\frac{\overline{\mathcal{L}}_1\big(
\overline{\mathcal{L}}_1(\kaux)\big)}{
\overline{\mathcal{L}}_1(\kaux)}\,
\big[
\zeta
\wedge
\overline{\kappa}
\big]
\bigg\{
\Big(
\frac{\overline{\sf c}}{{\sf c}}\,
\zeta
\Big)
\wedge
\Big(
\frac{1}{\overline{\sf c}}\,
\overline{\kappa}
\Big)
\bigg\}
\\
&
\ \ \ \ \
+
\frac{{\sf c}}{\overline{\sf c}}\,
\frac{\overline{\mathcal{L}}_1(\overline{\kaux})}{
\mathcal{L}_1(\overline{\kaux})}\,
\big[
\zeta\wedge\overline{\kappa}
\big]\,
\bigg\{
\Big(
\frac{\overline{\sf c}}{{\sf c}}\,
\zeta
\Big)
\wedge
\Big(
-
\frac{{\sf c}\overline{\sf e}}{\overline{\sf c}\overline{\sf c}}\,
\overline{\kappa}
\Big)
\bigg\}
\\
&
\,=\,
-\,
\frac{1}{\overline{\sf c}}\,
\frac{\overline{\mathcal{L}}_1\big(
\overline{\mathcal{L}}_1(\kaux)\big)}{
\overline{\mathcal{L}}_1(\kaux)}\,
-
\frac{{\sf c}\overline{\sf e}}{\overline{\sf c}\overline{\sf c}}\,
\frac{\overline{\mathcal{L}}_1(\overline{\kaux})}{
\mathcal{L}_1(\overline{\kaux})}\,
\endaligned
\]
Adding, we get $Z^8$.
\endproof

Observing that necessarily $-a_5 = 0$ from $\big[ \kappa \wedge 
\overline{\zeta} \big] \big\{ d\kappa \big\}$, we realize that
some other invariant 
relations between torsion coefficients appear:
\[
\aligned
R^{2\prime}
-
K^{5\prime}
\,=\,
R^2
-
K^5,
\\
\overline{R}^{2\prime}
+
Z^{9\prime}
\,=\,
\overline{R}^2
+
Z^9,
\endaligned
\]
that could potentially bring normalizations of
some group parameters, but will not, 
as it will come out that they are identically satisfied.
However, knowing them will be very useful later, hence we state a 
supplementary

\begin{Assertion}
Three other torsion coefficients have 
the common explicit expression:
\[
R^2
\,=\,
K^5
\,=\,
-\,
\overline{Z}^9
\,\,=\,\,
-\,
\frac{\overline{\sf c}}{{\sf c}}\,
\frac{\mathcal{L}_1(\kaux)}{\overline{\mathcal{L}}_1(\kaux)}.
\]
\end{Assertion}

\proof
Our technique gives:
\[
\aligned
R^2
&
\,=\,
\big[
\rho\wedge\zeta
\big]
\big\{
{\sf c}\overline{\sf c}\,
d\rho_0
\big\}
\\
&
\,=\,
0
-
{\sf c}\overline{\sf c}\,
\frac{\mathcal{L}_1(\kaux)}{\overline{\mathcal{L}}_1(\kaux)}\,
\frac{1}{{\sf c}\overline{\sf c}}\,
\frac{\overline{\sf c}}{{\sf c}}
+
0
+
0
+
0.
\endaligned
\]
Next:
\[
\aligned
K^5
&
\,=\,
\big[
\zeta\wedge\kappa
\big]
\big\{
{\sf b}\,d\rho_0
+
{\sf c}\,d\kappa_0
\big\}
\\
&
\,=\,
0
+
\big[
\zeta\wedge\kappa
\big]
\big\{
{\sf c}\,d\kappa_0
\big\}
\\
&
\,=\,
0
-
\frac{\mathcal{L}_1(\kaux)}{\overline{\mathcal{L}}_1(\kaux)}\,
{\sf c}\,
\frac{1}{{\sf c}}\,
\frac{\overline{\sf c}}{{\sf c}}
+
0.
\endaligned
\]
Lastly:
\begin{align}
Z^9
&
\,=\,
\big[
\overline{\kappa}\wedge\overline{\zeta}
\big]
\Big\{
{\sf d}\,d\rho_0
+
{\sf e}\,d\kappa_0
+
\frac{{\sf c}}{\overline{\sf c}}\,
d\zeta_0'
\Big\}
\notag
\\
&
\,=\,
0
+
0
+
\big[
\overline{\kappa}\wedge\overline{\zeta}
\big]
\Big\{
\frac{{\sf c}}{\overline{\sf c}}\,
d\zeta_0'
\Big\}
\notag
\\
&
\,=\,
0
+
0
+
0
+
\frac{{\sf c}}{\overline{\sf c}}\,
\frac{\overline{\mathcal{L}}_1(\overline{\kaux})}{
\mathcal{L}_1(\overline{\kaux})}\,
\frac{\overline{\sf c}}{{\sf c}}\,
\frac{{\sf c}}{\overline{\sf c}}.
\qedhere
\end{align}
\endproof

Coming back to 
Lemma~{\ref{Lemma-bar-R-1-K-6-Z-8}},
we can now compute in details, emphasizing one annihilation,
the expression of the interesting invariant torsion
combination:
\[
\aligned
\overline{R}^1
-
2\,K^6
+
Z^8
&
\,=\,
\frac{\overline{\Paux}}{\overline{\sf c}}
+
\zero{
\frac{{\sf c}\overline{\sf e}}{\overline{\sf c}\overline{\sf c}}\,
\frac{\overline{\mathcal{L}}_1(\overline{\kaux})}{
\mathcal{L}_1(\overline{\kaux})}}
-
\isqrt\,
\frac{{\sf b}}{{\sf c}\overline{\sf c}}
\,-
\\
&
\ \ \ \ \
-\,
2\,\isqrt\,
\frac{{\sf b}}{{\sf c}\overline{\sf c}}
+
2\,\frac{{\sf e}}{{\sf c}}
\,+
\\
&
\ \ \ \ \
+
\frac{{\sf e}}{{\sf c}}
-
\frac{1}{\overline{\sf c}}\,
\frac{\overline{\mathcal{L}}_1\big(
\overline{\mathcal{L}}_1(\kaux)\big)}{
\overline{\mathcal{L}}_1(\kaux)}
-
\zero{
\frac{{\sf c}\overline{\sf e}}{\overline{\sf c}\overline{\sf c}}\,
\frac{\overline{\mathcal{L}}_1(\overline{\kaux})}{
\mathcal{L}_1(\overline{\kaux})}}
\\
&
\,=\,
-\,
3\,\isqrt\,
\frac{{\sf b}}{{\sf c}\overline{\sf c}}
+
3\,\frac{{\sf e}}{{\sf c}}
-
\frac{1}{\overline{\sf c}}\,
\left(
\frac{\overline{\mathcal{L}}_1\big(
\overline{\mathcal{L}}_1(\kaux)\big)}{
\overline{\mathcal{L}}_1(\kaux)}
-
\overline{\Paux}
\right).
\endaligned
\]
Since the group parameter ${\sf b} \in \C$ is not on the diagonal,
there is no restriction for it to be nonzero, hence we can 
normalize it by requiring that:
\[
0
\,=\,
\overline{R}^1
-
2\,K^6
+
Z^8,
\]
and this produces the announced normalization:
\leqnomode\usetagform{default}
\begin{align}
\label{normalization-b}
\boxed{\,
{\sf b}
\,:=\,
-\,\isqrt\,
\overline{\sf c}{\sf e}
+
\frac{\isqrt}{3}\,
{\sf c}
\left(
\frac{\overline{\mathcal{L}}_1\big(
\overline{\mathcal{L}}_1(\kaux)\big)}{
\overline{\mathcal{L}}_1(\kaux)}
-
\overline{\Paux}
\right).\,}
\end{align}

For convenience, let us abbreviate:
\[
\Baux_0
\,:=\,
\frac{\overline{\mathcal{L}}_1\big(
\overline{\mathcal{L}}_1(\kaux)\big)}{
\overline{\mathcal{L}}_1(\kaux)}
-
\overline{\Paux},
\]
which is function on $M$, as its lower index ${}_0$ points out, 
so that:
\[
\boxed{\,
{\sf b}
\,:=\,
-\,\isqrt\,
\overline{\sf c}{\sf e}
+
\frac{\isqrt}{3}\,
{\sf c}\,
\Baux_0.\,}
\]

After this normalization, the lifted coframe becomes:
\[
\left(\!
\begin{array}{c}
\rho
\\
\kappa
\\
\zeta
\end{array}
\!\right)
\,:=\,
\left(\!
\begin{array}{ccc}
{\sf c}\overline{\sf c} & 0 & 0
\\
-\isqrt\,\overline{\sf c}{\sf e}
+
\frac{\isqrt}{3}\,{\sf c}\Baux_0
& {\sf c} & 0
\\
{\sf d} & {\sf e} & \frac{{\sf c}}{\overline{\sf c}}
\end{array}
\!\right)
\left(\!
\begin{array}{c}
\rho_0
\\
\kappa_0
\\
\zeta_0'
\end{array}
\!\right).
\]
Consequently, we can transform\big/rewrite in a natural way:
\[
\aligned
\kappa
&
\,=\,
\Big(
-\isqrt\,\overline{\sf c}{\sf e}
+
\frac{\isqrt}{3}\,
{\sf c}\,
\Baux_0
\Big)\,
\rho_0
+
{\sf c}\,\kappa_0
\\
&
\,=\,
\big(
-\isqrt\,\overline{\sf c}{\sf e}
\big)\,\rho_0
+
{\sf c}\,
\Big(
\underbrace{
\kappa_0
+
\frac{\isqrt}{3}\,
\Baux_0\,\rho_0}_{
=:\,\,\kappa_0'}
\Big),
\endaligned
\]
and this conducts us to change of initial coframe on $M$:
\[
\big\{
\rho_0,\,
\kappa_0,\,
\zeta_0',
\overline{\kappa}_0,\,
\overline{\zeta}_0'
\big\}
\ \ \ \ \
\leadsto
\ \ \ \ \
\big\{
\rho_0,\,
\kappa_0',\,
\zeta_0',\,
\overline{\kappa}_0',
\overline{\zeta}_0'
\big\},
\]
by introducing:
\leqnomode\usetagform{default}
\begin{align}
\label{def-kappa-prime-0}
\kappa_0'
&
\,:=\,
\kappa_0
+
\frac{\isqrt}{3}\,
\Baux_0\,\rho_0.
\end{align}
It follows that:
\[
\aligned
\zeta
\,=\,
{\sf d}\,\rho_0
+
{\sf e}\,\kappa_0
+
\frac{{\sf c}}{\overline{\sf c}}\,
\zeta_0'
&
\,=\,
{\sf d}\,\rho_0
+
{\sf e}\,
\Big(
\kappa_0'
-
\frac{\isqrt}{3}\,
\Baux_0\,
\rho_0
\Big)
+
\frac{{\sf c}}{\overline{\sf c}}\,
\zeta_0'
\\
&
\,=\,
\Big(
\underbrace{
{\sf d}
-
\frac{\isqrt}{3}\,
{\sf e}\,
\Baux_0}_{
=:\,\,{\sf d}'}
\Big)\,
\rho_0
+
{\sf e}\,\kappa_0'
+
\frac{{\sf c}}{\overline{\sf c}}\,
\zeta_0'.
\endaligned
\]
Before, ${\sf d} \in \C$ was a parameter representing some unknown
function. Introducing the new unknown\big/parameter:
\[
{\sf d}'
\,:=\,
{\sf d}
-
\frac{\isqrt}{3}\,
{\sf e},
\]
we come to a new $G$-structure of real dimension $6$ 
parametrized by ${\sf c}, {\sf e} \in \C^\ast$ and ${\sf d}' \in \C$
whose lifted coframe writes:
\[
\left(\!
\begin{array}{c}
\rho
\\
\kappa
\\
\zeta
\end{array}
\!\right)
\,:=\,
\left(\!
\begin{array}{ccc}
{\sf c}\overline{\sf c} & 0 & 0
\\
-\isqrt\,\overline{\sf c}{\sf e} & {\sf c} & 0
\\
{\sf d}' & {\sf e} & \frac{{\sf c}}{\overline{\sf c}}
\end{array}
\!\right)
\left(\!
\begin{array}{c}
\rho_0
\\
\kappa_0'
\\
\zeta_0'
\end{array}
\!\right).
\]
We will write again ${\sf d}$ instead of ${\sf d}'$.

\Section{\bf Darboux-Cartan Structure of the Coframe
$\big\{ \rho_0, \kappa_0', \zeta_0', \overline{\kappa}_0', 
\overline{\zeta}_0' \big\}$}
\label{D-C-structure-kappa-0-prime-zeta-0-prime}
\HEAD{{\ref{D-C-structure-kappa-0-prime-zeta-0-prime}}.~{\sf 
Darboux-Cartan Structure of the Coframe
$\big\{ \rho_0, \kappa_0', \zeta_0', \overline{\kappa}_0', 
\overline{\zeta}_0' \big\}$}
}{
Wei Guo {\sc Foo} (Beijing) and Joël {\sc Merker} (Orsay)}

Before continuing, we must compute the Darboux-Cartan structure 
of this new initial coframe 
$\big\{ \rho_0, \kappa_0', \zeta_0', \overline{\kappa}_0', 
\overline{\zeta}_0'\big\}$, for which absolutely
no details were provided 
in~{\cite{Pocchiola-2013, Merker-Pocchiola-2018}}.
Here, we offer complete explanations.

Abstractly, the structure in question will have the shape:
\leqnomode\usetagform{default}
\begin{align}
\label{inexplicit-0-structure-3-loop}
d\rho_0
&
\,=\,
\Raux_0^{1\prime}\,
\rho_0\wedge\kappa_0'
+
\Raux_0^{2\prime}\,
\rho_0\wedge\zeta_0'
+
\overline{\Raux}_0^{1\prime}\,
\rho_0\wedge\overline{\kappa}_0'
+
\overline{\Raux}_0^{2\prime}\,
\rho_0\wedge\overline{\zeta}_0'
+
\isqrt\,
\kappa_0'\wedge\overline{\kappa}_0',
\notag
\\
\ \ \ \ \ \
d\kappa_0'
&
\,=\,
\Kaux_0^{1\prime}\,
\rho_0\wedge\kappa_0'
+
\Kaux_0^{2\prime}\,
\rho_0\wedge\zeta_0'
+
\Kaux_0^{3\prime}\,
\rho_0\wedge\overline{\kappa}_0'
\,+
\\
&
\ \ \ \ \ \ \ \ \ \ \ \ \ \ \ \ \ \ \ \
+
\Kaux_0^{5\prime}\,
\kappa_0'\wedge\zeta_0'
+
\Kaux_0^{6\prime}\,
\kappa_0'\wedge\overline{\kappa}_0'
+
\zeta_0'\wedge\overline{\kappa}_0',
\notag
\\
d\zeta_0'
&
\,=\,
\Zaux_0^{2\prime}\,
\rho_0\wedge\zeta_0'
+
\Zaux_0^{5\prime}\,
\kappa_0'\wedge\zeta_0'
+
\Zaux_0^{8\prime}\,
\zeta_0'\wedge\overline{\kappa}_0'
+
\Zaux_0^{9\prime}\,
\zeta_0'\wedge\overline{\zeta}_0'.
\notag
\end{align}
Our goal is to compute {\em explicitly} all these coefficients,
and the answer is stated as follows.

\begin{Proposition}
\label{Proposition-3-loop-structure-0-coframe}
The Darboux-Cartan structure for the initial coframe 
$\big\{ \rho_0, \kappa_0', \zeta_0', \overline{\kappa}_0', 
\overline{\zeta}_0' \big\}$ expands as:
\[
\aligned
d\rho_0
&
\,=\,
\left(
\frac{1}{3}\,
\frac{\mathcal{L}_1\big(\mathcal{L}_1(\overline{\kaux})\big)}{
\mathcal{L}_1(\overline{\kaux})}
+
\frac{2}{3}\,
\Paux
\right)
\rho_0\wedge\kappa_0'
-
\frac{\mathcal{L}_1(\kaux)}{\overline{\mathcal{L}}_1(\kaux)}\,
\rho_0\wedge\zeta_0'
\,+
\\
&
\ \ \ \ \
+
\left(
\frac{1}{3}\,
\frac{\overline{\mathcal{L}}_1\big(\overline{\mathcal{L}}_1
(\kaux)\big)}{\overline{\mathcal{L}}_1(\kaux)}
+
\frac{2}{3}\,\overline{\Paux}
\right)
\rho_0\wedge\overline{\kappa}_0'
-
\frac{\overline{\mathcal{L}}_1(\overline{\kaux})}{
\mathcal{L}_1(\overline{\kaux})}\,
\rho_0\wedge\overline{\zeta}_0'
+
\isqrt\,
\kappa_0'\wedge\overline{\kappa}_0',
\endaligned
\]
\[
\aligned
d\kappa_0'
&
\,=\,
\left(
-\,\frac{\isqrt}{3}\,
\frac{\mathcal{L}_1\big(\overline{\mathcal{L}}_1\big(
\overline{\mathcal{L}}_1(\kaux)\big)\big)}{
\overline{\mathcal{L}}_1(\kaux)}
+
\frac{\isqrt}{9}\,
\frac{\mathcal{L}_1\big(\mathcal{L}_1(\overline{\kaux})\big)\,\,
\overline{\mathcal{L}}_1\big(\overline{\mathcal{L}}_1(\kaux)\big)}{
\mathcal{L}_1(\overline{\kaux})\,\,\overline{\mathcal{L}}_1
(\kaux)}
\,+
\right.
\\
&
\ \ \ \ \ \ \ \ \ \
+
\left.
\frac{\isqrt}{3}\,
\frac{\mathcal{L}_1\big(\overline{\mathcal{L}}_1(\kaux)\big)\,\,
\overline{\mathcal{L}}_1\big(\overline{\mathcal{L}}_1(\kaux)\big)}{
\overline{\mathcal{L}}_1(\kaux)^2}
-
\frac{\isqrt}{9}\,
\frac{\mathcal{L}_1\big(\mathcal{L}_1(\overline{\kaux})\big)}{
\mathcal{L}_1(\overline{\kaux})}\,
\overline{\Paux}
\,+
\right.
\\
&
\ \ \ \ \ \ \ \ \ \
\left.
+
\frac{2\,\isqrt}{9}\,
\frac{\overline{\mathcal{L}}_1\big(\overline{\mathcal{L}}_1
(\kaux)\big)}{\overline{\mathcal{L}}_1(\kaux)}\,
\Paux
+
\frac{\isqrt}{3}\,
\mathcal{L}_1\big(\overline{\Paux}\big)
-
\frac{2\,\isqrt}{9}\,
\Paux\,\overline{\Paux}
\right)
\rho_0\wedge\kappa_0'
\,+
\\
&
\ \ \ \ \
+
\left(
-\,\frac{\isqrt}{3}\,
\frac{\mathcal{K}\big(\overline{\mathcal{L}}_1\big(
\overline{\mathcal{L}}_1(\kaux)\big)\big)}{
\overline{\mathcal{L}}_1(\kaux)^2}
+
\frac{\isqrt}{3}\,
\frac{\mathcal{K}\big(\overline{\mathcal{L}}_1(\kaux)\big)\,\,
\overline{\mathcal{L}}_1\big(\overline{\mathcal{L}}_1(\kaux)\big)}{
\overline{\mathcal{L}}_1(\kaux)^3}
\,-
\right.
\\
&
\ \ \ \ \ \ \ \ \ \ \ \ \ \ \ \ \ \ \ \ \ \ \ \ \ \ \ \ \ \ \ \ \ \ 
\ \ \ \ \ \ \ \ \ \ \ \ \ 
\left.
\,-
\frac{\isqrt}{3}\,
\frac{\mathcal{L}_1\big(\mathcal{L}_1(\overline{\kaux})\big)}{
\mathcal{L}_1(\overline{\kaux})}
-
\frac{\isqrt}{3}\,
\frac{\overline{\mathcal{L}}_1\big(
\mathcal{L}_1(\kaux)\big)}{
\overline{\mathcal{L}}_1(\kaux)}
-
\frac{2}{3}\,
\frac{\mathcal{T}(\kaux)}{\overline{\mathcal{L}}_1(\kaux)}
\right)\,
\rho_0\wedge\zeta_0'
\,+
\\
&
\ \ \ \ \ 
+
\left(
-\,\frac{\isqrt}{3}\,
\frac{\overline{\mathcal{L}}_1\big(\overline{\mathcal{L}}_1\big(
\overline{\mathcal{L}}_1(\kaux)\big)\big)}{
\overline{\mathcal{L}}_1(\kaux)}
+
\frac{4\,\isqrt}{9}\,
\frac{\overline{\mathcal{L}}_1\big(\overline{\mathcal{L}}_1
(\kaux)\big)^2}{\overline{\mathcal{L}}_1(\kaux)^2}
\,+
\right.
\\
&
\ \ \ \ \ \ \ \ \ \ \ \ \ \ \ \ \ \ \ \ \ \ \ \ \ \ \ \ \ \ 
\left.
+
\frac{\isqrt}{9}\,
\frac{\overline{\mathcal{L}}_1\big(\overline{\mathcal{L}}_1
(\kaux)\big)}{\overline{\mathcal{L}}_1(\kaux)}\,
\overline{\Paux}
+
\frac{\isqrt}{3}\,
\overline{\mathcal{L}}_1\big(\overline{\Paux}\big)
-
\frac{2\,\isqrt}{9}\,
\overline{\Paux}\,\overline{\Paux}
\right)\,
\rho_0\wedge\overline{\kappa}_0'
\,+
\\
&
\ \ \ \ \
+
0\,
\rho_0\wedge\overline{\zeta}_0'
\,-
\\
&
\ \ \ \ \
-\,
\frac{\mathcal{L}_1(\kaux)}{\overline{\mathcal{L}}_1(\kaux)}\,
\kappa_0'\wedge\zeta_0'
+
\left(
-\frac{1}{3}\,
\frac{\overline{\mathcal{L}}_1\big(
\overline{\mathcal{L}}_1(\kaux)\big)}{
\overline{\mathcal{L}}_1(\kaux)}
+
\frac{1}{3}\,\overline{\Paux}
\right)
\kappa_0'\wedge\overline{\kappa}_0'
+
\zeta_0'\wedge\overline{\kappa}_0',
\endaligned
\]
\[
\aligned
d\zeta_0'
&
\,=\,
\left(
\frac{\isqrt}{3}\,
\frac{\mathcal{L}_1\big(\mathcal{L}_1(\overline{\kaux})\big)\,\,
\overline{\mathcal{L}}_1\big(\overline{\mathcal{L}}_1(\kaux)\big)}{
\mathcal{L}_1(\overline{\kaux})\,\,\overline{\mathcal{L}}_1(\kaux)}
-
\frac{\isqrt}{3}\,
\frac{\mathcal{L}_1\big(\overline{\mathcal{L}}_1(\kaux)\big)\,\,
\overline{\mathcal{L}}_1\big(\overline{\mathcal{L}}_1(\kaux)\big)}{
\overline{\mathcal{L}}_1(\kaux)^2}
\,-
\right.
\\
&
\ \ \ \ \ \ \ \ \ \
\left.
-\,\frac{\isqrt}{3}\,
\frac{\overline{\mathcal{L}}_1\big(
\overline{\mathcal{L}}_1(\kaux)\big)}{
\overline{\mathcal{L}}_1(\kaux)}\,
\Paux
+
\frac{\isqrt}{3}\,
\frac{\mathcal{L}_1\big(\overline{\mathcal{L}}_1(\kaux)\big)}{
\overline{\mathcal{L}}_1(\kaux)}\,
\overline{\Paux}
+
\frac{\mathcal{T}\big(\overline{\mathcal{L}}_1(\kaux)\big)}{
\overline{\mathcal{L}}_1(\kaux)}
\right)
\rho_0\wedge\zeta_0'
\,+
\\
&
\ \ \ \ \
+
\frac{\mathcal{L}_1\big(\overline{\mathcal{L}}_1(\kaux)\big)}{
\overline{\mathcal{L}}_1(\kaux)}\,
\kappa_0'\wedge\zeta_0'
-
\frac{\overline{\mathcal{L}}_1\big(
\overline{\mathcal{L}}_1(\kaux)\big)}{
\overline{\mathcal{L}}_1(\kaux)}\,
\zeta_0'\wedge\overline{\kappa}_0'
+
\frac{\overline{\mathcal{L}}_1(\overline{\kaux})}{
\mathcal{L}_1(\overline{\kaux})}\,
\zeta_0'\wedge\overline{\zeta}_0'.
\endaligned
\]
\end{Proposition}

Observe from these explicit expressions that:
\[
2\,\Kaux_0^{6\prime}
\,=\,
\overline{\Raux}_0^{1\prime}
+
\Zaux_0^{8\prime}
\ \ \ \ \ \ \ \ \ \ \ \ \ \ \ \ \ \
\text{and}
\ \ \ \ \ \ \ \ \ \ \ \ \ \ \ \ \ \
\Raux_0^{2\prime}
\,=\,
\Kaux_0^{5\prime}.
\]

\proof
We treat first $d\rho_0$ and $d\zeta_0'$, which are easier than
$d\kappa_0'$.

Observing from~({\ref{def-kappa-prime-0}}), that:
\[
\rho_0\wedge\kappa_0
\,=\,
\rho_0\wedge\kappa_0'
\ \ \ \ \ \ \ \ \ \ \ \ \ \ \ \ \ \
\text{and}
\ \ \ \ \ \ \ \ \ \ \ \ \ \ \ \ \ \
\rho_0\wedge\overline{\kappa}_0
\,=\,
\rho_0\wedge\overline{\kappa}_0',
\]
it comes by replacement 
in~({\ref{2-loop-structure-0-coframe}}):
\[
\aligned
d\rho_0
&
\,=\,
\Paux\,
\rho_0\wedge\kappa_0'
-
\frac{\mathcal{L}_1(\kaux)}{\overline{\mathcal{L}}_1(\kaux)}\,
\rho_0\wedge\zeta_0'
+
\overline{\Paux}\,
\rho_0\wedge\overline{\kappa}_0
-
\frac{\overline{\mathcal{L}}_1(\overline{\kaux})}{
\mathcal{L}_1(\overline{\kaux})}\,
\rho_0\wedge\overline{\zeta}_0'
\,+
\\
&
\ \ \ \ \
+\,
\isqrt
\left(
\kappa_0'
-
\frac{\isqrt}{3}\,
\bigg(
\frac{\overline{\mathcal{L}}_1\big(
\overline{\mathcal{L}}_1(\kaux)\big)}{\overline{\mathcal{L}}_1(\kaux)}
-
\overline{\Paux}
\bigg)\,
\rho_0
\right)
\wedge
\left(
\overline{\kappa}_0'
+
\frac{\isqrt}{3}\,
\bigg(
\frac{\mathcal{L}_1\big(\mathcal{L}_1(\overline{\kaux})\big)}{
\mathcal{L}_1(\overline{\kaux})}
\bigg)\,
\rho_0
\right),
\endaligned
\]
and a plain expansion yields the stated expression of $d\rho_0$.

\smallskip

Next, again from~({\ref{def-kappa-prime-0}}), it comes
by replacement in~({\ref{2-loop-structure-0-coframe}}):
\[
\aligned
d\zeta_0'
&
\,=\,
\frac{\mathcal{T}\big(\overline{\mathcal{L}}_1(\kaux)\big)}{
\overline{\mathcal{L}}_1(\kaux)}\,
\rho_0\wedge\zeta_0'
+
\frac{\mathcal{L}_1\big(\overline{\mathcal{L}}_1(\kaux)\big)}{
\overline{\mathcal{L}}_1(\kaux)}\,
\left(
\kappa_0'
-
\frac{\isqrt}{3}\,
\bigg(
\frac{\overline{\mathcal{L}}_1\big(
\overline{\mathcal{L}}_1(\kaux)\big)}{
\overline{\mathcal{L}}_1(\kaux)}
-
\overline{\Paux}
\bigg)\,\rho_0
\right)
\wedge\zeta_0'
\,-
\\
&
\ \ \ \ \
-\,
\frac{\overline{\mathcal{L}}_1\big(\overline{\mathcal{L}}_1
(\kaux)\big)}{\overline{\mathcal{L}}_1(\kaux)}\,
\zeta_0'
\wedge
\left(
\overline{\kappa}_0'
+
\frac{\isqrt}{3}\,
\bigg(
\frac{\mathcal{L}_1\big(\mathcal{L}_1(\overline{\kaux})\big)}{
\mathcal{L}_1(\overline{\kaux})}-
\Paux
\bigg)\,\rho_0
\right)
+
\frac{\overline{\mathcal{L}}_1(\overline{\kaux})}{
\mathcal{L}_1(\overline{\kaux})}\,
\zeta_0'\wedge\overline{\zeta}_0',
\endaligned
\]
and visually\,\,---\,\,no pen needed\,\,---, we obtain the stated 
result.

\smallskip

To treat $d\kappa_0'$, we start from:
\[
\kappa_0'
\,=\,
\kappa_0
+
\frac{\isqrt}{3}\,
\Baux_0\,\rho_0,
\]
and we exterior differentiate:
\leqnomode\usetagform{default}
\begin{align}
\label{d-kappa-prime-0}
d\kappa_0'
\,=\,
d\kappa_0
+
\frac{\isqrt}{3}\,
d\Baux_0
\wedge
\rho_0
+
\frac{\isqrt}{3}\,
\Baux_0\,d\rho_0.
\end{align}
As a preliminary, we need to know $d\Baux_0$. Let us recall that:
\[
\Baux_0
\,=\,
\frac{\overline{\mathcal{L}}_1\big(
\overline{\mathcal{L}}_1(\kaux)\big)}{
\overline{\mathcal{L}}_1(\kaux)}
-
\overline{\Paux}
\ \ \ \ \ \ \ \ \ \ \ \ \ \ \ \ \ \
\text{whence}
\ \ \ \ \ \ \ \ \ \ \ \ \ \ \ \ \ \
\overline{\Baux}_0
\,=\,
\frac{\mathcal{L}_1\big(\mathcal{L}_1(\overline{\kaux})\big)}{
\mathcal{L}_1(\overline{\kaux})}
-
\Paux.
\]
A plain application of Lemma~{\ref{Lemma-d-G-0}} provides 
this exterior differential:
\[
\aligned
d\bigg(
\frac{\overline{\mathcal{L}}_1\big(
\overline{\mathcal{L}}_1(\kaux)\big)}{
\overline{\mathcal{L}}_1(\kaux)}
-
\overline{\Paux}
\bigg)
&
\,=\,
\left(
\frac{\mathcal{T}\big(\overline{\mathcal{L}}_1\big(
\overline{\mathcal{L}}_1(\kaux)\big)\big)}{
\overline{\mathcal{L}}_1(\kaux)}
-
\frac{\mathcal{T}\big(\overline{\mathcal{L}}_1(\kaux)\big)\,\,
\overline{\mathcal{L}}_1\big(
\overline{\mathcal{L}}_1(\kaux)\big)}{
\overline{\mathcal{L}}_1(\kaux)^2}
-
\mathcal{T}\big(\overline{\Paux}\big)
\right)
\rho_0
\,+
\\
&
\ \ \ \ \
+
\left(
\frac{\mathcal{L}_1\big(\overline{\mathcal{L}}_1\big(
\overline{\mathcal{L}}_1(\kaux)\big)\big)}{
\overline{\mathcal{L}}_1(\kaux)}
-
\frac{\mathcal{L}_1\big(\overline{\mathcal{L}}_1(\kaux)\big)\,\,
\overline{\mathcal{L}}_1\big(\overline{\mathcal{L}}_1(\kaux)\big)}{
\overline{\mathcal{L}}_1(\kaux)^2}
-
\mathcal{L}_1\big(\overline{\Paux}\big)
\right)
\kappa_0
\,+
\\
&
\ \ \ \ \
+
\left(
\frac{\mathcal{K}\big(\overline{\mathcal{L}}_1\big(
\overline{\mathcal{L}}_1(\kaux)\big)\big)}{
\overline{\mathcal{L}}_1(\kaux)}
-
\frac{\mathcal{K}\big(\overline{\mathcal{L}}_1(\kaux)\big)\,\,
\overline{\mathcal{L}}_1\big(\overline{\mathcal{L}}_1(\kaux)\big)}{
\overline{\mathcal{L}}_1(\kaux)^2}
-
\mathcal{K}\big(\overline{\Paux}\big)
\right)
\zeta_0
\,+
\\
&
\ \ \ \ \
+
\left(
\frac{\overline{\mathcal{L}}_1\big(\overline{\mathcal{L}}_1\big(
\overline{\mathcal{L}}_1(\kaux)\big)\big)}{
\overline{\mathcal{L}}_1(\kaux)}
-
\frac{\overline{\mathcal{L}}_1\big(
\overline{\mathcal{L}}_1(\kaux)\big)^2}{
\overline{\mathcal{L}}_1(\kaux)^2}
-
\overline{\mathcal{L}}_1\big(\overline{\Paux}\big)
\right)
\overline{\kappa}_0
\\
&
\ \ \ \ \
+
\left(
\frac{\overline{\mathcal{K}}\big(\overline{\mathcal{L}}_1\big(
\overline{\mathcal{L}}_1(\kaux)\big)\big)}{
\overline{\mathcal{L}}_1(\kaux)}
-
\frac{\overline{\mathcal{K}}\big(
\overline{\mathcal{L}}_1(\kaux)\big)\,\,
\overline{\mathcal{L}}_1\big(\overline{\mathcal{L}}_1(\kaux)\big)}{
\overline{\mathcal{L}}_1(\kaux)^2}
-
\overline{\mathcal{K}}\big(\overline{\Paux}\big)
\right)
\overline{\zeta}_0,
\endaligned
\]
an expression that we will abbreviate as:
\[
d\Baux_0
\,=\,
\Uaux_0\,\rho_0
+
\Vaux_0\,\kappa_0
+
\Waux_0\,\zeta_0
+
\Xaux_0\,\overline{\kappa}_0
+
\Yaux_0\,\overline{\zeta}_0.
\]

\begin{Assertion}
After simplifications:
\[
\Yaux_0
\,=\,
-\,
\frac{\overline{\mathcal{L}}_1(\overline{\kaux})\,\,
\overline{\mathcal{L}}_1\big(\overline{\mathcal{L}}_1(\kaux)\big)}{
\overline{\mathcal{L}}_1(\kaux)}
+
\overline{\mathcal{L}}_1(\overline{\kaux})\,
\overline{\Paux}.
\]
\end{Assertion}

\proof
In the first two terms of $\Yaux_0$, we replace 
from Lemma~{\ref{Lemma-K-bar-penetrates-L1-bar-k}}:
\[
\aligned
\overline{\mathcal{K}}
\big(\overline{\mathcal{L}}_1
\big(\overline{\mathcal{L}}_1(\kaux)\big)\big)
&
\,=\,
-\,2\,
\overline{\mathcal{L}}_1
\big(\overline{\kaux}\big)\,\,
\overline{\mathcal{L}}_1\big(
\overline{\mathcal{L}}_1(\kaux)\big)
-
\overline{\mathcal{L}}_1
\big(\overline{\mathcal{L}}_1
\big(\overline{\kaux}\big)\big)\,\,
\overline{\mathcal{L}}_1(\kaux),
\\
\overline{\mathcal{K}}
\big(\overline{\mathcal{L}}_1(\kaux)\big)
&
\,=\,
-\,
\overline{\mathcal{L}}_1\big(\overline{\kaux}\big)\,
\overline{\mathcal{L}}_1(\kaux),
\endaligned
\]
and in the third term of $\Yaux_0$, we replace from
Lemma~{\ref{Lemma-K-bar-k-K-bar-P}}:
\[
\overline{\mathcal{K}}\big(\overline{\Paux}\big)
\,=\,
-\,\overline{\Paux}\,
\overline{\mathcal{L}}_1\big(\overline{\kaux}\big)
-
\overline{\mathcal{L}}_1\big(
\overline{\mathcal{L}}_1\big(\overline{\kaux}\big)\big),
\]
which yields the result after one (underlined) pair cancellation:
\begin{align}
\Yaux_0
&
\,=\,
-\,
\frac{2\,\overline{\mathcal{L}}_1(\overline{\kaux})\,\,
\overline{\mathcal{L}}_1\big(\overline{\mathcal{L}}_1(\kaux)\big)}{
\overline{\mathcal{L}}_1(\kaux)}
-
\zero{
\overline{\mathcal{L}}_1\big(\overline{\mathcal{L}}_1
(\overline{\kaux})\big)}
+
\frac{\overline{\mathcal{L}}_1(\overline{\kaux})\,\,
\overline{\mathcal{L}}_1\big(\overline{\mathcal{L}}_1(\kaux)\big)}{
\overline{\mathcal{L}}_1(\kaux)}
\,+
\notag
\\
&
\ \ \ \ \ \ \ \ \ \ \ \ \ \ \ \ \ \ \ \ \ \ \ \ \ \
\ \ \ \ \ \ \ \ \ \ \ \ \ \ \ \ \ \ \ \ \ \ \ \ \ \
+
\overline{\Paux}\,
\overline{\mathcal{L}}_1\big(\overline{\kaux}\big)
+
\zero{
\overline{\mathcal{L}}_1\big(\overline{\mathcal{L}}_1
(\overline{\kaux})\big)}.
\qedhere
\end{align}
\endproof

Temporarily, let us work with the abbreviations $\Uaux_0$, $\Vaux_0$,
$\Waux_0$, $\Xaux_0$, $\Yaux_0$.  So, using the previous structure
formulas~({\ref{2-loop-structure-0-coframe}}) in which, directly we
replace:
\[
\zeta_0
\,=\,
\frac{\zeta_0'}{\overline{\mathcal{L}}_1(\kaux)},
\]
let us add line-by-line all three terms
of~({\ref{d-kappa-prime-0}}):
\[
\!\!\!\!\!\!\!\!\!\!\!\!\!\!\!
\small
\aligned
d\kappa_0'
&
\,=\,
-\,
\frac{\mathcal{T}(\kaux)}{\overline{\mathcal{L}}_1(\kaux)}\,
\rho_0\wedge\zeta_0'
-
\frac{\mathcal{L}_1(\kaux)}{\overline{\mathcal{L}}_1(\kaux)}\,
\kappa_0\wedge\zeta_0'
+
\zeta_0'\wedge\overline{\kappa}_0'
\,+
\\
&
\ \ \ \ \
+
\frac{\isqrt}{3}\,
\zero{
\Uaux_0\,\rho_0\wedge\rho_0}
+
\frac{\isqrt}{3}\,
\Vaux_0\,
\kappa_0\wedge\rho_0
+
\frac{\isqrt}{3}\,
\Waux_0\,
\frac{\zeta_0'}{\overline{\mathcal{L}}_1(\kaux)}
\wedge\rho_0
+
\frac{\isqrt}{3}\,
\Xaux_0\,
\overline{\kappa}_0\wedge\rho_0
+
\frac{\isqrt}{3}\,
\Yaux_0\,
\frac{\overline{\zeta}_0'}{\mathcal{L}_1(\overline{\kaux})}\,
\wedge\rho_0
\,+
\\
&
\ \ \ \ \ 
+
\frac{\isqrt}{3}\,
\Baux_0\,\Paux\,
\rho_0\wedge\kappa_0
-
\frac{\isqrt}{3}\,
\Baux_0\,
\frac{\mathcal{L}_1(\kaux)}{\overline{\mathcal{L}}_1(\kaux)}\,
\rho_0\wedge
\zeta_0'
+
\frac{\isqrt}{3}\,
\Baux_0\,\overline{\Paux}\,
\rho_0\wedge\overline{\kappa}_0
-
\frac{\isqrt}{3}\,
\Baux_0\,
\frac{\overline{\mathcal{L}}_1(\overline{\kaux})}{
\mathcal{L}_1(\overline{\kaux})}\,
\rho_0\wedge\overline{\zeta}_0'
-
\frac{1}{3}\,
\Baux_0\,
\kappa_0\wedge\overline{\kappa}_0,
\endaligned
\]
hence after collecting coefficients of basic $2$-forms, we get:
\[
\small
\aligned
d\kappa_0'
&
\,=\,
\rho_0\wedge\zeta_0'
\left[
-\,\frac{\mathcal{T}(\kaux)}{\overline{\mathcal{L}}_1(\kaux)}
-
\frac{\isqrt}{3}\,
\frac{\Waux_0}{\overline{\mathcal{L}}_1(\kaux)}
-
\frac{\isqrt}{3}\,
\Baux_0\,
\frac{\mathcal{L}_1(\kaux)}{\overline{\mathcal{L}}_1(\kaux)}
\right]
+
\rho_0\wedge\kappa_0
\left[
-\,\frac{\isqrt}{3}\,
\Vaux_0
+
\frac{\isqrt}{3}\,
\Baux_0\,\Paux
\right]
\,+
\\
&
\ \ \ \ \
+
\rho_0\wedge\overline{\kappa}_0
\left[
-\,\frac{\isqrt}{3}\,\Xaux_0
+
\frac{\isqrt}{3}\,\Baux_0\,\overline{\Paux}
\right]
+
\rho_0\wedge\overline{\zeta}_0'
\left[
-\,\frac{\isqrt}{3}\,
\frac{\Yaux_0}{\mathcal{L}_1(\overline{\kaux})}
-
\frac{\isqrt}{3}\,\Baux_0\,
\frac{\overline{\mathcal{L}}_1(\overline{\kaux})}{
\mathcal{L}_1(\overline{\kaux})}
\right]
\,+
\\
&
\ \ \ \ \
+
\kappa_0\wedge\zeta_0'
\left[
-\,\frac{\mathcal{L}_1(\kaux)}{\overline{\mathcal{L}}_1(\kaux)}
\right]
+
\kappa_0\wedge\overline{\kappa}_0
\left[
-\,\frac{1}{3}\,\Baux_0
\right]
+
\zeta_0'\wedge\overline{\kappa}_0.
\endaligned
\]
Next, replace everywhere:
\[
\kappa_0
\,=\,
\kappa_0'
-
\frac{\isqrt}{3}\,
\Baux_0\,\rho_0.
\] 
Then using again $\kappa_0 \wedge \rho_0 = \kappa_0' \wedge \rho_0$,
only the last line changes, as it becomes:
\[
\small
\aligned
\!\!\!\!\!\!\!\!\!\!\!\!\!\!\!
\Big(
\kappa_0'
-
\frac{\isqrt}{3}\,
\Baux_0\,\rho_0
\Big)
\wedge
\zeta_0'
\left[
-\,\frac{\mathcal{L}_1(\kaux)}{\overline{\mathcal{L}}_1(\kaux)}
\right]
+
\Big(
\kappa_0'
-
\frac{\isqrt}{3}\,\Baux_0\,\rho_0
\Big)
\wedge
\Big(
\overline{\kappa}_0'
+
\frac{\isqrt}{3}\,
\overline{\Baux}_0\,
\rho_0
\Big)\,
\left[
-\,\frac{1}{3}\,\Baux_0
\right]
+
\zeta_0'
\wedge
\Big(
\overline{\kappa}_0'
+
\frac{\isqrt}{3}\,
\overline{\Baux}_0\,\rho_0
\Big).
\endaligned
\]
Expanding and collecting visually\,\,---\,\,no pen needed\,\,---,
we get:
\[
\aligned
d\kappa_0'
&
\,=\,
\rho_0\wedge\zeta_0'
\left[
-\,\frac{\mathcal{T}(\kaux)}{\overline{\mathcal{L}}_1(\kaux)}
-
\frac{\isqrt}{3}\,
\frac{\Waux_0}{\overline{\mathcal{L}}_1(\kaux)}
-
\zero{
\frac{\isqrt}{3}\,
\Baux_0\,
\frac{\mathcal{L}_1(\kaux)}{\overline{\mathcal{L}}_1(\kaux)}}
+
\zero{
\frac{\isqrt}{3}\,
\Baux_0\,
\frac{\mathcal{L}_1(\kaux)}{\overline{\mathcal{L}}_1(\kaux)}}
-
\frac{\isqrt}{3}\,
\overline{\Baux}_0
\right]
\,+
\\
&
\ \ \ \ \
+
\rho_0\wedge\kappa_0'
\left[
-\,\frac{\isqrt}{3}\,
\Vaux_0
+
\frac{\isqrt}{3}\,\Baux_0\,
\Paux
+
\frac{\isqrt}{9}\,\Baux_0\,\overline{\Baux}_0
\right]
\,+
\\
&
\ \ \ \ \
+
\rho_0\wedge\overline{\kappa}_0'
\left[
-\,\frac{\isqrt}{3}\,\Xaux_0
+
\frac{\isqrt}{3}\,\Baux_0\,\overline{\Paux}
+
\frac{\isqrt}{9}\,
\Baux_0\,\Baux_0
\right]
\,+
\\
&
\ \ \ \ \
+
\rho_0\wedge\overline{\zeta}_0'
\left[
\zero{
-\,\frac{\isqrt}{3}\,
\frac{\Yaux_0}{\mathcal{L}_1(\overline{\kaux})}
-
\frac{\isqrt}{3}\,
\Baux_0\,
\frac{\overline{\mathcal{L}}_1(\overline{\kaux})}{
\mathcal{L}_1(\overline{\kaux})}}
\right]
\,+
\\
&
\ \ \ \ \
\kappa_0'\wedge\zeta_0'
\left[
-\,
\frac{\mathcal{L}_1(\kaux)}{\overline{\mathcal{L}}_1(\kaux)}
\right]
+
\kappa_0'\wedge\overline{\kappa}_0'
\left[
-\,\frac{1}{3}\,\Baux_0
\right]
+
\zeta_0'\wedge\overline{\kappa}_0'.
\endaligned
\]

To finish, we must yet
replace $\Vaux_0$, $\Waux_0$, $\Xaux_0$, $\Yaux_0$
by their complete values, and we will realize, as
indicated by anticipation above, that the
coefficient of $\rho_0 \wedge \overline{\zeta}_0'$ vanishes
identically. 

\smallskip

Firstly, a replacement followed by a visual expansion finalizes:
\[
\aligned
\!\!\!\!\!\!\!\!\!\!
\big[
\rho_0\wedge\kappa_0'
\big]
\big\{
d\kappa_0'
\big\}
&
\,=\,
-\,
\frac{\isqrt}{3}\,
\frac{\mathcal{L}_1\big(
\overline{\mathcal{L}}_1\big(
\overline{\mathcal{L}}_1(\kaux)\big)\big)}{
\overline{\mathcal{L}}_1(\kaux)}
+
\frac{\isqrt}{3}\,
\frac{\mathcal{L}_1\big(\overline{\mathcal{L}}_1(\kaux)\big)\,\,
\overline{\mathcal{L}}_1\big(\overline{\mathcal{L}}_1(\kaux)\big)}{
\overline{\mathcal{L}}_1(\kaux)^2}
+
\frac{\isqrt}{3}\,
\mathcal{L}_1\big(\overline{\Paux}\big)
\,+
\\
&
\ \ \ \ \
+
\frac{\isqrt}{3}\,
\frac{\overline{\mathcal{L}}_1\big(
\overline{\mathcal{L}}_1(\kaux)\big)}{
\overline{\mathcal{L}}_1}\,
\Paux
-
\frac{\isqrt}{3}\,
\Paux\,\overline{\Paux}
+
\frac{\isqrt}{9}\,
\left(
\frac{\overline{\mathcal{L}}_1\big(
\overline{\mathcal{L}}_1(\kaux)\big)}{
\overline{\mathcal{L}}_1(\kaux)}
-
\overline{\Paux}
\right)
\left(
\frac{\mathcal{L}_1\big(\mathcal{L}_1(\overline{\kaux})\big)}{
\mathcal{L}_1(\overline{\overline{\kaux}})}
-
\Paux
\right).
\endaligned
\]

Secondly:
\[
\aligned
\big[
\rho_0\wedge\zeta_0'
\big]
\big\{
d\kappa_0'
\big\}
&
\,=\,
-\,
\frac{\mathcal{T}(\kaux)}{\overline{\mathcal{L}}_1(\kaux)}
-
\frac{\isqrt}{3}\,
\frac{\mathcal{K}\big(\overline{\mathcal{L}}_1\big(
\overline{\mathcal{L}}_1(\kaux)\big)\big)}{
\overline{\mathcal{L}}_1(\kaux)^2}
+
\frac{\isqrt}{3}\,
\frac{\mathcal{K}\big(\overline{\mathcal{L}}_1(\kaux)\big)\,\,
\overline{\mathcal{L}}_1\big(\overline{\mathcal{L}}_1
(\kaux)\big)}{
\overline{\mathcal{L}}_1(\kaux)^3}
\,+
\\
&
\ \ \ \ \ \ \ \ \ \ \ \ \ \ \
+
\frac{\isqrt}{3}\,
\frac{\boxed{\mathcal{K}(\overline{\Paux})}}{
\overline{\mathcal{L}}_1(\kaux)}
-
\frac{\isqrt}{3}\,
\frac{\mathcal{L}_1\big(\mathcal{L}_1(\overline{\kaux})\big)}{
\mathcal{L}_1(\overline{\kaux})}
+
\frac{\isqrt}{3}\,\Paux,
\endaligned
\]
but here, we must still replace the boxed term using 
Lemma~{\ref{Lemma-K-bar-k-K-bar-P}}:
\[
\aligned
\big[
\rho_0\wedge\zeta_0'
\big]
\big\{
d\kappa_0'
\big\}
&
\,=\,
-\,
\frac{\mathcal{T}(\kaux)}{\overline{\mathcal{L}}_1(\kaux)}
-
\frac{\isqrt}{3}\,
\frac{\mathcal{K}\big(\overline{\mathcal{L}}_1\big(
\overline{\mathcal{L}}_1(\kaux)\big)\big)}{
\overline{\mathcal{L}}_1(\kaux)^2}
+
\frac{\isqrt}{3}\,
\frac{\mathcal{K}\big(\overline{\mathcal{L}}_1(\kaux)\big)\,\,
\overline{\mathcal{L}}_1\big(\overline{\mathcal{L}}_1
(\kaux)\big)}{
\overline{\mathcal{L}}_1(\kaux)^3}
\,+
\\
&
\ \ \ \ \
-
\zero{\frac{\isqrt}{3}\,\Paux}
-
\frac{\isqrt}{3}\,
\frac{\overline{\mathcal{L}}_1\big(\mathcal{L}_1(\kaux)\big)}{
\overline{\mathcal{L}}_1(\kaux)}
+
\frac{1}{3}\,
\frac{\mathcal{T}(\kaux)}{\overline{\mathcal{L}}_1(\kaux)}
-
\frac{\isqrt}{3}\,
\frac{\mathcal{L}_1\big(\mathcal{L}_1(\overline{\kaux})\big)}{
\mathcal{L}_1(\overline{\kaux})}
+
\zero{\frac{\isqrt}{3}\,\Paux}.
\endaligned
\]
A pair cancellation 
makes the obtained expression match precisely with what 
Proposition~{\ref{Proposition-3-loop-structure-0-coframe}} stated,
after some permutation of terms.

\smallskip

The third replacement conducts directly to the stated result:
\[
\aligned
\big[
\rho_0\wedge\overline{\kappa}_0'
\big]
\big\{
d\kappa_0'
\big\}
&
\,=\,
-\,\frac{\isqrt}{3}\,
\frac{\overline{\mathcal{L}}_1\big(
\overline{\mathcal{L}}_1\big(
\overline{\mathcal{L}}_1(\kaux)\big)\big)}{
\overline{\mathcal{L}}_1(\kaux)}
+
\frac{\isqrt}{3}\,
\frac{\overline{\mathcal{L}}_1\big(
\overline{\mathcal{L}}_1(\kaux)\big)^2}{
\overline{\mathcal{L}}_1(\kaux)^2}
+
\frac{\isqrt}{3}\,
\overline{\mathcal{L}}_1\big(\overline{\Paux}\big)
\,+
\\
&
\ \ \ \ \ 
+
\frac{\isqrt}{3}\,
\frac{\overline{\mathcal{L}}_1\big(
\overline{\mathcal{L}}_1(\kaux)\big)}{
\overline{\mathcal{L}}_1(\kaux)}\,
\overline{\Paux}
-
\frac{\isqrt}{3}\,
\overline{\Paux}\,\overline{\Paux}
\,+
\\
&
\ \ \ \ \
+
\frac{\isqrt}{9}\,
\frac{\overline{\mathcal{L}}_1\big(
\overline{\mathcal{L}}_1(\kaux)\big)^2}{
\overline{\mathcal{L}}_1(\kaux)^2}
-
\frac{2\,\isqrt}{9}\,
\frac{\overline{\mathcal{L}}_1\big(
\overline{\mathcal{L}}_1(\kaux)\big)}{
\overline{\mathcal{L}}_1(\kaux)}\,
\overline{\Paux}
+
\frac{\isqrt}{9}\,
\overline{\Paux}\,\overline{\Paux},
\endaligned
\]
while the fourth (last) brings an identically zero result:
\begin{align}
\big[
\rho_0\wedge\overline{\zeta}_0'
\big]
\big\{
d\kappa_0'
\big\}
&
\,=\,
\zero{
\frac{\isqrt}{3}\,
\frac{\overline{\mathcal{L}}_1(\overline{\kaux})\,\,
\overline{\mathcal{L}}_1\big(\overline{\mathcal{L}}_1(\kaux)\big)}{
\mathcal{L}_1(\overline{\kaux})\,\,
\overline{\mathcal{L}}_1(\kaux)}}
-
\zerozero{
\frac{\isqrt}{3}\,
\frac{\overline{\mathcal{L}}_1(\overline{\kaux})}{
\mathcal{L}_1(\overline{\kaux})}\,
\overline{\Paux}}
-
\zero{
\frac{\isqrt}{3}\,
\frac{
\overline{\mathcal{L}}_1\big(
\overline{\mathcal{L}}_1(\kaux)\big)}{
\overline{\mathcal{L}}_1(\kaux)}\,
\frac{\overline{\mathcal{L}}_1(\overline{\kaux})}{
\mathcal{L}_1(\overline{\kaux})}
}
\,+
\\
&
\ \ \ \ \ 
+
\zerozero{
\frac{\isqrt}{3}\,
\overline{\Paux}\,
\frac{\overline{\mathcal{L}}_1(\overline{\kaux})}{
\mathcal{L}_1(\overline{\kaux})}}.
\qedhere
\end{align}
\endproof

\Section{\bf Third Loop: Reduction of the Group Parameter ${\sf d}$}
\label{third-loop-reduction-d}
\HEAD{{\ref{third-loop-reduction-d}}.~{\sf Third Loop: 
Reduction of the Group Parameter ${\sf d}$}
}{
Wei Guo {\sc Foo} (Beijing) and Joël {\sc Merker} (Orsay)}

After normalization of the group parameter ${\sf b}$
from~({\ref{normalization-b}}),
we have a new reduced group $G^6$ of real dimension $6$,
and the lifted coframe is:
\leqnomode\usetagform{default}
\begin{align}
\label{3-loop-initial}
\left(\!
\begin{array}{c}
\rho
\\
\kappa
\\
\zeta
\end{array}
\!\right)
\,:=\,
\left(\!
\begin{array}{ccc}
{\sf c}\overline{\sf c} & 0 & 0
\\
-\isqrt\,\overline{\sf c}{\sf e} & {\sf c} & 0
\\
{\sf d} & {\sf e} & \frac{{\sf c}}{\overline{\sf c}}
\end{array}
\!\right)
\left(\!
\begin{array}{c}
\rho_0
\\
\kappa_0'
\\
\zeta_0'
\end{array}
\!\right)
\ \ \ \ \ \ \ \ \ \ \ \ \ \ \ \ \ \
\Longleftrightarrow
\ \ \ \ \ \ \ \ \ \ \ \ \ \ \ \ \ \
\left\{
\aligned
\rho
&
\,:=\,
{\sf c}\overline{\sf c}\,
\rho_0,
\\
\kappa
&
\,:=\,
-\,\isqrt\,\overline{\sf c}{\sf e}\,
\rho_0
+
{\sf c}\,\kappa_0',
\\
\zeta
&
\,:=\,
{\sf d}\,\rho_0
+
{\sf e}\,\kappa_0'
+
\frac{{\sf c}}{\overline{\sf c}}\,
\zeta_0',
\endaligned\right.
\end{align}
with inverse formulas:
\leqnomode\usetagform{default}
\begin{align}
\label{3-loop-initial-inverted}
\rho_0
&
\,=\,
\frac{1}{{\sf c}\overline{\sf c}}\,
\rho,
\notag
\\
\kappa_0'
&
\,=\,
\isqrt\,\frac{{\sf e}}{{\sf c}{\sf c}}\,
\rho
+
\frac{1}{{\sf c}}\,
\kappa,
\\
\zeta_0'
&
\,=\,
\Big(
-\,\isqrt\,
\frac{\overline{\sf c}{\sf e}{\sf e}}{
{\sf c}{\sf c}{\sf c}}
-
\frac{{\sf d}}{{\sf c}{\sf c}}
\Big)\,
\rho
-
\frac{\overline{\sf c}{\sf e}}{{\sf c}{\sf c}}\,
\kappa
+
\frac{\overline{\sf c}}{\sf c}\,
\zeta.
\notag
\end{align}
The Maurer-Cartan matrix becomes:
\[
\aligned
dg\cdot g^{-1}
&
\,=\,
\left(\!
\begin{array}{ccc}
\overline{\sf c}\,d{\sf c}+{\sf c}d\overline{\sf c} & 0 & 0
\\
-\isqrt\,{\sf e}d\overline{\sf c}
-\isqrt\,\overline{\sf c}d{\sf e} & d{\sf c} & 0
\\
d{\sf d} & d{\sf e} & \frac{d{\sf c}}{\overline{\sf c}}
- \frac{{\sf c}\,d\overline{\sf c}}{\overline{\sf c}\overline{\sf c}}
\end{array}
\!\right)
\left(\!
\begin{array}{ccc}
\frac{1}{{\sf c}\overline{\sf c}} & 0 & 0
\\
\isqrt\,\frac{{\sf e}}{{\sf c}{\sf c}} & \frac{1}{{\sf c}} & 0
\\
-\isqrt\,\frac{\overline{\sf c}{\sf e}{\sf e}}{
{\sf c}{\sf c}{\sf c}}-\frac{{\sf d}}{{\sf c}{\sf c}}
& -\frac{\overline{\sf c}{\sf e}}{{\sf c}{\sf c}} & 
\frac{\overline{\sf c}}{{\sf c}}
\end{array}
\!\right)
\\
&
\,=:\,
\left(\!
\begin{array}{ccc}
\alpha+\overline{\alpha} & 0 & 0
\\
\beta & \alpha & 0
\\
\gamma & \isqrt\,\beta & \alpha-\overline{\alpha}
\end{array}
\!\right),
\endaligned
\]
in terms of the group-invariant $1$-forms:
\[
\aligned
\alpha
&
\,:=\,
\frac{d{\sf c}}{{\sf c}},
\\
\beta
&
\,:=\,
\isqrt\,
\frac{{\sf e}\,d{\sf c}}{{\sf c}{\sf c}}
-
\isqrt\,\frac{{\sf e}\,d\overline{\sf c}}{{\sf c}\overline{\sf c}}
-
\isqrt\,\frac{d{\sf e}}{{\sf c}},
\\
\gamma
&
\,:=\,
\Big(
\frac{{\sf c}{\sf d}+\isqrt\,\overline{\sf c}{\sf e}{\sf e}}{
{\sf c}{\sf c}\overline{\sf c}}
\Big)\,
\bigg(
-\frac{d{\sf c}}{\sf c}
+
\frac{d\overline{\sf c}}{\overline{\sf c}}
\bigg)
+
\frac{d{\sf d}}{{\sf c}\overline{\sf c}}
+
\isqrt\,\frac{{\sf e}\,d{\sf e}}{{\sf c}{\sf c}}.
\endaligned
\]

Now, if we exterior-differentiate the lifted coframe on the product
manifold equipped with coordinates:
\[
\big( 
z_1,z_2,\overline{z}_1,\overline{z}_2,v
\big)
\times
\big(
{\sf c},\overline{\sf c},
{\sf d},\overline{\sf d},{\sf e},\overline{\sf e}
\big)
\,\,\in\,\,
M^5
\times
G^6,
\]
after some computations, we may come to structure equations of the
abstract shape:
\[
\aligned
d\rho
&
\,=\,
\big(\alpha+\overline{\alpha}\big)
\wedge
\rho
\,+
\\
&
\ \ \ \ \ 
+
R^1\,
\rho\wedge\kappa
+
R^2\,
\rho\wedge\zeta
+
\overline{R}^1\,
\rho\wedge\overline{\kappa}
+
\overline{R}^2\,
\rho\wedge\overline{\zeta}
+
\isqrt\,
\kappa\wedge\overline{\kappa},
\\
d\kappa
&
\,=\,
\beta\wedge\rho
+
\alpha\wedge\kappa
\,+
\\
&
\ \ \ \ \ 
+
K^1\,
\rho\wedge\kappa
+
K^2\,\rho\wedge\zeta
+
\boxed{K^3}\,
\rho\wedge\overline{\kappa}
+
K^4\,
\rho\wedge\overline{\zeta}
\,+
\\
&
\ \ \ \ \
+
K^5\,
\kappa\wedge\zeta
+
K^6\,
\kappa\wedge\overline{\kappa}
+
\zeta\wedge\overline{\kappa},
\\
d\zeta
&
\,=\,
\gamma\wedge\rho
+
\isqrt\,\beta\wedge\kappa
+
\big(
\alpha-\overline{\alpha}
\big)
\wedge\zeta
\,+
\\
&
\ \ \ \ \
+
Z^1\,
\rho\wedge\kappa
+
Z^2\,
\rho\wedge\zeta
+
Z^3\,
\rho\wedge\overline{\kappa}
+
Z^4\,
\rho\wedge\overline{\zeta}
\,+
\\
&
\ \ \ \ \
+
Z^5\,
\kappa\wedge\zeta
+
\boxed{Z^6}\,
\kappa\wedge\overline{\kappa}
+
Z^7\,
\kappa\wedge\overline{\zeta}
+
Z^8\,
\zeta\wedge\overline{\kappa}
+
Z^9\,
\zeta\wedge\overline{\zeta}.
\endaligned
\]

Before really computing explicitly some of these torsion
coefficients, let us examine what are the absorption equations.
For this, we replace:
\[
\aligned
\alpha
&
\,=:\,
\alpha'
+
a_1\,\rho
+
a_2\,\kappa
+
a_3\,\zeta
+
a_4\,\overline{\kappa}
+
a_5\,\overline{\zeta},
\\
\beta
&
\,=:\,
\beta'
+
b_1\,\rho
+
b_2\,\kappa
+
b_3\,\zeta
+
b_4\,\overline{\kappa}
+
b_5\,\overline{\zeta},
\\
\gamma
&
\,=:\,
\gamma'
+
c_1\,\rho
+
c_2\,\kappa
+
c_3\,\zeta
+
c_4\,\overline{\kappa}
+
c_5\,\overline{\zeta}.
\endaligned
\]
A moment of reflection convinces that the result for
$d\rho$ is the same as in the proof of
Lemma~{\ref{Lemma-R-1-2-K-6-Z-8}}:
\[
\aligned
d\rho
&
\,=\,
\big(\alpha'+\overline{\alpha}'\big)
\wedge
\rho
\,+
\\
&
\ \ \ \ \
+
\rho\wedge\kappa
\Big(
R^1
-
a_2
-
\overline{a}_4
\Big)
+
\rho\wedge\zeta
\Big(
R^2
-
a_3
-
\overline{a}_5
\Big)
+
\rho\wedge\overline{\kappa}
\Big(
\overline{R}^1
-
a_4
-
\overline{a}_2
\Big)
\,+
\\
&
\ \ \ \ \ 
+
\rho\wedge\overline{\zeta}
\Big(
\overline{R}^2
-
a_5
-
\overline{a}_3
\Big)
+
\isqrt\,
\kappa\wedge\overline{\kappa}.
\endaligned
\]
Similarly, $d\kappa$ is unchanged:
\[
\aligned
d\kappa
&
\,=\,
\beta'\wedge\rho
+
\alpha'\wedge\kappa
\,+
\\
&
\ \ \ \ \ 
+
\rho\wedge\kappa\,
\Big(
K^1
+
a_1
-
b_2
\Big)
+
\rho\wedge\zeta\,
\Big(
K^2
-
b_3
\Big)
+
\rho\wedge\overline{\kappa}\,
\Big(
K^3
-
b_4
\Big)
+
\rho\wedge\overline{\zeta}\,
\Big(
K^4
-
b_5
\Big)
\,+
\\
&
\ \ \ \ \
+
\kappa\wedge\zeta\,
\Big(
K^5
-
a_3
\Big)
+
\kappa\wedge\overline{\kappa}\,
\Big(
K^6
-
a_4
\Big)
+
\kappa\wedge\overline{\zeta}\,
\big(
-a_5
\big)
+
\zeta\wedge\overline{\kappa}.
\endaligned
\]

However, for $d\zeta$, we have to compute:
\[
\!\!\!\!\!\!\!\!\!\!\!\!\!\!\!
\small
\aligned
\gamma\wedge\rho
+
\isqrt\,
\beta\wedge\kappa
+
\big(\alpha-\overline{\alpha}\big)
\wedge\zeta
&
\,=\,
\gamma'\wedge\rho
+
0
+
c_2\,
\kappa\wedge\rho
+
c_3\,
\zeta\wedge\rho
+
c_4\,
\overline{\kappa}\wedge\rho
+
c_5\,
\overline{\zeta}\wedge\rho
\,+
\\
&
\ \ \ \ \
+
\isqrt\,
\beta'\wedge\kappa
+
\isqrt\,
b_1\,
\rho\wedge\kappa
+
0
+
\isqrt\,b_3\,
\zeta\wedge\kappa
+
\isqrt\,b_4\,
\overline{\kappa}\wedge\kappa
+
\isqrt\,b_5\,
\overline{\zeta}\wedge\kappa
\,+
\\
&
\ \ \ \ \
+
\alpha'\wedge\zeta
+
a_1\,
\rho\wedge\zeta
+
a_2\,
\kappa\wedge\zeta
+
0
+
a_4\,
\overline{\kappa}\wedge\zeta
+
a_5\,
\overline{\zeta}\wedge\zeta
\,-
\\
&
\ \ \ \ \ \
-\,
\overline{\alpha}'\wedge\zeta
-
\overline{a}_1\,
\rho\wedge\zeta
-
\overline{a}_2\,
\overline{\kappa}\wedge\zeta
-
\overline{a}_3\,
\overline{\zeta}\wedge\zeta
-
\overline{a}_4\,
\kappa\wedge\zeta
-
0,
\endaligned
\]
and we get:
\[
\aligned
d\zeta
&
\,=\,
\gamma'\wedge\rho
+
\isqrt\,\beta'\wedge\kappa
+
\big(\alpha'-\overline{\alpha}'\big)
\wedge\zeta
\,+
\\
&
\ \ \ \ \
+
\rho\wedge\kappa
\Big(
Z^1
+
\isqrt\,b_1
-
c_2
\Big)
+
\rho\wedge\zeta\,
\Big(
Z^2
-
c_3
+
a_1
-
\overline{a}_1
\Big)
+
\rho\wedge\overline{\kappa}\,
\Big(
Z^3
-
c_4
\Big)
+
\rho\wedge\overline{\zeta}\,
\Big(
Z^4
-
c_5
\Big)
\,+
\\
&
\ \ \ \ \
+
\kappa\wedge\zeta\,
\Big(
Z^5
-
\isqrt\,b_3
+
a_2
-
\overline{a}_4
\Big)
+
\kappa\wedge\overline{\kappa}\,
\Big(
Z^6
-
\isqrt\,b_4
\Big)
+
\kappa\wedge\overline{\zeta}\,
\Big(
Z^7
-
\isqrt\,b_5
\Big)
\,+
\\
&
\ \ \ \ \
+
\zeta\wedge\overline{\kappa}\,
\Big(
Z^8
-
a_4
+
\overline{a}_2
\Big)
+
\zeta\wedge\overline{\zeta}\,
\Big(
Z^9
-
a_5
+
\overline{a}_3
\Big).
\endaligned
\]

\begin{Lemma}
Here is an essential linear combination of torsion terms:
\[
\isqrt\,K^3
-
Z^6.
\]
\end{Lemma}

\proof
Indeed:
\[
\aligned
K^{3\prime}
&
\,=\,
K^3
-
b_4,
\notag
\\
Z^{6\prime}
&
\,=\,
Z^6
-
\isqrt\,
b_4,
\endaligned
\]
whence:
\[
\isqrt\,K^{3\prime}
-
Z^{6\prime}
\,=\,
\isqrt\,K^3
-
Z^6.
\qedhere
\]
\endproof

\begin{Proposition}
Their explicit expressions are:
\[
\aligned
K^3
&
\,=\,
-\,
\frac{{\sf d}}{{\sf c}\overline{\sf c}}
+
\frac{{\sf e}}{{\sf c}\overline{\sf c}}
\left(
-2\isqrt\,
\frac{\overline{\mathcal{L}}_1\big(
\overline{\mathcal{L}}_1(\kaux)\big)}{
\overline{\mathcal{L}}_1(\kaux)}
-
\frac{\isqrt}{3}\,
\overline{\Paux}
\right)
-
\isqrt\,
\frac{{\sf e}\overline{\sf e}}{\overline{\sf c}\overline{\sf c}}\,
\frac{\overline{\mathcal{L}}_1(\kaux)}{
\mathcal{L}_1(\overline{\kaux})}
\,+
\\
&
\ \ \ \ \ \
+
\frac{1}{{\sf c}\overline{\sf c}}\,
\left(
-\,\frac{\isqrt}{3}\,
\frac{\overline{\mathcal{L}}_1\big(\overline{\mathcal{L}}_1\big(
\overline{\mathcal{L}}_1(\kaux)\big)\big)}{
\overline{\mathcal{L}}_1(\kaux)}
+
\frac{4\,\isqrt}{9}\,
\frac{\overline{\mathcal{L}}_1\big(\overline{\mathcal{L}}_1
(\kaux)\big)^2}{\overline{\mathcal{L}}_1(\kaux)^2}
\,+
\right.
\\
&
\ \ \ \ \ \ \ \ \ \ \ \ \ \ \ \ \ \ \ \ \ \ \ \ \ \ \ \ \ \ 
\left.
+
\frac{\isqrt}{9}\,
\frac{\overline{\mathcal{L}}_1\big(\overline{\mathcal{L}}_1
(\kaux)\big)}{\overline{\mathcal{L}}_1(\kaux)}\,
\overline{\Paux}
+
\frac{\isqrt}{3}\,
\overline{\mathcal{L}}_1\big(\overline{\Paux}\big)
-
\frac{2\,\isqrt}{9}\,
\overline{\Paux}\,\overline{\Paux}
\right),
\\
Z^6
&
\,=\,
\isqrt\,
\frac{{\sf d}}{{\sf c}\overline{\sf c}}
-
\frac{{\sf e}{\sf e}}{{\sf c}{\sf c}}
+
\frac{{\sf e}}{{\sf c}\overline{\sf c}}
\left(
\frac{1}{3}\,
\overline{\Paux}
+
\frac{2}{3}\,
\frac{\overline{\mathcal{L}}_1\big(
\overline{\mathcal{L}}_1(\kaux)\big)}{
\overline{\mathcal{L}}_1(\kaux)}
\right)
+
\frac{\overline{\sf e}\overline{\sf e}}{
\overline{\sf c}\overline{\sf c}}\,
\frac{\overline{\mathcal{L}}_1(\overline{\kaux})}{
\mathcal{L}_1(\overline{\kaux})}.
\endaligned
\]
\end{Proposition}

\proof
We start by differentiating~({\ref{3-loop-initial}}),
finalizing directly the Maurer-Cartan part, thanks
to the Maurer-Cartan matrix shown above, and setting aside
$d\rho$ for the moment:
\[
\aligned
d\kappa
&
\,=\,
\beta\wedge\rho
+
\alpha\wedge\kappa
\,+
\\
&
\ \ \ \ \ 
-\,
\isqrt\,\overline{\sf c}{\sf e}\,
d\rho_0
+
{\sf c}\,
d\kappa_0',
\\
d\zeta
&
\,=\,
\gamma\wedge\rho
+
\isqrt\,\beta\wedge\kappa
+
\big(\alpha-\overline{\alpha}\big)
\wedge\zeta
\,+
\\
&
\ \ \ \ \
+
{\sf d}\,d\rho_0
+
{\sf e}\,d\kappa_0'
+
\frac{{\sf c}}{\overline{\sf c}}\,
d\zeta_0'.
\endaligned
\]

So we have to compute first:
\[
\aligned
K^3
&
\,=\,
\big[
\rho\wedge\overline{\kappa}
\big]
\big\{
d\kappa
\big\}
\\
&
\,=\,
-\,\isqrt\,
\overline{\sf c}{\sf e}\,
\big[
\rho\wedge\overline{\kappa}
\big]
\big\{
d\rho_0
\big\}
+
{\sf c}\,
\big[
\rho\wedge\overline{\kappa}
\big]
\big\{
d\kappa_0'
\big\}.
\endaligned
\]
The first term is, by~({\ref{inexplicit-0-structure-3-loop}}),
using the inversion 
formulas~({\ref{3-loop-initial-inverted}}):
\[
\!\!\!\!\!\!\!\!\!\!\!\!\!\!\!
\aligned
\big[
\rho\wedge\overline{\kappa}
\big]
\big\{
d\rho_0
\big\}
&
\,=\,
\big[
\rho\wedge\overline{\kappa}
\big]
\bigg\{
0
+
0
+
\overline{\Raux}_0^{1\prime}\,
\Big(
\frac{1}{{\sf c}\overline{\sf c}}
\Big)
\wedge
\Big(
\frac{1}{\overline{\sf c}}\,
\overline{\kappa}
\Big)
+
\overline{\Raux}_0^{2\prime}\,
\Big(
\frac{1}{{\sf c}\overline{\sf c}}
\Big)
\wedge
\Big(
-\frac{{\sf c}\overline{\sf e}}{\overline{\sf c}\overline{\sf c}}\,
\overline{\kappa}
\Big)
+
\isqrt\,
\Big(
\isqrt\,
\frac{{\sf e}}{{\sf c}{\sf c}}\,
\rho
\Big)
\wedge
\Big(
\frac{1}{\overline{\sf c}}\,
\overline{\kappa}
\Big)
\bigg\}
\\
&
\,=\,
\frac{1}{{\sf c}\overline{\sf c}\overline{\sf c}}\,
\overline{\Raux}_0^{1\prime}
-
\frac{\overline{\sf e}}{\overline{\sf c}\overline{\sf c}
\overline{\sf c}}\,
\overline{\Raux}_0^{2\prime}
-
\frac{{\sf e}}{{\sf c}{\sf c}\overline{\sf c}}.
\endaligned
\]

Similarly:
\[
\aligned
\big[
\rho\wedge\overline{\kappa}
\big]
\big\{
d\kappa_0'
\big\}
&
\,=\,
\big[
\rho\wedge\overline{\kappa}
\big]
\bigg\{
0
+
0
+
\Kaux_0^{3\prime}\,
\Big(
\frac{1}{{\sf c}\overline{\sf c}}\,\rho
\Big)
\wedge
\Big(
\frac{1}{\overline{\sf c}}\,\overline{\kappa}
\Big)
\,+
\\
&
\ \ \ \ \ \ \ \ \ \ \ \ \ \ \ \ \ \ \ \ \ \ \ \ \ \ \ \ 
+
0
+
\Kaux_0^{6\prime}\,
\Big(
\isqrt\,
\frac{{\sf e}}{{\sf c}{\sf c}}\,\rho
\Big)
\wedge
\Big(
\frac{1}{\overline{\sf c}}\,
\overline{\kappa}
\Big)
+
\bigg(
\Big(
-\isqrt\,\frac{\overline{\sf c}{\sf e}{\sf e}}{
{\sf c}{\sf c}{\sf c}}
-
\frac{{\sf d}}{{\sf c}{\sf c}}
\Big)\,\rho
\bigg)
\wedge
\Big(
\frac{1}{\overline{\sf c}}\,\overline{\kappa}
\Big)
\bigg\}
\\
&
\,=\,
\frac{1}{{\sf c}\overline{\sf c}\overline{\sf c}}\,
\Kaux_0^{3\prime}
+
\isqrt\,
\frac{{\sf e}}{{\sf c}{\sf c}\overline{\sf c}}\,
\Kaux_0^{6\prime}
-
\isqrt\,
\frac{{\sf e}{\sf e}}{{\sf c}{\sf c}{\sf c}}
-
\frac{{\sf d}}{{\sf c}{\sf c}\overline{\sf c}}.
\endaligned
\]
Hence:
\[
\!\!\!\!\!\!\!\!\!\!\!\!\!\!\!
\aligned
\Kaux_0^{3\prime}
&
\,=\,
-\,\isqrt\,
\frac{{\sf e}}{{\sf c}\overline{\sf c}}\,
\overline{\Raux}_0^{1\prime}
+
\isqrt\,
\frac{{\sf e}\overline{\sf e}}{\overline{\sf c}\overline{\sf c}}\,
\overline{\Raux}_0^{2\prime}
+
\zero{
\isqrt\,\frac{{\sf e}{\sf e}}{{\sf c}{\sf c}}}
+
\frac{1}{\overline{\sf c}\overline{\sf c}}\,
\Kaux_0^{3\prime}
+
\isqrt\,
\frac{{\sf e}}{{\sf c}\overline{\sf c}}\,
\Kaux_0^{6\prime}
-
\zero{
\isqrt\,
\frac{{\sf e}{\sf e}}{{\sf c}{\sf c}}}
-
\frac{{\sf d}}{{\sf c}\overline{\sf c}}
\\
&
\,=\,
-\,
\frac{{\sf d}}{{\sf c}\overline{\sf c}}
+
\frac{{\sf e}}{{\sf c}\overline{\sf c}}
\left(
-\frac{\isqrt}{3}\,
\frac{\overline{\mathcal{L}}_1\big(
\overline{\mathcal{L}}_1\big(\kaux)\big)}{
\overline{\mathcal{L}}_1(\kaux)}
-
\frac{2\,\isqrt}{3}\,
\overline{\Paux}
-
\frac{\isqrt}{3}\,
\overline{\mathcal{L}}_1\big(
\overline{\mathcal{L}}_1(\kaux)\big)
+
\frac{\isqrt}{3}\,
\overline{\paux}
\right)
-
\isqrt\,
\frac{{\sf e}\overline{\sf e}}{\overline{\sf c}\overline{\sf c}}\,
\frac{\overline{\mathcal{L}}_1(\overline{\kaux})}{
\mathcal{L}_1(\overline{\kaux})}
+
\frac{1}{\overline{\sf c}\overline{\sf c}}\,
\Kaux_0^{3\prime}.
\endaligned
\]
Replacing this last term $\Kaux_0^{3\prime}$ by 
its value from
Proposition~{\ref{Proposition-3-loop-structure-0-coframe}},
we reach the stated explicit expression of $K^3$.

\smallskip

Next:
\[
\aligned
Z^6
&
\,=\,
\big[
\kappa\wedge\overline{\kappa}
\big]
\big\{
d\zeta
\big\}
\\
&
\,=\,
{\sf d}\,
\big[
\kappa\wedge\overline{\kappa}
\big]
\big\{
d\rho_0
\big\}
+
{\sf e}\,
\big[
\kappa\wedge\overline{\kappa}
\big]
\big\{
d\kappa_0'
\big\}
+
\frac{{\sf c}}{\overline{\sf c}}\,
\big[
\kappa\wedge\overline{\kappa}
\big]
\big\{
d\zeta_0'
\big\}.
\endaligned
\]
Separately:
\[
\aligned
\big[
\kappa\wedge\overline{\kappa}
\big]
\big\{
{\sf d}\,d\rho_0
\big\}
&
\,=\,
0+0+0+0+
{\sf d}\,\isqrt\,
\frac{1}{{\sf c}\overline{\sf c}}
\,\,=\,\,
\isqrt\,\frac{{\sf d}}{{\sf c}\overline{\sf c}},
\\
\big[
\kappa\wedge\overline{\kappa}
\big]
\big\{
{\sf e}\,d\kappa_0'
\big\}
&
\,=\,
0+0+0+0
+
{\sf e}\,\Kaux_0^{6\prime}\,\frac{1}{{\sf c}\overline{\sf c}}
+
{\sf e}\,
\Big(
-\frac{\overline{\sf c}{\sf e}}{{\sf c}{\sf c}}
\Big)\,
\frac{1}{\overline{\sf c}}
\,\,=\,\,
\frac{{\sf e}}{{\sf c}\overline{\sf c}}\,
\Kaux_0^{6\prime}
-
\frac{{\sf e}{\sf e}}{{\sf c}{\sf c}},
\\
\big[
\kappa\wedge\overline{\kappa}
\big]
\Big\{
\frac{{\sf c}}{\overline{\sf c}}\,d\zeta_0'
\Big\}
&
\,=\,
0+0
+
\frac{{\sf c}}{\overline{\sf c}}\,
\Zaux_0^{8\prime}\,
\Big(
-\frac{\overline{\sf c}{\sf e}}{{\sf c}{\sf c}}
\Big)\,
\Big(
\frac{1}{\overline{\sf c}}
\Big)
+
\frac{{\sf c}}{\overline{\sf c}}\,
\Zaux_0^{9\prime}\,
\Big(
-\frac{\overline{\sf c}{\sf e}}{{\sf c}{\sf c}}
\Big)\,
\Big(
-\frac{{\sf c}\overline{\sf e}}{\overline{\sf c}\overline{\sf c}}
\Big)
\\
&
\,=\,
-\,\frac{{\sf e}}{{\sf c}\overline{\sf c}}\,
\Zaux_0^{8\prime}
+
\frac{{\sf e}\overline{\sf e}}{\overline{\sf c}\overline{\sf c}}\,
\Zaux_0^{9\prime},
\endaligned
\]
hence summing and inserting the explicit expressions from
Proposition~{\ref{Proposition-3-loop-structure-0-coframe}},
we conclude:
\begin{align}
Z^6
&
\,=\,
\isqrt\,
\frac{{\sf d}}{{\sf c}\overline{\sf c}}
+
\frac{{\sf e}}{{\sf c}\overline{\sf c}}\,
\Kaux_0^{6\prime}
-
\frac{{\sf e}{\sf e}}{{\sf c}{\sf c}}
-
\frac{{\sf e}}{{\sf c}\overline{\sf c}}\,
\Zaux_0^{8\prime}
+
\frac{{\sf e}\overline{\sf e}}{\overline{\sf c}\overline{\sf c}}\,
\Zaux_0^{9\prime}
\notag
\\
&
\,=\,
\isqrt\,
\frac{{\sf d}}{{\sf c}\overline{\sf c}}
-
\frac{{\sf e}{\sf e}}{{\sf c}{\sf c}}
+
\frac{{\sf e}}{{\sf c}\overline{\sf c}}
\left(
\frac{1}{3}\,
\overline{\Paux}
+
\frac{2}{3}\,
\frac{\overline{\mathcal{L}}_1\big(
\overline{\mathcal{L}}_1(\kaux)\big)}{
\overline{\mathcal{L}}_1(\kaux)}
\right)
+
\frac{\overline{\sf e}\overline{\sf e}}{
\overline{\sf c}\overline{\sf c}}\,
\frac{\overline{\mathcal{L}}_1(\overline{\kaux})}{
\mathcal{L}_1(\overline{\kaux})}.
\qedhere
\end{align}
\endproof

Once we have reached the explicit expressions of both
$K^3$ and $Z^6$, when we perform the essential
combination $\isqrt\, K^3 - Z^6$, we see that 
both the coefficients of $\frac{{\sf e}}{{\sf c}\overline{\sf c}}$
and of $\frac{{\sf e}\overline{\sf e}}{\overline{\sf c}
\overline{\sf c}}$ disappear, and it remains:
\[
\aligned
\!\!\!\!\!\!\!\!\!\!\!\!\!\!\!\!\!\!\!\!
\isqrt\,
K^3
-
Z^6
&
\,=\,
-\,2\isqrt\,
\frac{{\sf d}}{{\sf c}\overline{\sf c}}
+
\frac{{\sf e}{\sf e}}{{\sf c}{\sf c}}
+
\isqrt\,
\frac{1}{{\sf c}\overline{\sf c}}\,
\Kaux_0^{3\prime}
\\
&
\,=\,
-\,2\isqrt\,
\frac{{\sf d}}{{\sf c}\overline{\sf c}}
+
\frac{{\sf e}{\sf e}}{{\sf c}{\sf c}}
\,+
\\
&
\ \ \ \ \ 
+
\frac{1}{{\sf c}\overline{\sf c}}\,
\Bigg(
\underbrace{
\frac{1}{3}\,
\frac{\overline{\mathcal{L}}_1\big(
\overline{\mathcal{L}}_1\big(
\overline{\mathcal{L}}_1(\kaux)\big)\big)}{
\overline{\mathcal{L}}_1(\kaux)}
-
\frac{4}{9}\,
\frac{\overline{\mathcal{L}}_1\big(
\overline{\mathcal{L}}_1(\kaux)\big)^2}{
\overline{\mathcal{L}}_1(\kaux)^2}
-
\frac{1}{9}\,
\frac{\overline{\mathcal{L}}_1\big(
\overline{\mathcal{L}}_1(\kaux)\big)}{
\overline{\mathcal{L}}_1(\kaux)}\,
\overline{\Paux}
-
\frac{1}{3}\,
\overline{\mathcal{L}}_1\big(\overline{\Paux}\big)
+
\frac{2}{9}\,
\overline{\Paux}\,
\overline{\Paux}}_{=:\,\,-\,2\,\Haux_0}
\Bigg).
\endaligned
\]
We introduce, as is underbraced:
\[
\Haux_0
\,:=\,
-\,\frac{1}{6}\,
\frac{\overline{\mathcal{L}}_1\big(
\overline{\mathcal{L}}_1\big(
\overline{\mathcal{L}}_1(\kaux)\big)\big)}{
\overline{\mathcal{L}}_1(\kaux)}
+
\frac{2}{9}\,
\frac{\overline{\mathcal{L}}_1\big(
\overline{\mathcal{L}}_1(\kaux)\big)^2}{
\overline{\mathcal{L}}_1(\kaux)^2}
+
\frac{1}{18}\,
\frac{\overline{\mathcal{L}}_1\big(
\overline{\mathcal{L}}_1(\kaux)\big)}{
\overline{\mathcal{L}}_1(\kaux)}\,
\overline{\Paux}
+
\frac{1}{6}\,
\overline{\mathcal{L}}_1\big(\overline{\Paux}\big)
-
\frac{1}{9}\,
\overline{\Paux}\,
\overline{\Paux},
\]
a function
which coincides with Pocchiola's function $H$. Then
by means of the invariant condition:
\[
0
\,=\,
\isqrt\,K^3
-
Z^6,
\]
we reach a convenient normalization of the group
parameter:
\[
\!\!\!\!\!\!\!\!\!\!\!\!\!\!\!\!\!\!\!\!\!\!\!\!\!
\small
\aligned
{\sf d}
\,:=\,
&\,
-\,\frac{\isqrt}{2}\,
\frac{\overline{\sf c}{\sf e}{\sf e}}{{\sf c}}
+
\isqrt\,
\frac{{\sf c}}{\overline{\sf c}}\,
\Haux_0
\\
\,=\,
&\,
-\,\frac{\isqrt}{2}\,
\frac{\overline{\sf c}{\sf e}{\sf e}}{{\sf c}}
+
\isqrt\,
\frac{{\sf c}}{\overline{\sf c}}
\left(
-\,\frac{1}{6}\,
\frac{\overline{\mathcal{L}}_1\big(
\overline{\mathcal{L}}_1\big(
\overline{\mathcal{L}}_1(\kaux)\big)\big)}{
\overline{\mathcal{L}}_1(\kaux)}
+
\frac{2}{9}\,
\frac{\overline{\mathcal{L}}_1\big(
\overline{\mathcal{L}}_1(\kaux)\big)^2}{
\overline{\mathcal{L}}_1(\kaux)^2}
+
\frac{1}{18}\,
\frac{\overline{\mathcal{L}}_1\big(
\overline{\mathcal{L}}_1(\kaux)\big)}{
\overline{\mathcal{L}}_1(\kaux)}\,
\overline{\Paux}
+
\frac{1}{6}\,
\overline{\mathcal{L}}_1\big(\overline{\Paux}\big)
-
\frac{1}{9}\,
\overline{\Paux}\,
\overline{\Paux}
\right).
\endaligned
\]

Before we really perform this normalization of the
group parameter ${\sf d}$, let us point out that
some other invariant 
relations between torsion coefficients appear.
In fact, we see above that:
\[
\aligned
\isqrt\,K^{4\prime}
&
\,=\,
\isqrt\,
K^4
-
\isqrt\,
b_5,
\\
Z^{7\prime}
&
\,=\,
Z^7
-
\isqrt\,
b_5,
\endaligned
\]
whence:
\[
\isqrt\,K^{4\prime}
-
Z^{7\prime}
\,=\,
\isqrt\,K^4
-
Z^7.
\]
However, the next lemma shows that no group parameter
can be normalized so.

\begin{Lemma}
Their explicit expressions are:
\[
\isqrt\,K^4
\,=\,
Z^7
\,=\,
-\,
\frac{{\sf e}}{{\sf c}}\,
\frac{\overline{\mathcal{L}}_1(\overline{\kaux})}{
\mathcal{L}_1(\overline{\kaux})}.
\]
\end{Lemma}

\proof
Indeed, 
by~({\ref{inexplicit-0-structure-3-loop}}),
replacing $\overline{\Raux}_0^{2\prime}$ from
Proposition~{\ref{Proposition-3-loop-structure-0-coframe}},
we can compute using~({\ref{3-loop-initial-inverted}}):
\[
\aligned
K^4
&
\,=\,
\big[
\rho\wedge\overline{\zeta}
\big]
\big\{
-\isqrt\,\overline{\sf c}{\sf e}\,
d\rho_0
+
{\sf c}\,
d\kappa_0'
\big\}
\\
&
\,=\,
-\,\isqrt\,\overline{\sf c}{\sf e}\,
\big[
\rho\wedge\overline{\zeta}
\big]
\big\{
d\rho_0
\big\}
+
{\sf c}\,
\big[
\rho\wedge\overline{\zeta}
\big]
\big\{
d\kappa_0'
\big\}
\\
&
\,=\,
-\,\isqrt\,
\overline{\sf c}{\sf e}\,
\bigg(
0
+
0
+
0
+
\overline{\Raux}_0^{2\prime}\,
\Big(
\frac{1}{{\sf c}\overline{\sf c}}
\Big)
\Big(
\frac{{\sf c}}{\overline{\sf c}}
\Big)
\bigg)
+
{\sf c}
\cdot 
0
\\
&
\,=\,
-\,\isqrt\,
\frac{{\sf e}}{\overline{\sf c}}\,
\left(
-\,
\frac{\overline{\mathcal{L}}_1(\overline{\kaux})}{
\overline{\mathcal{L}}_1(\kaux)}
\right),
\endaligned
\]
and similarly:
\begin{align}
Z^7
&
\,=\,
\big[
\kappa\wedge\overline{\zeta}
\big]
\Big\{
{\sf d}\,d\rho_0
+
{\sf e}\,d\kappa_0'
+
\frac{{\sf c}}{\overline{\sf c}}\,
d\zeta_0'
\Big\}
\notag
\\
&
\,=\,
{\sf d}\,
\big[
\kappa\wedge\overline{\zeta}
\big]
\big\{
d\rho_0
\big\}
+
{\sf e}\,
\big[
\kappa\wedge\overline{\zeta}
\big]
\big\{
d\kappa_0'
\big\}
+
\frac{{\sf c}}{\overline{\sf c}}\,
\big[
\kappa\wedge\overline{\zeta}
\big]
\big\{
d\zeta_0'
\big\}
\notag
\\
&
\,=\,
0+0
+
\frac{{\sf c}}{\overline{\sf c}}\,
\bigg(
0+0+0
+
\Zaux_0^{9\prime}\,
\Big(
-\frac{\overline{\sf c}{\sf e}}{{\sf c}{\sf c}}
\Big)
\Big(
\frac{{\sf c}}{\overline{\sf c}}
\Big)
\bigg)
\notag
\\
&
\,=\,
-\,
\frac{{\sf e}}{\overline{\sf c}}\,
\frac{\overline{\mathcal{L}}_1(\overline{\kaux})}{
\overline{\mathcal{L}}_1(\kaux)}.
\qedhere
\end{align}
\endproof

Another invariant torsion combination is the following.

\begin{Lemma}
\label{Lemma-K-2-Z-5-bar-Z-8}
Here is an essential linear combination of torsion terms:
\[
-\,
\isqrt\,
K^2
+
Z^5
-
\overline{Z}^8.
\]
\end{Lemma}

\proof
A glance at what precedes shows:
\[
\aligned
K^{2\prime}
&
\,=\,
K^2
-
b_3,
\\
Z^{5\prime}
&
\,=\,
Z^5
-
\isqrt\,b_3
+
a_2
-
\overline{a}_4,
\\
Z^{8\prime}
&
\,=\,
Z^8
-
a_4
+
\overline{a}_2,
\endaligned
\]
whence indeed:
\[
-\,\isqrt\,
K^{2\prime}
+
Z^{5\prime}
-
\overline{Z}^{8\prime}
\,=\,
-\,
\isqrt\,K^2
+
Z^5
-
\overline{Z}^8.
\qedhere
\]
\endproof

\begin{Lemma}
Their explicit expressions are:
\[
\aligned
K^2
&
\,=\,
\isqrt\,
\frac{\overline{\sf e}}{{\sf c}}
+
\frac{1}{{\sf c}}
\left(
-\,\frac{\isqrt}{3}\,
\frac{\mathcal{K}\big(\overline{\mathcal{L}}_1\big(
\overline{\mathcal{L}}_1(\kaux)\big)\big)}{
\overline{\mathcal{L}}_1(\kaux)^2}
+
\frac{\isqrt}{3}\,
\frac{\mathcal{K}\big(\overline{\mathcal{L}}_1(\kaux)\big)\,\,
\overline{\mathcal{L}}_1\big(\overline{\mathcal{L}}_1(\kaux)\big)}{
\overline{\mathcal{L}}_1(\kaux)^3}
\,-
\right.
\\
&
\left.
\ \ \ \ \ \ \ \ \ \ \ \ \ \ \ \ \ \ \ \ \ \ \ \
\ \ \ \ \ \ \ \ \ \ \ \ \ \ \ \ \ \ \ \ \ \ \ \
-
\frac{\isqrt}{3}\,
\frac{\mathcal{L}_1\big(\mathcal{L}_1(\overline{\kaux})\big)}{
\mathcal{L}_1(\overline{\kaux})}
-
\frac{\isqrt}{3}\,
\frac{\overline{\mathcal{L}}_1\big(\mathcal{L}_1(\kaux)\big)}{
\overline{\mathcal{L}}_1(\kaux)}
-
\frac{2}{3}\,
\frac{\mathcal{T}(\kaux)}{\overline{\mathcal{L}}_1(\kaux)}
\right),
\\
Z^5
&
\,=\,
\frac{1}{{\sf c}}\,
\frac{\mathcal{L}_1\big(\overline{\mathcal{L}}_1(\kaux)\big)}{
\overline{\mathcal{L}}_1(\kaux)}
-
\frac{\overline{\sf c}{\sf e}}{{\sf c}{\sf c}}\,
\frac{\mathcal{L}_1(\overline{\kaux})}{
\overline{\mathcal{L}}_1(\kaux)},
\\
Z^8
&
\,=\,
\frac{{\sf e}}{{\sf c}}
-
\frac{1}{\overline{\sf c}}\,
\frac{\overline{\mathcal{L}}_1\big(
\overline{\mathcal{L}}_1(\kaux)\big)}{
\overline{\mathcal{L}}_1(\kaux)}
-
\frac{{\sf c}\overline{\sf e}}{
\overline{\sf c}\overline{\sf c}}\,
\frac{\overline{\mathcal{L}}_1(\kaux)}{
\mathcal{L}_1(\overline{\kaux})}.
\endaligned
\]
\end{Lemma}

\proof
Recall:
\[
\aligned
d\rho
&
\,=\,
\big(\alpha+\overline{\alpha}\big)
\wedge
\rho
\,+
\\
&
\ \ \ \ \ \ \ \ \ \ 
+
{\sf c}\overline{\sf c}\,
d\rho_0,
\\
d\kappa
&
\,=\,
\beta\wedge\rho
+
\alpha\wedge\kappa
\,-
\\
&
\ \ \ \ \ \ \ \ \ \ 
-\,
\isqrt\,\overline{\sf c}{\sf e}\,
d\rho_0
+
{\sf c}\,
d\kappa_0',
\\
d\zeta
&
\,=\,
\gamma\wedge\rho
+
\isqrt\,
\beta\wedge\kappa
+
\big(\alpha-\overline{\alpha}\big)
\wedge\rho
\,+
\\
&
\ \ \ \ \ \ \ \ \ \ 
+
{\sf d}\,
d\rho_0
+
{\sf e}\,
d\kappa_0'
+
\frac{{\sf c}}{\overline{\sf c}}\,
d\zeta_0',
\endaligned
\]
hence:
\[
K^2
\,=\,
\big[
\rho\wedge\zeta
\big]
\Big\{
-\,\isqrt\,
\overline{\sf c}{\sf e}\,
d\rho_0
+
{\sf c}\,
d\kappa_0'
\Big\}.
\]
Visually:
\[
\aligned
\big[
\rho\wedge\zeta
\big]
\big\{
d\rho_0
\big\}
&
\,=\,
\Raux_0^{2\prime}\,
\Big(
\frac{1}{{\sf c}\overline{\sf c}}
\Big)
\Big(
\frac{\overline{\sf c}}{{\sf c}}
\Big)
\,\,=\,\,
\frac{1}{{\sf c}{\sf c}}\,
\Raux_0^{2\prime},
\\
\big[
\rho\wedge\zeta
\big]
\big\{
d\kappa_0'
\big\}
&
\,=\,
\Kaux_0^{2\prime}\,
\Big(
\frac{1}{{\sf c}\overline{\sf c}}
\Big)
\Big(
\frac{\overline{\sf c}}{{\sf c}}
\Big)
+
\Kaux_0^{5\prime}\,
\Big(
\isqrt\,
\frac{{\sf e}}{{\sf c}{\sf c}}
\Big)
\Big(
\frac{\overline{\sf c}}{{\sf c}}
\Big)
-
\Big(
\frac{\overline{\sf c}}{{\sf c}}
\Big)
\Big(
-\isqrt\,
\frac{\overline{\sf e}}{\overline{\sf c}\overline{\sf c}}
\Big)
\\
&
\,=\,
\frac{1}{{\sf c}{\sf c}}\,
\Kaux_0^{2\prime}
+
\isqrt\,
\frac{\overline{\sf c}{\sf e}}{
{\sf c}{\sf c}{\sf c}}\,
\Kaux_0^{5\prime}
+
\isqrt\,
\frac{\overline{\sf e}}{{\sf c}\overline{\sf c}},
\endaligned
\]
hence:
\[
\aligned
K^2
&
\,=\,
-\,\isqrt\,
\frac{\overline{\sf c}{\sf e}}{{\sf c}{\sf c}}\,
\Raux_0^{2\prime}
+
\frac{1}{{\sf c}}\,
\Kaux_0^{2\prime}
+
\isqrt\,
\frac{\overline{\sf c}{\sf e}}{{\sf c}{\sf c}}\,
\Kaux_0^{5\prime}
+
\isqrt\,\frac{\overline{\sf e}}{\overline{\sf c}}
\\
&
\,=\,
\isqrt\,
\frac{\overline{\sf e}}{\overline{\sf c}}
+
\frac{1}{{\sf c}}
\left(
-\,\frac{\isqrt}{3}\,
\frac{\mathcal{K}\big(\overline{\mathcal{L}}_1\big(
\overline{\mathcal{L}}_1(\kaux)\big)\big)}{
\overline{\mathcal{L}}_1(\kaux)^2}
+
\frac{\isqrt}{3}\,
\frac{\mathcal{K}\big(\overline{\mathcal{L}}_1(\kaux)\big)\,\,
\overline{\mathcal{L}}_1\big(\overline{\mathcal{L}}_1(\kaux)\big)
}{
\overline{\mathcal{L}}_1(\kaux)^3}
\,-
\right.
\\
&
\left.
\ \ \ \ \ \ \ \ \ \ \ \ \ \ \ \ \ \ \ \ \ \ \ \ \ \ \ \ \ \ \ \ \
\ \ \ \ \ \ \ \ \ \ \ \ \ \ \ \
-\,
\frac{\isqrt}{3}\,
\frac{\mathcal{L}_1\big(\mathcal{L}_1(\overline{\kaux})\big)}{
\mathcal{L}_1(\overline{\kaux})}
-
\frac{\isqrt}{3}\,
\frac{\overline{\mathcal{L}}_1\big(\mathcal{L}_1(\kaux)\big)}{
\overline{\mathcal{L}}_1(\kaux)}
-
\frac{2}{3}\,
\frac{\mathcal{T}(\kaux)}{\overline{\mathcal{L}}_1(\kaux)}
\right).
\endaligned
\]

Next, treat:
\[
Z^5
\,=\,
\big[
\kappa\wedge\zeta
\big]
\Big\{
{\sf d}\,d\rho_0
+
{\sf e}\,d\kappa_0'
+
\frac{{\sf c}}{\overline{\sf c}}\,
d\zeta_0'
\Big\},
\]
using:
\[
\aligned
\big[
\kappa\wedge\zeta
\big]
\big\{
d\rho_0
\big\}
&
\,=\,
0,
\\
\big[
\kappa\wedge\zeta
\big]
\big\{
d\kappa_0'
\big\}
&
\,=\,
\Kaux_0^{5\prime}\,
\Big(
\frac{1}{{\sf c}}
\Big)
\Big(
\frac{\overline{\sf c}}{{\sf c}}
\Big)
\,\,=\,\,
\frac{\overline{\sf c}}{{\sf c}{\sf c}}\,
\Kaux_0^{5\prime},
\\
\big[
\kappa\wedge\zeta
\big]
\big\{
d\zeta_0'
\big\}
&
\,=\,
\Zaux_0^{5\prime}\,
\Big(
\frac{1}{{\sf c}}
\Big)
\Big(
\frac{\overline{\sf c}}{{\sf c}}
\Big)
\,\,=\,\,
\frac{\overline{\sf c}}{{\sf c}{\sf c}}\,
\Zaux_0^{5\prime},
\endaligned
\]
so:
\[
Z^5
\,=\,
\frac{\overline{\sf c}{\sf e}}{{\sf c}{\sf c}}
\left(
-\,\frac{\mathcal{L}_1(\kaux)}{
\overline{\mathcal{L}}_1(\kaux)}
\right)
+
\frac{1}{{\sf c}}
\left(
\frac{\mathcal{L}_1\big(\overline{\mathcal{L}}_1(\kaux)\big)}{
\overline{\mathcal{L}}_1(\kaux)}
\right).
\]

Lastly treat:
\[
Z^8
\,=\,
\big[
\zeta\wedge\overline{\kappa}
\big]
\Big\{
{\sf d}\,
d\rho_0
+
{\sf e}\,
d\kappa_0'
+
\frac{{\sf c}}{\overline{\sf c}}\,
d\zeta_0'
\Big\},
\]
using:
\[
\aligned
\big[
\zeta\wedge\overline{\kappa}
\big]
\big\{
d\rho_0
\big\}
&
\,=\,
0,
\\
\big[
\zeta\wedge\overline{\kappa}
\big]
\big\{
d\kappa_0'
\big\}
&
\,=\,
\frac{\overline{\sf c}}{{\sf c}}\,
\frac{1}{\overline{\sf c}}
\,=\,
\frac{1}{{\sf c}},
\\
\big[
\zeta\wedge\overline{\kappa}
\big]
\big\{
d\zeta_0'
\big\}
&
\,=\,
\Zaux_0^{8\prime}\,
\frac{\overline{\sf c}}{{\sf c}}\,
\frac{1}{\overline{\sf c}}
+
\Zaux_0^{9\prime}\,
\frac{\overline{\sf c}}{{\sf c}}\,
\Big(
-\frac{{\sf c}\overline{\sf e}}{\overline{\sf c}\overline{\sf c}}
\Big)
\,\,=\,\,
\frac{1}{{\sf c}}\,
\Zaux_0^{8\prime}
-
\frac{\overline{\sf e}}{\overline{\sf c}}\,
\Zaux_0^{9\prime},
\endaligned
\]
which concludes:
\begin{align}
Z^8
&
\,=\,
\frac{{\sf e}}{{\sf c}}
+
\frac{1}{\overline{\sf c}}\,
\Zaux_0^{8\prime}
-
\frac{{\sf c}\overline{\sf e}}{\overline{\sf c}\overline{\sf c}}\,
\Zaux_0^{9\prime}
\notag
\\
&
\,=\,
\frac{{\sf e}}{{\sf c}}
-
\frac{1}{\overline{\sf c}}\,
\frac{\overline{\mathcal{L}}_1\big(\overline{\mathcal{L}}_1
(\kaux)\big)}{
\overline{\mathcal{L}}_1(\kaux)}
-
\frac{{\sf c}\overline{\sf e}}{\overline{\sf c}\overline{\sf c}}\,
\frac{\overline{\mathcal{L}}_1(\kaux)}{
\mathcal{L}_1(\overline{\kaux})}.
\qedhere
\end{align}
\endproof

Thanks to these explicit expressions, we can compute
the essential linear combination of torsion terms, emphasizing
two important annihilations by pairs:
\[
\aligned
-\,\isqrt\,K^2
+
Z^5
-
\overline{Z}^8
&
\,=\,
\zero{\frac{\overline{\sf e}}{{\sf c}}}
+
\frac{1}{{\sf c}}
\left(
-\,\frac{1}{3}\,
\frac{\mathcal{K}\big(\overline{\mathcal{L}}_1\big(
\overline{\mathcal{L}}_1(\kaux)\big)\big)}{
\overline{\mathcal{L}}_1(\kaux)^2}
+
\frac{1}{3}\,
\frac{\mathcal{K}\big(\overline{\mathcal{L}}_1(\kaux)\big)\,\,
\overline{\mathcal{L}}_1\big(\overline{\mathcal{L}}_1(\kaux)\big)}{
\overline{\mathcal{L}}_1(\kaux)^3}
\,-
\right.
\\
&
\left.
\ \ \ \ \ \ \ \ \ \ \ \ \ \ \ \ \ \ \ \ \ \ \ \
\ \ \ \ \ \ \ \ \ \ \ \ \
-
\frac{1}{3}\,
\frac{\mathcal{L}_1\big(\mathcal{L}_1(\overline{\kaux})\big)}{
\mathcal{L}_1(\overline{\kaux})}
-
\frac{1}{3}\,
\frac{\overline{\mathcal{L}}_1\big(\mathcal{L}_1(\kaux)\big)}{
\overline{\mathcal{L}}_1(\kaux)}
+
\frac{2\,\isqrt}{3}\,
\frac{\mathcal{T}(\kaux)}{\overline{\mathcal{L}}_1(\kaux)}
\right)
\\
&
\ \ \ \ \
+
\frac{1}{{\sf c}}\,
\frac{\mathcal{L}_1\big(\overline{\mathcal{L}}_1(\kaux)\big)}{
\overline{\mathcal{L}}_1(\kaux)}
-
\zerozero{\frac{\overline{\sf c}{\sf e}}{{\sf c}{\sf c}}\,
\frac{\mathcal{L}_1(\overline{\kaux})}{
\overline{\mathcal{L}}_1(\kaux)}}
\,
-
\\
&
\ \ \ \ \ \
-\,
\zero{\frac{\overline{\sf e}}{\overline{\sf c}}}
+
\frac{1}{{\sf c}}\,
\frac{\mathcal{L}_1\big(\mathcal{L}_1(\overline{\kaux})\big)}{
\mathcal{L}_1(\overline{\kaux})}
+
\zerozero{\frac{\overline{\sf c}{\sf e}}{{\sf c}{\sf c}}\,
\frac{\mathcal{L}_1(\overline{\kaux})}{
\overline{\mathcal{L}}_1(\kaux)}}.
\endaligned
\]
Also, in order to match exactly with 
Pocchiola's function $W$ introduced
in~{\cite{Pocchiola-2013, Merker-Pocchiola-2018}}, we decompose
the last term of the second line as:
\[
\frac{2\,\isqrt}{3}\,
\frac{\mathcal{T}(\kaux)}{\overline{\mathcal{L}}_1(\kaux)}
\,=\,
-\,\frac{1}{3}\,
\frac{\mathcal{L}_1\big(\overline{\mathcal{L}}_1(\kaux)\big)}{
\overline{\mathcal{L}}_1(\kaux)}
+
\frac{1}{3}\,
\frac{\overline{\mathcal{L}}_1\big(\mathcal{L}_1(\kaux)\big)}{
\overline{\mathcal{L}}_1(\kaux)}
+
\frac{\isqrt}{3}\,
\frac{\mathcal{T}(\kaux)}{\overline{\mathcal{L}}_1(\kaux)},
\]
so that a third pair of terms disappears, and after 
reorganization\,\,---\,\,no pen needed\,\,---,
the result is:
\[
\aligned
-\,\isqrt\,K^2
+
Z^5
-
\overline{Z}^8
\,=\,
&\,
\frac{1}{{\sf c}}
\left(
-\,\frac{1}{3}\,
\frac{\mathcal{K}\big(\overline{\mathcal{L}}_1\big(
\overline{\mathcal{L}}_1(\kaux)\big)\big)}{
\overline{\mathcal{L}}_1(\kaux)^2}
+
\frac{1}{3}\,
\frac{\mathcal{K}\big(\overline{\mathcal{L}}_1(\kaux)\big)\,\,
\overline{\mathcal{L}}_1\big(\overline{\mathcal{L}}_1(\kaux)\big)}{
\overline{\mathcal{L}}_1(\kaux)^3}
\,+
\right.
\\
&
\ \ \ \ \ \ \ \ \ \ \ \ \ \ \ \ \ \ \ \ \ \ \ \
\ \ \ \ \ \ \ \ \ \ \ \ \
\left.
+
\frac{2}{3}\,
\frac{\mathcal{L}_1\big(\mathcal{L}_1(\overline{\kaux})\big)}{
\mathcal{L}_1(\overline{\kaux})}
+
\frac{2}{3}\,
\frac{\mathcal{L}_1\big(\overline{\mathcal{L}}_1(\kaux)\big)}{
\overline{\mathcal{L}}_1(\kaux)}
+
\frac{\isqrt}{3}\,
\frac{\mathcal{T}(\kaux)}{\overline{\mathcal{L}}_1(\kaux)}
\right)
\\
\,=:\,
&\,
\frac{1}{{\sf c}}\,
\Waux_0,
\endaligned
\]
and this defines a new horizontal function $\Waux_0$, equal
to Pocchiola's function $W$.

For now, we will not use the potential normalization ${\sf c} =
\Waux_0$ on the open subset of $M^5 \subset \C^3$ on which:
\[
0
\,\neq\,
\Waux_0
\big(
z_1,z_2,\overline{z}_1,\overline{z}_2,v
\big),
\]
if nonempty\,\,---\,\,a hypothesis must be set up\,\,---,
but we will deal with this discussion later. 
In fact, before proceeding, we state a technical 
differential relation useful later, whose proof can 
be skipped in a first reading.

\begin{Lemma}
\label{Lemma-K-bar-H-0}
One has:
\[
\overline{\mathcal{K}}
\big(
\Haux_0
\big)
\,=\,
-\,2\,
\overline{\mathcal{L}}_1
\big(\overline{\kaux}\big)\,
\Haux_0.
\]
\end{Lemma}

\proof
Apply the derivation $\overline{\mathcal{K}}$ to $\Haux_0$:
\[
\aligned
\overline{\mathcal{K}}
\big(\Haux_0\big)
&
\,=\,
-\,\frac{1}{6}\,
\frac{\overline{\mathcal{K}}\big(
\overline{\mathcal{L}}_1\big(
\overline{\mathcal{L}}_1\big(
\overline{\mathcal{L}}_1(\kaux)\big)\big)\big)}{
\overline{\mathcal{L}}_1(\kaux)}
+
\frac{1}{6}\,
\frac{\overline{\mathcal{K}}\big(
\overline{\mathcal{L}}_1(\kaux)\big)\,\,
\overline{\mathcal{L}}_1\big(
\overline{\mathcal{L}}_1\big(
\overline{\mathcal{L}}_1(\kaux)\big)\big)}{
\overline{\mathcal{L}}_1(\kaux)^2}
\,+
\\
&
\ \ \ \ \ 
+
\frac{4}{9}\,
\frac{\overline{\mathcal{K}}\big(
\overline{\mathcal{L}}_1\big(
\overline{\mathcal{L}}_1(\kaux)\big)\big)\,\,
\overline{\mathcal{L}}_1\big(
\overline{\mathcal{L}}_1(\kaux)\big)}{
\overline{\mathcal{L}}_1(\kaux)^2}
-
\frac{4}{9}\,
\frac{\overline{\mathcal{K}}\big(
\overline{\mathcal{L}}_1(\kaux)\big)\,\,
\overline{\mathcal{L}}_1\big(
\overline{\mathcal{L}}_1(\kaux)\big)^2}{
\overline{\mathcal{L}}_1(\kaux)^3}
\,+
\\
&
\ \ \ \ \ \
+
\frac{1}{18}\,
\frac{\overline{\mathcal{K}}\big(
\overline{\mathcal{L}}_1\big(
\overline{\mathcal{L}}_1(\kaux)\big)\big)}{
\overline{\mathcal{L}}_1(\kaux)}\,
\overline{\Paux}
+
\frac{1}{18}\,
\frac{\overline{\mathcal{L}}_1\big(
\overline{\mathcal{L}}_1(\kaux)\big)\,\,
\overline{\mathcal{K}}\big(\overline{\Paux}\big)}{
\overline{\mathcal{L}}_1(\kaux)}
\,-
\\
&
\ \ \ \ \
-\,
\frac{1}{18}\,
\frac{\overline{\mathcal{K}}\big(
\overline{\mathcal{L}}_1(\kaux)\big)\,\,
\overline{\mathcal{L}}_1\big(
\overline{\mathcal{L}}_1(\kaux)\big)\,\,
\overline{\Paux}}{
\overline{\mathcal{L}}_1(\kaux)^2}
+
\frac{1}{6}\,
\overline{\mathcal{K}}\big(
\overline{\mathcal{L}}_1
\big(\overline{\Paux}\big)\big)
-
\frac{2}{9}\,
\overline{\Paux}\,
\overline{\mathcal{K}}\big(\overline{\Paux}\big),
\endaligned
\]
perform replacements using
Lemmas~{\ref{Lemma-K-bar-penetrates-L1-bar-k}} 
and~{\ref{Lemma-K-bar-k-K-bar-P}}:
\[
\!\!\!\!\!\!\!\!\!\!\!\!\!\!\!\!\!\!\!\!
\footnotesize
\aligned
\overline{\mathcal{K}}
\big(\Haux_0\big)
&
\,=\,
\frac{1}{2}\,
\frac{\overline{\mathcal{L}}_1(\overline{\kaux})\,\,
\overline{\mathcal{L}}_1\big(
\overline{\mathcal{L}}_1\big(
\overline{\mathcal{L}}_1(\kaux)\big)\big)}{
\overline{\mathcal{L}}_1(\kaux)}
+
\zero{
\frac{1}{2}\,
\frac{\overline{\mathcal{L}}_1\big(
\overline{\mathcal{L}}_1(\overline{\kaux})\big)\,\,
\overline{\mathcal{L}}_1\big(\overline{\mathcal{L}}_1(\kaux)\big)}{
\overline{\mathcal{L}}_1(\kaux)}}
+
\frac{1}{6}\,
\overline{\mathcal{L}}_1\big(
\overline{\mathcal{L}}_1\big(
\overline{\mathcal{L}}_1\big(\overline{\kaux}\big)\big)\big)
-
\frac{1}{6}\,
\frac{\overline{\mathcal{L}}_1(\overline{\kaux})\,\,
\overline{\mathcal{L}}_1\big(
\overline{\mathcal{L}}_1\big(
\overline{\mathcal{L}}_1(\kaux)\big)\big)}{
\overline{\mathcal{L}}_1(\kaux)}
\,-
\\
&
\ \ \ \ \
-\,
\frac{8}{9}\,
\frac{\overline{\mathcal{L}}_1\big(
\overline{\mathcal{L}}_1(\kaux)\big)^2\,\,
\overline{\mathcal{L}}_1(\overline{\kaux})}{
\overline{\mathcal{L}}_1(\kaux)^2}
-
\zero{
\frac{4}{9}\,
\frac{\overline{\mathcal{L}}_1\big(
\overline{\mathcal{L}}_1(\overline{\kaux})\big)\,\,
\overline{\mathcal{L}}_1\big(
\overline{\mathcal{L}}_1(\kaux)\big)}{
\overline{\mathcal{L}}_1(\kaux)}}
+
\frac{4}{9}\,
\frac{\overline{\mathcal{L}}_1(\overline{\kaux})\,\,
\overline{\mathcal{L}}_1\big(
\overline{\mathcal{L}}_1(\kaux)\big)^2}{
\overline{\mathcal{L}}_1(\kaux)^2}
\,-
\\
&
\ \ \ \ \
-\,
\frac{1}{9}\,
\frac{\overline{\mathcal{L}}_1\big(
\overline{\mathcal{L}}_1(\kaux)\big)\,\,
\overline{\mathcal{L}}_1(\overline{\kaux})\,\,
\overline{\Paux}}{
\overline{\mathcal{L}}_1(\kaux)}
-
\frac{1}{18}\,
\overline{\mathcal{L}}_1\big(
\overline{\mathcal{L}}_1\big(
\overline{\kaux}\big)\big)\,
\overline{\Paux}
-
\zerozero{
\frac{1}{18}\,
\frac{\overline{\mathcal{L}}_1\big(
\overline{\mathcal{L}}_1(\kaux)\big)\,\,
\overline{\mathcal{L}}_1(\overline{\kaux})\,\,
\overline{\Paux}}{
\overline{\mathcal{L}}_1(\kaux)}}
-
\zero{
\frac{1}{18}\,
\frac{\overline{\mathcal{L}}_1\big(
\overline{\mathcal{L}}_1(\kaux)\big)\,\,
\overline{\mathcal{L}}_1\big(
\overline{\mathcal{L}}_1(\overline{\kaux})\big)}{
\overline{\mathcal{L}}_1(\kaux)}}
\,+
\\
&
\ \ \ \ \
+
\zerozero{
\frac{1}{18}\,
\frac{\overline{\mathcal{L}}_1(\overline{\kaux})\,\,
\overline{\mathcal{L}}_1\big(
\overline{\mathcal{L}}_1(\kaux)\big)\,\,
\overline{\Paux}}{
\overline{\mathcal{L}}_1(\kaux)}}
+
\frac{1}{6}\,
\overline{\mathcal{K}}\big(
\overline{\mathcal{L}}_1\big(
\overline{\Paux}\big)\big)
+
\frac{2}{9}\,
\overline{\Paux}\,
\overline{\Paux}\,
\overline{\mathcal{L}}_1
\big(\overline{\kaux}\big)
+
\frac{2}{9}\,
\overline{\Paux}\,
\overline{\mathcal{L}}_1\big(
\overline{\mathcal{L}}_1\big(
\overline{\kaux}\big)\big)
\endaligned
\]
and observe some (underlined) cancellations to get an expression 
in which the last three terms must yet be transformed:
\[
\aligned
\overline{\mathcal{K}}
\big(
\Haux_0
\big)
&
\,=\,
\frac{1}{3}\,
\frac{\overline{\mathcal{L}}_1(\overline{\kaux})\,\,
\overline{\mathcal{L}}_1\big(
\overline{\mathcal{L}}_1\big(
\overline{\mathcal{L}}_1(\kaux)\big)\big)}{
\overline{\mathcal{L}}_1(\kaux)}
-
\frac{4}{9}\,
\frac{\overline{\mathcal{L}}_1(\overline{\kaux})\,\,
\overline{\mathcal{L}}_1\big(
\overline{\mathcal{L}}_1(\kaux)\big)^2}{
\overline{\mathcal{L}}_1(\kaux)^2}
-
\frac{1}{9}\,
\frac{\overline{\mathcal{L}}_1(\overline{\kaux})\,\,
\overline{\mathcal{L}}_1\big(
\overline{\mathcal{L}}_1(\kaux)\big)\,\,
\overline{\Paux}}{
\overline{\mathcal{L}}_1(\kaux)}
\,+
\\
&
\ \ \ \ \
+
\frac{2}{9}\,
\overline{\Paux}\,
\overline{\Paux}\,
\overline{\mathcal{L}}_1(\overline{\kaux})
+
\frac{1}{6}\,
\overline{\mathcal{L}}_1\big(
\overline{\mathcal{L}}_1\big(
\overline{\mathcal{L}}_1\big(
\overline{\kaux}\big)\big)\big)
+
\frac{1}{6}\,
\overline{\mathcal{L}}_1\big(
\overline{\mathcal{L}}_1
(\kaux)\big)\,
\overline{\Paux}
+
\frac{1}{6}\,
\overline{\mathcal{K}}\big(
\overline{\mathcal{L}}_1\big(
\overline{\Paux}\big)\big).
\endaligned
\]

\begin{Lemma}
One has:
\[
\overline{\mathcal{L}}_1\big(
\overline{\mathcal{L}}_1\big(
\overline{\mathcal{L}}_1\big(
\overline{\kaux}\big)\big)\big)
+
\overline{\mathcal{L}}_1\big(
\overline{\mathcal{L}}_1
(\kaux)\big)\,
\overline{\Paux}
+
\overline{\mathcal{K}}\big(
\overline{\mathcal{L}}_1\big(
\overline{\Paux}\big)\big)
\,\,=\,\,
-\,
2\,
\overline{\mathcal{L}}_1
\big(\overline{\kaux}\big)\,
\overline{\mathcal{L}}_1\big(
\overline{\Paux}\big).
\]
\end{Lemma}

\proof
Apply the vector field $\overline{\mathcal{L}}_1$ to
Lemma~{\ref{Lemma-K-bar-k-K-bar-P}}:
\[
\overline{\mathcal{L}}_1\big(
\overline{\mathcal{K}}\big(
\overline{\Paux}\big)\big)
\,\,=\,\,
-\,
\overline{\mathcal{L}}_1
\big(\overline{\Paux}\big)\,
\overline{\mathcal{L}}_1\big(\overline{\kaux}\big)
-
\overline{\Paux}\,
\overline{\mathcal{L}}_1\big(
\overline{\mathcal{L}}_1\big(
\overline{\kaux}\big)\big)
-
\overline{\mathcal{L}}_1\big(
\overline{\mathcal{L}}_1\big(
\overline{\mathcal{L}}_1\big(
\overline{\kaux}\big)\big)\big).
\]
On the other hand, apply the Lie bracket $\big[
\overline{\mathcal{L}}_1, \overline{\mathcal{K}} \big]
(\centersmallbullet)$ to the function $\overline{\Paux}$,
using the concerned known commutation relation
shown in 
Section~{\ref{local-geometry-hypersurfaces-C-3}}:
\[
\overline{\mathcal{L}}_1\big(
\overline{\mathcal{K}}\big(
\overline{\Paux}\big)\big)
-
\overline{\mathcal{K}}\big(
\overline{\mathcal{L}}_1\big(
\overline{\Paux}\big)\big)
\,\,=\,\,
\big[
\overline{\mathcal{L}}_1,
\overline{\mathcal{K}}
\big]
\big(
\overline{\Paux}
\big)
\,\,=\,\,
\overline{\mathcal{L}}_1
\big(\overline{\kaux}\big)\,
\overline{\mathcal{L}}_1
\big(\overline{\Paux}\big),
\]
and replace the first term 
$\overline{\mathcal{L}}_1\big(
\overline{\mathcal{K}}\big(
\overline{\Paux}\big)\big)$ by its value above to 
get the result.
\endproof

Consequently, after this transformation, we see that
$\overline{\mathcal{K}} \big( \Haux_0 \big)$ 
is a multiple of $\overline{\mathcal{L}}_1
(\overline{\kaux})$ in which we recognize
$-2\, \Haux_0$ as stated:
\begin{align}
\overline{\mathcal{K}}
\big(\Haux_0\big)
&
\,=\,
\overline{\mathcal{L}}_1
\big(\overline{\kaux}\big)
\left(
\frac{1}{3}\,
\frac{\overline{\mathcal{L}}_1\big(
\overline{\mathcal{L}}_1\big(
\overline{\mathcal{L}}_1(\kaux)\big)\big)}{
\overline{\mathcal{L}}_1(\kaux)}
-
\frac{4}{9}\,
\frac{
\overline{\mathcal{L}}_1\big(\overline{\mathcal{L}}_1(\kaux)\big)^2}{
\overline{\mathcal{L}}_1(\kaux)^2}
-
\frac{1}{9}\,
\frac{\overline{\mathcal{L}}_1\big(\overline{\mathcal{L}}_1(\kaux)\big)\,
\overline{\Paux}}{
\overline{\mathcal{L}}_1(\kaux)}
\,+
\right.
\notag
\\
&
\left.
\ \ \ \ \ \ \ \ \ \ \ \ \ \ \ \ \ \ \ \ \ \ \ \ \ \ \ \ \ \ \ \ \ \ \ 
\ \ \ \ \ \ \ \ \ \ \ \ \ \ \ \ \ \ \ \ \ \ \ \ \ \ \ \ \ \ \ \ \ \ \
\ \ \ \ \ \ \  
-
\frac{1}{3}\,
\overline{\mathcal{L}_1}\big(\overline{\Paux}\big)
+
\frac{2}{9}\,
\overline{\Paux}^2
\right).
\qedhere
\end{align}
\endproof

As we already observed, the essential (invariant)
torsion $\isqrt\, K^3 - Z^6$ can be set
$0$ to normalize the group parameter ${\sf d}$ as:
\[
{\sf d}
\,:=\,
-\,\frac{\isqrt}{2}\,
\frac{\overline{\sf c}{\sf e}{\sf e}}{{\sf c}}
+
\isqrt\,
\frac{{\sf c}}{\overline{\sf c}}\,
\Haux_0,
\]
whence inserting in~({\ref{3-loop-initial}}):
\[
\left(\!
\begin{array}{c}
\rho
\\
\kappa
\\
\zeta
\end{array}
\!\right)
\,:=\,
\left(\!
\begin{array}{ccc}
{\sf c}\overline{\sf c} & 0 & 0
\\
-\isqrt\,\overline{\sf c}{\sf e} & {\sf c} & 0
\\
-\frac{\isqrt}{2}\,
\frac{\overline{\sf c}{\sf e}{\sf e}}{{\sf c}}
+
\isqrt\,\frac{{\sf c}}{\overline{\sf c}}\,
\Haux_0
& 
{\sf e} & \frac{{\sf c}}{\overline{\sf c}}
\end{array}
\!\right)
\left(\!
\begin{array}{c}
\rho_0
\\
\kappa_0'
\\
\zeta_0'
\end{array}
\!\right).
\]
Thus, we are naturally led to change the initial coframe on $M$:
\[
\big\{
\rho_0,\,
\kappa_0',\,
\zeta_0',\,
\overline{\kappa}_0',\,
\overline{\zeta}_0'
\big\}
\ \ \ \ \
\leadsto
\ \ \ \ \
\big\{
\rho_0,\,
\kappa_0',\,
\zeta_0'',\,
\overline{\kappa}_0',\,
\overline{\zeta}_0''
\big\},
\]
by introducing the new $1$-form:
\[
\zeta_0''
\,:=\,
\zeta_0'
+
\isqrt\,
\Haux_0\,
\rho_0,
\]
so that a new, reduced by two real dimensions,
$G$-structure, appears:
\[
\left(\!
\begin{array}{c}
\rho
\\
\kappa
\\
\zeta
\end{array}
\!\right)
\,:=\,
\left(\!
\begin{array}{ccc}
{\sf c}\overline{\sf c} & 0 & 0
\\
-\isqrt\,\overline{\sf c}{\sf e} & {\sf c} & 0
\\
-\frac{\isqrt}{2}\,
\frac{\overline{\sf c}{\sf e}{\sf e}}{{\sf c}}
& 
{\sf e} & \frac{{\sf c}}{\overline{\sf c}}
\end{array}
\!\right)
\left(\!
\begin{array}{c}
\rho_0
\\
\kappa_0'
\\
\zeta_0''
\end{array}
\!\right),
\]
which is justified by the computation\big/reorganization:
\[
\aligned
\zeta
&
\,=\,
\Big(
-\frac{\isqrt}{2}\,
\frac{\overline{\sf c}{\sf e}{\sf e}}{{\sf c}}
+
\isqrt\,
\frac{{\sf c}}{\overline{\sf c}}\,
\Haux_0
\Big)\,\rho_0
+
{\sf e}\,\kappa_0'
+
\frac{{\sf c}}{\overline{\sf c}}\,
\zeta_0'
\\
&
\,=\,
-\,\frac{\isqrt}{2}\,
\frac{\overline{\sf c}{\sf e}{\sf e}}{{\sf c}}\,
\rho_0
+
{\sf e}\,\kappa_0'
+
\frac{{\sf c}}{\overline{\sf c}}\,
\Big(
\underbrace{
\zeta_0'
+
\isqrt\,\Haux_0\,\rho_0}_{=:\,\,\zeta_0''}
\Big).
\endaligned
\]
Back to previous expressions, this last coframe writes out as:
\[
\aligned
\rho_0
&
\,:=\,
\frac{1}{\ell}\,
\Big(
dv
-
\Aaux^1\,dz_1
-
\Aaux^2\,dz_2
-
\overline{\Aaux}^1\,
d\overline{z}_1
-
\overline{\Aaux}^2\,
d\overline{z}_2
\Big),
\\
\kappa_0'
&
\,:=\,
dz_1
-
\kaux\,dz_2
+
\frac{\isqrt}{3}\,
\Baux_0\,\rho_0,
\\
\zeta_0''
&
\,:=\,
\overline{\mathcal{L}}_1(\kaux)\,
dz_2
+
\isqrt\,\Haux_0\,\rho_0.
\endaligned
\]

\Section{\bf Darboux-Cartan Structure of the Coframe
$\big\{ \rho_0, \kappa_0', \zeta_0'', \overline{\kappa}_0', 
\overline{\zeta}_0'' \big\}$}
\label{D-C-structure-kappa-0-prime-zeta-0-prime-prime}
\HEAD{{\ref{D-C-structure-kappa-0-prime-zeta-0-prime-prime}}.~{\sf 
Darboux-Cartan Structure of the Coframe
$\big\{ \rho_0, \kappa_0', \zeta_0'', \overline{\kappa}_0', 
\overline{\zeta}_0'' \big\}$}
}{
Wei Guo {\sc Foo} (Beijing) and Joël {\sc Merker} (Orsay)}

The present change of initial coframe expresses as:
\[
\zeta_0''
\,:=\,
\zeta_0'
+
\isqrt\,
\Haux_0\,\rho_0
\ \ \ \ \ \ \ \ \ \ \ \ \ \ \ \ \ \
\Longleftrightarrow
\ \ \ \ \ \ \ \ \ \ \ \ \ \ \ \ \ \
\zeta_0'
\,=\,
\zeta_0''
-
\isqrt\,
\Haux_0.
\]
The exterior differentiation of $\zeta_0''$ comprises $3$ terms
that we shall compute soon:
\[
d\zeta_0''
\,=\,
d\zeta_0'
+
\isqrt\,d\Haux_0
\wedge\rho_0
+
\isqrt\,\Haux_0\,
d\rho_0.
\]

Back to the previous structure equations written in the
abbreviated form~({\ref{inexplicit-0-structure-3-loop}}),
we may start by replacing $\zeta_0'$ in $d\rho_0$, 
while observing that:
\[
\rho_0\wedge\zeta_0'
\,=\,
\rho_0\wedge\zeta_0''
\ \ \ \ \ \ \ \ \ \ \ \ \ \ \ \ \ \
\text{and}
\ \ \ \ \ \ \ \ \ \ \ \ \ \ \ \ \ \
\rho_0\wedge\overline{\zeta}_0'
\,=\,
\rho_0\wedge\overline{\zeta}_0'',
\]
we come to unchanged coefficients for:
\[
d\rho_0
\,=\,
\Raux_0^{1\prime}\,
\rho_0\wedge\kappa_0'
+
\Raux_0^{2\prime}\,
\rho_0\wedge\zeta_0''
+
\overline{\Raux}_0^{1\prime}\,
\rho_0\wedge\overline{\kappa}_0'
+
\overline{\Raux}_0^{2\prime}\,
\rho_0\wedge\overline{\zeta}_0''
+
\isqrt\,
\kappa_0'\wedge\overline{\kappa}_0',
\]
hence without computation, the third term is:
\[
\isqrt\,\Haux_0\,
d\rho_0
\,=\,
\isqrt\,\Haux_0\,
\Raux_0^{1\prime}\,
\rho_0\wedge\kappa_0'
+
\isqrt\,\Haux_0\,
\Raux_0^{2\prime}\,
\rho_0\wedge\zeta_0'
+
\isqrt\,\Haux_0\,
\overline{\Raux}_0^{1\prime}\,
\rho_0\wedge\overline{\kappa}_0'
+
\isqrt\,\Haux_0\,
\overline{\Raux}_0^{2\prime}\,
\rho_0\wedge\overline{\zeta}_0'
-
\Haux_0\,
\kappa_0\wedge\overline{\kappa}_0'.
\]

Next, we do the same replacement of $\zeta_0'$ in:
\[
\aligned
d\kappa_0'
&
\,=\,
\Kaux_0^{1\prime}\,
\rho_0\wedge\kappa_0'
+
\Kaux_0^{2\prime}\,
\rho_0
\wedge
\Big(
\zeta_0''
-
\isqrt\,\Haux_0\,\rho_0
\Big)
+
\Kaux_0^{3\prime}\,
\rho_0\wedge\overline{\kappa}_0'
\,+
\\
&
\ \ \ \ \
+
\Kaux_0^{5\prime}\,
\kappa_0'
\wedge
\Big(
\zeta_0''
-
\isqrt\,\Haux_0\,\rho_0
\Big)
+
\Kaux_0^{6\prime}\,
\kappa_0'
\wedge
\overline{\kappa}_0'
+
\Big(
\zeta_0''
-
\isqrt\,\Haux_0\,\rho_0
\Big)
\wedge
\overline{\kappa}_0',
\endaligned
\]
hence:
\[
\aligned
d\kappa_0'
&
\,=\,
\Big(
\underbrace{
\Kaux_0^{1\prime}
+
\isqrt\,
\Kaux_0^{5\prime}\,
\Haux_0}_{=:\,\,\Kaux_0^{1\prime\prime}}
\Big)\,
\rho_0\wedge\kappa_0'
+
\Kaux_0^{2\prime}\,
\rho_0\wedge\zeta_0''
+
\Big(
\underbrace{
\Kaux_0^{3\prime}
-
\isqrt\,\Haux_0}_{=:\,\,\Kaux_0^{3\prime}}
\Big)\,
\rho_0\wedge\overline{\kappa}_0'
\,+
\\
&
\ \ \ \ \
+
\Kaux_0^{5\prime}\,
\kappa_0'
\wedge
\zeta_0''
+
\Kaux_0^{6\prime}\,
\kappa_0'
\wedge
\overline{\kappa}_0'
+
\zeta_0''
\wedge
\overline{\kappa}_0'.
\endaligned
\]

Similarly, do the same for:
\[
\aligned
d\zeta_0'
&
\,=\,
\Zaux_0^{2\prime}\,
\rho_0
\wedge
\Big(
\zeta_0''
-
\isqrt\,
\Haux_0\,\rho_0
\Big)
+
\Zaux_0^{5\prime}\,
\kappa_0'
\wedge
\Big(
\overline{\zeta}_0''
+
\isqrt\,\Haux_0\,
\rho_0
\Big)
\,+
\\
&
\ \ \ \ \ 
+
\Zaux_0^{8\prime}\,
\Big(
\zeta_0''
-
\isqrt\,\Haux_0\,
\rho_0
\Big)
\wedge
\overline{\kappa}_0'
+
\Zaux_0^{9\prime}\,
\Big(
\zeta_0''
-
\isqrt\,\Haux_0\,\rho_0
\Big)
\wedge
\Big(
\overline{\zeta}_0''
+
\isqrt\,\overline{\Haux}_0\,\rho_0
\Big),
\endaligned
\]
hence:
\[
\aligned
d\zeta_0'
&
\,=\,
\isqrt\,\Zaux_0^{5\prime}\,
\Haux_0\,
\rho_0\wedge\kappa_0'
+
\Big(
\Zaux_0^{2\prime}
-
\isqrt\,\Zaux_0^{9\prime}\,
\overline{\Haux}_0
\Big)\,
\rho_0\wedge\zeta_0''
-
\isqrt\,\Zaux_0^{8\prime}\,
\Haux_0\,
\rho_0\wedge\overline{\kappa}_0'
\,-
\\
&
\ \ \ \ \
-\,\isqrt\,
\Zaux_0^{9\prime}\,
\Haux_0\,
\rho_0\wedge\overline{\zeta}_0''
+
\Zaux_0^{5\prime}\,
\kappa_0'\wedge\zeta_0''
+
\Zaux_0^{8\prime}\,
\zeta_0''\wedge\overline{\kappa}_0'
+
\Zaux_0^{9\prime}\,
\zeta_0''\wedge\overline{\zeta}_0''.
\endaligned
\]

Next, we have to compute the second term in $d\zeta_0''$, and using:
\[
d\Haux_0
\,=\,
\mathcal{T}\big(\Haux_0\big)\,
\rho_0
+
\mathcal{L}_1\big(\Haux_0\big)\,
\kappa_0
+
\mathcal{K}\big(\Haux_0\big)\,
\zeta_0
+
\overline{\mathcal{L}}_1\big(\Haux_0\big)\,
\overline{\kappa}_0
+
\overline{\mathcal{K}}\big(\Haux_0\big)\,
\overline{\zeta}_0,
\]
it comes:
\[
\!\!\!\!\!\!\!\!\!\!\!\!\!\!\!
\aligned
d\Haux_0
\wedge
\rho_0
&
\,=\,
0
-
\mathcal{L}_1\big(\Haux_0\big)\,
\rho_0\wedge\kappa_0
-
\mathcal{K}\big(\Haux_0\big)\,
\rho_0\wedge\zeta_0
-
\overline{\mathcal{L}}_1\big(\Haux_0\big)\,
\rho_0\wedge\overline{\kappa}_0
-
\overline{\mathcal{K}}\big(\Haux_0\big)\,
\rho_0\wedge\overline{\zeta}_0
\\
&
\,=\,
-\,\mathcal{L}_1\big(\Haux_0\big)\,
\rho_0\wedge
\Big(
\kappa_0'
-
\frac{\isqrt}{3}\,
\Baux_0\,\rho_0
\Big)
-
\mathcal{K}\big(\Haux_0\big)\,
\rho_0\wedge\frac{\zeta_0'}{\overline{\mathcal{L}}_1(\kaux)}
-
\overline{\mathcal{L}}_1\big(\Haux_0\big)\,
\rho_0\wedge
\Big(
\overline{\kappa}_0'
+
\frac{\isqrt}{3}\,
\overline{\Baux}_0\,\rho_0
\Big)
\,-
\\
&
\ \ \ \ \ \ \ \ \ \ \ \ \ \ \ \ \ \ \ \ \ \ \ \ \ \
\ \ \ \ \ \ \ \ \ \ \ \ \ \ \ \ \ \ \ \ \ \ \ \ \ \
\ \ \ \ \ \ \ \ \ \ \ \ \ \ \ \ \ \ \ \ \ \ \ \ \ \
\ \ \ \ \ \ \ \ \ \ \ \ \ \ \ \ \ \ \ \ \ \ \ \ \ \
-\,
\overline{\mathcal{K}}\big(\Haux_0\big)\,
\rho_0
\wedge
\frac{\overline{\zeta}_0'}{\mathcal{L}_1(\overline{\kaux})}
\\
&
\,=\,
-\,\mathcal{L}_1\big(\Haux_0\big)\,
\rho_0\wedge\kappa_0'
-
\frac{\mathcal{K}(\Haux_0)}{\overline{\mathcal{L}}_1(\kaux)}\,
\rho_0\wedge\zeta_0'
-
\overline{\mathcal{L}}_1\big(\Haux_0\big)\,
\rho_0\wedge\overline{\kappa}_0'
-
\frac{\overline{\mathcal{K}}(\Haux_0)}{\mathcal{L}_1(\overline{\kaux})}\,
\rho_0\wedge\overline{\zeta}_0',
\endaligned
\]
hence multipliying by $\isqrt$, we get the expression of the second
term:
\[
\isqrt\,
d\Haux_0
\wedge\rho_0
\,=\,
-\,\isqrt\,
\mathcal{L}_1\big(\Haux_0\big)\,
\rho_0\wedge\kappa_0'
-
\isqrt\,
\frac{\mathcal{K}(\Haux_0)}{\overline{\mathcal{L}}_1(\kaux)}\,
\rho_0\wedge\zeta_0''
-
\isqrt\,
\overline{\mathcal{L}}_1\big(\Haux_0\big)\,
\rho_0\wedge\overline{\kappa}_0'
-
\isqrt\,
\frac{\overline{\mathcal{K}}(\Haux_0)}{\mathcal{L}_1(\overline{\kaux})}\,
\rho_0\wedge\overline{\zeta}_0''.
\]

Summing and collecting the three computed terms yields:
\[
\!\!\!\!\!\!\!\!\!\!\!\!\!\!\!
\footnotesize
\aligned
d\zeta_0''
&
\,=\,
\rho_0\wedge\kappa_0'
\Big[
\underbrace{
\isqrt\,\Zaux_0^{5\prime}\,\Haux_0
-
\isqrt\,\mathcal{L}_1\big(\Haux_0\big)
+
\isqrt\,\Haux_0\,\Raux_0^{1\prime}}_{=:\,\,\Zaux_0^{1\prime\prime}}
\Big]
+
\rho_0\wedge\zeta_0''
\Big[
\underbrace{
\Zaux_0^{2\prime}
-
\isqrt\,\Zaux_0^{9\prime}\,
\overline{\Haux}_0
-
\isqrt\,
\frac{\mathcal{K}(\Haux_0)}{\overline{\mathcal{L}}_1(\kaux)}
+
\isqrt\,\Haux_0\,
\Raux_0^{2\prime}}_{=:\,\,\Zaux_0^{2\prime\prime}}
\Big]
\,+
\\
&
\ \ \ \ \
+
\rho_0\wedge\overline{\kappa}_0'
\Big[
\underbrace{
-\,\isqrt\,\Zaux_0^{8\prime}\,\Haux_0
-
\isqrt\,\overline{\mathcal{L}}_1\big(\Haux_0\big)
+
\isqrt\,\Haux_0\,
\overline{\Raux}_0^{1\prime}}_{=:\,\,\Zaux_0^{3\prime\prime}}
\Big]
+
\rho_0\wedge\overline{\zeta}_0''
\Big[
\zero{
-\,\isqrt\,\Zaux_0^{9\prime}\,
\Haux_0
-
\isqrt\,\frac{\overline{\mathcal{K}}(\Haux_0)}{
\mathcal{L}_1(\overline{\kaux})}
+
\isqrt\,
\Haux_0\,
\overline{\Raux}_0^{2\prime}}
\Big]
\,+
\\
&
\ \ \ \ \
+
\Zaux_0^{5\prime}\,
\kappa_0'\wedge\zeta_0''
+
\kappa_0'\wedge\overline{\kappa}_0'
\Big[
\underbrace{
-\,\Haux_0}_{=:\,\,\Zaux_0^{6\prime\prime}}
\Big]
+
\Zaux_0^{8\prime}\,
\zeta_0''\wedge\overline{\kappa}_0'
+
\Zaux_0^{9\prime}\,
\zeta_0''\wedge\overline{\zeta}_0''.
\endaligned
\]

\begin{Lemma}
\label{Lemma-Z-0-4-prime-prime}
One has the identical vanishing of the coefficient
of $\rho_0 \wedge \overline{\zeta}_0''$ in $d \zeta_0''$:
\[
\aligned
\Zaux_0^{4\prime\prime}
\,:=\,
&\,
-\,\isqrt\,
\Zaux_0^{9\prime}\,
\Haux_0
-
\isqrt\,
\frac{\overline{\mathcal{K}}(\Haux_0)}{
\mathcal{L}_1(\overline{\kaux})}
+
\isqrt\,
\Haux_0\,
\overline{\Raux}_0^{2\prime}
\\
\,\equiv\,
&\,
0.
\endaligned
\]
\end{Lemma}

\proof
This is equivalent to:
\[
\overline{\mathcal{K}}
\big(\Haux_0\big)
\overset{\text{\bf ?}}{\,\equiv\,}
\mathcal{L}_1\big(\overline{\kaux}\big)\,
\Haux_0\,
\Big(
-\,\Zaux_0^{9\prime}
+
\overline{\Raux}_0^{2\prime}
\Big),
\]
and after a replacement using
Proposition~{\ref{Proposition-3-loop-structure-0-coframe}}, to:
\[
\overline{\mathcal{K}}
\big(\Haux_0\big)
\overset{\text{\bf ?}}{\,\equiv\,}
\mathcal{L}_1\big(\overline{\kaux}\big)\,
\Haux_0\,
\left(
-\,
\frac{\overline{\mathcal{L}}_1(\overline{\kaux})}{
\mathcal{L}_1(\overline{\kaux})}
-
\frac{\overline{\mathcal{L}}_1(\overline{\kaux})}{
\mathcal{L}_1(\overline{\kaux})}
\right),
\]
an identity which was already seen by
Lemma~{\ref{Lemma-K-bar-H-0}}.
\endproof

In summary:
\[
\aligned
d\rho_0
&
\,=\,
\Raux_0^{1\prime}\,
\rho_0\wedge\kappa_0'
+
\Raux_0^{2\prime}\,
\rho_0\wedge\zeta_0''
+
\overline{\Raux}_0^{1\prime}\,
\rho_0\wedge\overline{\kappa}_0'
+
\overline{\Raux}_0^{2\prime}\,
\rho_0\wedge\overline{\zeta}_0''
+
\isqrt\,
\kappa_0'\wedge\overline{\kappa}_0',
\\
d\kappa_0'
&
\,=\,
\Kaux_0^{1\prime\prime}\,
\rho_0\wedge\kappa_0'
+
\Kaux_0^{2\prime}\,
\rho_0\wedge\zeta_0''
+
\Kaux_0^{3\prime\prime}\,
\rho_0\wedge\overline{\kappa}_0'
\,+
\\
&
\ \ \ \ \
+
\Kaux_0^{5\prime}\,
\kappa_0'\wedge\zeta_0''
+
\Kaux_0^{6\prime}\,
\kappa_0'\wedge\overline{\kappa}_0'
+
\zeta_0''\wedge\overline{\kappa}_0',
\\
d\zeta_0''
&
\,=\,
\Zaux_0^{1\prime\prime}\,
\rho_0\wedge\kappa_0'
+
\Zaux_0^{2\prime\prime}\,
\rho_0\wedge\zeta_0''
+
\Zaux_0^{3\prime\prime}\,
\rho_0\wedge\overline{\kappa}_0'
\,+
\\
&
\ \ \ \ \
+
\Zaux_0^{5\prime}\,
\kappa_0'\wedge\zeta_0''
+
\Zaux_0^{6\prime\prime}\,
\kappa_0'\wedge\overline{\kappa}_0'
+
\Zaux_0^{8\prime}\,
\zeta_0''\wedge\overline{\kappa}_0'
+
\Zaux_0^{9\prime}\,
\zeta_0''\wedge\overline{\zeta}_0''.
\endaligned
\]
Notice that new coefficients 
$\Zaux_0^{2\prime\prime}$, 
$\Zaux_0^{3\prime\prime}$,
$\Zaux_0^{4\prime\prime}$
appear in $d\zeta_0''$, which were absent
in $d\zeta_0'$, as they are coming from the second
term $\isqrt\, d\Haux_0 \wedge\rho_0$.

\Section{\bf Absorption and apparition of two $1$-forms $\pi^1$,
$\pi^2$}
\label{absorption-pi-1-pi-2}
\HEAD{{\ref{absorption-pi-1-pi-2}}.~{\sf Absorption and 
apparition of two $1$-forms $\pi^1$,
$\pi^2$}
}{
Wei Guo {\sc Foo} (Beijing) and Joël {\sc Merker} (Orsay)}

With the $4$-dimensional group parametrized by $\big({\sf c},
\overline{\sf c}, {\sf e}, \overline{\sf e} \big)$,
the lifted coframe writes:
\[
\left(\!
\begin{array}{c}
\rho
\\
\kappa
\\
\zeta
\end{array}
\!\right)
\,:=\,
\left(\!
\begin{array}{ccc}
{\sf c}\overline{\sf c} & 0 & 0
\\
-\isqrt\,\overline{\sf c}{\sf e} & {\sf c} & 0
\\
-\frac{\isqrt}{2}\,
\frac{\overline{\sf c}{\sf e}{\sf e}}{{\sf c}}
& {\sf e} & \frac{{\sf c}}{\overline{\sf c}}
\end{array}
\!\right)
\left(\!
\begin{array}{c}
\rho_0
\\
\kappa_0'
\\
\zeta_0''
\end{array}
\!\right)
\ \ \ \ \ \ \ \ \ \ \ \ \ \ \ \ \ \
\Longleftrightarrow
\ \ \ \ \ \ \ \ \ \ \ \ \ \ \ \ \ \
\left\{
\aligned
\rho
&
\,:=\,
{\sf c}\overline{\sf c}\,
\rho_0,
\\
\kappa
&
\,:=\,
-\,\isqrt\,\overline{\sf c}{\sf e}\,
\rho_0
+
{\sf c}\,\kappa_0',
\\
\zeta
&
\,:=\,
-\,\frac{\isqrt}{2}\,
\frac{\overline{\sf c}{\sf e}{\sf e}}{{\sf c}}\,
\rho_0
+
{\sf e}\,\kappa_0'
+
\frac{{\sf c}}{\overline{\sf c}}\,
\zeta_0'',
\endaligned\right.
\]
with inverse formulas:
\leqnomode\usetagform{default}
\begin{align}
\label{4-loop-initial-inverted}
\rho_0
&
\,=\,
\frac{1}{{\sf c}\overline{\sf c}}\,
\rho,
\notag
\\
\kappa_0'
&
\,=\,
\isqrt\,\frac{{\sf e}}{{\sf c}{\sf c}}\,
\rho
+
\frac{1}{{\sf c}}\,
\kappa,
\\
\zeta_0''
&
\,=\,
-\,\frac{\isqrt}{2}\,
\frac{\overline{\sf c}{\sf e}{\sf e}}{
{\sf c}{\sf c}{\sf c}}\,
\rho
-
\frac{\overline{\sf c}{\sf e}}{{\sf c}{\sf c}}\,
\kappa
+
\frac{\overline{\sf c}}{\sf c}\,
\zeta.
\notag
\end{align}
The Maurer-Cartan matrix becomes:
\[
\aligned
dg\cdot g^{-1}
&
\,=\,
\left(\!
\begin{array}{ccc}
\overline{\sf c}\,d{\sf c}+{\sf c}d\overline{\sf c} & 0 & 0
\\
-\isqrt\,{\sf e}d\overline{\sf c}
-\isqrt\,\overline{\sf c}d{\sf e} & d{\sf c} & 0
\\
-\frac{\isqrt}{2}
\frac{{\sf e}{\sf e}\,d\overline{\sf c}}{{\sf c}}
-\isqrt\,
\frac{\overline{\sf c}{\sf e}\,d{\sf e}}{{\sf c}}
+
\frac{\isqrt}{2}
\frac{\overline{\sf c}{\sf e}{\sf e}\,d{\sf c}}{
{\sf c}{\sf c}}
& d{\sf e} & \frac{d{\sf c}}{\overline{\sf c}}
-\frac{{\sf c}\,d\overline{\sf c}}{\overline{\sf c}\overline{\sf c}}
\end{array}
\!\right)
\left(\!
\begin{array}{ccc}
\frac{1}{{\sf c}\overline{\sf c}} & 0 & 0
\\
\isqrt\,\frac{{\sf e}}{{\sf c}{\sf c}} & \frac{1}{{\sf c}} & 0
\\
-\frac{\isqrt}{2}
\frac{\overline{\sf c}{\sf e}{\sf e}}{
{\sf c}{\sf c}{\sf c}}
& -\frac{\overline{\sf c}{\sf e}}{{\sf c}{\sf c}} & 
\frac{\overline{\sf c}}{{\sf c}}
\end{array}
\!\right)
\\
&
\,=:\,
\left(\!
\begin{array}{ccc}
\alpha+\overline{\alpha} & 0 & 0
\\
\beta & \alpha & 0
\\
0 & \isqrt\,\beta & \alpha-\overline{\alpha}
\end{array}
\!\right),
\endaligned
\]
in terms of the group-invariant $1$-forms:
\[
\aligned
\alpha
&
\,:=\,
\frac{d{\sf c}}{{\sf c}},
\\
\beta
&
\,:=\,
\isqrt\,
\frac{{\sf e}\,d{\sf c}}{{\sf c}{\sf c}}
-
\isqrt\,\frac{{\sf e}\,d\overline{\sf c}}{{\sf c}\overline{\sf c}}
-
\isqrt\,\frac{d{\sf e}}{{\sf c}}.
\endaligned
\]

Now, if we exterior-differentiate the lifted coframe on the product
manifold equipped with coordinates:
\[
\big( 
z_1,z_2,\overline{z}_1,\overline{z}_2,v
\big)
\times
\big(
{\sf c},\overline{\sf c},{\sf e},\overline{\sf e}
\big)
\,\,\in\,\,
M^5
\times
G^4,
\]
after hard computations, we may come to structure equations of the
abstract shape:
\[
\aligned
d\rho
&
\,=\,
\big(\alpha+\overline{\alpha}\big)
\wedge
\rho
\,+
\\
&
\ \ \ \ \ 
+
R^1\,
\rho\wedge\kappa
+
R^2\,
\rho\wedge\zeta
+
\overline{R}^1\,
\rho\wedge\overline{\kappa}
+
\overline{R}^2\,
\rho\wedge\overline{\zeta}
+
\isqrt\,
\kappa\wedge\overline{\kappa},
\\
d\kappa
&
\,=\,
\beta\wedge\rho
+
\alpha\wedge\kappa
\,+
\\
&
\ \ \ \ \ 
+
K^1\,
\rho\wedge\kappa
+
K^2\,\rho\wedge\zeta
+
K^3\,
\rho\wedge\overline{\kappa}
+
K^4\,
\rho\wedge\overline{\zeta}
\,+
\\
&
\ \ \ \ \
+
K^5\,
\kappa\wedge\zeta
+
K^6\,
\kappa\wedge\overline{\kappa}
+
\zeta\wedge\overline{\kappa},
\\
d\zeta
&
\,=\,
\gamma\wedge\rho
+
\isqrt\,\beta\wedge\kappa
+
\big(
\alpha-\overline{\alpha}
\big)
\wedge\zeta
\,+
\\
&
\ \ \ \ \
+
Z^1\,
\rho\wedge\kappa
+
Z^2\,
\rho\wedge\zeta
+
Z^3\,
\rho\wedge\overline{\kappa}
+
Z^4\,
\rho\wedge\overline{\zeta}
\,+
\\
&
\ \ \ \ \
+
Z^5\,
\kappa\wedge\zeta
+
Z^6\,
\kappa\wedge\overline{\kappa}
+
Z^7\,
\kappa\wedge\overline{\zeta}
+
Z^8\,
\zeta\wedge\overline{\kappa}
+
Z^9\,
\zeta\wedge\overline{\zeta}.
\endaligned
\]

A moment of reflection convinces of the truth of 

\begin{Assertion}
\label{Assertion-relations-R-K-Z}
The relations coming from the normalizations
of the group parameters ${\sf f}$, ${\sf b}$, ${\sf c}$ are preserved:
\[
\aligned
1
&
\,=\,
\big[
\zeta\wedge\overline{\kappa}
\big]
\big\{
d\kappa
\big\},
\\
0
&
\,=\,
\overline{R}^1
-
2\,K^6
+
Z^8,
\\
0
&
\,=\,
\isqrt\,K^3
-
Z^6,
\endaligned
\]
as well as the auxiliary relations:
\reqnomode\usetagform{EngelLie}
\begin{align}
K^5
&
\,=\,
R^2,
\notag
\\
Z^7
&
\,=\,
\isqrt\,
K^4,
\notag
\\
Z^9
&
\,=\,
-\,\overline{R}^2.
\tag{\qed}
\end{align}
\end{Assertion}

Now, we want to {\sl absorb} as many as possible of 
these torsion coefficients. So we introduce {\sl modified
Maurer-Cartan forms}\,\,---\,\,mind notations:
\[
\aligned
\pi^1
&
\,:=\,
\alpha
-
a_1\,\rho
-
a_2\,\kappa
-
a_3\,\zeta
-
a_4\,\overline{\kappa}
-
a_5\,\overline{\zeta},
\\
\pi^2
&
\,:=\,
\beta
-
b_1\,\rho
-
b_2\,\kappa
-
b_3\,\zeta
-
b_4\,\overline{\kappa}
-
b_5\,\overline{\zeta},
\endaligned
\]
and we try to determine (fix) the unknown coefficients $a_i$, $b_i$.
By replacement, setting $c_i := 0$ in the formula seen above
for $d\zeta$, we obtain without pain:
\[
\aligned
d\rho
&
\,=\,
\big(\pi^1+\overline{\pi}^1\big)
\,+
\\
&
\ \ \ \ \
+
\rho\wedge\kappa
\Big(
R^1
-
a_2
-
\overline{a}_4
\Big)
+
\rho\wedge\zeta
\Big(
R^2
-
a_3
-
\overline{a}_5
\Big)
+
\rho\wedge\overline{\kappa}
\Big(
\overline{R}^1
-
a_4
-
\overline{a}_2
\Big)
\,+
\\
&
\ \ \ \ \ 
+
\rho\wedge\overline{\zeta}
\Big(
\overline{R}^2
-
a_5
-
\overline{a}_3
\Big)
+
\isqrt\,
\kappa\wedge\overline{\kappa},
\\
d\kappa
&
\,=\,
\pi^2\wedge\rho
+
\pi^1\wedge\kappa
\,+
\\
&
\ \ \ \ \ 
+
\rho\wedge\kappa\,
\Big(
K^1
+
a_1
-
b_2
\Big)
+
\rho\wedge\zeta\,
\Big(
K^2
-
b_3
\Big)
+
\rho\wedge\overline{\kappa}\,
\Big(
K^3
-
b_4
\Big)
+
\rho\wedge\overline{\zeta}\,
\Big(
K^4
-
b_5
\Big)
\,+
\\
&
\ \ \ \ \
+
\kappa\wedge\zeta\,
\Big(
K^5
-
a_3
\Big)
+
\kappa\wedge\overline{\kappa}\,
\Big(
K^6
-
a_4
\Big)
+
\kappa\wedge\overline{\zeta}\,
\big(
-a_5
\big)
+
\zeta\wedge\overline{\kappa},
\\
d\zeta
&
\,=\,
\isqrt\,\pi^2\wedge\kappa
+
\big(\pi^1-\overline{\pi}^1\big)
\wedge\zeta
\,+
\\
&
\ \ \ \ \
+
\rho\wedge\kappa
\Big(
Z^1
+
\isqrt\,b_1
\Big)
+
\rho\wedge\zeta\,
\Big(
Z^2
+
a_1
-
\overline{a}_1
\Big)
+
\rho\wedge\overline{\kappa}\,
\Big(
Z^3
\Big)
+
\rho\wedge\overline{\zeta}\,
\Big(
Z^4
\Big)
\,+
\\
&
\ \ \ \ \
+
\kappa\wedge\zeta\,
\Big(
Z^5
-
\isqrt\,b_3
+
a_2
-
\overline{a}_4
\Big)
+
\kappa\wedge\overline{\kappa}\,
\Big(
Z^6
-
\isqrt\,b_4
\Big)
+
\kappa\wedge\overline{\zeta}\,
\Big(
Z^7
-
\isqrt\,b_5
\Big)
\,+
\\
&
\ \ \ \ \
+
\zeta\wedge\overline{\kappa}\,
\Big(
Z^8
-
a_4
+
\overline{a}_2
\Big)
+
\zeta\wedge\overline{\zeta}\,
\Big(
Z^9
-
a_5
+
\overline{a}_3
\Big).
\endaligned
\]

Now, replacing from
Assertion~{\ref{Assertion-relations-R-K-Z}}:
\[
Z^8
\,:=\,
-\,\overline{R}^1
+
2\,K^6,
\ \ \ \ \ \ \ \ \ \ \ \ \ \
Z^6
\,:=\,
\isqrt\,K^3,
\ \ \ \ \ \ \ \ \ \ \ \ \ \
K^5
\,:=\,
R^2,
\ \ \ \ \ \ \ \ \ \ \ \ \ \
Z^7
\,:=\,
\isqrt\,K^4,
\ \ \ \ \ \ \ \ \ \ \ \ \ \
Z^9
\,:=\,
-\,\overline{R}^2,
\]
the absorption equations write out as:
\[
\aligned
a_2+\overline{a}_4
&
\,=\,
R^1,
\\
\underline{a_3+\overline{a}_5}
&
\,=\,
R^2,
\\
\\
\\
\\
\\
\\
\\
\\
\\
\endaligned
\ \ \ \ \ \ \ \ \ \ \ \ \ \ \ \ \ \ \ \ \ \ \ \ \
\aligned
-\,a_1+b_2
&
\,=\,
K^1,
\\
b_3
&
\,=\,
K^2,
\\
b_4
&
\,=\,
K^3,
\\
b_5
&
\,=\,
K^4,
\\
\underline{a_3}
&
\,=\,
R^2,
\\
a_4
&
\,=\,
K^6,
\\
\underline{-\,a_5}
&
\,=\,
0,
\\
\\
\\
\\
\endaligned
\ \ \ \ \ \ \ \ \ \ \ \ \ \ \ \ \ \ \ \ \ \ \ \ \
\aligned
\isqrt\,b_1
&
\,=\,
-\,Z^1,
\\
-\,a_1+\overline{a}_1
&
\,=\,
Z^2,
\\
0
&
\,=\,
\boxed{Z^3},
\\
0
&
\,=\,
\boxed{Z^4},
\\
-\,a_2+\overline{a}_4+\isqrt\,b_3
&
\,=\,
Z^5,
\\
\isqrt\,b_4
&
\,=\,
\isqrt\,K^3,
\\
\isqrt\,b_5
&
\,=\,
\isqrt\,K^4,
\\
-\,\overline{a}_2+a_4
&
\,=\,
-\,\overline{R}^1
+
2\,K^6,
\\
\underline{-\,\overline{a}_3+a_5}
&
\,=\,
-\,\overline{R}^2.
\endaligned
\]
The boxed $Z^3$ and $Z^4$ are clearly essential torsions,
since they cannot be annihilated by any choice of $a_i$, $b_i$.
We will compute them explicitly a bit later.

At the end of the second colon, $a_5 = 0$, whence at the ends
of the other two colons, we get $a_3 := R^2$, hence all the $4$
underlined equations drop. 

Also, unique assignments exist for:
\[
\ \ \ \ \ \ \ \ \ \ \ \ \ \ \ \ \ \ \ \ \ \ \ \ \
\ \ \ \ \ \ \ \ \ \ \ \ \ \ \ \ \ \ \
\aligned
b_3
&
\,:=\,
K^2,
\\
b_4
&
\,:=\,
K^3,
\\
b_5
&
\,:=\,
K^4,
\\
a_4
&
\,:=\,
K^6,
\endaligned
\ \ \ \ \ \ \ \ \ \ \ \ \ \ \ \ \ \ \ \ \ \ \ \ \
\aligned
b_1
&
\,:=\,
\isqrt\,Z^1,
\\
b_4
&
\,:=\,
K^3,
\\
b_5
&
\,:=\,
K^4,
\\
\\
\endaligned
\]
and it remains to solve:
\[
\aligned
a_2+\overline{K}^6
\
&
\overset{\ast}{\,=\,}
R^1,
\\
\\ 
\\
\endaligned
\ \ \ \ \ \ \ \ \ \ \ \ \ \ \ \ \ \ \ \ \ \ \ \ \
\aligned
-\,a_1+b_2
&
\,=\,
K^1,
\\
\\ 
\\
\endaligned
\ \ \ \ \ \ \ \ \ \ \ \ \ \ \ \ \ \ \ \ \ \ \ \ \
\aligned
-\,a_1+\overline{a}_1
&
\,=\,
Z^2,
\\
-\,a_2+\overline{K}^6+\isqrt\,K^2
&
\overset{\text{\bf ?}}{\,=\,}
Z^5,
\\
-\,\overline{a}_2+K^6
&
\overset{\ast}{\,=\,}
-\,\overline{R}^1
+
2\,K^6.
\endaligned
\]
Certainly:
\[
b_2
\,:=\,
K^1
+
a_1,
\]
and the two equations $\overset{\ast}{\,=\,}$ for $a_2$
are equivalent\,\,---\,\,this comes
from the normalization relation $0 = \overline{R}^1 - 2\, K^6
+ Z^8$ already taken account of\,\,---, yielding:
\[
a_2
\,:=\,
R^1
-
\overline{K}^6.
\]

However, the equation $\overset{\text{\bf ?}}{\,=\,}$
cannot be satisfied automatically, and this provides
an essential torsion combination:
\[
-\,R^1
+
\overline{K}^6
+
\overline{K}^6
+
\isqrt\,K^2
\,=\,
Z^5
\ \ \ \ \ \ \ \ \ \ \ \ \ \ \ \ \ \
\Longleftrightarrow
\ \ \ \ \ \ \ \ \ \ \ \ \ \ \ \ \ \
-\,\isqrt\,K^2
+
Z^5
-
\overline{Z}^8
\,=\,
0,
\]
which was already seen in
Lemma~{\ref{Lemma-K-2-Z-5-bar-Z-8}}.

The last remaining equation:
\[
-\,a_1
+
\overline{a}_1
\,=\,
Z^2,
\]
shows that one can annihilate $\Im\, Z^2$
by choosing:
\[
\Im\,a_1
\,:=\,
-\,\frac{1}{2}\,
\Im\,Z^2,
\]
and it still remains precisely {\em one} real degree
of freedom, a free variable that we will re-denote:
\[
{\sf t}
\,:=\,
\Re\,a_1.
\]
In summary, we have established a fundamental

\begin{Proposition}
\label{Proposition-intermediate-before-R-secondary}
With ${\sf t} \in \R$ being a free variable, by defining
the precise modified Maurer-Cartan forms:
\[
\aligned
\pi^1
&
\,:=\,
\alpha
-
\Big(
{\sf t}
-
\frac{\isqrt}{2}\,
\Im\,Z^2
\Big)\,
\rho
-
\Big(
R^1
-
\overline{K}^6
\Big)\,
\kappa
-
R^2\,
\zeta
-
K^6\,\overline{\kappa}
-
0,
\\
\pi^2
&
\,:=\,
\beta
-
\isqrt\,Z^1\,\rho
-
\Big(
{\sf t}
-
\frac{\isqrt}{2}\,
\Im\,Z^2
+
K^1
\Big)\,
\kappa
-
K^2\,\zeta
-
K^3\,\overline{\kappa}
-
K^4\,\overline{\zeta},
\endaligned
\]
it holds:
\reqnomode\usetagform{EngelLie}
\begin{align}
d\rho
&
\,=\,
\big(
\pi^1
+
\overline{\pi}^1
\big)
\wedge\rho
+
\isqrt\,\kappa\wedge\overline{\kappa},
\notag
\\
d\kappa
&
\,=\,
\pi^2\wedge\rho
+
\pi^1\wedge\kappa
+
\zeta\wedge\overline{\kappa},
\notag
\\
d\zeta
&
\,=\,
\big(\pi^1-\overline{\pi}^1\big)
\wedge\zeta
+
\isqrt\,\pi^2\wedge\kappa
\,+
\notag
\\
&
\ \ \ \ \
\big(\Re\,Z^2\big)\,
\rho\wedge\zeta
+
Z^3\,\rho\wedge\overline{\kappa}
+
Z^4\,\rho\wedge\overline{\zeta}
+
\Big(
Z^5
+
R^1
-
2\,\overline{K}^6
-
\isqrt\,K^2
\Big)\,
\kappa\wedge\zeta.
\tag{\qed}
\end{align}
\end{Proposition}

We yet have to compute the remaining $4$ essential torsion 
coefficients:
\[
\Re\,Z^2,
\ \ \ \ \ \ \ \ \ \ \ \ \ \ \ \ \ \
Z^3,
\ \ \ \ \ \ \ \ \ \ \ \ \ \ \ \ \ \
Z^4,
\ \ \ \ \ \ \ \ \ \ \ \ \ \ \ \ \ \
Z^5
+
R^1
-
2\,\overline{K}^6
-
\isqrt\,K^2.
\]
Fortunately, by anticipation,
we have already explored and finalized:
\[
\aligned
Z^5
+
R^1
-
2\,\overline{K}^6
-
\isqrt\,K^2
&
\,=\,
-\,\isqrt\,K^2
+
Z^5
-
\overline{Z}^8
\\
&
\,=\,
\frac{1}{\sf c}\,
\Waux_0.
\endaligned
\]

\begin{Assertion}
One torsion coefficient vanishes identically:
\[
0
\,\equiv\,
Z^4.
\]
\end{Assertion}

\proof
Recall:
\[
\aligned
Z^4
&
\,=\,
\big[
\rho\wedge\overline{\zeta}
\big]
\big\{
d\zeta
\big\}
\\
&
\,=\,
\big[
\rho\wedge\overline{\zeta}
\big]
\Big\{
-\,\frac{\isqrt}{2}\,
\frac{\overline{\sf c}{\sf e}{\sf e}}{{\sf c}}\,
d\rho_0
+
{\sf e}\,d\kappa_0'
+
\frac{{\sf c}}{\overline{\sf c}}\,
d\zeta_0''
\Big\}.
\endaligned
\]
Compute separately:
\[
\aligned
-\,\frac{\isqrt}{2}\,
\frac{\overline{\sf c}{\sf e}{\sf e}}{{\sf c}}\,
\big[
\rho\wedge\overline{\zeta}
\big]
\big\{
d\rho_0
\big\}
&
\,=\,
-\,\frac{\isqrt}{2}\,
\frac{\overline{\sf c}{\sf e}{\sf e}}{{\sf c}}\,
\overline{\Raux}_0^{2\prime}\,
\Big(
\frac{1}{{\sf c}\overline{\sf c}}
\Big)
\Big(
\frac{{\sf c}}{\overline{\sf c}}
\Big)
\,=\,
-\,\frac{\isqrt}{2}\,
\frac{{\sf e}{\sf e}}{{\sf c}\overline{\sf c}}\,
\overline{\Raux}_0^{2\prime},
\\
{\sf e}\,
\big[
\rho\wedge\overline{\zeta}
\big]
\big\{
d\kappa_0'
\big\}
&
\,=\,
0,
\\
\frac{{\sf c}}{\overline{\sf c}}\,
\big[
\rho\wedge\overline{\zeta}
\big]
\big\{
d\zeta_0''
\big\}
&
\,=\,
\frac{{\sf c}}{\overline{\sf c}}\,
\zero{
\Zaux_0^{4\prime\prime}}\,
\Big(
\frac{1}{{\sf c}\overline{\sf c}}
\Big)
\Big(
\frac{{\sf c}}{\overline{\sf c}}
\Big)
+
\frac{{\sf c}}{\overline{\sf c}}\,
\Zaux_0^{9\prime}\,
\Big(
-\frac{\isqrt}{2}\,
\frac{\overline{\sf c}{\sf e}{\sf e}}{{\sf c}{\sf c}{\sf c}}
\Big)
\Big(
\frac{{\sf c}}{\overline{\sf c}}
\Big)
\\
&
\,=\,
0
-
\frac{\isqrt}{2}\,
\frac{{\sf e}{\sf e}}{{\sf c}{\sf c}}\,
\Zaux_0^{9\prime},
\endaligned
\]
and since we have already seen in
Lemma~{\ref{Lemma-Z-0-4-prime-prime}} 
that $\Zaux_0^{4\prime\prime} \equiv 0$, in the proof
of which we have used $\overline{\Raux}_0^{2\prime} + 
\Zaux_0^{9\prime} \equiv 0$, the sum of these $3$ terms
is indeed zero, and we done.
\endproof

It remains to analyze $Z^3$ and $\Re\, Z^2$, a substantial
task to which the two next sections are devoted.
At least, we know that:
\[
\aligned
d\zeta
&
\,=\,
\big(
\pi^1
-
\overline{\pi}^1
\big)\wedge\rho
+
\isqrt\,
\kappa\wedge\overline{\kappa}
\,+
\\
&
\ \ \ \ \
+
\big(\Re\,Z^2\big)\,
\rho\wedge\zeta
+
Z^3\,\rho\wedge\overline{\kappa}
+
\frac{1}{{\sf c}}\,
\Waux_0\,
\kappa\wedge\zeta.
\endaligned
\]

\Section{\bf Computation of Pocchiola's invariant $\Jaux_0$}
\label{computation-Jaux-0}
\HEAD{{\ref{computation-Jaux-0}}.~{\sf Computation 
of Pocchiola's invariant $\Jaux_0$}
}{
Wei Guo {\sc Foo} (Beijing) and Joël {\sc Merker} (Orsay)}

We now determine:
\[
\aligned
Z^3
&
\,=\,
\big[
\rho\wedge\overline{\kappa}
\big]
\big\{
d\zeta
\big\}
\\
&
\,=\,
-\,\frac{\isqrt}{2}\,
\frac{\overline{\sf c}{\sf e}{\sf e}}{{\sf c}}\,
\big[
\rho\wedge\overline{\kappa}
\big]
\big\{
d\rho_0
\big\}
+
{\sf e}\,
\big[
\rho\wedge\overline{\kappa}
\big]
\big\{
d\kappa_0'
\big\}
+
\frac{{\sf c}}{\overline{\sf c}}\,
\big[
\rho\wedge\overline{\kappa}
\big]
\big\{
d\zeta_0''
\big\}
\\
&
\,=\,
-\,\frac{\isqrt}{2}\,
\frac{\overline{\sf c}{\sf e}{\sf e}}{{\sf c}}\,
\bigg[
\overline{\Raux}_0^{1\prime}\,
\Big(
\frac{1}{{\sf c}\overline{\sf c}}
\Big)
\Big(
\frac{1}{\overline{\sf c}}
\Big)
+
\overline{\Raux}_0^{2\prime}\,
\Big(
\frac{1}{{\sf c}\overline{\sf c}}
\Big)
\Big(
-
\frac{{\sf c}\overline{\sf e}}{\overline{\sf c}\overline{\sf c}}
\Big)
+
\zero{
\isqrt\,
\Big(
\isqrt\,\frac{{\sf e}}{{\sf c}{\sf c}}
\Big)
\Big(
\frac{1}{\overline{\sf c}}
\Big)}
\bigg]
\\
&
\ \ \ \ \
+
{\sf e}\,
\bigg[
\Kaux_0^{3\prime\prime}
\Big(
\frac{1}{{\sf c}\overline{\sf c}}
\Big)
\Big(
\frac{1}{\overline{\sf c}}
\Big)
+
\Kaux_0^{6\prime}
\Big(
\isqrt\,\frac{{\sf e}}{{\sf c}{\sf c}}
\Big)
\Big(
\frac{1}{\overline{\sf c}}
\Big)
+
\zero{
\Big(
-\frac{\isqrt}{2}\,
\frac{\overline{\sf c}{\sf e}{\sf e}}{{\sf c}{\sf c}{\sf c}}
\Big)
\Big(
\frac{1}{\overline{\sf c}}
\Big)}
\bigg]
\\
&
\ \ \ \ \
+
\frac{{\sf c}}{\overline{\sf c}}\,
\bigg[
\Zaux_0^{3\prime\prime}
\Big(
\frac{1}{{\sf c}\overline{\sf c}}
\Big)
\Big(
\frac{1}{\overline{\sf c}}
\Big)
+
\Zaux_0^{6\prime\prime}
\Big(
\isqrt\,
\frac{{\sf e}}{{\sf c}{\sf c}}
\Big)
\Big(
\frac{1}{\overline{\sf c}}
\Big)
+
\Zaux_0^{8\prime}
\Big(
-\frac{\isqrt}{2}\,
\frac{\overline{\sf c}{\sf e}{\sf e}}{{\sf c}{\sf c}{\sf c}}
\Big)
\Big(
\frac{1}{\overline{\sf c}}
\Big)
+
\Zaux_0^{9\prime}
\Big(
-\frac{\isqrt}{2}\,
\frac{\overline{\sf c}{\sf e}{\sf e}}{{\sf c}{\sf c}{\sf c}}
\Big)
\Big(
\frac{\overline{\sf c}\overline{\sf e}}{
\overline{\sf c}\overline{\sf c}}
\Big)
\bigg],
\endaligned
\]
hence after collecting:
\[
\aligned
Z^3
&
\,=\,
\frac{{\sf e}{\sf e}}{{\sf c}{\sf c}\overline{\sf c}}
\left[
-\frac{\isqrt}{2}\,
\overline{\Raux}_0^{1\prime}
+
\isqrt\,
\Kaux_0^{6\prime}
-
\frac{\isqrt}{2}\,
\Zaux_0^{8\prime}
\right]
+
\frac{{\sf e}{\sf e}\overline{\sf e}}{{\sf c}\overline{\sf c}
\overline{\sf c}}
\left[
\zero{
\frac{\isqrt}{2}\,
\overline{\Raux}_0^{2\prime}
+
\frac{\isqrt}{2}\,
\Zaux_0^{9\prime}}
\right]
\,+
\\
&
\ \ \ \ \
+
\frac{{\sf e}}{{\sf c}\overline{\sf c}\overline{\sf c}}
\Big[
\Kaux_0^{3\prime\prime}
+
\isqrt\,
\Zaux_0^{6\prime\prime}
\Big]
+
\frac{1}{
\overline{\sf c}\overline{\sf c}\overline{\sf c}}\,
\Zaux_0^{3\prime\prime}.
\endaligned
\]
As we already know, the second term vanishes, the third one as well:
\[
\Kaux_0^{3\prime\prime}
+
\isqrt\,
\Zaux_0^{6\prime\prime}
\,\,=\,\,
2\isqrt\,\Haux_0
-
\isqrt\,\Haux_0
-
\isqrt\,\Haux_0,
\]
and also the first one:
\[
\!\!\!\!\!\!\!\!\!\!\!\!\!\!\!\!\!\!\!\!\!\!\!\!\!
\aligned
-\,
\frac{\isqrt}{2}\,
\overline{\Raux}_0^{1\prime}
+
\isqrt\,
\Kaux_0^{6\prime}
-
\frac{\isqrt}{2}\,
\Zaux_0^{8\prime}
&
\,=\,
-\,\frac{\isqrt}{2}
\left(
\frac{1}{3}\,
\frac{\overline{\mathcal{L}}_1\big(
\overline{\mathcal{L}}_1(\kaux)\big)}{
\overline{\mathcal{L}}_1(\kaux)}
+
\frac{2}{3}\,
\overline{\Paux}
\right)
+
\isqrt
\left(
-
\frac{1}{3}\,
\frac{\overline{\mathcal{L}}_1\big(
\overline{\mathcal{L}}_1(\kaux)\big)}{
\overline{\mathcal{L}}_1(\kaux)}
+
\frac{1}{3}\,
\overline{\Paux}
\right)
-
\frac{\isqrt}{2}
\left(
-
\frac{\overline{\mathcal{L}}_1\big(
\overline{\mathcal{L}}_1(\kaux)\big)}{
\overline{\mathcal{L}}_1(\kaux)}
\right).
\endaligned
\]
It remains only one term:
\[
\aligned
Z^3
&
\,=\,
\frac{1}{\overline{\sf c}\overline{\sf c}\overline{\sf c}}\,
\Zaux_0^{3\prime\prime}
\\
&
\,=\,
\frac{1}{\overline{\sf c}\overline{\sf c}\overline{\sf c}}\,
\left(
-\,\isqrt\,\Zaux_0^{8\prime}\,\Haux_0
-
\isqrt\,\overline{\mathcal{L}}_1\big(\Haux_0\big)
+
\isqrt\,\Haux_0\,
\overline{\Raux}_0^{1\prime}
\right)
\\
&
\,=\,
\frac{\isqrt}{\overline{\sf c}\overline{\sf c}\overline{\sf c}}\,
\left(
\frac{\overline{\mathcal{L}}_1\big(\overline{\mathcal{L}}_1
(\kaux)\big)}{\overline{\mathcal{L}}_1(\kaux)}\,
\Haux_0
-
\overline{\mathcal{L}}_1\big(\Haux_0\big)
+
\frac{1}{3}\,
\frac{\overline{\mathcal{L}}_1\big(
\overline{\mathcal{L}}_1(\kaux)\big)}{
\overline{\mathcal{L}}_1(\kaux)}\,
\Haux_0
+
\frac{2}{3}\,
\Haux_0\,
\overline{\Paux}
\right)
\\
&
\,=\,
\frac{\isqrt}{\overline{\sf c}\overline{\sf c}\overline{\sf c}}\,
\Bigg(
\underbrace{
\frac{4}{3}\,
\frac{\overline{\mathcal{L}}_1\big(
\overline{\mathcal{L}}_1(\kaux)\big)}{
\overline{\mathcal{L}}_1(\kaux)}\,
\Haux_0
+
\frac{2}{3}\,
\overline{\Paux}_0\,
\Haux_0
-
\overline{\mathcal{L}}_1\big(\Haux_0\big)
}_{=:\,\,\overline{\Jaux}_0}
\Bigg).
\endaligned
\]

Then a direct expansion of the derivative
$\overline{\mathcal{L}}_1 \big( \Haux_0\big)$ which
uses neither Lemma~{\ref{Lemma-K-bar-k-K-bar-P}},
nor Lemma~{\ref{Lemma-K-bar-penetrates-L1-bar-k}}, provides (exercise)
exactly the same expression as the one of Pocchiola:
\[
\aligned
\overline{\Jaux}_0
&
\,=\,
\frac{1}{6}\,
\frac{\overline{\mathcal{L}}_1\big(
\overline{\mathcal{L}}_1\big(
\overline{\mathcal{L}}_1\big(
\overline{\mathcal{L}}_1(\kaux)\big)\big)\big)}{
\overline{\mathcal{L}}_1(\kaux)}
-
\frac{5}{6}\,
\frac{\overline{\mathcal{L}}_1\big(
\overline{\mathcal{L}}_1\big(
\overline{\mathcal{L}}_1(\kaux)\big)\big)\,\,
\overline{\mathcal{L}}_1\big(
\overline{\mathcal{L}}_1(\kaux)\big)
}{
\overline{\mathcal{L}}_1(\kaux)^2}
-
\frac{1}{6}\,
\frac{\overline{\mathcal{L}}_1\big(
\overline{\mathcal{L}}_1\big(
\overline{\mathcal{L}}_1(\kaux)\big)\big)
}{
\overline{\mathcal{L}}_1(\kaux)}\,
\overline{\Paux}
\,+
\\
&
\ \ \ \ \
+
\frac{20}{27}\,
\frac{\overline{\mathcal{L}}_1\big(\overline{\mathcal{L}}_1
(\kaux)\big)^3}{
\overline{\mathcal{L}}_1(\kaux)^3}
+
\frac{5}{18}\,
\frac{\overline{\mathcal{L}}_1\big(
\overline{\mathcal{L}}_1(\kaux)\big)^2}{
\overline{\mathcal{L}}_1(\kaux)^2}\,
\overline{\Paux}
+
\frac{1}{6}\,
\frac{\overline{\mathcal{L}}_1\big(
\overline{\mathcal{L}}_1(\kaux)\big)\,\,
\overline{\mathcal{L}}_1\big(\overline{\Paux}\big)}{
\overline{\mathcal{L}}_1(\kaux)}
-
\frac{1}{9}\,
\frac{\overline{\mathcal{L}}_1\big(
\overline{\mathcal{L}}_1(\kaux)\big)}{
\overline{\mathcal{L}}_1(\kaux)}\,\,
\overline{\Paux}\,\overline{\Paux}
\,-
\\
&
\ \ \ \ \
-
\frac{1}{6}\,
\overline{\mathcal{L}}_1\big(
\overline{\mathcal{L}}_1\big(
\overline{\Paux}\big)\big)
+
\frac{1}{3}\,
\overline{\mathcal{L}}_1\big(\overline{\Paux}\big)\,
\overline{\Paux}
-
\frac{2}{27}\,
\overline{\Paux}\,
\overline{\Paux}\,
\overline{\Paux}.
\endaligned
\]

\Section{\bf Computation of the derived invariant $R := \Re\, Z^2$}
\label{computation-R}
\HEAD{{\ref{computation-R}}.~{\sf Computation 
of the derived invariant $R := \Re\, Z^2$}
}{
Wei Guo {\sc Foo} (Beijing) and Joël {\sc Merker} (Orsay)}

Next, we determine:
\[
\aligned
Z^2
&
\,=\,
\big[
\rho\wedge\zeta
\big]
\big\{
d\zeta
\big\}
\\
&
\,=\,
-\,\frac{\isqrt}{2}\,
\frac{\overline{\sf c}{\sf e}{\sf e}}{{\sf c}}\,
\big[
\rho\wedge\zeta
\big]
\big\{
d\rho_0
\big\}
+
{\sf e}\,
\big[
\rho\wedge\zeta
\big]
\big\{
d\kappa_0'
\big\}
+
\frac{{\sf c}}{\overline{\sf c}}\,
\big[
\rho\wedge\zeta
\big]
\big\{
d\zeta_0''
\big\}
\\
&
\,=\,
-\,\frac{\isqrt}{2}\,
\frac{\overline{\sf c}{\sf e}{\sf e}}{{\sf c}}
\left[
\Raux_0^{2\prime}
\Big(
\frac{1}{{\sf c}\overline{\sf c}}
\Big)
\Big(
\frac{\overline{\sf c}}{{\sf c}}
\Big)
\right]
\,+
\\
&
\ \ \ \ \
+
{\sf e}
\left[
\Kaux_0^{2\prime}
\Big(
\frac{1}{{\sf c}\overline{\sf c}}
\Big)
\Big(
\frac{\overline{\sf c}}{{\sf c}}
\Big)
+
\Kaux_0^{5\prime}
\Big(
\isqrt\,
\frac{{\sf e}}{{\sf c}{\sf c}}
\Big)
\Big(
\frac{\overline{\sf c}}{{\sf c}}
\Big)
-
\Big(
\frac{\overline{\sf c}}{{\sf c}}
\Big)
\Big(
-\isqrt\,
\frac{\overline{\sf e}}{\overline{\sf c}\overline{\sf c}}
\Big)
\right]
\,+
\\
&
\ \ \ \ \
+
\frac{{\sf c}}{\overline{\sf c}}
\left[
\Zaux_0^{2\prime\prime}
\Big(
\frac{1}{{\sf c}\overline{\sf c}}
\Big)
\Big(
\frac{\overline{\sf c}}{{\sf c}}
\Big)
+
\Zaux_0^{5\prime}
\Big(
\isqrt\,
\frac{{\sf e}}{{\sf c}{\sf c}}
\Big)
\Big(
\frac{\overline{\sf c}}{{\sf c}}
\Big)
-
\Zaux_0^{8\prime}
\Big(
\frac{\overline{\sf c}}{{\sf c}}
\Big)
\Big(
-\isqrt\,
\frac{\overline{\sf e}}{\overline{\sf c}\overline{\sf c}}
\Big)
-
\Zaux_0^{9\prime}
\Big(
\frac{\overline{\sf c}}{{\sf c}}
\Big)
\Big(
\frac{\isqrt}{2}\,
\frac{{\sf c}\overline{\sf e}\overline{\sf e}}{
\overline{\sf c}\overline{\sf c}\overline{\sf c}}
\Big)
\right]
\endaligned
\]
hence after collecting:
\[
\aligned
Z^2
&
\,=\,
\isqrt\,
\frac{{\sf e}\overline{\sf e}}{{\sf c}\overline{\sf c}}\,
+
\frac{\overline{\sf c}{\sf e}{\sf e}}{{\sf c}{\sf c}{\sf c}}
\left(
-
\frac{\isqrt}{2}\,
\Raux_0^{2\prime}
+
\isqrt\,\Kaux_0^{5\prime}
\right)
+
\frac{{\sf c}\overline{\sf e}\overline{\sf e}}{
\overline{\sf c}\overline{\sf c}\overline{\sf c}}
\left(
-\frac{\isqrt}{2}\,
\Zaux_0^{9\prime}
\right)
\,+
\\
&
\ \ \ \ \
+
\frac{{\sf e}}{{\sf c}{\sf c}}
\Big(
\Kaux_0^{2\prime}
+
\isqrt\,
\Zaux_0^{5\prime}
\Big)
+
\frac{\overline{\sf e}}{\overline{\sf c}\overline{\sf c}}
\Big(
\isqrt\,
\Zaux_0^{8\prime}
\Big)
+
\frac{1}{{\sf c}\overline{\sf c}}\,
\Zaux_0^{2\prime\prime},
\endaligned
\]
that is to say:
\[
\aligned
Z^2
&
\,=\,
\isqrt\,
\frac{{\sf e}\overline{\sf e}}{{\sf c}\overline{\sf c}}\,
+
\frac{\overline{\sf c}{\sf e}{\sf e}}{{\sf c}{\sf c}{\sf c}}
\left(
\frac{\isqrt}{2}\,
\frac{\mathcal{L}_1(\kaux)}{\overline{\mathcal{L}}_1(\kaux)}
-
\isqrt\,
\frac{\mathcal{L}_1(\kaux)}{\overline{\mathcal{L}}_1(\kaux)}
\right)
+
\frac{{\sf c}\overline{\sf e}\overline{\sf e}}{
\overline{\sf c}\overline{\sf c}\overline{\sf c}}
\left(
-\frac{\isqrt}{2}\,
\frac{\overline{\mathcal{L}}_1(\overline{\kaux})}{
\mathcal{L}_1(\overline{\kaux})}
\right)
\,+
\\
&
\ \ \ \ \
+
\frac{{\sf e}}{{\sf c}{\sf c}}
\left(
-\,\frac{\isqrt}{3}\,
\frac{\mathcal{K}\big(\overline{\mathcal{L}}_1\big(
\overline{\mathcal{L}}_1(\kaux)\big)\big)}{
\overline{\mathcal{L}}_1(\kaux)^2}
+
\frac{\isqrt}{3}\,
\frac{\mathcal{K}\big(\overline{\mathcal{L}}_1(\kaux)\big)\,\,
\overline{\mathcal{L}}_1\big(\overline{\mathcal{L}}_1(\kaux)\big)}{
\overline{\mathcal{L}}_1(\kaux)^3}
-
\frac{\isqrt}{3}\,
\frac{\mathcal{L}_1\big(\mathcal{L}_1(\overline{\kaux})\big)}{
\mathcal{L}_1(\overline{\kaux})}
\,-
\right.
\\
&
\ \ \ \ \ \ \ \ \ \ \ \ \ \ \ \ \ \ \ \ \ \ \ \ \ \
\ \ \ \ \ \ \ \ \ \ \ \ \ \ \ \ \ \ \ \ \ \ \ \ \ \
\ \ \ \ \ \ \ \ \ \ \ \ \ \ \ \ \ \ \ \ \ \ \ \ \ \
\left.
-\,
\frac{\isqrt}{3}\,
\frac{\overline{\mathcal{L}}_1\big(
\mathcal{L}_1(\kaux)\big)}{
\overline{\mathcal{L}}_1(\kaux)}
-
\frac{2}{3}\,
\frac{\mathcal{T}(\kaux)}{\overline{\mathcal{L}}_1(\kaux)}
+
\isqrt\,
\frac{\mathcal{L}_1\big(\overline{\mathcal{L}}_1(\kaux)\big)}{
\overline{\mathcal{L}}_1(\kaux)}
\right)
+
\\
&
\ \ \ \ \
+
\frac{\overline{\sf e}}{\overline{\sf c}\overline{\sf c}}
\left(
-\isqrt\,
\frac{\overline{\mathcal{L}}_1\big(
\overline{\mathcal{L}}_1(\kaux)\big)}{
\overline{\mathcal{L}}_1(\kaux)}
\right)
+
\frac{1}{{\sf c}\overline{\sf c}}
\left(
\underbrace{
-\,
\isqrt\,\Zaux_0^{9\prime}\,
\overline{\Haux}_0
+
\isqrt\,\Haux_0\,
\Raux_0^{2\prime}}_{\text{\sf on hold}}
+
\Zaux_0^{2\prime}
-
\isqrt\,
\frac{\mathcal{K}(\Haux_0)}{\overline{\mathcal{L}}_1(\kaux)}
\right).
\endaligned
\]

Now, observe firstly that when we consider:
\[
2\,\Re\,Z^2
\,=\,
Z^2
+
\overline{Z}^2,
\]
the real part of the sum of the first three terms of $Z^2$:
\[
\isqrt\,
\frac{{\sf e}\overline{\sf e}}{{\sf c}\overline{\sf c}}
+
\frac{\overline{\sf c}{\sf e}{\sf e}}{{\sf c}{\sf c}{\sf c}}\,
\left(
-\,
\frac{\isqrt}{2}\,
\frac{\mathcal{L}_1(\kaux)}{\overline{\mathcal{L}}_1(\kaux)}
\right)
+
\frac{{\sf c}\overline{\sf e}\overline{\sf e}}{
\overline{\sf c}\overline{\sf c}\overline{\sf c}}
\left(
-\,\frac{\isqrt}{2}\,
\frac{\overline{\mathcal{L}}_1(\overline{\kaux})}{
\mathcal{L}_1(\overline{\kaux})}
\right)
\]
vanishes, visibly. 

Secondly, in the sum $Z^2 + 
\overline{Z}^2$, if the terms multiples of $\frac{{\sf e}}{
{\sf c}{\sf c}}$ are grouped together, we realize
that we recover $\Waux_0$ exactly:
\[
\aligned
\isqrt\,
\frac{{\sf e}}{{\sf c}{\sf c}}
&
\left(
-\,\frac{1}{3}\,
\frac{\mathcal{K}\big(\overline{\mathcal{L}}_1\big(
\overline{\mathcal{L}}_1(\kaux)\big)\big)}{
\overline{\mathcal{L}}_1(\kaux)^2}
+
\frac{1}{3}\,
\frac{\mathcal{K}\big(\overline{\mathcal{L}}_1(\kaux)\big)\,\,
\overline{\mathcal{L}}_1\big(\overline{\mathcal{L}}_1(\kaux)\big)}{
\overline{\mathcal{L}}_1(\kaux)^3}
-
\frac{1}{3}\,
\frac{\mathcal{L}_1\big(\mathcal{L}_1(\overline{\kaux})\big)}{
\mathcal{L}_1(\overline{\kaux})}
\,-
\right.
\\
&
\ \ \ \ \
\left.
-\,
\frac{1}{3}\,
\frac{\overline{\mathcal{L}}_1\big(
\mathcal{L}_1(\kaux)\big)}{
\overline{\mathcal{L}}_1(\kaux)}
+
\frac{2\isqrt}{3}\,
\frac{\mathcal{T}(\kaux)}{\overline{\mathcal{L}}_1(\kaux)}
+
\frac{\mathcal{L}_1\big(\overline{\mathcal{L}}_1(\kaux)\big)}{
\overline{\mathcal{L}}_1(\kaux)}
+
\frac{\mathcal{L}_1\big(\mathcal{L}_1(\overline{\kaux})\big)}{
\mathcal{L}_1(\overline{\kaux})}
\right)
\\
&
\,=\,
i\,
\frac{{\sf e}}{{\sf c}{\sf c}}\,
\left(
-\,\frac{1}{3}\,
\frac{\mathcal{K}\big(\overline{\mathcal{L}}_1\big(
\overline{\mathcal{L}}_1(\kaux)\big)\big)}{
\overline{\mathcal{L}}_1(\kaux)^2}
+
\frac{1}{3}\,
\frac{\mathcal{K}\big(\overline{\mathcal{L}}_1(\kaux)\big)\,\,
\overline{\mathcal{L}}_1\big(\overline{\mathcal{L}}_1(\kaux)\big)}{
\overline{\mathcal{L}}_1(\kaux)^3}
+
\frac{2}{3}\,
\frac{\mathcal{L}_1\big(\mathcal{L}_1(\overline{\kaux})\big)}{
\mathcal{L}_1(\overline{\kaux})}
\,+
\right.
\\
&
\ \ \ \ \ \ \ \ \ \ \ \ \ \ \
\left.
+
\frac{2}{3}\,
\frac{\mathcal{L}_1\big(
\overline{\mathcal{L}}_1(\kaux)\big)}{
\overline{\mathcal{L}}_1(\kaux)}
+
\frac{\isqrt}{3}\,
\frac{\mathcal{T}(\kaux)}{\overline{\mathcal{L}}_1(\kaux)}
\right)
\\
&
\,=\,
i\,\frac{{\sf e}}{{\sf c}{\sf c}}\,
\Waux_0,
\endaligned
\]
as we remember its explicit expression from
Section~{\ref{third-loop-reduction-d}}.

In addition thirdly, using the explicit expressions
from Proposition~{\ref{Proposition-3-loop-structure-0-coframe}}:
\[
\Raux_0^{2\prime}
\,=\,
-\,
\frac{\mathcal{L}_1(\kaux)}{\overline{\mathcal{L}}_1(\kaux)}
\ \ \ \ \ \ \ \ \ \ \ \ \ \ \ \ \ \
\text{and}
\ \ \ \ \ \ \ \ \ \ \ \ \ \ \ \ \ \
\Zaux_0^{9\prime}
\,=\,
\frac{\overline{\mathcal{L}}_1(\overline{\kaux})}{
\mathcal{L}_1(\overline{\kaux})},
\]
and the explicit expression of:
\[
\Haux_0
\,=\,
-\,\frac{1}{6}\,
\frac{\overline{\mathcal{L}}_1\big(
\overline{\mathcal{L}}_1\big(
\overline{\mathcal{L}}_1(\kaux)\big)\big)}{
\overline{\mathcal{L}}_1(\kaux)}
+
\frac{2}{9}\,
\frac{\overline{\mathcal{L}}_1\big(
\overline{\mathcal{L}}_1(\kaux)\big)^2}{
\overline{\mathcal{L}}_1(\kaux)^2}
+
\frac{1}{18}\,
\frac{\overline{\mathcal{L}}_1\big(
\overline{\mathcal{L}}_1(\kaux)\big)}{
\overline{\mathcal{L}}_1(\kaux)}\,
\overline{\Paux}
+
\frac{1}{6}\,
\overline{\mathcal{L}}_1\big(\overline{\Paux}\big)
-
\frac{1}{9}\,
\overline{\Paux}\,
\overline{\Paux},
\]
we verify by a direct computation the identical vanishing:
\[
0
\,\equiv\,
-\,\isqrt\,
\Zaux_0^{9\prime}\,\overline{\Haux}_0
+
\isqrt\,\Haux_0\,\Raux_0^{2\prime}
+
\overline{
-\,\isqrt\,\Zaux_0^{9\prime}\,\overline{\Haux}_0
+
\isqrt\,\Haux_0\,\Raux_0^{2\prime}},
\]
which means that the term `{\sl on hold}' underbraced above
disappears when taking $2\, \Re\, Z^2$, and we receive:
\[
2\,\Re\,Z^2
\,=\,
\isqrt\,\frac{{\sf e}}{{\sf c}{\sf c}}\,
\Waux_0
-
\isqrt\,\frac{\overline{\sf e}}{\overline{\sf c}\overline{\sf c}}\,
\overline{\Waux}_0
+
\frac{1}{{\sf c}\overline{\sf c}}
\left(
\Zaux_0^{2\prime}
-
\isqrt\,
\frac{\mathcal{K}(\Haux_0)}{\overline{\mathcal{L}}_1(\kaux)}
+
\overline{\Zaux}_0^{2\prime}
+
\isqrt\,
\frac{\overline{\mathcal{K}}(\overline{\Haux}_0)}{
\mathcal{L}_1(\overline{\kaux})}
\right).
\] 

Fourthly and lastly, by replacing:
\[
\Haux_0
\,=\,
-\,\frac{\isqrt}{2}\,
\Kaux_0^{3\prime},
\]
we get:
\leqnomode\usetagform{default}
\begin{align}
\label{2-Re-Z-2-intermediate}
2\,\Re\,Z^2
\,=\,
2\,\Re
\Bigg(
\isqrt\,\frac{{\sf e}}{{\sf c}{\sf c}}\,
\Waux_0
+
\frac{1}{{\sf c}\overline{\sf c}}\,
\bigg(
\Zaux_0^{2\prime}
\underbrace{
-
\frac{1}{2}\,
\frac{\mathcal{K}\big(\Kaux_0^{3\prime}\big)}{
\overline{\mathcal{L}}_1(\kaux)}}_{\text{\sf on hold}}
\bigg)
\Bigg).
\end{align}

A miraculous re-expression of $2\, \Re\, Z^2$ was discovered
by Pocchiola on his computer, and was shown
in~{\cite{Pocchiola-2013, Merker-Pocchiola-2018}}, 
but without any details of proof.

\begin{Lemma}
One has in fact:
\[
2\,\Re\,Z^2
\,=\,
2\,\Re\,
\left[
\isqrt\,
\frac{{\sf e}}{{\sf c}{\sf c}}\,
\Waux_0
+
\frac{1}{{\sf c}\overline{\sf c}}
\bigg(
-\,\frac{\isqrt}{2}\,
\overline{\mathcal{L}}_1\big(\Waux_0\big)
+
\frac{\isqrt}{2}\,
\bigg(
-\,\frac{1}{3}\,
\frac{\overline{\mathcal{L}}_1\big(
\overline{\mathcal{L}}_1(\kaux)\big)}{
\overline{\mathcal{L}}_1(\kaux)}
+
\frac{1}{3}\,
\overline{\Paux}
\bigg)\,
\Waux_0
\bigg)
\right].
\]
\end{Lemma}

This expression shows that $\Re(Z^{2})$ depends on the first jet of
$\Waux_0$, that it vanishes when $\Waux_0 = 0$, and therefore,
$\Re\, Z^2$ is {\em not}
a primary invariant.  We provide details of proof, with no
computer help.

\proof
To transform the term `{\sl on hold}' above, we need
a technical lemma, whose proof, to be done afterwards,
uses mainly the Poincaré relation $d \circ d = 0$ applied
to the structure equations~({\ref{inexplicit-0-structure-3-loop}}). 

\begin{Lemma}
\label{Lemma-two-identities}
The following two identities hold identically:
\leqnomode\usetagform{default}
\begin{align}
\label{identity-K-K-0-3}
\frac{\mathcal{K}\big(\Kaux_0^{3\prime}\big)}{
\overline{\mathcal{L}}_1(\kaux)}
&
\,=\,
\overline{\mathcal{L}}_1\big(\Kaux_0^{2\prime}\big)
-
\Kaux_0^{2\prime}\,
\Kaux_0^{6\prime}
-
\Kaux_0^{1\prime}
+
\overline{\Kaux}_0^{1\prime}
+
\Zaux_0^{2\prime},
\\
\label{identity-L-1-bar-Z-0-5-8}
\overline{\mathcal{L}}_1\big(\Zaux_0^{5\prime}\big)
+
\mathcal{L}_1
\big(\Zaux_0^{8\prime}\big)
&
\,=\,
\Zaux_0^{5\prime}\,
\Kaux_0^{6\prime}
+
\Zaux_0^{8\prime}\,
\overline{\Kaux}_0^{6\prime}
+
\isqrt\,
\Zaux_0^{2\prime}.
\end{align}
\end{Lemma}

Admitting these identities temporarily, let us prove
the proposition. 
In order to replace the term `{\sl on hold}' 
in~({\ref{2-Re-Z-2-intermediate}}) above, 
let us multiply by $-\frac{1}{2}$ the first 
identity~({\ref{identity-K-K-0-3}}), and
take $2\, \Re ( \centersmallbullet)$:
\[
2\,\Re\,
\bigg(
-\,\frac{1}{2}\,
\frac{\mathcal{K}\big(\Kaux_0^{3\prime}\big)}{
\overline{\mathcal{L}}_1(\kaux)}
\bigg)
\,=\,
2\,\Re\,
\bigg(
-\,
\frac{1}{2}\,
\boxed{
\overline{\mathcal{L}}_1
\big(\Kaux_0^{2\prime}\big)}
+
\frac{1}{2}\,
\Kaux_0^{2\prime}\,
\Kaux_0^{6\prime}
+
0
-
\frac{1}{2}\,
\Zaux_0^{2\prime}
\bigg).
\]
We yet have to transform the boxed term.
To this aim, we conjugate the second 
identity~({\ref{identity-L-1-bar-Z-0-5-8}}):
\[
\mathcal{L}_1
\big(\overline{\Zaux}_0^{5\prime}\big)
+
\overline{\mathcal{L}}_1
\big(\overline{\Zaux}_0^{8\prime}\big)
\,=\,
\overline{\Zaux}_0^{5\prime}\,
\overline{\Kaux}_0^{6\prime}
+
\overline{\Zaux}_0^{8\prime}\,
\Kaux_0^{6\prime}
-
\isqrt\,
\overline{\Zaux}_0^{2\prime},
\]
and to this identity multiplied by $\isqrt$, we 
subtract~({\ref{identity-K-K-0-3}}) also multiplied by $\isqrt$,
to get:
\[
-\,
\isqrt\,
\overline{\mathcal{L}}_1
\Big(
\Zaux_0^{5\prime}
-
\overline{\Zaux}_0^{8\prime}
\Big)
+
\isqrt\,
\mathcal{L}_1
\Big(
\overline{\Zaux}_0^{5\prime}
-
\Zaux_0^{8\prime}
\Big)
\,\,=\,\,
-\,\isqrt\,
\Kaux_0^{6\prime}\,
\Big(
\Zaux_0^{5\prime}
-
\overline{\Zaux}_0^{8\prime}
\Big)
+
\isqrt\,
\overline{\Kaux}_0^{6\prime}\,
\Big(
\overline{\Zaux}_0^{5\prime}
-
\Zaux_0^{8\prime}
\Big)
+
\Zaux_0^{2\prime}
+
\overline{\Zaux}_0^{2\prime}.
\]
But here, remembering that, by definition of $\Waux_0$:
\[
\Zaux_0^{5\prime}
-
\overline{\Zaux}_0^{8\prime}
\,=\,
\Waux_0
+
\isqrt\,
\Kaux_0^{2\prime},
\]
we can replace to get:
\[
-\,\isqrt\,
\overline{\mathcal{L}}_1
\big(\Waux_0\big)
+
\overline{\mathcal{L}}_1
\big(\Kaux_0^{2\prime}\big)
+
\isqrt\,\mathcal{L}_1
\big(\overline{\Waux}_0\big)
+
\mathcal{L}_1
\big(\overline{\Kaux}_0^{2\prime}\big)
\,\,=\,\,
-\,\isqrt\,
\Kaux_0^{6\prime}\,
\Waux_0
+
\Kaux_0^{6\prime}\,
\Kaux_0^{2\prime}
+
\isqrt\,\overline{\Kaux}_0^{6\prime}\,
\overline{\Waux}_0
+
\overline{\Kaux}_0^{6\prime}\,
\overline{\Kaux}_0^{2\prime}
+
\Zaux_0^{2\prime}
+
\overline{\Zaux}_0^{2\prime},
\]
that is to say for the mentioned boxed term:
\[
2\,\Re\,
\Big(
\overline{\mathcal{L}}_1
\big(\Kaux_0^{2\prime}\big)
\Big)
\,=\,
2\,\Re\,
\bigg(
\isqrt\,\overline{\mathcal{L}}_1
\big(\Waux_0\big)
-
\isqrt\,\Kaux_0^{6\prime}\,\Waux_0
+
\Kaux_0^{2\prime}\,
\Kaux_0^{6\prime}
+
\Zaux_0^{2\prime}
\bigg).
\]
Multiplying this result by $-\, \frac{1}{2}$, 
and replacing above yields:
\[
2\,\Re\,
\bigg(
-\,\frac{1}{2}\,
\frac{\mathcal{K}\big(\Kaux_0^{3\prime}\big)}{
\overline{\mathcal{L}}_1(\kaux)}
\bigg)
\,\,=\,\,
2\,\Re\,
\bigg(
-\,\frac{\isqrt}{2}\,
\overline{\mathcal{L}}_1\big(\Waux_0\big)
+
\frac{\isqrt}{2}\,
\Kaux_0^{6\prime}\,
\Waux_0
-
\zero{
\frac{1}{2}\,
\Kaux_0^{2\prime}\,
\Kaux_0^{6\prime}}
-
\frac{1}{2}\,
\Zaux_0^{2\prime}
+
\zero{
\frac{1}{2}\,
\Kaux_0^{2\prime}\,
\Kaux_0^{6\prime}}
-
\frac{1}{2}\,
\Zaux_0^{2\prime}
\bigg),
\]
and a final replacement 
in~({\ref{2-Re-Z-2-intermediate}})
concludes, if one remembers that:
\[
\Kaux_0^{6\prime}
\,=\,
-\,\frac{1}{3}\,
\frac{\overline{\mathcal{L}}_1\big(
\overline{\mathcal{L}}_1(\kaux)\big)}{
\overline{\mathcal{L}}_1(\kaux)}
+
\frac{1}{3}\,
\overline{\Paux}.
\qedhere
\]
\endproof

\proof[Proof of Lemma~{\ref{Lemma-two-identities}}]
To treat the first 
identity~({\ref{identity-K-K-0-3}}), apply the exterior
differentiation operator $d$ to the structure
equation for $d\kappa_0'$ from~({\ref{inexplicit-0-structure-3-loop}}):
\begin{equation*}
\begin{aligned}
0 &= 
d^{2}\kappa_{0}'\\
&=
d\Kaux_{0}^{1\prime}\wedge\rho_{0}\wedge\kappa_{0}'
+
\Kaux_{0}^{1\prime}\ d\rho_{0}\wedge\kappa_{0}'
-
\Kaux_{0}^{1\prime}\ \rho_{0}\wedge d\kappa_{0}'\\
&+
d\Kaux_{0}^{2\prime}\wedge\rho_{0}\wedge\zeta_{0}'
+\Kaux_{0}^{2\prime}\  d\rho_{0}\wedge\zeta_{0}'
-\Kaux_{0}^{2\prime}\ \rho_{0}\wedge d\zeta_{0}'\\
&+
\underbrace{d\Kaux_{0}^{3\prime}\wedge\rho_{0}
\wedge\overline{\kappa}_{0}'}_{\text{needed}}
+\Kaux_{0}^{3\prime}\ d\rho_{0}\wedge\overline{\kappa}_{0}'
-\Kaux_{0}^{3\prime}\ \rho_{0}\wedge d\overline{\kappa}_{0}'\\
&+
d\Kaux_{0}^{5\prime}\wedge\kappa_{0}'\wedge\zeta_{0}'
+\Kaux_{0}^{5\prime}\ d\kappa_{0}'\wedge\zeta_{0}'
-\Kaux_{0}^{5\prime}\ \kappa_{0}'\wedge d\zeta_{0}'\\
&+
d\Kaux_{0}^{6\prime}\wedge\kappa_{0}'\wedge\overline{\kappa}_{0}'
+\Kaux_{0}^{6\prime}\ d\kappa_{0}'\wedge\overline{\kappa}_{0}'
-\Kaux_{0}^{6\prime}\ \kappa_{0}'\wedge d\overline{\kappa}_{0}'\\
&
\ \ \ \ \ \ \ \ \ \ \ \ \ \ \ \ \ \ \ \ \ \ \ \ \ \ \ \ \ 
+
d\zeta_{0}'\wedge\overline{\kappa}_{0}'
-\zeta_{0}'\wedge d\overline{\kappa}_{0}'.
\end{aligned}
\end{equation*}
Because we are dealing with $\mathcal{K}\big(
\Kaux_{0}^{3\prime} \big)$, we
can wedge throughout with 
$\kappa_{0}' \wedge \overline{\zeta}_{0}'$ to obtain
$\mathcal{K} \big( \Kaux_{0}^{3\prime} \big) \big/
\overline{\mathcal{L}}_{1}(\kaux)$
from the term marked `{\sl needed}', and we get:
\begin{equation*}
\begin{aligned}
0
&=
0
\hspace{4.2cm}+
0
\hspace{3.7cm}-
\Kaux_{0}^{1\prime}\ 
\rho_{0}\wedge d\kappa_{0}'\wedge\kappa_{0}'\wedge\overline{\zeta}_{0}'\\
&+
d\Kaux_{0}^{2\prime}\wedge\rho_{0}\wedge\zeta_{0}'\wedge\kappa_{0}'\wedge\overline{\zeta}_{0}'
\hspace{0.2cm}+
\Kaux_{0}^{2\prime}\ d\rho_{0}\wedge\zeta_{0}'\wedge\kappa_{0}'\wedge\overline{\zeta}_{0}'
\hspace{0.17cm}-
\Kaux_{0}^{2\prime}\ \rho_{0}\wedge d\zeta_{0}'\wedge\kappa_{0}'\wedge\overline{\zeta}_{0}'\\
&+
d\Kaux_{0}^{3\prime}\wedge\rho_{0}\wedge\overline{\kappa}_{0}'\wedge\kappa_{0}'\wedge\overline{\zeta}_{0}'
\hspace{0.15cm}+
\Kaux_{0}^{3\prime}\ d\rho_{0}\wedge\overline{\kappa}_{0}'\wedge\kappa_{0}'\wedge\overline{\zeta}_{0}'
\hspace{0.15cm}-
\Kaux_{0}^{3\prime}\ \rho_{0}\wedge d\overline{\kappa}_{0}'\wedge\kappa_{0}'\wedge\overline{\zeta}_{0}'\\
&+ 0
\hspace{4.3cm}+ \Kaux_{0}^{5\prime}\  d\kappa_{0}'\wedge\zeta_{0}'\wedge\kappa_{0}'\wedge\overline{\zeta}_{0}'
\hspace{0.3cm}-
0\\
&+
0
\hspace{4.3cm}+
\Kaux_{0}^{6\prime}\ d\kappa_{0}'\wedge\overline{\kappa}_{0}'\wedge\kappa_{0}'\wedge\overline{\zeta}_{0}'
\hspace{0.2cm}-
0\\
&+
0
\ \ \ \ \ \ \ \ \ \ \ \ \ \ \ \ \ \ \ \ \ \ \ \ \ \ \ \ \ \ \ \ \ \ 
\ \ \ \ \ \ \
+
d\zeta_0'\wedge\overline{\kappa}_{0}'\wedge\kappa_{0}'\wedge\overline{\zeta}_{0}'
\hspace{1.10cm}
-\zeta_{0}'\wedge d\overline{\kappa}_{0}'\wedge\kappa_{0}'\wedge\overline{\zeta}_{0}'.
\end{aligned}
\end{equation*}

In the left column,
observe that two exterior differentials appear, $d\Kaux_0^{2\prime}$,
$d\Kaux_0^{3\prime}$.
Already in 
Section~{\ref{D-C-structure-kappa-0-prime-zeta-0-prime-prime}},
we have implicitly used the following
companion of Lemma~{\ref{Lemma-d-G-0}}.

\begin{Lemma}
\label{Lemma-d-G-kappa-prime-zeta-prime}
The exterior differential of 
any function $\Gaux = \Gaux \big(z_1, z_2, \overline{z}_1,
\overline{z}_2, v \big)$ on $M$ expresses as:
\[
d\Gaux
\,=\,
\bigg(
\mathcal{T}\big(\Gaux\big)
-
\frac{i}{3}\,\Baux_0
+
\frac{i}{3}\,
\overline{\Baux}_0
\bigg)\,
\rho_0
+
\mathcal{L}_1\big(\Gaux\big)\,
\kappa_0'
+
\frac{\mathcal{K}(\Gaux)}{\overline{\mathcal{L}}_1(\kaux)}\,
\zeta_0'
+
\overline{\mathcal{L}}_1\big(\Gaux\big)\,
\overline{\kappa}_0'
+
\frac{\mathcal{K}(\Gaux)}{\mathcal{L}_1(\overline{\kaux})}\,
\overline{\zeta}_0'.
\]
\end{Lemma}

\proof
Replacing $\kappa_0$ by $\kappa_0' - \frac{i}{3}\,
\Baux_0\, \rho_0$ from~({\ref{def-kappa-prime-0}}), 
and $\zeta_0$ by $\frac{\zeta_0'}{\overline{\mathcal{L}}_1
(\kaux)}$ from~({\ref{def-zeta-prime-0}}), we indeed obtain:
\begin{align}
d\Gaux
&
\,=\,
\mathcal{T}\big(\Gaux\big)\,
\rho_0
+
\mathcal{L}_1\big(\Gaux\big)\,
\kappa_0
+
\mathcal{K}\big(\Gaux\big)\,
\zeta_0
+
\overline{\mathcal{L}}_1\big(\Gaux\big)\,
\overline{\kappa}_0
+
\overline{\mathcal{K}}\big(\Gaux\big)\,
\overline{\zeta}_0
\notag
\\
&
\,=\,
\mathcal{T}\big(\Gaux\big)\,
\rho_0
+
\mathcal{L}_1\big(\Gaux\big)\,
\Big(
\kappa_0'
-
\frac{i}{3}\,
\Baux_0\,
\rho_0
\Big)
+
\mathcal{K}\big(\Gaux\big)\,
\frac{\zeta_0'}{\overline{\mathcal{L}}_1(\kaux)}
\,+
\notag
\\
&
\ \ \ \ \ \ \ \ \ \ \ \ \ \ \ \ \ \ \ \ \ \ \ \ \ \ \ \ \ \ \
+
\overline{\mathcal{L}}_1\big(\Gaux\big)\,
\Big(
\overline{\kappa}_0'
+
\frac{i}{3}\,\overline{\Baux}_0\,\rho_0
\Big)\,
+
\overline{\mathcal{K}}\big(\Gaux\big)\,
\frac{\overline{\zeta}_0'}{\mathcal{L}_1(\overline{\kaux})}.
\qedhere
\end{align}
\endproof

Using this lemma for $d\Kaux_0^{2\prime}$, $d\Kaux_0^{3\prime}$, and
replacing also $d\rho_0$, $d\kappa_0'$, $d\zeta_0'$,
$d\overline{\kappa}_0'$, $d\overline{\zeta}_0'$ by means
of~({\ref{inexplicit-0-structure-3-loop}}), we have:
\[
\!\!\!\!\!\!\!\!\!\!\!\!\!\!\!
\aligned
0
&
\,=\,
0
\ \ \ \ \ \ \ \ \ \ \ \ \ \ \ \ \ \ \ \ \ \ \ \ \ \ \ \ \ \ \ \ \ \ 
\ \ \ \ \ \ \ \ \ \ \ \ \ \ \ \ \ \ \
+
0
\ \ \ \ \ \ \ \ \ \ \ \ \ \ \ \ \ \ \ \ \ \ \ \ \ \ \ \ \ \ \ \ \ \ 
\ \ \ \ \ \ \ \ \ \ \ 
-
\Kaux_0^{1\prime}\,
\rho_0\wedge\zeta_0'\wedge\overline{\kappa}_0'\wedge
\kappa_0'\wedge\overline{\zeta}_0'
\,+
\\
&
\ \ \ \ \
+
\overline{\mathcal{L}}_1\big(\Kaux_0^{2\prime}\big)\,
\overline{\kappa}_0'\wedge\rho_0\wedge\zeta_0'
\wedge\kappa_0'\wedge\overline{\zeta}_0'
+
\Kaux_0^{2\prime}\,
\overline{\Raux}_0^{1\prime}\,
\rho_0\wedge\overline{\kappa}_0'\wedge\zeta_0'
\wedge\kappa_0'\wedge\overline{\zeta}_0'
-
\Kaux_0^{2\prime}\,\rho_0\wedge
\Zaux_0^{8\prime}\,\zeta_0'\wedge\overline{\kappa}_0'
\wedge\overline{\zeta}_0'
\,+
\\
&
\ \ \ \ \ 
+
\frac{\mathcal{K}(\Kaux_0^{3\prime})}{
\overline{\mathcal{L}}_1(\kaux)}\,
\zeta_0'\wedge\rho_0\wedge\overline{\kappa}_0'\wedge\kappa_0'
\wedge\overline{\zeta}_0'
+
\Kaux_0^{3\prime}\,\Raux_0^{2\prime}\,
\rho_0\wedge\zeta_0'\wedge\overline{\kappa}_0'\wedge\kappa_0'
\wedge\overline{\zeta}_0'
-
0
\ \ \ \ \ \ \ \ \ \ \ \ \ \ \ \ \ \ \ \ \ \ \ \ \ \ \ \ \ \ \ \ \ \ 
\ \ \ \ \ \ 
+
\\
&
\ \ \ \ \
+
0
\ \ \ \ \ \ \ \ \ \ \ \ \ \ \ \ \ \ \ \ \ \ \ \ \ \ \ \ \ \ \ \ \ \ 
\ \ \ \ \ \ \ \ \ \ \ \ \ \ \ 
+
\Kaux_0^{5\prime}\,\Kaux_0^{3\prime}\,
\rho_0\wedge\overline{\kappa}_0'\wedge\zeta_0'\wedge\kappa_0'
\wedge\overline{\zeta}_0'
-
0
\ \ \ \ \ \ \ \ \ \ \ \ \ \ \ \ \ \ \ \ \ \ \ \ \ \ \ \ \ \ \ \ \ \ 
\ \ \ \ \ \
+
\\
&
\ \ \ \ \
+
0
\ \ \ \ \ \ \ \ \ \ \ \ \ \ \ \ \ \ \ \ \ \ \ \ \ \ \ \ \ \ \ \ \ \ 
\ \ \ \ \ \ \ \ \ \ \ \ \ \ \
+
\Kaux_0^{6\prime}\,\Kaux_0^{2\prime}\,
\rho_0\wedge\zeta_0'\wedge\overline{\kappa}_0'\wedge\kappa_0'
\wedge\overline{\zeta}_0'
-
0
\ \ \ \ \ \ \ \ \ \ \ \ \ \ \ \ \ \ \ \ \ \ \ \ \ \ \ \ \ \ \ \ \ \ 
\ \ \ \ \ \
+
\\
&
\ \ \ \ \
+
0
\ \ \ \ \ \ \ \ \ \ \ \ \ \ \ \ \ \ \ \ \ \ \ \ \ \ \ \ \ \ \ \ \ \ 
\ \ \ \ \ \ \ \ \ \ \ \ \ \ \
+
\Zaux_0^{2\prime}\,
\rho_0\wedge\zeta_0'\wedge\overline{\kappa}_0'\wedge\kappa_0'
\wedge\overline{\zeta}_0'
\ \ \ \ \ \ \
-
\zeta_0'\wedge\overline{\Kaux}_0^{1\prime}\,
\wedge\rho_0\wedge\overline{\kappa}_0'
\wedge\kappa_0'\wedge\overline{\zeta}_0',
\endaligned
\]
hence caring about signs when factoring by the naturally
appearing $5$-form:
\[
\aligned
0
\,=\,
\rho_0\wedge\kappa_0'\wedge\zeta_0'
\wedge\overline{\kappa}_0'\wedge\overline{\zeta}_0'\,
\Big(
&
0
\ \ \ \ \ \ \ \ \ \ \ \ \ \ \ \ \ \
+
0
\ \ \ \ \ \ \ \ \ 
-
\Kaux_0^{1\prime}
\ \ \ \ \ \
+
\\
&
+
\overline{\mathcal{L}}_1\big(\Kaux_0^{2\prime}\big)
-
\Kaux_0^{2\prime}\,
\overline{\Raux}_0^{1\prime}
-
\Kaux_0^{2\prime}\,
\Zaux_0^{8\prime}
-
\\
&
-\,
\frac{\mathcal{K}(\Kaux_0^{3\prime})}{
\overline{\mathcal{L}}_1(\kaux)}
+
\Kaux_0^{3\prime}\,\Raux_0^{2\prime}
-
0
\ \ \ \ \ \ \ \ \
+
\\
&
+
0
\ \ \ \ \ \ \ \ \ \ \ \ \ 
-
\Kaux_0^{5\prime}\,\Kaux_0^{3\prime}
-
0
\ \ \ \ \ \ \ \ \
+
\\
&
+
0
\ \ \ \ \ \ \ \ \ \ \ \ \
+
\Kaux_0^{6\prime}\,\Kaux_0^{2\prime}
-
0
\ \ \ \ \ \ \ \ \
+
\\
&
+
0
\ \ \ \ \ \ \ \ \ \ \ \ \
+
\Zaux_0^{2\prime}
\ \ \ \ \ \ \ 
+
\overline{\Kaux}_0^{1\prime}
\ \ \ \ \ \ \ \ \ \ \ \ \ \ \ 
\Big),
\endaligned
\]
whence we arrive at the announced first 
identity~({\ref{identity-K-K-0-3}}) 
by remembering some useful relations:
\[
\frac{\mathcal{K}\big(\Kaux_0^{3\prime}\big)}{
\overline{\mathcal{L}}_1(\kaux)}
\,=\,
\overline{\mathcal{L}}_1
\big(\Kaux_0^{2\prime}\big)
+
\Kaux_0^{2\prime}\,
\Kaux_0^{6\prime}
-
\Kaux_0^{2\prime}\,
\Big(
\underbrace{
\overline{\Raux}_0^{1\prime}
+
\Zaux_0^{8\prime}}_{=\,\,2\,\Kaux_0^{6\prime}}
\Big)
+
\Kaux_0^{3\prime}\,
\Big(
\underbrace{
\Raux_0^{2\prime}
-
\Kaux_0^{5\prime}}_{=\,\,0!}
\Big)
-
\Kaux_0^{1\prime}
+
\overline{\Kaux}_0^{1\prime}
+
\Zaux_0^{2\prime}.
\]

For the second 
identity~({\ref{identity-L-1-bar-Z-0-5-8}}), we proceed
similarly, applying the exterior differentiation
operator $d$ to the structure equation for
$d\zeta_0'$ from~({\ref{inexplicit-0-structure-3-loop}}):
\begin{equation*}
\begin{aligned}
0 &= d^{2}\zeta_{0}'\\
&=
\underbrace{d(\Zaux_{0}^{2\prime})\wedge\rho_{0}\wedge\zeta_{0}'}_{\text{don't want}}
+
\Zaux_{0}^{2\prime}\ d\rho_{0}\wedge\zeta_{0}'
-
\Zaux_{0}^{2\prime}\ \rho_{0}\wedge d\zeta_{0}'\\
&+
\underbrace{d(\Zaux_{0}^{5\prime})\wedge\kappa_{0}'\wedge\zeta_{0}'}_{\text{want}}
+
\Zaux_{0}^{5\prime}\ d\kappa_{0}'\wedge\zeta_{0}'
-
\Zaux_{0}^{5\prime}\ \kappa_{0}'\wedge d\zeta_{0}'\\
&+
\underbrace{d(\Zaux_{0}^{8\prime})\wedge\zeta_{0}'\wedge\overline{\kappa}_{0}'}_{\text{want}}
+
\Zaux_{0}^{8\prime}\ d\zeta_{0}'\wedge\overline{\kappa}_{0}'
-
\Zaux_{0}^{8\prime}\ \zeta_{0}'\wedge d\overline{\kappa}_{0}'\\
&+
\underbrace{d(\Zaux_{0}^{9\prime})\wedge\zeta_{0}'\wedge\overline{\zeta}_{0}'}_{\text{don't want}}
+
\Zaux_{0}^{9\prime}\ d\zeta_{0}'\wedge\overline{\zeta}_{0}'
-
\Zaux_{0}^{9\prime}\ \zeta_{0}'\wedge d\overline{\zeta}_{0}'.
\end{aligned}
\end{equation*}
Observe that the desired identity involves the derivatives of
$\Zaux_{0}^{5\prime}$ and $\Zaux_{0}^{8\prime}$. Hence we may
conserve those terms marked `{\sl want}' by wedging with the
appropriate 2-form $\rho_{0}\wedge\overline{\zeta}_{0}'$:
\begin{equation*}
\begin{aligned}
0 
&
=
0 
\hspace{4cm}+
\Zaux_{0}^{2\prime}\ d\rho_{0}\wedge\zeta_{0}'\wedge\rho_{0}\wedge\overline{\zeta}_{0}'
+
0\\
&+
d\Zaux_{0}^{5\prime}\wedge\kappa_{0}'\wedge\zeta_{0}'\wedge\rho_{0}\wedge\overline{\zeta}_{0}'
+
\Zaux_{0}^{5\prime}\ d\kappa_{0}'\wedge\zeta_{0}'\wedge\rho_{0}\wedge\overline{\zeta}_{0}'
-
\Zaux_{0}^{5\prime}\ \kappa_{0}'\wedge d\zeta_{0}'\wedge\rho_{0}\wedge\overline{\zeta}_{0}'\\
&+
d\Zaux_{0}^{8\prime}\wedge\zeta_{0}'\wedge\overline{\kappa}_{0}'\wedge\rho_{0}\wedge\overline{\zeta}_{0}'
+
\Zaux_{0}^{8\prime}\ d\zeta_{0}'\wedge\overline{\kappa}_{0}'\wedge\rho_{0}\wedge\overline{\zeta}_{0}'
-
\Zaux_{0}^{8\prime}\ \zeta_{0}'\wedge d\overline{\kappa}_{0}'\wedge\rho_{0}\wedge\overline{\zeta}_{0}'\\
&+
0
\hspace{4cm}+
0
\hspace{3.6cm}-
\Zaux_{0}^{9\prime}\ \zeta_{0}'\wedge d\overline{\zeta}_{0}'\wedge\rho_{0}\wedge\overline{\zeta}_{0}'.
\end{aligned}
\end{equation*}
Using Lemma~{\ref{Lemma-d-G-kappa-prime-zeta-prime}}
for $d\Zaux_0^{5\prime}$,
$d\Zaux_0^{8\prime}$, 
and replacing also
$d\rho_0$, $d\kappa_0'$, $d\zeta_0'$, 
$d\overline{\kappa}_0'$,
$d\overline{\zeta}_0'$ by means 
of~({\ref{inexplicit-0-structure-3-loop}}), we have:
\[
\!\!\!\!\!\!\!\!\!\!\!\!\!\!\!
\aligned
0
&
\,=\,
0
\ \ \ \ \ \ \ \ \ \ \ \ \ \ \ \ \ \ \ \ \ \ \ \ \ \ \ \ \ \ \ \
\ \ \ \ \ \ \ \ \ \ \ \ \ \ \ \ \ \ \ \ \
+
\Zaux_0^{2\prime}\,\isqrt\,
\kappa_0'\wedge\overline{\kappa}_0'
\wedge\zeta_0'\wedge\rho_0\wedge\overline{\zeta}_0'
\ \ \ 
+
0
\,+
\\
&
\ \ \ \ \
+
\overline{\mathcal{L}}_1\big(\Zaux_0^{5\prime}\big)\,
\overline{\kappa}_0'\wedge\kappa_0'\wedge\zeta_0'
\wedge\rho_0\wedge\overline{\zeta}_0'
+
\Zaux_0^{5\prime}\,
\Kaux_0^{6\prime}\,
\kappa_0'\wedge\overline{\kappa}_0'\wedge\zeta_0'
\wedge\rho_0\wedge\overline{\zeta}_0'
+
\Zaux_0^{5\prime}\,\kappa_0'\wedge
\Zaux_0^{8\prime}\,\zeta_0'\wedge
\overline{\kappa}_0'\wedge\rho_0\wedge\overline{\zeta}_0'
\,+
\\
&
\ \ \ \ \
+
\mathcal{L}_1\big(\Zaux_0^{8\prime}\big)\,
\kappa_0'\wedge\zeta_0'\wedge\overline{\kappa}_0'\wedge\rho_0
\wedge\overline{\zeta}_0'
\ 
+
\Zaux_0^{8\prime}\,
\Zaux_0^{5\prime}\,
\kappa_0'\wedge\zeta_0'\wedge\overline{\kappa}_0'\wedge\rho_0
\wedge\overline{\zeta}_0'
+
\Zaux_0^{8\prime}\,\zeta_0'\wedge\overline{\Kaux}_0^{6\prime}
\wedge\overline{\kappa}_0'\wedge\kappa_0'\wedge\rho_0\wedge
\overline{\zeta}_0'
\,+
\\
&
\ \ \ \ \
+
0
\ \ \ \ \ \ \ \ \ \ \ \ \ \ \ \ \ \ \ \ \ \ \ \ \ \ \ \ \ \ \ \
\ \ \ \ \ \ \ \ \ \ \ \ \ \ \ \ \ \ 
+
0
\ \ \ \ \ \ \ \ \ \ \ \ \ \ \ \ \ \ \ \ \ \ \ \ \ \ \ \ \ \ \ \
\ \ \ \ \ \ \ \ \ \ \ \ \ 
-
0,
\endaligned
\]
hence caring about signs when factoring by the naturally
appearing $5$-form,
we arrive at the announced second
identity~({\ref{identity-L-1-bar-Z-0-5-8}}):
\begin{align}
0
\,=\,
\rho_0\wedge\kappa_0'\wedge\zeta_0'
\wedge\overline{\kappa}_0'\wedge\overline{\zeta}_0'\,
\Big(
&
0
\ \ \ \ \ \ \ \ \ \ \ \ \ \ \ \ \ \
+
\isqrt\,\Zaux_0^{2\prime} 
\ \ \ \
+
0
\ \ \ \ \ \ \ \ \ \ \
-
\notag
\\
&
-\,
\overline{\mathcal{L}}_1\big(\Zaux_0^{5\prime}\big)
+
\Zaux_0^{5\prime}\,\Kaux_0^{6\prime}
+
\zero{
\Zaux_0^{5\prime}\,\Zaux_0^{8\prime}}
\,-
\notag
\\
&
\,-
\mathcal{L}_1\big(\Zaux_0^{8\prime}\big)
-
\zero{
\Zaux_0^{8\prime}\,\Zaux_0^{5\prime}}
+
\Zaux_0^{8\prime}\,
\overline{\Kaux}_0^{6\prime}
\,+
\notag
\\
&
+
0
\ \ \ \ \ \ \ \ \ \ \ \ \
+
0
\ \ \ \ \ \ \ \ \ \ 
-
0
\ \ \ \ \ \ \ \ \ \ \ \ \ \ \ 
\Big).
\qedhere
\end{align}
\endproof

\Section{\bf Summarized structure equations}
\label{summarized-structure-equations}
\HEAD{{\ref{summarized-structure-equations}}.~{\sf Summarized 
structure equations}
}{
Wei Guo {\sc Foo} (Beijing) and Joël {\sc Merker} (Orsay)}

All this work conducted us to finalize the statement
of Proposition~{\ref{Proposition-intermediate-before-R-secondary}},
but before, 
let us make an ample summary.

After normalizations of the group parameters
${\sf f}$, ${\sf b}$, ${\sf d}$, the equivalence
problem for $2$-nondegenerate (constant) Levi rank $1$
$\mathcal{C}^\omega$ or $\mathcal{C}^\infty$ real hypersurfaces
$M^5 \subset \C^3$ conducts to a $4$-dimensional $G$-structure:
\[
\left(\!
\begin{array}{ccc}
{\sf c}\overline{\sf c} & 0 & 0
\\
-\isqrt\,\overline{\sf c}{\sf e} & {\sf c} & 0
\\
-\frac{\isqrt}{2}\,
\frac{\overline{\sf c}{\sf e}{\sf e}}{{\sf c}}
& {\sf e} & \frac{{\sf c}}{\overline{\sf c}}
\end{array}
\!\right),
\]
where ${\sf c} \in \C^\ast$ and ${\sf e} \in \C$, 
with Maurer-Cartan forms (conjutates are not written):
\[
\aligned
\alpha
&
\,:=\,
\frac{d{\sf c}}{{\sf c}},
\\
\beta
&
\,:=\,
\isqrt\,
\frac{{\sf e}\,d{\sf c}}{{\sf c}{\sf c}}
-
\isqrt\,\frac{{\sf e}\,d\overline{\sf c}}{{\sf c}\overline{\sf c}}
-
\isqrt\,\frac{d{\sf e}}{{\sf c}}.
\endaligned
\]

Furthermore, 
$2$ fundamental
primary differential invariants occur:
\[
J
\,=\,
\frac{\isqrt}{\overline{\sf c}\overline{\sf c}\overline{\sf c}}\,
\overline{\Jaux}_0
\ \ \ \ \ \ \ \ \ \ \ \ \ \ \ \ \ \
\text{and}
\ \ \ \ \ \ \ \ \ \ \ \ \ \ \ \ \ \
W
\,=\,
\frac{1}{{\sf c}}\,
\Waux_0,
\]
where $\Jaux_0$ and $\Waux_0$ are explicit functions
on $M$, together with $1$ secondary invariant:
\[
\aligned
R
\,:=\,
&\,
\Re\,
Z^2
\\
\,=\,
&\,
\Re\,
\left[
\isqrt\,
\frac{{\sf e}}{{\sf c}{\sf c}}\,
\Waux_0
+
\frac{1}{{\sf c}\overline{\sf c}}
\bigg(
-\,\frac{\isqrt}{2}\,
\overline{\mathcal{L}}_1\big(\Waux_0\big)
+
\frac{\isqrt}{2}\,
\bigg(
-\,\frac{1}{3}\,
\frac{\overline{\mathcal{L}}_1\big(
\overline{\mathcal{L}}_1(\kaux)\big)}{
\overline{\mathcal{L}}_1(\kaux)}
+
\frac{1}{3}\,
\overline{\Paux}
\bigg)\,
\Waux_0
\bigg)
\right].
\endaligned
\]

On the $10$-dimensional manifold $M^5 \times G^4 \times \R$
equipped with coordinates:
\[
\big(
z_1,z_2,\overline{z}_1,\overline{z}_2,v
\big)
\times
\big(
{\sf c},\overline{\sf c},
{\sf e},\overline{\sf e}
\big)
\times
({\sf t}),
\]
there are two modified-prolonged Maurer-Cartan forms:\[
\aligned
\pi^1
&
\,:=\,
\alpha
-
\Big(
{\sf t}
-
\frac{\isqrt}{2}\,
\Im\,Z^2
\Big)\,
\rho
-
\Big(
R^1
-
\overline{K}^6
\Big)\,
\kappa
-
R^2\,
\zeta
-
K^6\,\overline{\kappa}
-
0,
\\
\pi^2
&
\,:=\,
\beta
-
\isqrt\,Z^1\,\rho
-
\Big(
{\sf t}
-
\frac{\isqrt}{2}\,
\Im\,Z^2
+
K^1
\Big)\,
\kappa
-
K^2\,\zeta
-
K^3\,\overline{\kappa}
-
K^4\,\overline{\zeta},
\endaligned
\]
where $R^i$, $K^i$, $Z^i$ are explicit functions
on $M^5 \times G^4$. 

\begin{Theorem}
\label{Theorem-finalization-absorption-d-rho-d-kappa-d-zeta}
After finalization of absorption, 
the structure equations read:
\reqnomode\usetagform{EngelLie}
\begin{align}
d\rho
&
\,=\,
\big(
\pi^1
+
\overline{\pi}^1
\big)
\wedge\rho
+
\isqrt\,\kappa\wedge\overline{\kappa},
\notag
\\
d\kappa
&
\,=\,
\pi^2\wedge\rho
+
\pi^1\wedge\kappa
+
\zeta\wedge\overline{\kappa},
\notag
\\
d\zeta
&
\,=\,
\big(\pi^1-\overline{\pi}^1\big)
\wedge\zeta
+
\isqrt\,\pi^2\wedge\kappa
\,+
\notag
\\
&
\ \ \ \ \
+
R\,
\rho\wedge\zeta
+
J\,\rho\wedge\overline{\kappa}
+
W\,
\kappa\wedge\zeta.
\tag{\qed}
\end{align}
\end{Theorem}

\Section{\bf The final $\{e\}$-structure}
\label{final-e-structure}
\HEAD{{\ref{final-e-structure}}.~{\sf The final $\{e\}$-structure}
}{
Wei Guo {\sc Foo} (Beijing) and Joël {\sc Merker} (Orsay)}

Let $\Omega_1$ and $\Omega_2$ be the two $2$-forms defined by:
\[
\aligned
\Omega_1
&
\,:=\,
d\pi^1
-
\isqrt\,
\kappa\wedge\overline{\pi}^2
-
\zeta\wedge\overline{\zeta},
\\
\Omega_2
&
\,:=\,
d\pi^2
-
\pi^2\wedge\overline{\pi}^1
-
\zeta\wedge\overline{\pi}^2.
\endaligned
\]
When the two fundamental invariants $\Jaux_0 \equiv 0 \equiv \Waux_0$
vanish identically, since we know that:
\[
\aligned
R
&
\,=\,
\Re\,
\left[
\isqrt\,
\frac{{\sf e}}{{\sf c}{\sf c}}\,
\Waux_0
+
\frac{1}{{\sf c}\overline{\sf c}}
\bigg(
-\,\frac{\isqrt}{2}\,
\overline{\mathcal{L}}_1\big(\Waux_0\big)
+
\frac{\isqrt}{2}\,
\bigg(
-\,\frac{1}{3}\,
\frac{\overline{\mathcal{L}}_1\big(
\overline{\mathcal{L}}_1(\kaux)\big)}{
\overline{\mathcal{L}}_1(\kaux)}
+
\frac{1}{3}\,
\overline{\Paux}
\bigg)\,
\Waux_0
\bigg)
\right],
\\
J
&
\,=\,
\frac{\isqrt}{\overline{\sf c}\overline{\sf c}\overline{\sf c}}\,
\overline{\Jaux}_0,
\\
W
&
\,=\,
\frac{1}{{\sf c}}\,
\Waux_0,
\endaligned
\]
it comes:
\[
0
\,\equiv\,
R
\,\equiv\,
J
\,\equiv\,
W.
\]

Independently, the 
addendum
to~{\cite{Merker-Pocchiola-2018}},
reproduced below on p.~{\pageref{addendum-Merker-Pocchiola-JGA}},
shows that in the
case where all invariants vanish, these auxiliary $2$-forms
$\Omega_1$ and $\Omega_2$ satisfy:
\[
\aligned
\big(
\Omega_1+\overline{\Omega}_1
\big)
\wedge
\rho
&
\,=\,
0,
\\
\Omega_2\wedge\rho
+
\Omega_1\wedge\kappa
&
\,=\,
0,
\\
\big(
\Omega_1-\overline{\Omega}_1
\big)
\wedge\zeta
+
\isqrt\,
\Omega_2\wedge\kappa
&
\,=\,
0.
\endaligned
\]
In general, the right-hand sides of these structure equations
are not necessarily zero, and they depend on the invariants
$R$, $J$, $W$.

\begin{Proposition}
\label{Proposition-3-structure-Omega}
The two $2$-forms $\Omega_1$ and $\Omega_2$ satisfy:
\leqnomode\usetagform{default}
\begin{align}
\big(
\Omega_1
+
\overline{\Omega}_1
\big)
\wedge\rho
&
\,=\,
0,
\label{1-Omega-R-J-W}
\\
\Omega_2
\wedge\rho
+
\Omega_1
\wedge\kappa
&
\,=\,
-\,R\,
\rho\wedge\zeta\wedge\overline{\kappa}
-
W\,
\kappa\wedge\zeta\wedge\overline{\kappa},
\label{2-Omega-R-J-W}
\\
\isqrt\,
\Omega_2\wedge\kappa
+
\big(
\Omega_1-\overline{\Omega}_1
\big)
\wedge\zeta
&
\,=\,
-\,dR
\wedge\rho\wedge\zeta
-
R\,\big(
\pi^1+\overline{\pi}^1
\big)\wedge
\rho\wedge\zeta
-
\isqrt\,R\,
\pi^2\wedge\rho\wedge\kappa
+
\isqrt\,R\,\kappa\wedge\zeta\wedge\overline{\zeta}
\,-
\label{3-Omega-R-J-W}
\\
&
\ \ \ \ \
-
dJ\wedge\rho\wedge\overline{\kappa}
-
3\,J\,
\overline{\pi}^1\wedge\rho\wedge\overline{\kappa}
-
J\,\rho\wedge\kappa\wedge\overline{\zeta}
\,-
\notag
\\
&
\ \ \ \ \
-
dW\wedge\kappa\wedge\zeta
-
W\,\pi^2\wedge\rho\wedge\zeta
-
W\,\pi^1\wedge\kappa\wedge\zeta
-
WJ\,\rho\wedge\kappa\wedge\overline{\kappa}.
\notag
\end{align}
\end{Proposition}

\proof
These relations come from Poincaré's identities:
\[
0
\,\equiv\,
d\circ d\rho
\,\equiv\,
d\circ d\kappa
\,\equiv\,
d\circ d\zeta,
\]
applied to the finalized structure equations of
Theorem~{\ref{Theorem-finalization-absorption-d-rho-d-kappa-d-zeta}},
in which $d\rho$, $d\kappa$, $d\zeta$ should be replaced again using
Theorem~{\ref{Theorem-finalization-absorption-d-rho-d-kappa-d-zeta}},
followed by a reorganization of the obtained $3$-forms.

For the first line~({\ref{1-Omega-R-J-W}}):
\[
\aligned
0
&
\,=\,
d\circ d\rho
\\
&
\,=\,
\big(
d\pi^1
+
d\overline{\pi}^1
\big)
\wedge\rho
-
\big(
\pi^1
+
\overline{\pi}^1
\big)
\wedge d\rho
+
\isqrt\,d\kappa\wedge\overline{\kappa}
-
\isqrt\,\kappa\wedge d\overline{\kappa}
\\
&
\,=\,
\big(
d\pi^1
+
d\overline{\pi}^1
\big)
\wedge\rho
-
\big(
\pi^1
+
\overline{\pi}^1
\big)
\wedge
\Big(
\big(
\zero{
\pi^1
+
\overline{\pi}^1}
\big)
\wedge\rho
+
\isqrt\,\kappa\wedge\overline{\kappa}
\Big)
\,+
\\
&
\ \ \ \ \
+
\isqrt\,
\Big(
\pi^2\wedge\rho
+
\pi^1\wedge\rho
+
\zero{
\zeta\wedge\overline{\kappa}}
\Big)
\wedge
\overline{\kappa}
-
\isqrt\,
\kappa
\wedge
\Big(
\overline{\pi}^2\wedge\rho
+
\overline{\pi}^1\wedge\overline{\kappa}
+
\zero{
\overline{\zeta}\wedge\kappa}
\Big).
\endaligned
\]
Afer simplification, this becomes:
\[
0
\,=\,
\Big(
d\pi^1
-
\isqrt\,\kappa\wedge\overline{\pi}^2
\Big)
\wedge\rho
+
\Big(
d\overline{\pi}^1
+
\isqrt\,\overline{\kappa}\wedge\pi^2
\Big)
\wedge\rho,
\]
and after insertion of twice $-\zeta \wedge \overline{\zeta}$
which is purely imaginary\,\,---\,\,hence disappears\,\,---,
we obtain~({\ref{1-Omega-R-J-W}}):
\[
\aligned
0
&
\,=\,
\Big(
d\pi^1
-
\isqrt\,\kappa\wedge\overline{\pi}^2
-
\zeta\wedge\overline{\zeta}
\Big)
\wedge\rho
+
\Big(
d\overline{\pi}^1
+
\isqrt\,\overline{\kappa}\wedge\pi^2
-
\overline{\zeta}\wedge\zeta
\Big)
\wedge\rho
\\
&
\,=\,
\Omega_1
\wedge\rho
+
\overline{\Omega}_1\wedge\rho.
\endaligned
\]

For~({\ref{2-Omega-R-J-W}}), we proceed analogously, starting from
the second structure equation of
Theorem~{\ref{Theorem-finalization-absorption-d-rho-d-kappa-d-zeta}}:
\[
\aligned
0
&
\,=\,
d\circ d\kappa
\\
&
\,=\,
d\pi^2\wedge\rho
-
\pi^2\wedge d\rho
+
d\pi^1\wedge\kappa
-
\pi^1\wedge d\kappa
+
d\zeta\wedge\overline{\kappa}
-
\zeta\wedge d\overline{\kappa}
\\
&
\,=\,
d\pi^2\wedge\rho
-
\pi^2\wedge
\Big(
\big(\pi^1+\overline{\pi}^1\big)
\wedge\rho
+
\isqrt\,\kappa\wedge\overline{\kappa}
\Big)
+
d\pi^1\wedge\kappa
-
\pi^1\wedge
\Big(
\pi^2\wedge\rho
+
\zeta\wedge\overline{\kappa}
\Big)
\,+
\\
&
\ \ \ \ \
+
\Big(
\big(
\pi^1
-
\overline{\pi}^1
\big)
\wedge\zeta
+
\isqrt\,
\pi^2\wedge\kappa
+
R\,\rho\wedge\zeta
+
W\,\kappa\wedge\zeta
\Big)
\wedge\overline{\kappa}
-
\zeta\wedge
\Big(
\overline{\pi}^2\wedge\rho
+
\overline{\pi}^1\wedge\overline{\kappa}
+
\overline{\zeta}\wedge\kappa
\Big).
\endaligned
\]
After four annihilations by pairs and a reorganization, this becomes:
\[
\!\!\!\!\!\!\!\!\!\!\!\!\!\!\!\!\!\!\!\!\!\!\!\!\!
\small
\aligned
0
&
\,=\,
d\pi^2\wedge\rho
-
\underline{
\pi^2\wedge\pi^1\wedge\rho}_1
-
\pi^2\wedge\overline{\pi}^1\wedge\rho
-
\underline{
\isqrt\,\pi^2\wedge\kappa\wedge\overline{\kappa}}_2
+
d\pi^1\wedge\kappa
-
\underline{
\pi^1\wedge\pi^2\wedge\rho}_1
-
\underline{
\pi^1\wedge\zeta\wedge\overline{\kappa}}_3
\,+
\\
&
\ \ \ \ \
+
\underline{
\pi^1\wedge\zeta\wedge\overline{\kappa}}_3
-
\underline{
\overline{\pi}^1\wedge\zeta\wedge\overline{\kappa}}_4
+
\underline{
\isqrt\,\pi^2\wedge\kappa\wedge\overline{\kappa}}_2
+
R\,\rho\wedge\zeta\wedge\overline{\kappa}
+
W\,\kappa\wedge\zeta\wedge\overline{\kappa}
-
\zeta\wedge\overline{\pi}^2\wedge\rho
-
\underline{
\zeta\wedge\overline{\pi}^1\wedge\overline{\kappa}}_4
-
\zeta\wedge\overline{\zeta}\wedge\kappa
\\
&
\,=\,
\Big(
d\pi^2
-
\pi^2\wedge\overline{\pi}^1
-
\zeta\wedge\overline{\pi}^2
\Big)
\wedge\rho
+
\big(
d\pi^1
-
\zeta\wedge\overline{\zeta}
\big)
\wedge\kappa
\,+
\\
&
\ \ \ \ \
+
R\,\rho\wedge\zeta\wedge\overline{\kappa}
+
W\,\kappa\wedge\zeta\wedge\overline{\kappa},
\endaligned
\]
which is~({\ref{2-Omega-R-J-W}}), 
since we can insert $\big( -\,\isqrt\, \kappa \wedge 
\overline{\pi}^2 \big) \wedge \kappa = 0$.

Lastly:
\[
\aligned
0
&
\,=\,
d\circ d\zeta
\\
&
\,=\,
\isqrt\,d\pi^2\wedge\kappa
-
\isqrt\,\pi^2\wedge d\kappa
+
d\pi^1\wedge\zeta
-
\pi^1\wedge d\zeta
-
d\overline{\pi}^1\wedge\zeta
+
\overline{\pi}^1\wedge d\zeta
\,+
\\
&
\ \ \ \ \
+
dR\wedge\rho\wedge\zeta
+
R\,d\rho\wedge\zeta
-
R\,\rho\wedge d\zeta
\,+
\\
&
\ \ \ \ \
+
dJ\wedge\rho\wedge\overline{\kappa}
+
J\,d\rho\wedge\overline{\kappa}
-
J\,\rho\wedge d\overline{\kappa}
\,+
\\
&
\ \ \ \ \
+
dW\wedge\kappa\wedge\zeta
+
W\,d\kappa\wedge\zeta
-
W\,\kappa\wedge d\zeta,
\endaligned
\]
whence by replacements:
\[
\!\!\!\!\!\!\!\!\!\!\!\!\!
\small
\aligned
0
&
\,=\,
\isqrt\,d\pi^2\wedge\kappa
-
\isqrt\,\pi^2\wedge
\big(
\pi^1\wedge\kappa
+
\zeta\wedge\overline{\kappa}
\big)
+
d\pi^1\wedge\zeta
-
\pi^1\wedge
\Big(
\isqrt\,\pi^2\wedge\kappa
-
\overline{\pi}^1\wedge\zeta
+
R\,\rho\wedge\zeta
+
J\,\rho\wedge\overline{\kappa}
+
W\,\kappa\wedge\zeta
\Big)
\,-
\\
&
\ \ \ \ \
-\,
d\overline{\pi}^1\wedge\zeta
+
\overline{\pi}^1\wedge
\Big(
\isqrt\,\pi^2\wedge\kappa
+
\pi^1\wedge\zeta
+
R\,\rho\wedge\zeta
+
J\,\rho\wedge\overline{\kappa}
+
W\,\kappa\wedge\zeta
\Big)
\,+
\\
&
\ \ \ \ \
+
dR\wedge\rho\wedge\zeta
+
R\,
\Big(
\big(
\pi^1
+
\overline{\pi}^1
\big)\wedge\rho
+
\isqrt\,\kappa\wedge\overline{\kappa}
\Big)\wedge\zeta
-
R\,\rho\wedge
\Big(
\isqrt\,\pi^2\wedge\kappa
+
\big(
\pi^1
-
\overline{\pi}^1
\big)\wedge\zeta
+
W\,\kappa\wedge\zeta
\Big)
\,+
\\
&
\ \ \ \ \
+
dJ\wedge\rho\wedge\overline{\kappa}
+
J\,\big(
\pi^1
+
\overline{\pi}^1
\big)\wedge\rho\wedge\overline{\kappa}
-
J\,\rho\wedge
\big(
\overline{\pi}^1\wedge\overline{\kappa}
+
\overline{\zeta}\wedge\kappa
\big)
\,+
\\
&
\ \ \ \ \ 
+
dW\wedge\kappa\wedge\zeta
+
W\,
\big(
\pi^2\wedge\rho
+
\pi^1\wedge\kappa
\big)\wedge\zeta
-
W\,
\kappa\wedge
\Big(
\big(
\pi^1-\overline{\pi}^1
\big)
\wedge\zeta
+
R\,\rho\wedge\zeta
+
J\,\rho\wedge\overline{\kappa}
\Big).
\endaligned
\]
Let us expand this and underline the eight annihilating pairs:
\[
\!\!\!\!\!\!\!\!\!\!\!\!\!\!\!\!\!\!\!\!\!\!\!\!\!
\scriptsize
\aligned
0
&
\,=\,
\isqrt\,d\pi^2\wedge\kappa
-
\underline{
\isqrt\,\pi^2\wedge\pi^1\wedge\kappa}_1
-
\isqrt\,\pi^2\wedge\zeta\wedge\overline{\kappa}
+
d\pi^1\wedge\zeta
-
\underline{
\isqrt\,\pi^1\wedge\pi^2\wedge\kappa}_1
+
\underline{
\pi^1\wedge\overline{\pi}^1\wedge\zeta}_2
-
\underline{
R\,\pi^1\wedge\rho\wedge\zeta}_3
-
\underline{
J\,\pi^1\wedge\rho\wedge\overline{\kappa}}_6
-
\underline{
W\,\pi^1\wedge\kappa\wedge\zeta}_7
\,-
\\
&
\ \ \ \ \
-\,
d\overline{\pi}^1\wedge\zeta
+
\isqrt\,\overline{\pi}^1\wedge\pi^2\wedge\kappa
+
\underline{
\overline{\pi}^1\wedge\pi^1\wedge\zeta}_2
+
\underline{
R\,\overline{\pi}^1\wedge\rho\wedge\zeta}_4
+
J\,\overline{\pi}^1\wedge\rho\wedge\overline{\kappa}
+
\underline{
W\,\overline{\pi}^1\wedge\kappa\wedge\zeta}_8
\,+
\\
&
\ \ \ \ \
+
dR\wedge\rho\wedge\zeta
+
\underline{
R\,\pi^1\wedge\rho\wedge\zeta}_3
+
R\,\overline{\pi}^1\wedge\rho\wedge\zeta
+
\isqrt\,R\,\kappa\wedge\overline{\kappa}\wedge\zeta
-
\isqrt\,R\,\rho\wedge\pi^2\wedge\kappa
-
R\,\rho\wedge\pi^1\wedge\zeta
+
\underline{
R\,\rho\wedge\overline{\pi}^1\wedge\zeta}_4
-
\underline{
RW\,\rho\wedge\kappa\wedge\zeta}_5
\,+
\\
&
\ \ \ \ \
+
dJ\wedge\rho\wedge\overline{\kappa}
+
\underline{
J\,\pi^1\wedge\rho\wedge\overline{\kappa}}_6
+
J\,\overline{\pi}^1\wedge\rho\wedge\overline{\kappa}
-
J\,\rho\wedge\overline{\pi}^1\wedge\overline{\kappa}
-
J\,\rho\wedge\overline{\zeta}\wedge\kappa
\,+
\\
&
\ \ \ \ \
+
dW\wedge\kappa\wedge\zeta
+
W\,\pi^2\wedge\rho\wedge\zeta
+
\underline{
W\,\pi^1\wedge\kappa\wedge\zeta}_7
-
W\,\kappa\wedge\pi^1\wedge\zeta
+
\underline{
W\,\kappa\wedge\overline{\pi}^1\wedge\zeta}_8
-
\underline{
WR\,\kappa\wedge\rho\wedge\zeta}_5
-
WJ\,\kappa\wedge\rho\wedge\overline{\kappa}.
\endaligned
\]
After simplification and reorganization:
\[
\aligned
0
&
\,=\,
\isqrt\,
\Big(
d\pi^2
-
\pi^2
\wedge\overline{\pi}^1
\Big)
\wedge
\kappa
+
\Big(
d\pi^1
-
d\overline{\pi}^1
-
\isqrt\,\overline{\kappa}\wedge\pi^2
\Big)
\wedge\zeta
\,+
\\
&
\ \ \ \ \
+
dR\wedge\rho\wedge\zeta
+
R\,\overline{\pi}^1\wedge\rho\wedge\zeta
-
\isqrt\,R\,\kappa\wedge\zeta\wedge\overline{\kappa}
+
\isqrt\,R\,\pi^2\wedge\rho\wedge\kappa
+
R\,\pi^1\wedge\rho\wedge\zeta
\,+
\\
&
\ \ \ \ \
+
dJ\wedge\rho\wedge\overline{\kappa}
+
3\,J\,\overline{\pi}^1\wedge\rho\wedge\overline{\kappa}
+
J\,\rho\wedge\kappa\wedge\overline{\zeta}
\,+
\\
&
\ \ \ \ \
+
dW\wedge\kappa\wedge\zeta
+
W\,\pi^2\wedge\rho\wedge\zeta
+
W\,\pi^1\wedge\kappa\wedge\zeta
+
WJ\,\rho\wedge\kappa\wedge\overline{\kappa}.
\endaligned
\]
To reach~({\ref{3-Omega-R-J-W}}) completely, 
only the first line must yet be transformed, 
and it suffices to insert into it two terms
which cancel together:
\[
\isqrt\,
\Big(
d\pi^2
-
\pi^2
\wedge\overline{\pi}^1
-
\zero{
\zeta\wedge\overline{\pi}^2}
\Big)
\wedge
\kappa
+
\Big(
d\pi^1
-
\zero{
\isqrt\,
\kappa\wedge\overline{\pi}^2}
-
d\overline{\pi}^1
-
\isqrt\,\overline{\kappa}\wedge\pi^2
\Big)
\wedge\zeta.
\qedhere
\]
\endproof

Remind that all present considerations hold on the $9$-dimensional 
manifold $M^5 \times G^4$ equipped with the coordinates:
\[
\big(
z_1,z_2,\overline{z}_1,\overline{z}_2,v
\big)
\times
\big(
{\sf c},{\sf e},\overline{\sf c},\overline{\sf e}
\big),
\]
the supplementary real variable ${\sf t} \in \R$ being considered
as a parameter until it becomes a variable
at the very end of the process for an $\{e\}$-structure
on the $10$-dimensional manifold $M^5 \times G^4 \times \R$.
In order to build up such an $\{e\}$-structure, the goal now
is to fully determine the two $2$-forms
$\Omega_1$, $\Omega_2$, and precisely, to 
determine how they express in terms of the coframe:
\[
\big\{
\pi^1,\pi^2,\overline{\pi}^1,\overline{\pi}^2,
\rho,\kappa,\zeta,\overline{\kappa},\overline{\zeta}
\big\}.
\]

To begin with, 
suppose that there are two ways of solving 
for $\big\{ \Omega_1, \Omega_2 \big\}$ the structure
equations of 
Proposition~{\ref{Proposition-3-structure-Omega}}, leading to another
set of solutions $\big\{ \Omega_1', \Omega_2' \big\}$. Then their
differences $\Gamma_1 := \Omega_1' - \Omega_1$ and $\Gamma_2 :=
\Omega_2' - \Omega_2$ must necessarily satisfy the homogeneous
equations:
\[
\aligned
\big(
\Gamma_1+\overline{\Gamma}_1
\big)
\wedge
\rho
&
\,=\,
0,
\\
\Gamma_2\wedge\rho
+
\Gamma_1\wedge\kappa
&
\,=\,
0,
\\
\isqrt\,
\Gamma_2\wedge\kappa
+
\big(
\Gamma_1-\overline{\Gamma}_1
\big)
\wedge
\zeta
&
\,=\,
0.
\endaligned
\]
The addendum to
the article~{\cite{Merker-Pocchiola-2018}},
reproduced on
p.~{\pageref{addendum-Merker-Pocchiola-JGA}},
provides
a detailed proof of the elementary

\begin{Proposition}
\label{Proposition-proved-addendum-MP-2018}
The general solution $\big\{ \Gamma_1, \Gamma_2 \big\}$
to these homogeneous equations is given by:
\[
\Gamma_1
\,:=\,
\Lambda
\wedge\rho,
\ \ \ \ \ \ \ \ \ \ \ \ \ \ \ \ \ \ \ \ \ \ \ \ \ \
\Gamma_2
\,:=\,
\Lambda\wedge\kappa
+
h\,\rho\wedge\kappa,
\]
where $\Lambda$ is a real $1$-form and $h$ is purely imaginary
function.\qed
\end{Proposition}

This means that the two sets of solutions are related to each
other by:
\[
\Omega_1'
\,=\,
\Omega_1
+
\Lambda\wedge\rho,
\ \ \ \ \ \ \ \ \ \ \ \ \ \ \ \ \ \ \ \ \ \ \ \ \ \
\Omega_2'
\,=\,
\Omega_2
+
\Lambda\wedge\kappa
+
h\,\rho\wedge\kappa.
\]

Due to this flexibility represented by $\Lambda$, $h$,
it will be necessary to prolong the structure
equations by adding this real $1$-form:
\[
\Lambda
\,=\,
d{\sf t}
+
\cdots,
\]
the remainder terms being very complicated, while the function $h$
{\em could} be some new invariant. However, it will be later shown
that $h$ expresses in terms of the 3\textsuperscript{rd}-order jets
of $W$ and $J$, thus eliminating the possibility
of appearance of new primary CR invariants. On the other hand,
the existence of $\Lambda$ can be explained by an application (not
detailed here) of Cartan's test, due to the fact that there is one
degree of real-valued indeterminancy during the fourth absorption.

It therefore suffices to find a particular set of solution $\Omega_1$
and $\Omega_2$, and then to parametrize the solution space by
means of
$\Lambda$, $h$.  We will adopt the following strategy. First, we
will find the simplest forms for $\Omega_1$ and $\Omega_2$ restrained
by the first two equations~({\ref{1-Omega-R-J-W}}),
({\ref{2-Omega-R-J-W}}) of the starting
Proposition~{\ref{Proposition-3-structure-Omega}}.  Then we will
simplify these 
$2$-forms by means of Cartan's lemma to eliminate as many
unknown variables as possible using the third, more subtle,
equation~({\ref{3-Omega-R-J-W}}).  At the end of the elimination,
those remaining unknowns which cannot be computed due to the lack of
information turn out to behave like $\Lambda$ and $h$, and hence we
will terminate the process of solving for solutions.

In $M^5 \times G^4$, it will be useful to adopt the following 
notations for the covariant derivatives:
\leqnomode\usetagform{default}
\begin{small}
\begin{align}
\label{covariant-derivatives-dR-dJ-dW}
dR
&
\,=\,
R_{\pi^1}\,\pi^1
+
R_{\pi^2}\,\pi^2
+
R_{\overline{\pi}^1}\,\overline{\pi}^1
+
R_{\overline{\pi}^2}\,\overline{\pi}^2
+
R_\rho\,\rho
+
R_\kappa\,\kappa
+
R_\zeta\,\zeta
+
R_{\overline{\kappa}}\,\overline{\kappa}
+
R_{\overline{\zeta}}\,\overline{\zeta},
\notag
\\
dJ
&
\,=\,
J_{\pi^1}\,\pi^1
+
J_{\pi^2}\,\pi^2
+
J_{\overline{\pi}^1}\,\overline{\pi}^1
+
J_{\overline{\pi}^2}\,\overline{\pi}^2
+
J_\rho\,\rho
+
J_\kappa\,\kappa
+
J_\zeta\,\zeta
+
J_{\overline{\kappa}}\,\overline{\kappa}
+
J_{\overline{\zeta}}\,\overline{\zeta},
\\
\ \ \ \ 
dW
&
\,=\,
W_{\pi^1}\,\pi^1
+
W_{\pi^2}\,\pi^2
+
W_{\overline{\pi}^1}\,\overline{\pi}^1
+
W_{\overline{\pi}^2}\,\overline{\pi}^2
+
W_\rho\,\rho
+
W_\kappa\,\kappa
+
W_\zeta\,\zeta
+
W_{\overline{\kappa}}\,\overline{\kappa}
+
W_{\overline{\zeta}}\,\overline{\zeta}.
\notag
\end{align}
\end{small}\noindent
Some of these coefficients will be revealed during the course
of solving the structure equations. We first turn ourselves
to finding the simplest form of $\Omega_1$, $\Omega_2$
satisfying {\em only} the first two 
equations~({\ref{1-Omega-R-J-W}}), ({\ref{2-Omega-R-J-W}}).

\begin{Proposition}
\label{Proposition-representations-Omega-1-2}
There exists a real-valued function $p$ and two
differential $1$-forms $\Pi$, $\Psi$ such that:
\[
\aligned
\Omega_1
&
\,=\,
\Pi\wedge\rho
+
p\,\kappa\wedge\overline{\kappa}
-
\overline{W}\,\kappa\wedge\overline{\zeta}
-
W\,\zeta\wedge\overline{\kappa},
\\
\Omega_2
&
\,=\,
\Psi\wedge\rho
+
\Pi\wedge\kappa
-
R\,\zeta\wedge\overline{\kappa}.
\endaligned
\]
\end{Proposition}

\proof 
We can rearrange the terms in~({\ref{2-Omega-R-J-W}}):
\leqnomode\usetagform{default}
\begin{align}
\label{Omega-kappa-rho}
0
\,=\,
\big(
\Omega_1
+
W\,\zeta\wedge\overline{\kappa}
\big)
\wedge\kappa
+
\big(
\Omega_2
+
R\,\zeta\wedge\overline{\kappa}
\big)
\wedge\rho,
\end{align}
in order that an application of the Cartan Lemma yield functions
$\Delta$, $\Theta$, $\Pi''$, $\Psi$ so that:
\[
\aligned
\Omega_1
+
W\,\zeta\wedge\overline{\kappa}
&
\,=\,
\Delta\wedge\kappa
+
\Theta\wedge\rho,
\\
\Omega_2
+
R\,\zeta\wedge\overline{\kappa}
&
\,=\,
\Pi''\wedge\kappa
+
\Psi\wedge\rho,
\endaligned
\]
with a double prime on $\Pi''$ meaning that we will soon modify
it two times. 

In fact, substituting these representations back 
into~({\ref{Omega-kappa-rho}}), 
we see that there are constraints on $\Theta$
and $\Pi''$:
\[
\aligned
0
&
\,=\,
\big(
\zero{
\Delta\wedge\kappa}
+
\Theta\wedge\rho
\big)
\wedge\kappa
+
\big(
\Pi''\wedge\kappa
+
\zero{
\Psi\wedge\rho}
\big)\wedge\rho
\\
&
\,=\,
\big(
\Theta
-
\Pi''
\big)
\wedge\rho\wedge\kappa.
\endaligned
\]
By the Cartan Lemma again, this implies the existence of 
two functions
$a$, $b$ so that $\Theta$ and $\Pi''$ are related to each other
by:
\[
\Theta
\,=\,
\Pi''
+
a\,\rho
+
b\,\kappa.
\]
Next, putting this into the expression of $\Omega_1$, while
letting $\Pi' := \Pi'' + b\, \kappa$, it follows that:
\[
\aligned
\Omega_1
&
\,=\,
\Delta\wedge\kappa
+
\Theta\wedge\rho
-
W\,\zeta\wedge\overline{\kappa}
\\
&
\,=\,
\Delta\wedge\kappa
+
\big(
\Pi''
+
\zero{
a\,\rho}
+
b\,\kappa
\big)
\wedge\rho
-
W\,\zeta\wedge\overline{\kappa}
\\
&
\,=\,
\Delta\wedge\kappa
+
\Pi'\wedge\rho
-
W\,\zeta\wedge\overline{\kappa},
\endaligned
\]
while $\Omega_2$ becomes:
\[
\aligned
\Omega_2
&
\,=\,
\Pi''\wedge\kappa
+
\Psi\wedge\rho
-
R\,\zeta\wedge\overline{\kappa}
\\
&
\,=\,
\big(
\Pi''+b\,\kappa
\big)
\wedge\kappa
+
\Psi\wedge\rho
-
R\,\zeta\wedge\overline{\kappa}
\\
&
\,=\,
\Pi'\wedge\kappa
+
\Psi\wedge\rho
-
R\,\zeta\wedge\overline{\kappa}.
\endaligned
\]

The next observation is that $\Delta$ can be further simplified.
Indeed, let us replace $\Omega_1$ in~({\ref{1-Omega-R-J-W}}):
\[
\aligned
0
&
\,=\,
\big(
\Omega_1+\overline{\Omega}_1
\big)
\wedge\rho
\\
&
\,=\,
\Delta\wedge\kappa\wedge\rho
-
W\,\zeta\wedge\overline{\kappa}\wedge\rho
+
\overline{\Delta}\wedge\overline{\kappa}\wedge\rho
-
\overline{W}\,\overline{\zeta}\wedge\kappa\wedge\rho.
\endaligned
\]
Then decomposing $\Delta$ as a linear combination along the
coframe:
\[
\Delta
\,=\,
d_1\,\pi^1
+
d_2\,\pi^2
+
d_3\,\overline{\pi}^1
+
d_4\,\overline{\pi}^2
+
d_5\,\rho
+
d_6\,\kappa
+
d_7\,\zeta
+
d_8\,\overline{\kappa}
+
d_9\,\overline{\zeta},
\]
we obtain the following values for these coefficients:
\[
d_1
\,=\,
d_2
\,=\,
d_3
\,=\,
d_4
\,=\,
0,
\ \ \ \ \
d_8
\,=\,
\overline{d}_8,
\ \ \ \ \
d_9
\,=\,
\overline{W},
\]
except for $d_5$ and $d_6$ which on which no constraint is deduced 
so, and hence:
\[
\Delta
\,=\,
d_5\,\rho
+
d_6\,\kappa
+
d_8\,\overline{\kappa}
+
\overline{W}\,\overline{\zeta}.
\]

Finally, if we write $p := -\, d_8$ and if we set $\Pi := \Pi' -
d_5\, \kappa$, we obtain by reorganization:
\[
\aligned
\Omega_1
&
\,=\,
\Delta\wedge\kappa
+
\Pi'\wedge\rho
-
W\,\zeta\wedge\overline{\kappa}
\\
&
\,=\,
\Big(
d_5\,\rho
+
\zero{
d_6\,\kappa}
+
d_8\,\overline{\kappa}
+
\overline{W}\,\overline{\zeta}
\Big)
\wedge\kappa
+
\Pi'\wedge\rho
-
W\,\zeta\wedge\overline{\kappa}
\\
&
\,=\,
-\,d_8\,\kappa\wedge\overline{\kappa}
+
\big(
\Pi'
-
d_5\,\kappa
\big)
\wedge\rho
-
\overline{W}\,\kappa\wedge\overline{\zeta}
-
W\,\zeta\wedge\overline{\kappa}
\\
&
\,=\,
p\,\kappa\wedge\overline{\kappa}
+
\Pi\wedge\rho
-
\overline{W}\,\kappa\wedge\overline{\zeta}
-
W\,\zeta\wedge\overline{\kappa},
\endaligned
\]
and moreover:
\begin{align}
\Omega_2
&
\,=\,
\Psi\wedge\rho
+
\Pi'\wedge\kappa
-
R\,\zeta\wedge\overline{\kappa}
\notag
\\
&
\,=\,
\Psi\wedge\rho
+
\big(
\Pi'
-
d_5\,\kappa
\big)
\wedge\kappa
-
R\,\zeta\wedge\overline{\kappa}
\notag
\\
&
\,=\,
\Psi\wedge\rho
+
\Pi\wedge\kappa
-
R\,\zeta\wedge\overline{\kappa}.
\qedhere
\end{align}
\endproof

Now, using the representations of $\Omega_1$ and of $\Omega_2$
offered by this 
Proposition~{\ref{Proposition-representations-Omega-1-2}},
we can therefore rewrite the third (still not taken account of)
equation~({\ref{3-Omega-R-J-W}}) as:
\leqnomode\usetagform{default}
\begin{align}
\label{start-dR-dJ-dW}
\isqrt\,
\Psi\wedge\rho\wedge\kappa
&
-
\isqrt\,R\,\zeta\wedge\overline{\kappa}\wedge\kappa
+
\big(
\Pi
-
\overline{\Pi}
\big)
\wedge\rho\wedge\zeta
+
2\,p\,
\kappa\wedge\overline{\kappa}\wedge\zeta
-
2\,\overline{W}\,\kappa\wedge\overline{\zeta}\wedge\zeta
\,=\,
\notag
\\
&
\,=\,
-\,dR\wedge\rho\wedge\zeta
-
R\,\big(\pi^1+\overline{\pi}^1\big)
\wedge\rho\wedge\zeta
-
\isqrt\,R\,\pi^2\wedge\rho\wedge\kappa
+
\isqrt\,R\,\kappa\wedge\zeta\wedge\overline{\kappa}
\,-
\\
&
\ \ \ \ \
-\,
dJ\wedge\rho\wedge\overline{\kappa}
-
3J\,\overline{\pi}^1\wedge\rho\wedge\overline{\kappa}
-
J\,\rho\wedge\kappa\wedge\overline{\zeta}
\,-
\notag
\\
&
\ \ \ \ \
-\,
dW\wedge\kappa\wedge\zeta
-
W\,\pi^2\wedge\rho\wedge\zeta
-
W\,\pi^1\wedge\kappa\wedge\zeta
-
WJ\,\rho\wedge\kappa\wedge\overline{\kappa}.
\notag
\end{align}

But before we commence with analyzing this equation (a long task),
we make a side remark. As we can rewrite:
\[
\aligned
\Omega_1
&
\,=\,
{\textstyle{\frac{1}{2}}}\,
\big(
\Pi+\overline{\Pi}
\big)
\wedge\rho
+
{\textstyle{\frac{1}{2}}}\,
\big(
\Pi-\overline{\Pi}
\big)
\wedge\rho
+
p\,\kappa\wedge\overline{\kappa}
-
\overline{W}\,
\kappa\wedge\zeta
-
W\,\zeta\wedge\overline{\kappa},
\\
\Omega_2
&
\,=\,
\Psi\wedge\rho
+
{\textstyle{\frac{1}{2}}}\,
\big(
\Pi+\overline{\Pi}
\big)
\wedge\kappa
+
{\textstyle{\frac{1}{2}}}\,
\big(
\Pi-\overline{\Pi}
\big)
\wedge\kappa
-
R\,\zeta\wedge\overline{\kappa},
\endaligned
\]
we remark that 
Proposition~{\ref{Proposition-proved-addendum-MP-2018}}
already tells us that the real part 
$\frac{1}{2}\,\big(\Pi + \overline{\Pi} \big)$ of $\Pi$
is {\em a priori} not fully determined, as can be formulated by an

\begin{Observation}
For an arbitrary real $1$-form $\Lambda$, the $2$-forms:
\[
\Omega_1'
\,:=\,
\Omega_1
+
\Lambda\wedge\rho
\ \ \ \ \ \ \ \ \ \ \ \ \ \ \ \ \ \
\text{and}
\ \ \ \ \ \ \ \ \ \ \ \ \ \ \ \ \ \
\Omega_2'
\,:=\,
\Omega_2
+
\Lambda\wedge\kappa
\]
still satisfy the structure equations of
Proposition~{\ref{Proposition-3-structure-Omega}}.
\end{Observation}

\proof
For the sake of completeness, let us detail the arguments.
The first equation~({\ref{1-Omega-R-J-W}}) is clear:
\[
\big(
\Omega_1'+\overline{\Omega}_1'
\big)
\wedge\rho
\,=\,
\Big(
\Omega_1
+
\zero{
\Lambda\wedge\rho}
+
\overline{\Omega}_1
+
\zero{
\Lambda\wedge\rho}
\Big)
\wedge\rho
\,=\,
\big(
\Omega_1
+
\overline{\Omega}_1
\big)
\wedge
\rho.
\]

The second equation~({\ref{2-Omega-R-J-W}}) also:
\[
\aligned
\Omega_2'
\wedge\rho
+
\Omega_1'\wedge\kappa
&
\,=\,
\big(
\Omega_2+\Lambda\wedge\kappa
\big)
\wedge\rho
+
\big(
\Omega_1
+
\Lambda\wedge\rho
\big)
\wedge\kappa
\\
&
\,=\,
\Omega_2\wedge\rho
+
\zero{
\Lambda\wedge\kappa\wedge\rho}
+
\Omega_1\wedge\kappa
+
\zero{
\Lambda\wedge\rho\wedge\kappa}
\,=\,
\Omega_2\wedge\rho
+
\Omega_1\wedge\kappa,
\endaligned
\]
and the third one as well:
\begin{align}
\isqrt\,\Omega_2'\wedge\kappa
+
\big(
\Omega_1'
-
\overline{\Omega}_1'
\big)\wedge\zeta
&
\,=\,
\isqrt\,
\big(
\Omega_2
+
\zero{
\Lambda
\wedge\kappa}
\big)
\wedge\kappa
+
\Big(
\Omega_1
+
\zero{
\Lambda\wedge\rho}
-
\overline{\Omega}_1
-
\zero{
\Lambda\wedge\rho}
\Big)
\wedge\zeta
\notag
\\
&
\,=\,
\isqrt\,\Omega_2\wedge\kappa
+
\big(
\Omega_1
-
\overline{\Omega}_1
\big)
\wedge\zeta.
\qedhere
\end{align}
\endproof

Now, coming back to~({\ref{start-dR-dJ-dW}}),
we remember that we should insert the covariant derivatives
$dR$, $dJ$, $dW$ from~({\ref{covariant-derivatives-dR-dJ-dW}}),
and we will do this {\em in a progressive way}, not in one stroke.

Indeed, by wedging $(\centersmallbullet) \wedge \rho$
both sides of~({\ref{start-dR-dJ-dW}}), 
we get rid of $dJ$, $dR$ and it remains 
only:
\[
\aligned
-\,\isqrt\,R\,
&
\zeta\wedge\overline{\kappa}\wedge\kappa\wedge\rho
+
2p\,\kappa\wedge\overline{\kappa}\wedge\zeta\wedge\rho
-
2\overline{W}\,
\kappa\wedge\overline{\zeta}\wedge\zeta\wedge\rho
\,=\,
\\
&
\ \ \ \ \
\,=\,
\isqrt\,R\,
\kappa\wedge\zeta\wedge\overline{\kappa}\wedge\rho
-
dW\wedge\kappa\wedge\zeta\wedge\rho
-
W\,\pi^1\wedge\kappa\wedge\zeta\wedge\rho,
\endaligned
\]
that is to say after putting everything to the right:
\[
0
\,=\,
-\,dW\wedge\rho\wedge\kappa\wedge\zeta
-
\big(
2\,p
+
2\,\isqrt\,R
\big)\,
\rho\wedge\kappa\wedge\zeta\wedge\overline{\kappa}
+
2\overline{W}\,
p\wedge\kappa\wedge\zeta\wedge\overline{\zeta}
-
W\,\pi^1\wedge\rho\wedge\kappa\wedge\zeta.
\]

Thus, inserting the expansion of $dW$ 
from~({\ref{covariant-derivatives-dR-dJ-dW}}):
\[
\!\!\!\!\!\!\!\!\!\!\!\!\!\!\!\!\!\!\!\!
\small
\aligned
-\,dW\wedge\rho\wedge\kappa\wedge\zeta
\,=\,
-\,
W_{\pi^1}\,\pi^1\wedge\rho\wedge\kappa\wedge\zeta
-
W_{\pi^2}\,\pi^2\wedge\rho\wedge\kappa\wedge\zeta
-
W_{\overline{\pi}^1}\,
\overline{\pi}^1\wedge\rho\wedge\kappa\wedge\zeta
-
W_{\overline{\pi}^2}\,
\overline{\pi}^2\wedge\rho\wedge\kappa\wedge\zeta
\,-
\\
-\,
W_{\overline{\kappa}}\,
\overline{\kappa}\wedge\rho\wedge\kappa\wedge\zeta
-
W_{\overline{\zeta}}\,
\overline{\zeta}\wedge\rho\wedge\kappa\wedge\zeta,
\endaligned
\]
we get:
\[
\small
\aligned
0
&
\,=\,
-\,
\big(
W_{\pi^1}+W
\big)\,
\pi^1\wedge\rho\wedge\kappa\wedge\zeta
-
W_{\pi^2}\,
\pi^2\wedge\rho\wedge\kappa\wedge\zeta
-
W_{\overline{\pi}^1}\,
\overline{\pi}^1\wedge\rho\wedge\kappa\wedge\zeta
-
W_{\overline{\pi}^2}\,
\overline{\pi}^2\wedge\rho\wedge\kappa\wedge\zeta
\,-
\\
&
\ \ \ \ \
-\,
\big(
2\,p
+
2\,\isqrt\,R
-
W_{\overline{\kappa}}
\big)\,
\rho\wedge\kappa\wedge\zeta\wedge\overline{\kappa}
-
\big(
2\,\overline{W}
+
W_{\overline{\zeta}}
\big)\,
\rho\wedge\kappa\wedge\zeta\wedge\overline{\zeta},
\endaligned
\]
whence by identification of coefficients of these independent
$4$-forms:
\[
\aligned
W_{\pi^1}
&
\,=\,
-\,W,
\ \ \ \ \ \ \ \ \ \ \ \ \ \ \ \ \ \ \ \ 
W_{\pi^2}
\,=\,
0,
\ \ \ \ \ \ \ \ \ \ 
W_{\overline{\pi}^1}
\,=\,
0,
\ \ \ \ \ \ \ \ \ \ 
W_{\overline{\pi}^2}
\,=\,
0,
\\
W_{\overline{\kappa}}
&
\,=\,
2\,p
+
2\,\isqrt\,R,
\ \ \ \ \ \ \ \ \ \ 
W_{\overline{\zeta}}
\,=\,
-\,2\,\overline{W},
\endaligned
\]
while no condition is imposed so on $W_\rho$, $W_\kappa$,
$W_\zeta$, and thus:
\[
dW
\,=\,
-\,W\,\pi^1
+
W_\rho\,\rho
+
W_\kappa\,\kappa
+
W_\zeta\,\zeta
+
\big(
2\,p
+
2\,\isqrt\,R
\big)\,
\overline{\kappa}
-
2\,\overline{W}\,
\overline{\zeta}.
\]

Next, putting this expression of $dW$ back 
into~({\ref{start-dR-dJ-dW}}) 
allows us to eliminate $p$ so that we can focus
only on $\Pi - \overline{\Pi}$ and $\Psi$, which we
place on the left:
\[
\!\!\!\!\!\!\!\!\!\!\!\!\!\!\!\!\!\!\!\!\!\!\!\!\!
\small
\aligned
\isqrt\,\Psi\wedge\rho\wedge\kappa
+
\big(\Pi-\overline{\Pi}\big)\wedge\rho\wedge\zeta
&
\,=\,
\underline{
\isqrt\,R\,
\zeta\wedge\overline{\kappa}\wedge\kappa}_1
-
\underline{
2\,p\,
\kappa\wedge\overline{\kappa}\wedge\zeta}_2
-
\underline{
2\,\overline{W}\,
\kappa\wedge\overline{\zeta}\wedge\zeta}_3
\,-
\\
&
\ \ \ \ \
-\,
dR\wedge\rho\wedge\zeta
-
R\,\big(\pi^1+\overline{\pi}^1\big)
\wedge\rho\wedge\zeta
-
\isqrt\,R\,
\pi^2\wedge\rho\wedge\kappa
+
\underline{
\isqrt\,R\,
\kappa\wedge\zeta\wedge\overline{\kappa}}_1
\,-
\\
&
\ \ \ \ \
-\,
dJ\wedge\rho\wedge\overline{\kappa}
-
3J\,\overline{\pi}^1\wedge\rho\wedge\overline{\kappa}
-
J\,\rho\wedge\kappa\wedge\overline{\zeta}
\,+
\\
&
\ \ \ \ \
+
\underline{
W\,\pi^1\wedge\kappa\wedge\zeta}_4
-
W_\rho\,\rho\wedge\kappa\wedge\zeta
-
\big(
\underline{
2\,p}_2
+
\underline{
2\,\isqrt\,R}_1
\big)\,
\overline{\kappa}\wedge\kappa\wedge\zeta
+
\underline{
2\,\overline{W}\,
\overline{\zeta}\wedge\kappa\wedge\zeta}_3
\,-
\\
&
\ \ \ \ \
-\,
W\,\pi^2\wedge\rho\wedge\zeta
-
\underline{
W\,\pi^1\wedge\kappa\wedge\zeta}_4
-
WJ\,\rho\wedge\kappa\wedge\overline{\kappa}.
\endaligned
\]
Here, four simplifications by pairs are underlined, in which 
we observe that $p$ eliminates itself, and 
if we collect at first the terms divisible by $\rho \wedge \kappa$,
we get:
\[
\aligned
\isqrt\,\Psi\wedge\rho\wedge\kappa
+
\big(
\Pi-\overline{\Pi}
\big)
\wedge\rho\wedge\zeta
&
\,=\,
\Big(
-\,\isqrt\,R\,\pi^2
-
J\,\overline{\zeta}
-
W_\rho\,\zeta
-
WJ\,\overline{\kappa}
\Big)
\wedge\rho\wedge\kappa
\,-
\\
&
\ \ \ \ \
-\,
dR\wedge\rho\wedge\zeta
-
R\,
\big(
\pi^1+\overline{\pi}^1
\big)
\wedge\rho\wedge\zeta
\,-
\\
&
\ \ \ \ \
-\,
dJ\wedge\rho\wedge\overline{\kappa}
-
3J\,\overline{\pi}^1\wedge\rho\wedge\overline{\kappa}
-
W\,\pi^2\wedge\rho\wedge\zeta.
\endaligned
\]
By introducing the modified $1$-form:
\[
\Psi'
\,:=\,
\Psi
-
\isqrt\,
\Big(
\isqrt\,R\,\pi^2
+
J\,\overline{\zeta}
+
W_\rho\,\zeta
+
WJ\,\overline{\kappa}
\Big),
\]
the equation becomes:
\leqnomode\usetagform{default}
\begin{align}
\isqrt\,
\Psi'\wedge\rho\wedge\kappa
+
\big(
\Pi-\overline{\Pi}
\big)
\wedge\rho\wedge\zeta
&
\,=\,
-\,dR\wedge\rho\wedge\zeta
-
R\,\big(\pi^1+\overline{\pi}^1\big)
\wedge\rho\wedge\zeta
\,-
\notag
\\
\label{i-Psi-prime-dJ}
&
\ \ \ \ \
-\,
dJ\wedge\rho\wedge\overline{\kappa}
-
3J\,\overline{\pi}^1\wedge\rho\wedge\overline{\kappa}
-
W\,\pi^2\wedge\rho\wedge\zeta.
\end{align}
Now, let us wedge $(\centersmallbullet) \wedge \kappa \wedge\zeta$
all this to make $\Psi$ and $\Pi - \overline{\Pi}$ disappear,
replacing simultaneously:
\[
dJ
\,=\,
J_{\pi^1}\,\pi^1
+
J_{\pi^2}\,\pi^2
+
J_{\overline{\pi}^1}\,\overline{\pi}^1
+
J_{\overline{\pi}^2}\,\overline{\pi}^2
+
J_\rho\,\rho
+
J_\kappa\,\kappa
+
J_\zeta\,\zeta
+
J_{\overline{\kappa}}\,\overline{\kappa}
+
J_{\overline{\zeta}}\,\overline{\zeta},
\]
to obtain:
\[
\!\!\!\!\!\!\!\!\!\!\!\!\!
\small
\aligned
0
&
\,=\,
-\,
J_{\pi^1}\,
\pi^1\wedge\rho\wedge\overline{\kappa}\wedge\kappa\wedge\zeta
-
J_{\pi^2}\,
\pi^2\wedge\rho\wedge\overline{\kappa}\wedge\kappa\wedge\zeta
-
J_{\overline{\pi}^1}\,
\overline{\pi}^1\wedge\rho\wedge\overline{\kappa}
\wedge\kappa\wedge\zeta
-
J_{\overline{\pi}^2}\,
\overline{\pi}^2\wedge\rho\wedge\overline{\kappa}
\wedge\kappa\wedge\zeta
\,-
\\
&
\ \ \ \ \ \ \ \ \ \ \ \ \ \ \ \ \ \ \ \ \ \ \ \ \ \
\ \ \ \ \ \ \ \ \ \ \ \ \ \ \ \ \ \ \ \ \ \ \ \ \ \
\ \ \ \ \ \ \ \ \ \ \ \ \ \ \ \ \ \ \ \ \ \ \ \ \ \
\ \ \ \ \ \ \ \ \ \ \ \ \ \ 
-\,
J_{\overline{\zeta}}\,
\overline{\zeta}\wedge\rho\wedge\overline{\kappa}
\wedge\kappa\wedge\zeta
-
3J\,
\overline{\pi}^1\wedge\rho\wedge\overline{\kappa}
\wedge\kappa\wedge\zeta
\\
&
\,=\,
-\,
J_{\pi^1}\,
\pi^1\wedge\rho\wedge\kappa\wedge\zeta\wedge\overline{\kappa}
-
J_{\pi^2}\,
\pi^2\wedge\rho\wedge\kappa\wedge\zeta\wedge\overline{\kappa}
-
\big(
J_{\overline{\pi}^1}
+
3\,J
\big)\,
\overline{\pi}^1
\wedge\rho\wedge\kappa\wedge\zeta\wedge\overline{\kappa}
-
J_{\overline{\pi}^2}\,
\overline{\pi}^2
\wedge\rho\wedge\kappa\wedge\zeta\wedge\overline{\kappa}
\,-
\\
&
\ \ \ \ \ \ \ \ \ \ \ \ \ \ \ \ \ \ \ \ \ \ \ \ \ \
\ \ \ \ \ \ \ \ \ \ \ \ \ \ \ \ \ \ \ \ \ \ \ \ \ \
\ \ \ \ \ \ \ \ \ \ \ \ \ \ \ \ \ \ \ \ \ \ \ \ \ \
\ \ \ \ \ \ \ \ \ \ \ \ \ \ 
-\,
J_{\overline{\zeta}}\,
\rho\wedge\kappa\wedge\zeta\wedge\overline{\kappa}
\wedge\overline{\zeta},
\endaligned
\]
and since these $5$-forms are linearly independent, we get by
identification:
\[
J_{\pi^1}
\,=\,
0,
\ \ \ \ \ \ \ \ \ \ 
J_{\pi^2}
\,=\,
0,
\ \ \ \ \ \ \ \ \ \ 
J_{\overline{\pi}^1}
\,=\,
-\,3J,
\ \ \ \ \ \ \ \ \ \
J_{\overline{\pi}^2}
\,=\,
0,
\ \ \ \ \ \ \ \ \ \
J_{\overline{\zeta}}
\,=\,
0,
\]
while no condition is imposed in this way on $J_\rho$, $J_\kappa$,
$J_\zeta$, $J_{\overline{\kappa}}$.
Consequently, the $1$-form $dJ$ contracts as:
\[
dJ
\,=\,
-\,3J\,\overline{\pi}^1
+
J_\rho\,\rho
+
J_\kappa\,\kappa
+
J_\zeta\,\zeta
+
J_{\overline{\kappa}}\,
\overline{\kappa},
\]
hence putting this expression back 
into~({\ref{i-Psi-prime-dJ}}), we obtain:
\[
\!\!\!\!\!\!\!\!\!\!\!\!\!\!\!\!\!\!\!\!\!\!
\small
\aligned
\isqrt\,\Psi'\wedge\rho\wedge\kappa
+
\big(
\Pi
-
\overline{\Pi}
\big)
\wedge\rho\wedge\zeta
&
\,=\,
-\,dR\wedge\rho\wedge\zeta
-
R\,\big(\pi^1+\overline{\pi}^1\big)
\wedge\rho\wedge\zeta
\,+
\\
&
\ \ \ \ \
+
\zero{
3J\,\overline{\pi}^1\wedge\rho\wedge\overline{\kappa}}
-
J_\kappa\,
\kappa\wedge\rho\wedge\overline{\kappa}
-
J_\zeta\,\zeta\wedge\rho\wedge\overline{\kappa}
-
\zero{
3J\,\overline{\pi}^1\wedge\rho\wedge\overline{\kappa}}
-
W\,\pi^2\wedge\rho\wedge\zeta.
\endaligned
\]
We can yet absorb in $\Psi'$ one term from the right-hand side
by introducing:
\[
\Psi''
\,:=\,
\Psi'
+
\isqrt\,
J_\kappa\,
\overline{\kappa},
\]
so that our equation becomes:
\[
\aligned
\isqrt\,\Psi''
\wedge\rho\wedge\kappa
+
\big(
\Pi
-
\overline{\Pi}
\big)
\wedge\rho\wedge\zeta
&
\,=\,
-\,dR
\wedge\rho\wedge\zeta
-
R\,
\big(
\pi^1
+
\overline{\pi}^1
\big)
\wedge\rho\wedge\zeta
\,+
\\
&
\ \ \ \ \
+
J_\zeta\,
\rho\wedge\zeta\wedge\overline{\kappa}
-
W\,\pi^2\wedge\rho\wedge\zeta.
\endaligned
\]

Now, observe that all terms except the first one 
$\isqrt\, \Psi'' \wedge\rho \wedge\kappa$ are multiple
of $\rho \wedge\zeta$. Consequently, wedging on both sides
by $(\centersmallbullet) \wedge \zeta$, we annihilate
everything except:
\[
\isqrt\,
\Psi''\wedge\rho\wedge\kappa\wedge\zeta
\,=\,
0.
\]
Thanks to the Cartan Lemma, there exist function $e$, $f$, $g$
so that:
\[
\Psi''
\,=\,
e\,\rho
+
f\,\kappa
+
g\,\zeta.
\]
For later use, we also observe in passing that:
\leqnomode\usetagform{default}
\begin{align}
\label{expression-Psi}
\Psi
&
\,=\,
\Psi'
+
\isqrt\,
W_\rho\,\zeta
+
\isqrt\,WJ\,
\overline{\kappa}
-
R\,\pi^2
+
\isqrt\,
J\,\overline{\zeta}
\\
&
\,=\,
\Psi''
-
\isqrt\,
J_\kappa\,
\overline{\kappa}
+
\isqrt\,
W_\rho\,\zeta
+
\isqrt\,WJ\,
\overline{\kappa}
-
R\,\pi^2
+
\isqrt\,
J\,\overline{\zeta}
\notag
\\
&
\,=\,
-\,R\,\pi^2
+
e\,\rho
+
f\,\kappa
+
\big(
\isqrt\,W_\rho
+
g
\big)\,
\zeta
+
\isqrt\,
\big(
WJ
-
J_\kappa
\big)\,
\overline{\kappa}
+
\isqrt\,J\,\overline{\zeta}.
\notag
\end{align}
Inserting this just above conducts to an identity:
\[
\aligned
\isqrt\,
g\,
\rho\wedge\kappa\wedge\zeta
+
\big(
\Pi
-
\overline{\Pi}
\big)
\wedge\rho\wedge\zeta
&
\,=\,
-\,
dR\wedge\rho\wedge\zeta
-
R\,\big(\pi^1+\overline{\pi}^1\big)
\wedge\rho\wedge\zeta
\,+
\\
&
\ \ \ \ \
+
J_\zeta\,
\rho\wedge\zeta\wedge\overline{\kappa}
-
W\,\pi^2\wedge\rho\wedge\zeta,
\endaligned
\]
in which {\em all} terms are now multiples of $\rho \wedge \zeta$.
Consequently, the Cartan Lemma implies the existence of functions
$r$ and $s$ such that:
\[
\Pi
-
\overline{\Pi}
\,=\,
\isqrt\,g\,\kappa
-
dR
-
R\,\pi^1
-
R\,\overline{\pi}^1
+
J_\zeta\,\overline{\kappa}
-
W\,\pi^2
+
r\,\rho
+
s\,\zeta.
\]
But here, we can take advantage of the fact that $\Pi - 
\overline{\Pi}$ is purely imaginary to obtain some information
about $g$, $r$, $s$. Indeed, conjugating:
\[
\overline{\Pi}
-
\Pi
\,=\,
-\,\isqrt\,\overline{g}\,\overline{\kappa}
-
dR
-
R\,\overline{\pi}^1
-
R\,\pi^1
-
\overline{J_\zeta}\,\kappa
-
\overline{W}\,\overline{\pi}^2
+
\overline{r}\,\rho
+
\overline{s}\,\overline{\zeta},
\]
and summing, we eliminate $\Pi - \overline{\Pi}$, hence
we are left after reorganization with:
\[
\aligned
0
&
\,=\,
-\,2\,dR
-
2R\,\pi^1
-
W\,\pi^2
-
2R\,\overline{\pi}^1
-
\overline{W}\,\overline{\pi}^2
\,+
\\
&
\ \ \ \ \
+
(r+\overline{r})\,
\rho
+
\big(
\isqrt\,g
+
\overline{J_\zeta}
\big)\,
\kappa
+
s\,\zeta
+
\big(
-\isqrt\,\overline{g}
+
J_\zeta
\big)\,
\overline{\kappa}
+
\overline{s}\,
\overline{\zeta}.
\endaligned
\]
Naturally, one has to use the expansion of $dR$
from~({\ref{covariant-derivatives-dR-dJ-dW}}) 
to continue the computation:
\[
\small
\aligned
0
&
\,=\,
-\,
\big(
2\,R_{\pi^1}+2\,R
\big)\,\pi^1
-
\big(
2\,R_{\pi^2}+W
\big)\,\pi^2
-
\big(
2\,R_{\overline{\pi}^1}+2\,R
\big)\,
\overline{\pi}^1
-
\big(
2\,R_{\overline{\pi}^2}+\overline{W}
\big)\,
\overline{\pi}^2
\,-
\\
&
\ \ \ \ \ 
-\,
\big(
2\,R_\rho-r-\overline{r}
\big)\,\rho
-
\big(
2\,R_\kappa-\isqrt\,g-\overline{J_\zeta}
\big)\,\kappa
-
\big(
2\,R_\zeta-s
\big)\,\zeta
-
\big(
2\,R_{\overline{\kappa}}+\isqrt\,\overline{g}-J_\zeta
\big)\,\overline{\kappa}
-
\big(
2\,R_{\overline{\zeta}}-\overline{s}
\big)\,\overline{\zeta}.
\endaligned
\]
An identification to zero of all the nine coefficients
of $\pi^1$, $\pi^2$, $\overline{\pi}^1$, $\overline{\pi}^2$,
$\rho$, $\kappa$, $\zeta$, $\overline{\kappa}$, $\overline{\zeta}$
gives:
\[
\!\!\!\!\!\!\!\!\!\!\!\!
\small
\aligned
R_{\pi^1}
&
\,=\,
-\,R,
\ \ \ \ \ \ \ \ \ \ \ \ \ \ \ \ \ 
R_{\pi^2}
\,=\,
-\,{\textstyle{\frac{1}{2}}}\,
W,
\ \ \ \ \ \ \ \ \ \ \ \ \ \ \ \ \ 
R_{\overline{\pi}^1}
\,=\,
-\,R,
\ \ \ \ \ \ \ \ \ \ 
R_{\overline{\pi}^2}
\,=\,
-\,{\textstyle{\frac{1}{2}}}\,
\overline{W},
\\
R_\rho
&
\,=\,
{\textstyle{\frac{1}{2}}}\,
\big(r+\overline{r}\big),
\ \ \ \ \ \ \ \ \ \ 
R_\kappa
\,=\,
{\textstyle{\frac{1}{2}}}\,
\big(\isqrt\,g+\overline{J_\zeta}\big),
\ \ \ \ \ \ \ \ \ \ 
R_\zeta
\,=\,
{\textstyle{\frac{1}{2}}}\,
s,
\ \ \ \ \ \ \ \ \ \ \ \ \ 
R_{\overline{\kappa}}
\,=\,
{\textstyle{\frac{1}{2}}}\,
\big(
-\,\isqrt\,\overline{g}+J_\zeta
\big),
\ \ \ \ \ \ \ \ \ \ 
R_{\overline{\zeta}}
\,=\,
{\textstyle{\frac{1}{2}}}\,
\overline{s},
\endaligned
\]
and so:
\[
dR
\,=\,
-\,R\,\pi^1
-
{\textstyle{\frac{1}{2}}}\,
W\,\pi^2
-
R\,\overline{\pi}^1
-
{\textstyle{\frac{1}{2}}}\,
\overline{W}\,\overline{\pi}^2
+
R_\rho\,\rho
+
R_\kappa\,\kappa
+
R_\zeta\,\zeta
+
R_{\overline{\kappa}}\,\overline{\kappa}
+
R_{\overline{\zeta}}\,\overline{\zeta}.
\]
Inserting this back into what precedes, we can therefore obtain
both:
\[
\Pi-\overline{\Pi}
= 
-{\textstyle{\frac{1}{2}}}W\pi^{2}
+{\textstyle{\frac{1}{2}}}\overline{W}\overline{\pi}^{2}
+R_{\zeta}\zeta
-R_{\overline{\zeta}}\overline{\zeta}
+(R_{\kappa}-\overline{J_{\zeta}})\kappa
-(R_{\overline{\kappa}}-J_{\zeta})\overline{\kappa}
+{\textstyle{\frac{1}{2}}}(g_{\rho}-\overline{g}_{\rho})\rho,
\]
and replacing $g = -\, 2\isqrt\, R_\kappa + \isqrt\,
\overline{J_\zeta}$ 
in~({\ref{expression-Psi}}):
\[
\Psi=
-R\pi^{2}+e\rho+f\kappa
+i(W_{\rho}-2R_{\kappa}+\overline{J_{\zeta}})\zeta
+i(WJ-J_{\kappa})\overline{\kappa}
+iJ\overline{\zeta}.
\]
Thus:
\begin{equation*}
\begin{aligned}
\Omega_{1}
&=p\kappa\wedge\overline{\kappa}
+\Pi\wedge\rho
+\overline{W}\overline{\zeta}\wedge\kappa
-W\zeta\wedge\overline{\kappa}\\
&=p\kappa\wedge\overline{\kappa}
+{\textstyle{\frac{1}{2}}}(\Pi-\overline{\Pi})\wedge\rho
+\overline{W}\overline{\zeta}\wedge\kappa
-W\zeta\wedge\overline{\kappa}
+
{\textstyle{\frac{1}{2}}}(\Pi+\overline{\Pi})\wedge\rho\\
&=
-{\textstyle{\frac{1}{4}}}W\pi^{2}\wedge\rho
+{\textstyle{\frac{1}{4}}}\overline{W}\overline{\pi}^{2}\wedge\rho
-{\textstyle{\frac{1}{2}}}(R_{\kappa}-\overline{J_{\zeta}})\rho\wedge\kappa
-{\textstyle{\frac{1}{2}}}R_{\zeta}\rho\wedge\zeta
+{\textstyle{\frac{1}{2}}}(R_{\overline{\kappa}}-J_{\zeta})\rho\wedge\overline{\kappa}\\
&\hspace{0.5cm}
+{\textstyle{\frac{1}{2}}}R_{\overline{\zeta}}\rho\wedge\overline{\zeta}
+({\textstyle{\frac{1}{2}}}W_{\overline{\kappa}}-iR)\kappa\wedge\overline{\kappa}
-\overline{W}\kappa\wedge\overline{\zeta}
-W\zeta\wedge\overline{\kappa}
+{\textstyle{\frac{1}{2}}}(\Pi+\overline{\Pi})\wedge\rho,
\end{aligned}
\end{equation*}
and:
\begin{equation*}
\begin{aligned}
\Omega_{2}
&=
-\,R\pi^{2}\wedge\rho
-{\textstyle{\frac{1}{4}}}W\pi^{2}\wedge\kappa
+{\textstyle{\frac{1}{4}}}\overline{W}\overline{\pi}^{2}\wedge\kappa
-i(W_{\rho}-2R_{\kappa}+\overline{J_{\zeta}})\rho\wedge\zeta\\
&\hspace{0.5cm}
-i(WJ-J_{\kappa})\rho\wedge\overline{\kappa}
-iJ\rho\wedge\overline{\zeta}
-{\textstyle{\frac{1}{2}}}R_{\zeta}\kappa\wedge\zeta
+{\textstyle{\frac{1}{2}}}(R_{\overline{\kappa}}-J_{\zeta})\kappa\wedge\overline{\kappa}
+{\textstyle{\frac{1}{2}}}R_{\overline{\zeta}}\kappa\wedge\overline{\zeta}\\
&\hspace{0.5cm}
-R\zeta\wedge\overline{\kappa}
+{\textstyle{\frac{1}{2}}}(\Pi+\overline{\Pi})\wedge\kappa
+({\textstyle{\frac{1}{2}}}(r-\overline{r})-f)\rho\wedge\kappa.
\end{aligned}
\end{equation*}
If we define:
\[
\Lambda:=
{\textstyle{\frac{1}{2}}}(\Pi+\overline{\Pi})+
\text{real part of }
\big({\textstyle{\frac{1}{2}}}(g_{\rho}-\overline{g_{\rho}})-d_{\kappa}
\big)\rho,
\]
and:
\[
h:= 
\text{imaginary part of }
\big({\textstyle{\frac{1}{2}}}(g_{\rho}-\overline{g_{\rho}})-d_{\kappa}
\big),
\]
we conclude that:
\begin{equation*}
\begin{aligned}
\Omega_{1}&=
-{\textstyle{\frac{1}{4}}}W\pi^{2}\wedge\rho
+{\textstyle{\frac{1}{4}}}\overline{W}\overline{\pi}^{2}\wedge\rho
-{\textstyle{\frac{1}{2}}}(R_{\kappa}-\overline{J_{\zeta}})\rho\wedge\kappa
-{\textstyle{\frac{1}{2}}}R_{\zeta}\rho\wedge\zeta
+{\textstyle{\frac{1}{2}}}(R_{\overline{\kappa}}-J_{\zeta})\rho\wedge\overline{\kappa}\\
&\hspace{0.5cm}
+{\textstyle{\frac{1}{2}}}R_{\overline{\zeta}}\rho\wedge\overline{\zeta}
+\big({\textstyle{\frac{1}{2}}}W_{\overline{\kappa}}-iR\big)\kappa\wedge\overline{\kappa}
-\overline{W}\kappa\wedge\overline{\zeta}
-W\zeta\wedge\overline{\kappa}
+\Lambda\wedge\rho,\\
\Omega_{2}
&=
-\,R\pi^{2}\wedge\rho
-{\textstyle{\frac{1}{4}}}W\pi^{2}\wedge\kappa
+{\textstyle{\frac{1}{4}}}\overline{W}\overline{\pi}^{2}\wedge\kappa
-i(W_{\rho}-2R_{\kappa}+\overline{J_{\zeta}})\rho\wedge\zeta\\
&\hspace{0.5cm}
-i(WJ-J_{\kappa})\rho\wedge\overline{\kappa}
-iJ\rho\wedge\overline{\zeta}
-{\textstyle{\frac{1}{2}}}R_{\zeta}\kappa\wedge\zeta
+{\textstyle{\frac{1}{2}}}(R_{\overline{\kappa}}-J_{\zeta})\kappa\wedge\overline{\kappa}
+{\textstyle{\frac{1}{2}}}R_{\overline{\zeta}}\kappa\wedge\overline{\zeta}\\
&\hspace{0.5cm}
-R\zeta\wedge\overline{\kappa}
+\Lambda\wedge\kappa+h\rho\wedge\kappa.
\end{aligned}
\end{equation*}
Notice that all coefficients of $2$-forms\,\,---\,\,except
only $h$\,\,---\,\,depend on $R$, $J$, $W$ and
their coframe derivatives.

We are now close to the termination towards an $\{e\}$-structure.
In summary, we have obtained the following structure equations:
\begin{equation*}
\begin{aligned}
d\rho &= \pi^{1}\wedge\rho+\overline{\pi}^{1}\wedge\rho+i\kappa\wedge\overline{\kappa},\\
d\kappa &= \pi^{1}\wedge\kappa+\pi^{2}\wedge\rho+\zeta\wedge\overline{\kappa},\\
d\zeta &= i\pi^{2}\wedge\kappa+\pi^{1}\wedge\zeta-\overline{\pi}^{1}\wedge\zeta
+W\kappa\wedge\zeta+R\rho\wedge\zeta+J\rho\wedge\overline{\kappa},\\
d\pi^{1} &=
\Lambda\wedge\rho
-{\textstyle{\frac{1}{4}}}W\pi^{2}\wedge\rho
+{\textstyle{\frac{1}{4}}}\overline{W}\overline{\pi}^{2}\wedge\rho
-i\overline{\pi}^{2}\wedge\kappa\\
&\hspace{0.5cm}
-{\textstyle{\frac{1}{2}}}(R_{\kappa}-\overline{J_{\zeta}})\rho\wedge\kappa
-{\textstyle{\frac{1}{2}}}R_{\zeta}\rho\wedge\zeta
+{\textstyle{\frac{1}{2}}}(R_{\overline{\kappa}}-J_{\zeta})\rho\wedge\overline{\kappa}
+{\textstyle{\frac{1}{2}}}R_{\overline{\zeta}}\rho\wedge\overline{\zeta}\\
&\hspace{0.5cm}
+\big({\textstyle{\frac{1}{2}}}W_{\overline{\kappa}}-iR\big)\kappa\wedge\overline{\kappa}
-\overline{W}\kappa\wedge\overline{\zeta}
-W\zeta\wedge\overline{\kappa}
+\zeta\wedge\overline{\zeta},\\
d\pi^{2} &=
\Lambda\wedge\kappa
+
\pi^{2}\wedge\overline{\pi}^{1}
-\overline{\pi}^{2}\wedge\zeta
-R\pi^{2}\wedge\rho
-{\textstyle{\frac{1}{4}}}W\pi^{2}\wedge\kappa
+{\textstyle{\frac{1}{4}}}\overline{W}\overline{\pi}^{2}\wedge\kappa\\
&\hspace{0.5cm}
+h\rho\wedge\kappa
-i(W_{\rho}-2R_{\kappa}+\overline{J_{\zeta}})\rho\wedge\zeta
-i(WJ-J_{\kappa})\rho\wedge\overline{\kappa}
-iJ\rho\wedge\overline{\zeta}\\
&\hspace{0.5cm}
-{\textstyle{\frac{1}{2}}}R_{\zeta}\kappa\wedge\zeta
+{\textstyle{\frac{1}{2}}}(R_{\overline{\kappa}}-J_{\zeta})\kappa\wedge\overline{\kappa}
+{\textstyle{\frac{1}{2}}}R_{\overline{\zeta}}\kappa\wedge\overline{\zeta}
-R\zeta\wedge\overline{\kappa}.
\end{aligned}
\end{equation*}

But at this stage, we cannot directly deduce from these 
equations an appropriate expression for $h$. 
For example, any attempt to isolate $h$ by wedging the equation
$d\pi^{2} = \cdots$ 
with any appropriate differential form will include a
component of Maurer-Cartan type.  This is to be expected, because
$h$ will soon be shown below
to depend on higher order jets of $R$, $J$, $W$,
while the torsions above only depend up to the
2\textsuperscript{nd}-order jets of these 
invariants. Therefore, an application of the exterior
differentiation on both sides of the equation $d\pi^{2} = \cdots$
appears necessary to 
reach an expression for $h$ from the Poincar\'{e} relation $d\circ
d=0$.

To facilitate the discussion, we set:
\begin{equation*}
\begin{aligned}
\widehat{\Omega}_{1} &=
-{\textstyle{\frac{1}{4}}}W\pi^{2}\wedge\rho
+
{\textstyle{\frac{1}{4}}}\overline{W}\overline{\pi}^{2}\wedge\rho
-
{\textstyle{\frac{1}{2}}}(R_{\kappa}-\overline{J_{\zeta}})\rho\wedge\kappa
-
{\textstyle{\frac{1}{2}}}R_{\zeta}\rho\wedge\zeta\\
&\hspace{0.5cm}
+
{\textstyle{\frac{1}{2}}}(R_{\overline{\kappa}}-J_{\zeta})\rho\wedge\overline{\kappa}
+
{\textstyle{\frac{1}{2}}}R_{\overline{\zeta}}\rho\wedge\overline{\zeta}
+
\bigg({\textstyle{\frac{1}{2}}}W_{\overline{\kappa}}-iR\bigg)\kappa\wedge\overline{\kappa}
-\overline{W}\kappa\wedge\overline{\zeta}
-W\zeta\wedge\overline{\kappa},\\
\widehat{\Omega}_{2}
&=
-R\pi^{2}\wedge\rho
-{\textstyle{\frac{1}{4}}}W\pi^{2}\wedge\kappa
+{\textstyle{\frac{1}{4}}}\overline{W}\overline{\pi}^{2}\wedge\kappa
-i(W_{\rho}-2R_{\kappa}+\overline{J_{\zeta}})\rho\wedge\zeta\\
&\hspace{0.5cm}
-i(WJ-J_{\kappa})\rho\wedge\overline{\kappa}
-iJ\rho\wedge\overline{\zeta}
-{\textstyle{\frac{1}{2}}}R_{\zeta}\kappa\wedge\zeta
+{\textstyle{\frac{1}{2}}}(R_{\overline{\kappa}}-J_{\zeta})\kappa\wedge\overline{\kappa}
+{\textstyle{\frac{1}{2}}}R_{\overline{\zeta}}\kappa\wedge\overline{\zeta}\\
&\hspace{0.5cm}
-R\zeta\wedge\overline{\kappa},
\end{aligned}
\end{equation*} 
so that: 
\begin{equation*}
\begin{aligned}
d\pi^{1} &= \Lambda\wedge\rho-i\overline{\pi}^{2}\wedge\kappa+\zeta\wedge\overline{\zeta}+\widehat{\Omega}_{1},\\
d\pi^{2} &= \Lambda\wedge\kappa +\pi^{2}\wedge\overline{\pi}^{1}-\overline{\pi}^{2}\wedge\zeta+\widehat{\Omega}_{2}+h\rho\wedge\kappa.
\end{aligned}
\end{equation*}

\begin{Proposition}
The function $h$ is a function of the
3\textsuperscript{rd}-order jets of $W$ and $J$.
\end{Proposition}

\begin{proof}
By applying exterior differentiation $d$ to the equation of
$d\pi^{2}$, while wedging on both sides with
$\kappa\wedge\pi^{1}\wedge\overline{\pi}^{1}\wedge
\pi^{2}\wedge\overline{\pi}^{2}$,
we obtain:
\begin{equation*}
\begin{aligned}
2h\ \rho\wedge\kappa\wedge\overline{\kappa}\wedge\zeta\wedge\pi^{1}\wedge\overline{\pi}^{1}\wedge\pi^{2}\wedge\overline{\pi}^{2}
&=
-\overline{\widehat{\Omega}}_{2}\wedge\kappa\wedge\zeta\wedge\pi^{1}\wedge\overline{\pi}^{1}\wedge\pi^{2}\wedge\overline{\pi}^{2}\\
&\hspace{0.5cm}-d\widehat{\Omega}_{2}\wedge\kappa\wedge\pi^{1}\wedge\overline{\pi}^{1}\wedge\pi^{2}\wedge\overline{\pi}^{2}.\qedhere
\end{aligned}
\end{equation*}
\end{proof}

At this point, let $\Phi$ be the auxiliary real $2$-form:
\[
\Phi:=d\Lambda-\Lambda\wedge\pi^{1}-\Lambda\wedge\overline{\pi}^{1}-i\pi^{2}\wedge\overline{\pi}^{2}.
\]
Again this comes from the consideration of the model case. 
The structure equations therefore become:
\begin{equation*}
\begin{aligned}
d\rho &= 
\pi^{1}\wedge\rho+\overline{\pi}^{1}\wedge\rho+i\kappa\wedge\overline{\kappa},\\
d\kappa &= \pi^{1}\wedge\kappa+\pi^{2}\wedge\rho+\zeta\wedge\overline{\kappa},\\
d\zeta &= i\pi^{2}\wedge\kappa+\pi^{1}\wedge\zeta-\overline{\pi}^{1}\wedge\zeta
+W\kappa\wedge\zeta+R\rho\wedge\zeta+J\rho\wedge\overline{\kappa},\\
d\pi^{1} &= \Lambda\wedge\rho-i\overline{\pi}^{2}\wedge\kappa+\zeta\wedge\overline{\zeta}+\widehat{\Omega}_{1},\\
d\pi^{2} &= \Lambda\wedge\kappa+\pi^{2}\wedge\overline{\pi}^{1}-\overline{\pi}^{2}\wedge\zeta+\widehat{\Omega}_{2}+h\rho\wedge\kappa,\\
d\Lambda &= \Lambda\wedge\pi^{1}+\Lambda\wedge\overline{\pi}^{1}+i\pi^{2}\wedge\overline{\pi}^{2}+\Phi.
\end{aligned}
\end{equation*}

\begin{Proposition}
The real $2$-form $\Phi$ is a function of
the 
4\textsuperscript{th}-order jets of $W$ and $J$.
\end{Proposition}

\begin{proof}
By taking exterior derivative of $d\pi^{1}$ and $d\pi^{2}$ again, this time using the expression of $d\Lambda$, we have:
\begin{equation*}
\begin{aligned}
\Phi\wedge\rho &=
i\overline{\widehat{\Omega}}_{2}\wedge\kappa
+ih\rho\wedge\kappa\wedge\overline{\kappa}
-W\kappa\wedge\zeta\wedge\overline{\zeta}
+\overline{W}\zeta\wedge\overline{\kappa}\wedge\overline{\zeta}
-2R\rho\wedge\zeta\wedge\overline{\zeta}\\
&\hspace{0.5cm}
-J\rho\wedge\overline{\kappa}\wedge\overline{\zeta}
+\overline{J}\rho\wedge\kappa\wedge\zeta
-d\widehat{\Omega}_{1},\\
\Phi\wedge\kappa &=
-\widehat{\Omega}_{2}\wedge\overline{\pi}^{1}
-h\rho\wedge\kappa\wedge\overline{\pi}^{1}
+\pi^{2}\wedge\overline{\widehat{\Omega}}_{1}
+\overline{\widehat{\Omega}}_{2}\wedge\zeta
-h\rho\wedge\overline{\kappa}\wedge\zeta
-W\overline{\pi}^{2}\wedge\kappa\wedge\zeta\\
&\hspace{0.5cm}
-R\overline{\pi}^{2}\wedge\rho\wedge\zeta
-J\overline{\pi}^{2}\wedge\rho\wedge\overline{\kappa}
-d\widehat{\Omega}_{2}
-d(h\rho\wedge\kappa).
\end{aligned}
\end{equation*}
Writing $\Phi$ as:
\[
\Phi=\widehat{\Omega}_{3}
+
u\rho\wedge\kappa,
\]
where $\widehat{\Omega}_{3}$ is the $2$-form not containing
$\rho\wedge\kappa$, then each of the coefficients in
$\widehat{\Omega}_{3}$ is a function of the
4\textsuperscript{th}-order jet of $W$ and $J$. Since $\Phi$ is real,
taking conjugate on both sides, we must have:
\[
\widehat{\Omega}_{3}
+
u\rho\wedge\kappa
=
\overline{\widehat{\Omega}}_{3}
+
\overline{u}\rho\wedge\overline{\kappa}.
\]
Therefore by inspection, $\overline{u}$ is also a function of the 
4\textsuperscript{th}-order
jets of $W$ and $J$, and therefore so is $u$. This finishes the proof.
\end{proof}

With this proposition, we have therefore fully constructed
an $\{e\}$-structure.


\newpage

\setcounter{section}{0}

$\:$

\bigskip\bigskip\bigskip\bigskip\bigskip

\begin{center}

{\large\bf Addendum to:}
\label{addendum-Merker-Pocchiola-JGA}

\medskip

{\large\bf Explicit absolute parallelism}

\smallskip

{\large\bf for $2$-nondegenerate real hypersurfaces}

\smallskip

{\large\bf $M^5 \subset \C^3$ of constant Levi rank $1$}

\medskip

Journal of Geometric Analysis

\smallskip

DOI 10.1007/s12220-018-9988-3

\bigskip\bigskip

Joël~{\sc Merker}\footnote{Laboratoire de Mathématiques d'Orsay,
Université Paris-Sud, CNRS, 
Université Paris-Saclay, 91405 Orsay Cedex, France.}
and Samuel~{\sc Pocchiola}

\end{center}\bigskip

\Section{\bf Introduction}
\label{introduction-addendum}
\HEAD{{\ref{introduction-addendum}}.~{\sf Introduction}
}{
and Joël {\sc Merker} (Orsay) and Samuel {\sc Pocchiola}}

The article~{\cite{Merker-Pocchiola-2018-addendum}}
and its preprint version~{\cite{Pocchiola-2013-addendum}}
(incorportating more calculations)
study $\mathcal{C}^\omega$ or $\mathcal{C}^\infty$
$2$-nondegenerate graphed hypersurfaces $M^5 \subset \C^3$:
\[
u
\,=\,
F\big(z_1,z_2,\overline{z}_1,\overline{z}_2,v\big)
\eqno
{\scriptstyle{(w\,=\,u\,+\,i\,v)}},
\]
whose Levi form has rank $1$ at every point.  With no details of
proof, the final Section~9
of~{\cite[pp.~42--43]{Merker-Pocchiola-2018-addendum}} states that when
Pocchiola's two fundamental differential invariants $J \equiv 0 \equiv
W$ of these CR structures vanish identically, a certain collection of
ten $1$-forms:
\[
\big\{
\rho,\,
\kappa,\,
\zeta,\,
\overline{\kappa},\,
\overline{\zeta},\,\,
\pi^1,\,
\pi^2,\,
\overline{\pi}^1,\,
\overline{\pi}^2,\,
\Lambda
\big\}
\eqno
{\scriptstyle{(\overline{\rho}\,=\,\rho,\,\,
\overline{\Lambda}\,=\,\Lambda)}},
\]
enjoy a Lie-Cartan structure
having {\em constant} coefficients 
(conjugate equations are unwritten):
\[
\aligned
d\rho
&
\,=\,
\pi^1\wedge\rho
+
\overline{\pi}^1\wedge\rho
+
i\,\kappa\wedge\overline{\kappa},
\\
d\kappa
&
\,=\,
\pi^1\wedge\kappa
+
\pi^2\wedge\rho
+
\zeta\wedge\overline{\kappa},
\\
d\zeta
&
\,=\,
i\,\pi^2\wedge\kappa
+
\pi^1\wedge\zeta
-
\overline{\pi}^1\wedge\zeta,
\\
d\pi^1
&
\,=\,
i\,\kappa\wedge\overline{\pi}^2
+
\zeta\wedge\overline{\zeta}
+
\Lambda\wedge\rho,
\\
d\pi^2
&
\,=\,
\pi^2\wedge\overline{\pi}^1
+
\zeta\wedge\overline{\pi}^2
+
\Lambda\wedge\kappa,
\\
d\Lambda
&
\,=\,
i\,\pi^2\wedge\overline{\pi}^2
+
\Lambda\wedge\pi^1
+
\Lambda\wedge\overline{\pi}^1.
\endaligned
\]
Furthermore, these equations are exactly equal to the Lie structure
(shown in~{\cite{Pocchiola-2014}}) of the model
light cone $M_{\sf LC}$:
\[
u
\,=\,
\frac{z_1\overline{z}_1
+
\frac{1}{2}\,z_1^2\overline{z}_2
+
\frac{1}{2}\,\overline{z}_1^2z_2}{
1-z_2\overline{z}_2},
\]
so that by general Cartan theory, 
hypersurfaces having $J \equiv 0 \equiv W$ are {\em all}
(locally) biholomorphic to $M_{\sf LC}$.

Importantly, Pocchiola discovered that {\em before} 
prolongation of the equivalence problem,
the structure equations are of the form:
\[
\aligned
d\rho
&
\,=\,
\big(
\pi^1
+
\overline{\pi}^1
\big)
\wedge\rho
+
\isqrt\,\kappa\wedge\overline{\kappa},
\notag
\\
d\kappa
&
\,=\,
\pi^2\wedge\rho
+
\pi^1\wedge\kappa
+
\zeta\wedge\overline{\kappa},
\notag
\\
d\zeta
&
\,=\,
\big(\pi^1-\overline{\pi}^1\big)
\wedge\zeta
+
\isqrt\,\pi^2\wedge\kappa
\,+
\notag
\\
&
\ \ \ \ \
+
R\,
\rho\wedge\zeta
+
{\textstyle{\frac{i}{\overline{\sf c}^3}}}\,
\overline{J}\,\rho\wedge\overline{\kappa}
+
{\textstyle{\frac{1}{{\sf c}}}}\,
W\,
\kappa\wedge\zeta.
\endaligned
\]
Furthermore, he showed that the real function
$R \equiv 0$ vanishes identically when $W \equiv 0$,
so that when $J \equiv 0 \equiv W$, these
equations reduce to constant coefficients:
\begin{align}
\label{structure-constant-coefficients}
d\rho
&
\,=\,
\big(\pi^1+\overline{\pi}^1\big)
\wedge\rho
+
i\,\kappa\wedge\overline{\kappa},
\notag
\\
d\kappa
&
\,=\,
\pi^2
\wedge\rho
+
\pi^1\wedge\kappa
+
\zeta\wedge\overline{\kappa},
\\
d\zeta
&
\,=\,
\big(
\pi^1
-
\overline{\pi}^1
\big)
\wedge\zeta
+
\isqrt\,
\pi^2\wedge\kappa.
\notag
\end{align}
At the end of~{\cite{Pocchiola-2013-addendum}},
he showed with elementary
reasonings that after prolongation, one obtains the last
$3$ structure equations above for $d\pi^1$, $d\pi^2$,
$d\Lambda$, 
in which {\em no} nonconstant structure
functions appear.

This phenomenon is in some sense
`counter-intuitive' to CR geometers, 
since for Levi nondegenerate CR structures,
and for the corresponding second order {\sc pde} systems,
{\em no}
curvatures appear after absorption before prolongation:
\[
\aligned
d\omega
&
\,=\,
\smallsum{\alpha}\,
\omega^\alpha
\wedge
\omega_\alpha
+
\omega\wedge\varphi,
\\
d\omega^\alpha
&
\,=\,
\smallsum{\beta}\,
\omega^\beta\wedge\varphi_\beta^\alpha
+
\omega\wedge\varphi^\alpha,
\\
d\omega_\alpha
&
\,=\,
\smallsum{\beta}\,
\varphi_\alpha^\beta
\wedge
\omega_\beta
+
\omega_\alpha\wedge\varphi
+
\omega\wedge\varphi_\alpha,
\endaligned
\]
while primary and secondary
invariants appear afterwards, {\em e.g.} 
like
$S_{\beta\rho}^{\alpha\sigma}$
and $R_{\beta\gamma}^\alpha$,
$T_\beta^{\alpha\gamma}$
in:
\[
\aligned
d\varphi_\beta^\alpha
&
\,=\,
{\textstyle{\frac{1}{2}}}\,
\delta_\beta^\alpha\,
\psi\wedge\omega
-
\varphi_\beta^\gamma
\wedge
\varphi_\gamma^\alpha
-
\varphi_\beta\wedge\omega^\alpha
-
\varphi^\alpha\wedge\omega_\beta
+
\delta_\beta^\alpha\,
\omega^\gamma\wedge\varphi_\gamma
\,+
\\
&
\ \ \ \ \
+
S_{\beta\rho}^{\alpha\sigma}\,
\omega^\rho
\wedge\omega_\sigma
+
R_{\beta\gamma}^\alpha\,
\omega^\gamma\wedge\omega
+
T_\beta^{\alpha\gamma}\,
\omega_\gamma\wedge\omega.
\endaligned
\]

This addendum is devoted to expose
the structure equations for $d\pi^1$, $d\pi^2$,
$d\Lambda$, the details of which were completely
skipped in Section~9 of~{\cite{Merker-Pocchiola-2018-addendum}},
hence asked to appear in print by some experts.
At the end of Section~{\ref{section-reduction-e-structure-constant}}, 
Theorem~{\ref{Theorem-reduction-e-structure-J-W-0}}
summarizes all constructions.

\smallskip\noindent{\bf Acknowledgments.} The authors thank
Wei Guo Foo (Beijing) who brought to their attention the presence of
a purely imaginary function $h = i\, H$
with $\overline{H} = H$ in Lemma~{\ref{Lemma-M-5-G-4-R}},
missed in Pocchiola's prepublication~{\cite{Pocchiola-2013-addendum}},
but anyway shown later to vanish in Lemma~{\ref{Lemma-H-vanishes}}.

\Section{\bf Summary of Pocchiola's constructions up to the assumption
$J \equiv 0 \equiv W$}
\label{summary-addendum}
\HEAD{{\ref{summary-addendum}}.~{\sf Summary of Pocchiola's 
constructions up to the assumption $J \equiv 0 \equiv W$}
}{
and Joël {\sc Merker} (Orsay) and Samuel {\sc Pocchiola}}

This brief technical section is independent of the rest
of the present addendum, and serves only as a (non-self-contained)
reminder. Once the mentioned
equations~({\ref{structure-constant-coefficients}})
will be reached at the end
of this section, all considerations will be self-contained.

In addition to $J$ and $W$, a certain derived invariant
also occurs in~{\cite{Pocchiola-2013-addendum,
Merker-Pocchiola-2018-addendum}}, 
which, after absorption of its real part in some
modified Maurer-Cartan form, reduces to:
\[
W_{\rho\zeta}^\zeta
\,=\,
\frac{1}{{\sf c}\overline{\sf c}}
\left(
-\frac{i}{6}\,
\bigg(
\frac{\overline{\mathcal{L}}_1\big(
\overline{\mathcal{L}}_1(k)\big)}{
\overline{\mathcal{L}}_1(k)}
-
\overline{P}
\bigg)\,
W
-
\frac{i}{2}\,
\overline{\mathcal{L}}_1(W)
\right)
+
i\,
\frac{{\sf e}}{{\sf c}{\sf c}}\,
W,
\]
whence $W_{\rho\zeta}^\zeta \equiv 0$ as soon as $W \equiv 0$.

In the product manifold $M^5 \times G^4$ corresponding to the
$4$-dimensional $G$-structure obtained in Section~6
of~{\cite{Merker-Pocchiola-2018-addendum}} after
several reductions and computations, 
which is equipped with coordinates:
\[
\big(
z_1,z_2,\overline{z}_1,\overline{z}_2,v
\big)
\times
\big(
{\sf c},\overline{\sf c},
{\sf e},\overline{\sf e}
\big),
\]
the Maurer-Cartan forms are:
\[
\delta^1
\,:=\,
\frac{d{\sf c}}{{\sf c}},
\ \ \ \ \ \ \ \ \ \ \ \ \ \ \ \ \ \ \ \ \ \ \ \ \ \
\delta^2
\,:=\,
i\,
\frac{{\sf e}\,d{\sf c}}{{\sf c}{\sf c}}
-
i\,
\frac{{\sf e}\,d\overline{\sf c}}{{\sf c}\overline{\sf c}}
-
i\,\frac{d{\sf e}}{{\sf c}},
\]
together with their conjugates $\overline{\delta}^1$, 
$\overline{\delta}^2$.
Furthermore, the modified Maurer-Cartan forms
introduced by Pocchiola are:
\[
\aligned
\widehat{\delta}^1
&
\,:=\,
\delta^1
+
{\textstyle{\frac{1}{2}}}\,
V_{\rho\zeta}^\zeta\,
\rho
-
\overline{V}_{\rho\kappa}^\rho\,
\kappa
-
V_{\rho\zeta}^\rho\,
\zeta
-
V_{\kappa\overline{\kappa}}^\kappa\,
\overline{\kappa},
\\
\widehat{\delta}^2
&
\,:=\,
\delta^2
-
V_{\rho\kappa}^\zeta\,
\rho
-
\big(
V_{\rho\kappa}^\kappa
-
{\textstyle{\frac{1}{2}}}\,
V_{\rho\zeta}^\zeta
\big)\,
\kappa
-
V_{\rho\zeta}^\kappa\,
\zeta
-
V_{\rho\overline{\kappa}}^{\kappa}\,
\overline{\kappa}
-
V_{\rho\overline{\zeta}}^\kappa\,
\overline{\zeta},
\endaligned
\]
and they satisfy the neat structure equations:
\[
\aligned
d\rho
&
\,=\,
\big(\widehat{\delta}^1+\overline{\widehat{\delta}}^1\big)
\wedge\rho
+
i\,\kappa\wedge\overline{\kappa},
\\
d\kappa
&
\,=\,
\widehat{\delta}^2
\wedge\rho
+
\widehat{\delta}^1\wedge\kappa
+
\zeta\wedge\overline{\kappa},
\\
d\zeta
&
\,=\,
\isqrt\,
\widehat{\delta}^2\wedge\kappa
+
\widehat{\delta}^1\wedge\zeta
-
\overline{\widehat{\delta}}^1\wedge\zeta.
\endaligned
\]
Lastly, with a single real free variable ${\sf t} = -\Re\,
w_\rho^1$ which parametrizes the $1$-dimensional prolongation
of this $G$-structure, introducing:
\[
\aligned
\pi^1
&
\,:=\,
\widehat{\delta}^1
+
{\sf t}\,\rho,
\\
\pi^2
&
\,:=\,
\widehat{\delta}^2
+
{\sf t}\,\kappa,
\endaligned
\]
the structure equations under the assumption $J \equiv 0 \equiv W$ 
become, as was already seen:
\begin{align}
\label{1-structure-d}
d\rho
&
\,=\,
\big(\pi^1+\overline{\pi}^1\big)
\wedge\rho
+
i\,\kappa\wedge\overline{\kappa},
\\
\label{2-structure-d}
d\kappa
&
\,=\,
\pi^2
\wedge\rho
+
\pi^1\wedge\kappa
+
\zeta\wedge\overline{\kappa},
\\
\label{3-structure-d}
d\zeta
&
\,=\,
\isqrt\,
\pi^2\wedge\kappa
+
\pi^1\wedge\zeta
-
\overline{\pi}^1\wedge\zeta,
\end{align}
together with the conjugate equations
$d\overline{\kappa} = \overline{\pi}^2 \wedge \rho + \cdots$
and $d\overline{\zeta} = -i\, \overline{\pi}^2 \wedge
\overline{\kappa} + \cdots$, while $d\overline\rho = d\rho$.

We now proceed to explain why the structure equations for 
$d\pi^1$ and $d\pi^2$ incorporate only {\em constant} 
coefficients, as well as the structure
equations for the exterior differential $d\Lambda$,
where $\Lambda = 
\overline{\Lambda}$ will be a final prolonged real $1$-form
completing an $\{e\}$-structure on $M^5 \times G^4 \times \R$.

\Section{\bf The $2$-forms $\Omega_1$ and $\Omega_2$}
\label{two-forms-Omega-1-2}
\HEAD{{\ref{two-forms-Omega-1-2}}.~{\sf The $2$-forms $\Omega_1$ 
and $\Omega_2$}
}{
and Joël {\sc Merker} (Orsay) and Samuel {\sc Pocchiola}}

On the manifold $M^5 \times G^4 \times \R$, we thus start 
from~({\ref{1-structure-d}}), ({\ref{2-structure-d}}),
({\ref{3-structure-d}}). All reasonings will be elementary,
and rely upon Poincaré's relation $d \circ d = 0$ and upon
several application of the famous
{\sl Cartan Lemma} which helps to bypass hard explicit computations
of Pocchiola's style.

\begin{Lemma}
The differential $2$-forms:
\[
\aligned
\Omega_1
&
\,:=\,
d\pi^1
-
\isqrt\,\kappa\wedge\overline{\pi}^2
-
\zeta\wedge\overline{\zeta},
\\
\Omega_2
&
\,:=\,
d\pi^2
-
\pi^2\wedge\overline{\pi}^1
-
\zeta\wedge\overline{\pi}^2,
\endaligned
\]
satisfy:
\begin{align}
0
&
\,=\,
\big(\Omega_1+\overline{\Omega}_1\big)
\wedge\rho,
\label{1-Omega}
\\
0
&
\,=\,
\Omega_2\wedge\rho
+
\Omega_1\wedge\kappa,
\label{2-Omega}
\\
0
&
\,=\,
\big(\Omega_1-\overline{\Omega}_1\big)
\wedge
\zeta
+
\isqrt\,
\Omega_2\wedge\kappa.
\label{3-Omega}
\end{align}
\end{Lemma}

\proof
Using the Poincaré relation $d \circ d = 0$, 
apply the exterior differentiation operator $d$ 
to~({\ref{1-structure-d}}), and replace $d\rho$,
$d\kappa$, $d\overline{\kappa}$ by means of~({\ref{1-structure-d}}),
({\ref{2-structure-d}}):
\[
\aligned
0
\,=\,
d\circ d\rho
&
\,=\,
\big(
d\pi^1+d\overline{\pi}^1
\big)
\wedge\rho
-
\big(\pi^1+\overline{\pi}^1\big)
\wedge
d\rho
+
i\,d\kappa\wedge\overline{\kappa}
-
i\,\kappa\wedge
d\overline{\kappa}
\\
&
\,=\,
\big(d\pi^1+d\overline{\pi}^1\big)
\wedge\rho
-
\big(\pi^1+\overline{\pi}^1\big)
\wedge
\Big(
\zero{
\big(\pi^1+\overline{\pi}^1\big)}
\wedge\rho
+
i\,\kappa\wedge\overline{\kappa}
\Big)
\,+
\\
&
\ \ \ \ \
+
i\,
\Big(
\pi^2\wedge\rho
+
\pi^1\wedge\kappa
+
\zero{
\zeta\wedge\overline{\kappa}}
\Big)
\wedge
\overline{\kappa}
-
i\,\kappa
\wedge
\Big(
\overline{\pi}^2\wedge\rho
+
\overline{\pi}^1\wedge\overline{\kappa}
+
\overline{\zeta}\wedge\kappa
\Big),
\endaligned
\]
which becomes after simplification:
\[
0
\,=\,
\Big(
d\pi^1
-
i\,\kappa\wedge\overline{\pi}^2
+
d\overline{\pi}^1
+
i\,\overline{\kappa}\wedge\pi^2
\Big)
\wedge
\rho,
\]
and proves~({\ref{1-Omega}}), since
$-\zeta \wedge \overline{\zeta}$ is purely imaginary.

For~({\ref{2-Omega}}), proceed similarly 
with~({\ref{2-structure-d}}):
\[
\aligned
0
\,=\,
d\circ d\kappa
&
\,=\,
d\pi^2\wedge\rho
-
\pi^2\wedge d\rho
+
d\pi^1\wedge\kappa
-
\pi^1\wedge d\kappa
+
d\zeta\wedge\overline{\kappa}
-
\zeta\wedge d\overline{\kappa}
\\
&
\,=\,
d\pi^2\wedge\rho
-
\pi^2\wedge
\Big(
\big(
\pi^1
+
\overline{\pi}^1
\big)
\wedge\rho
+
i\,\kappa\wedge\overline{\kappa}
\Big)
+
d\pi^1\wedge\kappa
-
\pi^1\wedge
\big(
\pi^2\wedge\rho
+
\zeta\wedge\overline{\kappa}
\big)
\,+
\\
&
\ \ \ \ \ 
+
\Big(
i\,\pi^2\wedge\kappa
+
\pi^1\wedge\zeta
-
\overline{\pi}^1\wedge\zeta
\Big)
\wedge\overline{\kappa}
-
\zeta\wedge
\Big(
\overline{\pi}^2\wedge\rho
+
\overline{\pi}^1
\wedge
\overline{\kappa}
+
\overline{\zeta}\wedge\kappa
\Big),
\endaligned
\]
and receive, after simplification and reorganization:
\[
0
\,=\,
\Big(
d\pi^2
-
\pi^2
\wedge\overline{\pi}^1
-
\zeta\wedge\overline{\pi}^2
\Big)
\wedge\rho
+
\big(
d\pi^1
-
\zeta\wedge\overline{\zeta}
\big)
\wedge\kappa,
\]
which is~({\ref{2-Omega}}), since $\big( - i \kappa \wedge
\overline{\pi}^2 \big) \wedge \kappa = 0$. 

Lastly, the Poincaré relation applied 
to~({\ref{3-structure-d}}):
\[
\aligned
0
\,=\,
d\circ d\zeta
&
\,=\,
i\,d\pi^2\wedge\kappa
-
i\,\pi^2\wedge d\kappa
+
d\pi^1\wedge\zeta
-
\pi^1\wedge d\zeta
-
d\overline{\pi}^1\wedge\zeta
+
\overline{\pi}^1\wedge d\zeta
\\
&
\,=\,
i\,d\pi^2\wedge\kappa
-
i\,\pi^2\wedge
\big(
\pi^1\wedge\kappa
+
\zeta\wedge\overline{\kappa}
\big)
+
d\pi^1\wedge\zeta
\,-
\\
&
\ \ \ \ \ 
-\,
\pi^1\wedge\big(
i\,\pi^2\wedge\kappa
-
\overline{\pi}^1\wedge\zeta
\big)
-
d\overline{\pi}^1\wedge\zeta
+
\overline{\pi}^1\wedge
\big(
i\,\pi^2\wedge\kappa
+
\pi^1\wedge\zeta
\big),
\endaligned
\]
contracts as and rewrites as~({\ref{3-Omega}}):
\begin{align}
0
&
\,=\,
\Big(
d\pi^1
-
d\overline{\pi}^1
-
i\,\overline{\kappa}
\wedge
\pi^2
\Big)
\wedge\zeta
+
i\,
\big(
d\pi^2
-
\pi^2\wedge\overline{\pi}^1
\big)
\wedge
\kappa
\notag
\\
&
\,=\,
\Big(
d\pi^1
-
i\,\underline{
\kappa\wedge\overline{\pi}^2}_A
-
d\overline{\pi}^1
-
i\,\overline{\kappa}
\wedge\pi^2
\Big)
\wedge\zeta
+
i\,
\Big(
d\pi^2
-
\pi^2\wedge\overline{\pi}^1
-
\underline{
\zeta\wedge\overline{\pi}^2}_A
\Big)
\wedge\kappa
\notag
\\
&
\,=\,
\big(
\Omega
-
\overline{\Omega}_1
\big)
\wedge\zeta
+
i\,\Omega_2\wedge\kappa.
\qedhere
\end{align}
\endproof

\begin{Lemma}
\label{Lemma-M-5-G-4-R}
On $M^5 \times G^4 \times \R$, the solutions in $2$-forms
$\Omega_1$, $\Omega_2$ of the three equations:
\begin{align}
0
&
\,=\,
\big(\Omega_1+\overline{\Omega}_1\big)
\wedge\rho,
\label{1-Omega-bis}
\\
0
&
\,=\,
\Omega_2\wedge\rho
+
\Omega_1\wedge\kappa,
\label{2-Omega-bis}
\\
0
&
\,=\,
\big(\Omega_1-\overline{\Omega}_1\big)
\wedge
\zeta
+
\isqrt\,
\Omega_2\wedge\kappa.
\label{3-Omega-bis}
\end{align}
are:
\[
\aligned
\Omega_1
&
\,=\,
\Lambda
\wedge
\rho,
\\
\Omega_2
&
\,=\,
\Lambda\wedge\kappa
+
h\,\rho\wedge\kappa,
\endaligned
\]
where $\overline{\Lambda} = \Lambda$ is a certain
{\em real} $1$-form, and where $h = - \overline{h}$
is a certain purely imaginary function.
\end{Lemma}

At this step, in a passage with no details of proof, 
Pocchiola missed $h$ in~{\cite{Pocchiola-2013-addendum}}, 
but a bit later we will show that $h = 0$ anyway
thanks to some other constraints.

\proof
Wedge~({\ref{2-Omega-bis}})$\wedge\kappa$:
\begin{align}
\label{Omega-2-rho-kappa}
0
\,=\,
\Omega_2\wedge\rho\wedge\kappa
+
0,
\end{align}
hence by the Cartan Lemma, there exist two $1$-forms
$\gamma$ and $\delta$ such that:
\begin{align}
\label{Omega-2-gamma-delta}
\Omega_2
\,=\,
\gamma\wedge\rho
+
\delta\wedge\kappa.
\end{align}

Next, wedge~({\ref{1-Omega-bis}})$\wedge\zeta$
and wedge~({\ref{3-Omega-bis}})$\wedge\rho$,
taking account of~({\ref{Omega-2-rho-kappa}}):
\[
\left.
\aligned
0
&
\,=\,
\big(
\Omega_1
+
\overline{\Omega}_1
\big)
\wedge\rho\wedge\zeta,
\\
0
&
\,=\,
-\,\big(\Omega_1-\overline{\Omega}_1\big)
\wedge\rho\wedge\zeta
+
0
\endaligned
\right\}
\ \ \ \ \ \ \ \ \ \ \ \ \ \ \ \ \ \
\Longrightarrow
\ \ \ \ \ \ \ \ \ \ \ \ \ \ \ \ \ \
0
\,=\,
\Omega_1
\wedge\rho\wedge\zeta,
\]
hence by the Cartan Lemma again, there exist two $1$-forms
$\alpha$ and $\beta$ such that:
\[
\Omega_1
\,=\,
\alpha\wedge\rho
+
\beta\wedge\zeta.
\]
Insert this $\Omega_1$ into~({\ref{2-Omega-bis}})$\wedge\rho$:
\[
0
\,=\,
0
+
\big(
\alpha\wedge\rho
+
\beta\wedge\zeta
\big)
\wedge\kappa\wedge\rho
\,\,=\,\,
-\,\beta
\wedge\rho\wedge\kappa\wedge\zeta,
\]
whence by the Cartan Lemma, there exist certain functions
$A$, $B$, $C$ such that:
\[
\beta
\,=\,
A\,\rho
+
B\,\kappa
+
C\,\zeta.
\]
We may assume $C = 0$, since this does not change
$\Omega_1$, and we assert that $B = 0$.
Indeed, replacing $\Omega_1$ in~({\ref{1-Omega-bis}})
shows $B = 0 = \overline{B}$:
\[
\aligned
0
&
\,=\,
\Big(
\zero{
\alpha\wedge\rho}
+
\beta\wedge\zeta
+
\zero{
\overline{\alpha}\wedge\rho}
+
\overline{\beta}\wedge\overline{\zeta}
\Big)
\wedge\rho
\\
&
\,=\,
\Big(
\zero{
A\,\rho\wedge\zeta}
+
B\,\kappa\wedge\zeta
+
0
+
\zero{
\overline{A}\,\rho\wedge\overline{\zeta}}
+
\overline{B}\,
\overline{\kappa}\wedge\overline{\zeta}
+
0
\Big)
\wedge\rho.
\endaligned
\]

Hence $\beta = A\, \rho$ and we get:
\[
\Omega_1
\,=\,
\big(\alpha-A\,\zeta\big)
\wedge
\rho
\,\,=:\,\,
\tau\wedge\rho,
\]
in terms of a new $1$-form $\tau$. But by inserting this
$\Omega_1$ into~({\ref{3-Omega-bis}})$\wedge\kappa$, we get:
\[
0
\,=\,
\big(
-\tau+\overline{\tau}
\big)
\wedge\rho\wedge\kappa\wedge\zeta
+
0,
\]
and the Cartan Lemma (again!) gives:
\[
\tau
-
\overline{\tau}
\,=\,
U\,\rho
+
V\,\kappa
+
W\,\zeta
\ \ \ \ \ \ \ \ \ \ \ \ \ \ \ \ \ \
\Longleftrightarrow
\ \ \ \ \ \ \ \ \ \ \ \ \ \ \ \ \ \
\overline{\tau}
-
\tau
\,=\,
\overline{U}\,\rho
+
\overline{V}\,\overline{\kappa}
+
\overline{W}\,\overline{\zeta},
\]
and since $\big\{ \rho, \kappa, \zeta, \overline{\kappa},
\overline{\zeta} \big\}$ are linearly independent
at every point, we see that $V = W = \overline{V} = \overline{W} = 0$
necessarily, while the function $U = - \overline{U}$ 
must be purely imaginary, and hence we can write:
\[
\tau
-
\overline{\tau}
\,=\,
U\,\rho
\,=\,
\big(
{\textstyle{\frac{U-\overline{U}}{2}}}
\big)\,
\rho.
\]
Consequently, if we introduce the $1$-form:
\[
\Lambda
\,:=\,
\tau
-
{\textstyle{\frac{1}{2}}}\,
U\,\rho
\,\,=\,\,
\overline{\Lambda},
\]
which is now {\em real}, we indeed obtain as announced:
\[
\Omega_1
\,=\,
\Lambda
\wedge
\rho.
\]

Next, inserting this into~({\ref{2-Omega-bis}}):
\[
0
\,=\,
\Omega_2
\wedge\rho
+
\Lambda\wedge\rho\wedge\kappa
\,\,=\,\,
\big(
\Omega_2
-
\Lambda\wedge\kappa
\big)
\wedge\rho,
\]
on application of the Cartan Lemma,
we receive a $1$-form $\gamma$ such that:
\[
\Omega_2
-
\Lambda\wedge\kappa
\,=\,
\gamma\wedge\rho,
\]
so that above in~({\ref{Omega-2-gamma-delta}}),
we had in fact $\delta = \Lambda$!
Next, insert all this into~({\ref{3-Omega-bis}}):
\[
0
\,=\,
\big(
\zero{
\Lambda\wedge\rho
-
\Lambda\wedge\rho}
\big)
\wedge\zeta
+
i\,\big(
\zero{
\Lambda\wedge\kappa}
+
\gamma\wedge\rho
\big)
\wedge\kappa
\,\,=\,\,
i\,
\gamma\wedge\rho\wedge\kappa.
\]
The Cartan Lemma gives two functions $G$ and $K$ with:
\[
\gamma
\,=\,
G\,\rho
+
K\,\kappa.
\]
But then:
\[
\aligned
\Omega_2
&
\,=\,
\Lambda\wedge\kappa
+
G\,
\zero{
\rho\wedge\rho}
-
\big(
\Re\,K
+
i\,\Im\,K
\big)\,
\rho\wedge\kappa
\\
&
\,=\,
\big(
\Lambda
-
\Re\,K
\cdot
\rho
\big)
\wedge\kappa
-
i\,\Im\,K\,
\rho\wedge\kappa,
\endaligned
\]
hence by redefining $\Lambda := \Lambda - \Re\, K \cdot \rho$
which leaves $\Omega_1$ untouched, and by setting
$h := - i\, \Im\, K$, we conclude this detailed proof.
\endproof

\Section{\bf Reduction to an $\{e\}$-structure with constant
coefficients}
\label{section-reduction-e-structure-constant}
\HEAD{{\ref{section-reduction-e-structure-constant}}.~{\sf Reduction 
to an $\{e\}$-structure with constant
coefficients}
}{
and Joël {\sc Merker} (Orsay) and Samuel {\sc Pocchiola}}

Changing the notation by setting $h =: i\, H$ with $H = \overline{H}$
being a real-valued function, we have therefore obtained:
\[
\aligned
d\pi^1
&
\,=\,
\Lambda\wedge\rho
-
i\,\overline{\pi}^2\wedge\kappa
+
\zeta\wedge\overline{\zeta},
\\
d\pi^2
&
\,=\,
\Lambda\wedge\kappa
+
\pi^2\wedge\overline{\pi}^1
-
\overline{\pi}^2\wedge\zeta
+
i\,H\,\rho\wedge\kappa.
\endaligned
\]

\begin{Lemma}
\label{Lemma-H-vanishes}
The function $H = -\, i\, h$ vanishes identically.
\end{Lemma}

\proof
Thanks to the Poincaré relation:
\[
\aligned
0
&
\,=\,
d\circ d\pi^1
\\
&
\,=\,
d\Lambda\wedge\rho
-
\Lambda\wedge d\rho
-
i\,d\overline{\pi}^2\wedge\kappa
+
i\,\overline{\pi}^2\wedge d\kappa
+
d\zeta\wedge\overline{\zeta}
-
\zeta\wedge d\overline{\zeta}
\\
&
\,=\,
d\Lambda\wedge\rho
-
\Lambda\wedge
\Big(
\big(\pi^1+\overline{\pi}^1\big)
\wedge
\rho
+
i\,
\kappa\wedge\overline{\kappa}
\Big)
\,-
\\
&
\ \ \ \ \
-\,
i\,\Big(
\Lambda\wedge\overline{\kappa}
+
\overline{\pi}^2\wedge\pi^1
-
\pi^2\wedge\overline{\zeta}
-
i\,H\,\rho\wedge\overline{\kappa}
\Big)
\wedge\kappa
\,+
\\
&
\ \ \ \ \
+
i\,\overline{\pi}^2
\wedge
\Big(
\pi^2\wedge\rho
+
\pi^1\wedge\kappa
+
\zeta\wedge\overline{\kappa}
\Big)
+
\Big(
i\,\pi^2\wedge\kappa
+
\pi^1\wedge\zeta
-
\overline{\pi}^1\wedge\zeta
\Big)
\wedge\overline{\zeta}
\,-
\\
&
\ \ \ \ \
-
\zeta
\wedge
\Big(
-i\,\overline{\pi}^2\wedge\overline{\kappa}
+
\overline{\pi}^1\wedge\overline{\zeta}
-
\pi^1\wedge\overline{\zeta}
\Big),
\endaligned
\]
we obtain:
\[
0
\,=\,
\Big(
\underbrace{
d\Lambda
-
\Lambda\wedge\pi^1
-
\Lambda\wedge\overline{\pi}^1
-
i\,\pi^2\wedge\overline{\pi}^2}_{=:\,\,\Phi}
\Big)
\wedge
\rho
+
H\,\rho\wedge\kappa\wedge\overline{\kappa}.
\]

Similarly:
\[
\aligned
0
&
\,=\,
d\circ d\pi^2
\\
&
\,=\,
d\Lambda\wedge\kappa
-
\Lambda\wedge d\kappa
+
d\pi^2\wedge\overline{\pi}^1
-
\pi^2\wedge d\overline{\pi}^1
-
d\overline{\pi}^2\wedge\zeta
+
\overline{\pi}^2\wedge d\zeta
\,+
\\
&
\ \ \ \ \
+
i\,dH\wedge\rho\wedge\kappa
+
i\,H\,d\rho\wedge\kappa
-
i\,H\,\rho\wedge d\kappa
\\
&
\,=\,
d\Lambda\wedge\kappa
-
\Lambda\wedge
\Big(
\pi^2\wedge\rho
+
\pi^1\wedge\kappa
+
\zeta\wedge\overline{\kappa}
\Big)
+
\Big(
\Lambda\wedge\kappa
-
\overline{\pi}^2\wedge\zeta
+
i\,H\,\rho\wedge\kappa
\Big)
\wedge\overline{\pi}^1
\,-
\\
&
\ \ \ \ \
-\,
\pi^2
\wedge
\Big(
\Lambda\wedge\rho
+
\overline{\zeta}\wedge\zeta
\Big)
-
\Big(
\Lambda\wedge\overline{\kappa}
+
\overline{\pi}^2\wedge\pi^1
-
\pi^2\wedge\overline{\zeta}
-
i\,H\,\rho\wedge\overline{\kappa}
\Big)
\wedge\zeta
\,+
\\
&
\ \ \ \ \ \ \ \ \ \ \ \ \ \ \ \ \ \ \ \ \ \ \ \ \ \ \ \ \ \ \
+
\overline{\pi}^2
\wedge
\Big(
i\,\pi^2\wedge\kappa
+
\pi^1\wedge\zeta
-
\overline{\pi}^1\wedge\zeta
\Big)
\\
&
\ \ \ \ \
+
i\,dH\wedge\rho\wedge\kappa
+
i\,H\,\big(\pi^1+\overline{\pi}^1\big)
\wedge\rho\wedge\kappa
-
i\,H\,\rho\wedge
\Big(
\pi^1\wedge\kappa
+
\zeta\wedge\overline{\kappa}
\Big),
\endaligned
\]
which simplifies and reorganizes as:
\[
\aligned
0
&
\,=\,
\Big(
\underbrace{
d\Lambda
-
\Lambda\wedge\pi^1
-
\Lambda\wedge\overline{\pi}^1
-
i\,\pi^2\wedge\overline{\pi}^2}_{=\,\,\Phi\,\,{\sf again!}}
\Big)
\wedge\kappa
\,+
\\
&
\ \ \ \ \
+
i\,dH\wedge\rho\wedge\kappa
+
2i\,H\,\rho\wedge\kappa\wedge\pi^1
+
2i\,H\,\rho\wedge\kappa\wedge\overline{\pi}^1
-
2i\,H\,\rho\wedge\zeta\wedge\overline{\kappa}.
\endaligned
\]

Then wedging this $(\centersmallbullet) \wedge\kappa$ 
gives:
\[
0
\,=\,
0
+
0
+
0
+
0
-
2i\,H\,\rho\wedge\zeta\wedge\overline{\kappa}\wedge\kappa,
\]
whence $H = 0$ necessarily.
\endproof

We thus have shown that the {\em real} $2$-form:
\[
\Phi
\,:=\,
d\Lambda
-
\Lambda\wedge\pi^1
-
\Lambda\wedge\overline{\pi}^1
-
i\,\pi^2\wedge\overline{\pi}^2
\]
satisfies:
\[
0
\,=\,
\Phi\wedge\rho
\,=\,
\Phi\wedge\kappa,
\]
hence by the Cartan Lemma, there exists a function $F$ such that:
\[
\Phi
\,=\,
F\,\rho\wedge\kappa.
\]
But since $\Phi = \overline{\Phi}$ is real, and since $\rho 
\wedge \kappa$ is linearly independent with its conjugate
$\rho \wedge \overline{\kappa}$, we necessarily have
$F = 0 = \overline{F}$, and therefore:
\[
d\Lambda
-
\Lambda\wedge\pi^1
-
\Lambda\wedge\overline{\pi}^1
-
i\,\pi^2\wedge\overline{\pi}^2
\,=\,
0.
\]

We summarized the result obtained under a synthetic form,
not formulated at the end
of~{\cite{Pocchiola-2013-addendum,
Merker-Pocchiola-2018-addendum}}.

\begin{Theorem}
\label{Theorem-reduction-e-structure-J-W-0}
After normalizations of the group parameters
${\sf f}$, ${\sf b}$, ${\sf d}$, the equivalence
problem for $2$-nondegenerate (constant) Levi rank $1$
$\mathcal{C}^\omega$ or $\mathcal{C}^\infty$ real hypersurfaces
$M^5 \subset \C^3$ conducts to $2$ fundamental
primary differential invariants 
denoted $J$ and $W$ by Pocchiola in~{\cite{Pocchiola-2013-addendum}}.

When both $J \equiv 0 \equiv W$ vanish identically,
on the $10$-dimensional manifold $M^5 \times G^4 \times \R$
equipped with coordinates:
\[
\big(
z_1,z_2,\overline{z}_1,\overline{z}_2,v
\big)
\times
\big(
{\sf c},\overline{\sf c},
{\sf e},\overline{\sf e}
\big)
\times
({\sf t}),
\]
the prolongation of the structure equations with constant
coefficients:
\[
\aligned
d\rho
&
\,=\,
\big(\pi^1+\overline{\pi}^1\big)
\wedge\rho
+
i\,\kappa\wedge\overline{\kappa},
\\
d\kappa
&
\,=\,
\pi^2
\wedge\rho
+
\pi^1\wedge\kappa
+
\zeta\wedge\overline{\kappa},
\\
d\zeta
&
\,=\,
\isqrt\,
\pi^2\wedge\kappa
+
\pi^1\wedge\zeta
-
\overline{\pi}^1\wedge\zeta,
\endaligned
\]
where $\pi^1$ and $\pi^2$ are modified-prolonged Maurer-Cartan
forms, conducts to structure equations also having constant
coefficients (no new curvature functions appear):
\[
\aligned
d\pi^1
&
\,=\,
\Lambda\wedge\rho
-
i\,\overline{\pi}^2\wedge\kappa
+
\zeta\wedge\overline{\zeta},
\\
d\pi^2
&
\,=\,
\Lambda\wedge\kappa
+
\pi^2\wedge\overline{\pi}^1
-
\overline{\pi}^2\wedge\zeta,
\endaligned
\]
where $\Lambda = \overline{\Lambda} = d{\sf t} + \cdots$
is a real $1$-form also having constant
coefficients structure:
\[
d\Lambda
=
\Lambda\wedge\pi^1
+
\Lambda\wedge\overline{\pi}^1
+
i\,\pi^2\wedge\overline{\pi}^2.
\]

Lastly, the $\{e\}$-structure
defined by the collection of these $10$ structure equations
(conjugates are implicit) coincides with the
Maurer-Cartan equations~{\cite{Pocchiola-2014}}
of the CR automorphism Lie group
of the model light cone $M_{\sf LC}$.\qed
\end{Theorem}

\linestop


\vfill\end{document}